\documentclass[reqno]{memo-l}

\newtheorem{thm}{{\bf Theorem}}[section]
\newtheorem{lem}[thm]{{\bf Lemma}}
\newtheorem{prop}[thm]{{\bf Proposition}}
\newtheorem{cor}[thm]{{\bf Corollary}}
\newtheorem{NN}[thm]{}
\theoremstyle{definition}\newtheorem{df}[thm]{{\bf Definition}}
\theoremstyle{definition}\newtheorem{rem}[thm]{{\bf Remark}}
\theoremstyle{definition}

\numberwithin{equation}{section}

\newcommand{\N}{\mathbb{N}}
\newcommand{\Z}{\mathbb{Z}}

\newcommand{\R}{\mathbb{R}}
\newcommand{\C}{\mathbb{C}}
\newcommand{\T}{\mathbb{T}}

\newcommand{\Aff}{\operatorname{Aff}}

\newcommand{\id}{\operatorname{id}}

\newcommand{\morp}{contractive completely positive linear map}

\newcommand{\hm}{homomorphism}
\newcommand{\dt}{\delta}
\newcommand{\ep}{\epsilon}
\newcommand{\andeqn}{\,\,\,{\text{and}}\,\,\,}
\newcommand{\rforal}{\,\,\,{\rm for\,\,\,all}\,\,\,}
\newcommand{\tforal}{\,\,\,\text{for\,\,\,all}\,\,\,}
\newcommand{\CA}{$C^*$-algebra}
\newcommand{\SCA}{$C^*$-subalgebra}

\newcommand{\af}{{\alpha}}
\newcommand{\bt}{{\beta}}
\newcommand{\di}{{\rm dist}}

\newcommand{\beq}{\begin{eqnarray}}
\newcommand{\eneq}{\end{eqnarray}}

\makeindex

\begin{document}

\frontmatter

\title{Approximate Homotopy of Homomorphisms from $C(X)$ into a Simple $C^*$-algebra}


\author{Huaxin Lin}
\curraddr{Department of Mathematics\\
University of Oregon\\
Eugene, OR 97403\\
USA}
\email{hlin@uoreogn.edu}
\thanks{}

\author{}
\address{}
\curraddr{}
\email{}
\thanks{}

\date{}



\begin{abstract}
Let $X$ be a  finite CW complex and let $h_1, h_2: C(X)\to A$ be two unital \hm s,
where $A$ is a unital \CA. We study the problem when $h_1$ and $h_2$ are approximately
homotopic.
We present a $K$-theoretical  necessary and
sufficient condition for them to be approximately homotopic under the assumption that $A$ is a unital separable simple
\CA\, of tracial rank zero, or $A$ is a unital purely infinite simple \CA. When
they are approximately homotopic, we also give a bound for  the
length of the homotopy. These  results are also extended to the
case that $h_1$ and $h_2$ are approximately multiplicative \morp
s.

Suppose that $h: C(X)\to A$ is a monomorphism and $u\in A$ is a
unitary (with $[u]=\{0\}$ in $K_1(A)$). We prove that, for any
$\ep>0,$ and any compact subset ${\mathcal F}\subset C(X),$ there
exists $\dt>0$ and a finite subset ${\mathcal G}\subset C(X)$
satisfying the following: if $\|[h(f), u]\|<\dt$ and
$\text{Bott}(h,u)=\{0\},$ then there exists a continuous
rectifiable path $\{u_t: t\in [0,1]\}$ such that
\beq\label{ab1}
u_0=u,\,\,\,u_1=1_A \andeqn \|[h(g),u_t]\|<\ep\,\rforal g\in {\mathcal  F}\andeqn t\in [0,1].
\eneq
Moreover,
\beq\label{ab2}
{\rm {Length}}(\{u_t\})\le 2\pi+\ep.
\eneq
We show that if $\text{dim}X\le 1,$ or $A$ is purely infinite
simple, then $\dt$ and ${\mathcal G}$ are universal (independent of
$A$ or $h$). In the case that ${\rm dim} X=1,$ this provides an
improvement of the so-called the Basic Homotopy Lemma of Bratteli,
Elliott, Evans and Kishimoto for the case that $A$ is mentioned
above.
 Moreover, we show that $\dt$ and ${\mathcal G}$ can not be universal whenever
 $\text{dim} X\ge 2.$ Nevertheless, we also found that  $\dt$ can be chosen to be  dependent on
 a measure distribution but independent of $A$ and $h.$
The above version of the so-called Basic Homotopy is also extended
to the case that $C(X)$ is replaced by an AH-algebra.

 We also present some general versions of so-called Super Homotopy
 Lemma.

\end{abstract}

\maketitle


\setcounter{page}{4}

\tableofcontents


\mainmatter
\chapter{Prelude}
\section{Introduction}

Let $A$ be a unital \CA\, and let $u\in A$ be a unitary which is
in the connected component $U_0(A)$ of the unitary group of $A$
containing the identity. Then there is a continuous path of
unitaries in $U_0(A)$ starting at $u$ and ending at $1.$ It is
known that the path can be made rectifiable. But, in general, the
length of the path has no bound. N. C. Phillips (\cite{Ph1})
proved that, in any unital  purely infinite simple \CA\, $A,$ if
$u\in U_0(A),$ then the length of the path from $u$ to the
identity can be chosen to be smaller than $\pi+\ep$ for any
$\ep>0.$ A more general result proved by the author shows that
this holds for any unital \CA s with real rank zero (\cite{LnSF}).

Let $X$ be a path connected finite CW complex. Fix a point $\xi
\in X.$ Let $A$ be a unital \CA. There is a trivial \hm\, $h_0:
C(X)\to A$ defined by $h_0(f)=f(\xi)1_A$ for $f\in C(X).$
 Suppose that $h: C(X)\to
A$ is a unital monomorphism. When can $h$  be homotopic to $h_0?$
When $h$ is homotopic to $h_0,$ how long could the length of the
homotopy be? This is one of the questions that motivates this
paper.

Let $h_1, h_2: C(X)\to A$ be two unital \hm s. A more general
question is when $h_1$ and $h_2$ are homotopic? When $A$ is
commutative, by the Gelfand transformation, it is a purely
topological homotopy question. We only consider noncommutative
cases. To the other end of noncommutativity, we only consider the
case that $A$ is a unital simple \CA.  To be possible and useful,
we actually consider approximate homotopy. So we study the problem
when $h_1$ and $h_2$ are approximately homotopic. For all
applications that we know, the length of the homotopy is extremely
important. So we also ask how long the homotopy is if $h_1$ and
$h_2$ are actually approximately homotopic.

Let $X$ be a path connected metric space. Fix a point $\xi_X.$ Let
$Y_X=X\setminus \xi_X.$ Let $x\in X$ and let $L(x, \xi_X)$ be the
infimum of the length of continuous paths from $x$ to $\xi_X.$
Define
$$
L(X, \xi_X)=\sup\{L(x,\xi_X): x\in X\}.
$$
We prove that, for any unital simple \CA\, of tracial rank zero,
or any unital purely infinite simple \CA\, $A,$ if $h: C(X)\to A$
is a unital monomorphism with
\beq\label{kk0}
[h|_{C_0(Y_X)}]=\{0\}\, \,\,\,{ \rm in}\,\,\, KL(C_0(Y_X), A),
\eneq
 then, for any $\ep>0$ and any compact subset ${\mathcal F}\subset C(X),$
  there is a \hm\,
$H: C(X)\to C([0,1],A)$ such that
$$
\pi_0\circ H\approx_{\ep} h\,\,\,{\rm on}\,\,\, {\mathcal F},\,\,\,
\pi_1\circ H=h_0\andeqn
$$
$$
\overline{\rm{Length}}(\{\pi_t\circ H\})\le L(X, \xi_X),
$$
where $h_0(f)=f(\xi_X)1_A$ for all $f\in C(X)$ and $\pi_t:
C([0,1],A)\to A$ is the point-evaluation at $t\in [0,1],$ and
where $\overline{\rm{Length}}(\{\pi_t\circ H\})$ is appropriately
defined. Note that $[h_0|_{C_0(Y_X)}]=\{0\}.$ Thus the condition
(\ref{kk0}) is necessary. Moreover, the estimate of length can not
be improved.

Suppose that $h_1, h_2: C(X)\to A$
are two unital \hm s. We show that  $h_1$ and $h_2$ are approximately
homotopic if and only if
$$
[h_1]=[h_2]\,\,\,\text{in}\,\,\,KL(C(X),A),
$$
under the assumption that $A$ is a unital separable simple \CA\, of tracial rank zero, or
$A$ is a unital purely infinite simple \CA.
Moreover, we show that the length of the homotopy can be bounded
by a universal constant.

Bratelli, Elliott, Evans and Kishimoto  (\cite{BEEK}) considered
the following homotopy question: Let $u$ and $v$ be two unitaries
such that $u$ almost commutes with $v.$ Suppose that $v\in
U_0(A).$ Is there a rectifiable continuous path $\{v_t: t\in
[0,1]\}$ with $v_0=v,$ $v_1=1_A$ such that the entire path almost
commutes with $u?$ They found that there is an additional obstacle
to prevent the existence of such path of unitaries. The additional
obstacle is the Bott element \text{bott}(u,v) associated with the
pair $u$ and $v.$ They proved what they called the Basic Homotopy
Lemma: For any $\ep>0$ there exists $\dt>0$ satisfying the
following: if $u, v$ are two unitaries in a unital separable
simple \CA\, with real rank zero and stable rank one, or in a
unital separable purely infinite simple \CA\, $A,$ if $v\in
U_0(A)$ and ${\rm sp}(u)$ is $\dt$-dense in $S^1$ except  possibly for a single  gap,
$$
\|uv-vu\|<\dt \andeqn \text{bott}_1(u,v)=0,
$$
then there exists a rectifiable continuous path of unitaries
$\{v_t: t\in [0,1]\}$ such that
$$
v_0=v,\,\,\, v_1=1_A\andeqn \|uv_t-v_tu\|<\ep
$$
for all $t\in [0,1].$ Moreover,
$$
{\rm{Length}}(\{v_t\})\le 4\pi+1.
$$

Bratteli, Elliott, Evans and Kishimoto were motivated by the study of
 classification of
purely infinite simple \CA s. The Basic Homotopy Lemma played an
important role in their work related to the classification of
purely infinite simple \CA s and that of \cite{ER}. The renewed
interest of this type of results is at least partly motivated by
the study of automorphism groups of  simple \CA s (see \cite{KK}).
It is also important in the  study of AF-embedding of crossed
products (\cite{Ma}).

We now replace the unitary $u$ in the Basic Homotopy Lemma by a
monomorphism $h: C(X)\to A.$  We first replace $S^1$ by a path
connected finite CW complex $X.$ A bott element
$\text{bott}_1(h,v)$ can be similarly defined. We proved that,
with the  assumption $A$ is a unital simple \CA\, of real rank
zero and stable rank one, or $A$ is a unital purely infinite
simple \CA\, the  Basic Homotopy Lemma holds for any compact
metric space with dimension no more than one. Moreover, in the
case that $K_1(A)=\{0\},$ the constant $\dt$ does not depend on
the spectrum of $h$ (so that the condition on the spectrum of $u$
in original Homotopy Lemma can be removed). The proof is shorter
than that of the original Homotopy Lemma of Bratteli, Elliott,
Evans and Kishimoto. Furthermore, we are able to cut the length of
homotopy by more than half (see \ref{dmyi} and \ref{Tpi}).

For more general compact metric space, the simple bott element has
to be replaced by a more general map $\text{Bott}(h,v).$ Even with
vanishing $\text{Bott}(h,v)$ and with $A$ having tracial rank zero,
we show that the same statement is false whenever ${\rm dim} X\ge
2.$ However, if we allow the constant $\dt$ not only depends on
$X$ and $\ep$ but also depends on a measure distribution, then the
similar homotopy result holds (see \ref{MT1}) for unital separable
simple \CA s with tracial rank zero. On the other hand, if $A$ is
assumed to be purely infinite simple, then there is no such
measure distribution. Therefore, for purely infinite simple \CA s,
the Basic Homotopy Lemma holds for any compact metric space with
shorter lengths. In fact  our estimates on the lengths is
$2\pi+\ep$ (for the case that $A$ is purely infinite simple as well as for the case
that $A$ is a unital separable simple \CA\, with tracial rank zero).

Several other homotopy results are also discussed. In particular,
a version of Super Homotopy Lemma (of Bratelli, Elliott, Evans and
Kishimoto) for  finite CW complex $X$ is also presented.

The presentation is organized as follows:

In section 2, we provide some conventions and a number of facts
which will be used later.

In section 3, we present the Basic Homotopy Lemma for $X$ being a
compact metric space with covering dimension no more than  one
under the assumption that $A$ is a unital simple \CA\, of real rank
zero and stable rank one, or $A$ is a unital purely infinite
simple \CA. The improvement is made not only on the bound of the
length but is also made so that the constant $\dt$ does not depend
on the spectrum of the \hm s.

In section 4 and section 5, we present some results which are
preparations for  later sections.

In section 6, 7 and 8, we prove a version of the Basic Homotopy
Lemma for general compact metric space under the assumption that
$A$ is a unital separable simple \CA\, of tracial rank zero. The
lengthy proof is due partly to the complexity caused by our
insistence that the constant $\dt$ should not be dependent on \hm
s or $A$ but only on a measure distribution.

In section 9, we show why the constant $\dt$ can not be made
universal as in the dimension 1 case. A hidden topological obstacle is revealed.
We show that the original
version of the Basic Homotopy Lemma fails whenever $X$ has dimension at
least two for   simple \CA s with real rank
zero and stable rank one.

In section 10, we present some familiar results about purely
infinite simple \CA s.

In section 11, we show that the Basic Homotopy Lemma holds for
general compact metric spaces under the assumption that $A$ is a
unital purely infinite simple \CA.

In section 12, we discuss the length of homotopy. A definition
related to Lipschitz functions is given there and some elementary facts are also given.

In section 13, we show two \hm s are approximately homotopic when they induce the same
$KL$ element under the assumption that $A$ is a unital separable simple \CA\, of tracial rank zero,
or $A$ is a unital purely infinite simple \CA. We also give an estimate on the bound of the length of the homotopy.

In section 14, we extend the results in section 13 to the maps which are not necessary homomorphisms.

In section 15, we present a version of the so-called Super Homotopy Lemma  for unital  purely infinite simple
\CA\, $A.$

In section 16, we show that same version of the Super Homotopy
Lemma is valid for unital separable simple \CA\, of tracial rank
zero.

In section 17, we show that the Basic Homotopy Lemma in section 8
and 11 is valid if we replace $C(X)$ by a unital AH-algebra.

In section 18, we end this paper with a few  concluding remarks.

\vspace{0.2in}

{\bf Acknowledgment}

The most of this research was done in the summer 2006 when the
author was in East China Normal University
 where he had a nice office and necessary computer
equipments. It is partially supported by Shanghai Priority
Academic Disciplines. The work is also partially supported by a
NSF grant (00355273).


\section{Conventions and some facts}

\begin{NN}

{\rm Let $A$ be a \CA. Using notation introduced by \cite{DL}, we
denote
$$
\underline{K}(A)=\oplus_{i=0,1}(K_i(A)\bigoplus \oplus_{n=1}
K_i(A, \Z/n\Z)).
$$
Let $m\ge 1$ be an integer. Denote by $C_m$ a commutative \CA\,
with $K_0(C_m)=\Z/m\Z$ and $K_1(C_m)=0.$
So $K_i(A,\Z/m\Z)=K_i(A\otimes C_m),$
$i=0,1.$

A theorem of Dadarlat and Loring (\cite{DL}) states that
$$
Hom_{\Lambda}(\underline{K}(A), \underline{K}(B))\cong KL(A,B),
$$
if $A$ satisfies the Universal Coefficient Theorem and $B$ is
$\sigma$-unital (see \cite{DL} for the definition of
$Hom_{\Lambda}(\underline{K}(A), \underline{K}(B))$). We will
identify these two objects.

Let $m\ge 1$ be an integer. Put
$$
F_m\underline{K}(A)=\oplus_{i=0,1}(K_i(A)\bigoplus
\oplus_{k|m}K_i(A, \Z/k\Z)).
$$

We will also identify these two objects. }
\end{NN}

\begin{NN}
{\rm Let $B_n$ be a sequence of \CA s. Denote by
$l^{\infty}(\{B_n\})$ the product of $\{B_n\},$ i.e., the \CA\, of
all bounded sequences $\{a_n: a_n\in B_n\}.$ Denote by
$c_0(\{B_n\})$ the direct sum of $\{B_n\},$ i.e, the \CA\, of all
sequences $\{a_n: a_n\in B_n\}$ for which
$\lim_{n\to\infty}\|a_n\|=0.$ Denote by
$q_{\infty}(\{B_n\})=l^{\infty}(\{B_n\})/c_0(\{B_n\})$ and by $q:
l^{\infty}(\{B_n\})\to q_{\infty}(\{B_n\})$ the quotient map.

}
\end{NN}

\begin{NN}\label{pKK}
{\rm Let $A$ be a \CA\, and let $B$ be another \CA. Let $\ep>0$
and ${\mathcal G}\subset A$ be a finite subset. We say that a \morp\,
$L:A\to B$ is $\dt$-${\mathcal G}$-multiplicative if
$$
\|L(ab)-L(a)L(b)\|<\dt\rforal a, b\in {\mathcal G}.
$$

Denote by ${\bf P}(A)$ the set of projections and unitaries in
$$
M_{\infty}({\tilde{A}})\cup\, \bigcup_{m\ge
1}M_{\infty}({\widetilde {A\otimes C_m}}).
$$

We also use $L$ for the map $L\otimes {\id_{M_k\otimes C_m}}:
A\otimes C_m\otimes M_k\to B\otimes C_m\otimes M_k,$ $k=1,2,...,.$
As in 6.1.1 of \cite{Lnbk}, for a fixed $p\in {\bf P}(A),$ if $L$
is $\dt$-${\mathcal G}$-multiplicative with sufficiently small
$\dt$ and sufficiently large ${\mathcal G},$ $L(p)$ is close to a
projection (with the norm of difference is less than $1/2$) which
will be denote by $[L(p)].$ Note if two projections are both close
to $L(p)$ within $1/2,$ they are equivalent.

If $L: A\to B$ is $\dt$-${\mathcal G}$-multiplicative, then there
is a finite subset ${\mathcal Q}\subset {\bf P}(A),$ such that
$[L](x)$ is well defined for $x\in \overline{\mathcal Q},$ where
$\overline{\mathcal Q}$ is the image of ${\mathcal Q}$ in
$\underline{K}(A),$ which means that if $p_1, p_2\in {\mathcal Q}$
and $[p_1]=[p_2],$ then $[L(p_1)]$ and $[L(p_2)]$ defines the same
element in $\underline{K}(B).$ Moreover, if $p_1, p_2, p_1\oplus
p_2\in {\mathcal Q},$ $[L(p_1\oplus p_2)]=[L(p_1)]+[L(p_2)]$  (see
0.6 of \cite{Lnalm} and 4.5.1 and 6.1.1 of \cite{Lnbk}). This
finite subset ${\mathcal Q}$ will be denoted by ${\mathcal
Q}_{\dt,{\mathcal G}}.$ Let ${\mathcal P}\subset
\underline{K}(A).$ We say $[L]|_{\mathcal P}$ is well defined, if
$\overline{{\mathcal Q}_{\dt, {\mathcal G}}}\supset {\mathcal P}.$
In what follows, whenever we write $[L]|_{\mathcal P},$ we mean
that $[L]|_{\mathcal P}$ is well defined (see also 2.4 of
\cite{D2} for further explanation).

}
\end{NN}

The following  proposition is  known and has been implicitly used
many times.

\begin{prop}\label{apK1}
Let $A$ be a  separable \CA\, for which $K_i(A)$ is finitely
generated (for $i=0,1$), and let ${\mathcal P}\subset
\underline{K}(A)$ be a finite subset. Then, there is $\dt>0$ and a
finite subset ${\mathcal G}\subset A$ satisfying the following: If
$B$ is a unital \CA\, and if $L: A \to B$ is a $\dt$-${\mathcal
G}$-multiplicative \morp, there is an element $\kappa\in
Hom_{\Lambda}(\underline{K}(A), \underline{K}(B))$ such that
\beq\label{apK1-1}
[L]|_{\mathcal P}=\kappa|_{\mathcal P}.
\eneq

Moreover,  there is a finite subset ${\mathcal P}_A\subset
\underline{K}(A)$ such that, if $[L]|_{{\mathcal P}_A}$ is well
defined, there is a unique $\kappa\in
Hom_{\Lambda}(\underline{K}(A), \underline{K}(B))$ such that {\rm
(\ref{apK1-1})} holds.
\end{prop}

\begin{proof}
Since  $K_i(A)$ is finitely generated ($i=0,1$), by 2. 11 of
\cite{DL},
$$
Hom_{\Lambda}(\underline{K}(A), \underline{K}(B))=
Hom_{\Lambda}(F_m\underline{K}(A), F_m\underline{K}(B))
$$
for some $m\ge 1.$ Thus it is clear that it suffices to show the
first part of the proposition. Suppose that the first part of the
lemma fails. One obtains a finite subset ${\mathcal P}\subset
\underline{K}(A),$  a sequence of $\sigma$-unital \CA s $B_n, $  a
sequence of positive numbers $\{d_n\}$ with
$\sum_{n=1}^{\infty}\dt_n<\infty,$ a finite subsets ${\mathcal
G}_n\subset A$ with $\cup_{n=1}^{\infty} {\mathcal G}_n$ is dense
in $A,$ and a sequence of $\dt_n$-${\mathcal G}_n$-multiplicative
\morp s $L_n: A\to B_n$ such that there exists no $\kappa\in
Hom_{\Lambda}(\underline{K}(A), \underline{K}(B))$ satisfying
(\ref{apK1-1}).

Define $\Phi: A\to l^{\infty}(\{B_n\otimes {\mathcal K}\})$ by
$\Phi(a)=\{L_n(a)\}$ for $a\in A$ and define ${\bar \Phi}: A\to
q_{\infty}(\{B_n\otimes {\mathcal K}\})$ by ${\bar \Phi}=\pi\circ
\Phi,$ where $\pi: l^{\infty}(\{B_n\otimes {\mathcal K}\})\to
q_{\infty}(\{B_n\otimes {\mathcal K}\})$ is the quotient map. Thus
we obtain an element $\af\in Hom_{\Lambda}(F_m\underline{K}(A),
F_m\underline{K}(q_{\infty}(\{B_n\otimes {\mathcal K}\})))$ such
that $[\bar \Phi]=\af.$ Since $K_i(A)$ is finitely generated
($i=0,1$), by 2.11 of \cite{DL}, there is an integer $m\ge 1$ such
that
$$
Hom_{\Lambda}(\underline{K}(A),
\underline{K}(q_{\infty}(\{B_n\otimes {\mathcal K}\})))\cong
Hom_{\Lambda}(F_m\underline{K}(A),
F_m\underline{K}(q_{\infty}(\{B_n\otimes {\mathcal K}\}))\andeqn
$$
$$
Hom_{\Lambda}((\underline{K}(A), \underline{K}(B_n))\cong
Hom_{\Lambda}(F_m\underline{K}(A), F_m\underline{K}(B_n)).
$$
By applying 7.2 of \cite{Lnbdf} and the proof of 7.5 of
\cite{Lnbdf}, for all larger $n,$ there is an element
$$\kappa_n\in
Hom_{\Lambda}(F_m\underline{K}(A), F_m\underline{K}(B_n))$$ such
that
$$
[L_n]|_{\mathcal P}=\kappa_n|_{\mathcal P}.
$$
This contradicts the assumption that the first part of the lemma
fails.

\end{proof}

The following is well known and follows immediately from the definition.
\begin{prop}\label{Khp-1}
Let $A$ be a unital amenable \CA. For any finite subset ${\mathcal P}\subset \underline{K}(A),$
there exists $\dt>0$ and a finite subset ${\mathcal G}\subset A$ satisfying the following:
for any pair of $\dt$-${\mathcal G}$-multiplicative \morp s $L_1, L_2: A\to B$ (for any unital \CA\, $B$),
$$
[L_1]|_{\mathcal P}=[L_2]|_{\mathcal P}
$$
provided that
$$
L_1\approx_{\dt}L_2\,\,\,\text{on}\,\,\, {\mathcal G}.
$$
\end{prop}

\begin{NN}

{\rm Let $B$ be a \CA\, and $C=C([0,1],B).$ Define $\pi_t: C\to B$
by $\pi_t(f)=f(t)$ for all $f\in C.$ This notation will be used
throughout this article.

}
\end{NN}

The following follows immediately from \ref{Khp-1} and will be
used frequently without further notice.

\begin{prop}\label{Khp}
Let $A$ and $B$ be two  unital \CA s and let $L: A\to B$ be a
\morp. Let ${\mathcal Q}\subset {\bf P}(A)$ be a finite subset.
Suppose that, for some small $\dt>0$ and a large finite subset
${\mathcal G},$ $L$ is $\dt$-${\mathcal G}$-multiplicative and ${\mathcal
Q}_{\dt,{\mathcal G}}\supset {\mathcal Q}.$ Put ${\mathcal P}=\overline{\mathcal
Q}$ in $\underline{K}(A).$ Suppose that $H: A\to C([0,1], B)$ is a
\morp\, such that $\pi_0\circ H=L$ and $\pi_t\circ H$ is
$\dt$-${\mathcal G}$-multiplicative for each $t\in [0,1].$ Then, for
each $t\in [0,1],$
$$
[\pi_t\circ H]|_{\mathcal P}=[L]|_{\mathcal P}.
$$

\end{prop}

The following follows immediately from 2.1 of \cite{LP}.

\begin{lem}\label{appn}
Let $B$ be a separable amenable \CA. For any $\ep>0$ and any
finite subset ${\mathcal F}_{0}\subset B$
there exists a finite subset ${\mathcal F}_1\subset B$ and $\dt>0$
satisfying the following: Suppose that $A$ is a unital \CA, $\phi:
B\to A$ is a unital \hm\,  and $u\in A$ is a unitary such that
\beq\label{appn1}
\|[\phi(a), u]\|<\dt\,\rforal\, a\in {\mathcal F}_1.
\eneq
Then there is an $\ep$-${\mathcal F}_0\otimes S$-multiplicative
\morp\,\\ $\psi: B\otimes C(S^1)\to A$ such that
\beq\label{appn2}
\|\phi(a)-\psi(a)\|<\ep\andeqn \|\psi(a\otimes g
)-\phi(a)g(u)\|<\ep
\eneq
for all $a\in {\mathcal F}_0$ and $g\in S,$ where $S=\{1_{C(S^1)},
z\}$ and $z\in C(S^1)$ is the standard unitary generator of
$C(S^1).$
\end{lem}

\begin{NN}

{\rm Let $A$ be a unital \CA. Denote by $U(A)$ the group of all
unitaries in $A.$ Denote by $U_0(A)$ the path connected component
of $U(A)$ containing $1_A.$

Denote by ${\rm Aut}(A)$ the group of automorphisms on $A.$ If
$u\in U(A),$ denote by ${\rm ad}\, u$ the inner automorphism
defined by ${\rm ad}\, u(a)=u^*au$ for all $a\in A.$

}

\end{NN}

\begin{df}\label{Dbot2}
{\rm Let
$A$ and $B$ be  two unital \CA s.  Let $h: A\to B$ be a \hm\, and
$v\in U(B)$ such that
$$
h(g)v=vh(g)\,\rforal\, g\in A.
$$
 Thus we
obtain a \hm\, ${\bar h}: A\otimes C(S^1)\to B$ by ${\bar
h}(f\otimes g)=h(f)g(v)$ for $f\in A$ and $g\in C(S^1).$ The
tensor product induces two injective \hm s:
\beq\label{dbot01}
\bt^{(0)}&:& K_0(A)\to K_1(A\otimes C(S^1))\\
 \bt^{(1)}&:&
K_1(A)\to K_0(A\otimes C(S^1)).
\eneq
The second one is the usual Bott map.
Note, in this way, one write $K_i(A\otimes C(S^1))=K_i(A)\oplus \bt^{(i-1)}(K_{i-1}(A)).$
We use $\widehat{\bt^{(i)}}: K_i(A\otimes C(S^1))\to \bt^{(i-1)}(K_{i-1}(A))$ for the  the projection
to $\bt^{(i-1)}(K_{i-1}(A)).$

For each integer $k\ge 2,$
one also obtains the following injective \hm s:
\beq\label{dbot02}
\bt^{(i)}_k: K_i(A, \Z/k\Z))\to K_{i-1}(A\otimes C(S^1), \Z/k\Z),
i=0,1.
\eneq
Thus we write
\beq\label{dbot02-1}
K_{i-1}(A\otimes C(S^1), \Z/k\Z)=K_{i-1}(A,\Z/k\Z)\oplus \bt^{(i)}_k(K_i(A, \Z/k\Z)),\,\,i=0,1.
\eneq
Denote by $\widehat{\bt^{(i)}_k}: K_{i}(A\otimes C(S^1),
\Z/k\Z)\to \bt^{(i-1)}_k(K_{i-1}(A,\Z/k\Z))$ similarly to that of
$\widehat{\bt^{(i)}}.,$ $i=1,2.$ If $x\in \underline{K}(A),$ we
use ${\boldsymbol{\beta}}(x)$ for $\bt^{(i)}(x)$ if $x\in K_i(A)$
and for $\bt^{(i)}_k(x)$ if $x\in K_i(A, \Z/k\Z).$ Thus we have a
map ${\boldsymbol{ \bt}}: \underline{K}(A)\to
\underline{K}(A\otimes C(S^1))$ as well as
$\widehat{\boldsymbol{\bt}}: \underline{K}(A\otimes C(S^1))\to
 {\boldsymbol{ \bt}}(\underline{K}(A)).$ Thus one may write
 $\underline{K}(A\otimes C(S^1))=\underline{K}(A)\oplus {\boldsymbol{ \bt}}( \underline{K}(A)).$

On the other hand ${\bar h}$ induces \hm s ${\bar h}_{*i,k}:
K_i(A\otimes C(S^1)), \Z/k\Z)\to K_i(B,\Z/k\Z),$ $k=0,2,...,$ and
$i=0,1.$

We use $\text{Bott}(h,v)$ for all  \hm s ${\bar h}_{*i,k}\circ
\bt^{(i)}_k.$ We write
$$
\text{Bott}(h,v)=0,
$$
if ${\bar h}_{*i,k}\circ \bt^{(i)}_k=0$ for all $k\ge 1$ and
$i=0,1.$

We will use $\text{bott}_1(h,v)$ for the \hm\, ${\bar
h}_{1,0}\circ \bt^{(1)}: K_1(A)\to K_0(B),$ and
$\text{bott}_0(h,u)$ for the \hm\, ${\bar h}_{0,0}\circ \bt^{(0)}:
K_0(A)\to K_1(B).$

Since $A$ is unital, if $\text{bott}_0(h,v)=0,$ then $[v]=0$ in
$K_1(B).$

In what follows, we will use $z$ for the standard generator of
$C(S^1)$ and we will often identify $S^1$ with the unit circle
without further explanation. With this identification $z$ is the
identity map from the circle to the circle.

Now let $A=C(S^1)$ and $u=h(z).$
 Then
$\text{bott}(u,v)=\text{bott}_1(h,u)([z]).$ Note that, if $[v]=0$
in $K_1(B),$ then $\text{bott}_0(h,v)=0.$ In this case,
$K_i(C(S^1))$ is free, thus ${\bar h}_{*i,k}=0,$ if $k\ge 2.$ In
particular, if $[v]=0$ in $K_1(B)$ and $\text{bott}_1(h,v)=0,$
then
$$
\text{Bott}(h,v)=0.
$$

Suppose that $\{v_n\}$ is a sequence of unitaries in $A$ such that
$$
\lim_{n\to\infty} \|h(a)v_n-v_nh(a)\|=0\,\rforal\, a\in A.
$$
Then we obtain a sequential asymptotic morphism $ \psi_n: A\otimes
S^1\to B$ such that
$$\lim_{n\to\infty}\|\psi_n(a\otimes
g)-h(a)g(v_n)\|=0\,\,\,\rforal a\in A \andeqn g\in C(S^1)
$$ (see
\ref{appn}).

Fix a finite subset ${\mathcal P}\subset \underline{K}(A),$ for sufficiently large $n,$
as in \ref{apK1}, $[\psi_n]|_{\bt(\mathcal P)}$ is well defined. We will denote this by
$$
\text{Bott}(h, v_n)|_{\mathcal P}.
$$
 In other words,  for a fixed finite subset ${\mathcal P}\subset \underline{K}(A),$
 there exists $\dt>0$ and a finite subset ${\mathcal
G}\subset A$ such that, if $v\in B$ is a unitary for which
$$
\|h(a)v-vh(a)\|<\dt\,\rforal\, a\in {\mathcal G},
$$
then $\text{Bott}(h,v)|_{{\mathcal P}}$ is well defined. In what follows, whenever we write
$\text{Bott}(h,v)|_{\mathcal P},$ we mean that $\dt$ is sufficiently small and ${\mathcal G}$ is sufficiently large
so it is well defined.

Now suppose that $A$ is also amenable and $K_i(A)$ is finitely
generated ($i=0,1$). For example, $A=C(X),$ where $X$ is a finite
CW complex. So, for all sufficiently large $n,$ $\psi_n$ defines
an element $[\psi_n]$ in $KL(A, B)$ (see \ref{apK1}).
 Therefore, for a fixed finite subset ${\mathcal P}_0\subset \underline{K}(A),$
 there exists $\dt_0>0$ and a finite subset ${\mathcal
G}_0\subset A$ such that, if $v\in B$ is a unitary for which
$$
\|h(a)v-vh(a)\|<\dt_0\,\rforal\, a\in {\mathcal G}_0,
$$
then $\text{Bott}(h,v)|_{{\mathcal P}_0}$ is well defined. Thus when
$K_i(A)$ is finitely generated,  $\text{Bott}(h,v)|_{{\mathcal P}_0}$
defines $\text{Bott}(h,v)$ for some sufficiently large finite
subset ${\mathcal P}_0.$  In what follows such ${\mathcal P}_0$ may be
denoted by ${\mathcal P}_A.$ Suppose that ${\mathcal P}\subset
\underline{K}(A)$ is a larger finite subset, and ${\mathcal G}\supset
{\mathcal G}_0$ and $0<\dt<\dt_0.$

{\it A fact that we be used in this paper is that,
$\text{Bott}(h,v)|_{\mathcal P}$ defines the same map
$\text{Bott}(h,v)$ as $\text{Bott}(h,v)|_{{\mathcal P}_0}$ defines,
if}
$$
\|h(a)v-vh(a)\|<\dt\rforal \, a\in {\mathcal G}.
$$

In what follows, in the case that $K_i(A)$ is finitely generated,
whenever we write $\text{Bott}(h,v),$ we always assume that $\dt$
is smaller than $\dt_0$  and ${\mathcal G}$ is  larger than ${\mathcal
G}_0$  so that $\text{Bott}(h,v)$ is well-defined.

}

\end{df}

\begin{NN}\label{botsmall}

{\rm
In the case that $A=C(S^1),$ there is a concrete way to visualize $\text{bott}_1(h,v).$
It is perhaps helpful to describe it here.
The map $\text{bott}_1(h,v)$ is determined by $\text{bott}_1(h, v)([z]),$ where $z$ is, again,
the identity map on the unit circle.

Denote $u=h(z)$ and
define
\beq\nonumber
f(e^{2\pi i t})=\begin{cases} 1-2t, &\text{if $0\le t\le 1/2$,}\\
  -1+2t, & \text{if $1/2<t\le 1$,}
  \end{cases}
  \eneq
  \beq\nonumber
g(e^{2\pi i t})=\begin{cases} (f(e^{2\pi i t})-f(e^{2\pi it})^2)^{1/2} &\text{if $0\le t\le 1/2$,}\\
 0 , & \text{if $1/2<t\le 1$ \,\,\,and}
  \end{cases}
  \eneq
\beq\nonumber
h(e^{2\pi i t})=\begin{cases}0 , &\text{if $0\le t\le 1/2$,}\\
 (f(e^{2\pi i t})-f(e^{2\pi it})^2)^{1/2} , & \text{if $1/2<t\le 1$,}
  \end{cases}
  \eneq
These are non-negative continuous functions defined on the unit
circle. Suppose that $uv=vu.$ Define
\beq\label{bot1}
e(u,v)=\begin{pmatrix} f(v) & g(v)+h(v)u^*\\
                                     g(v)+uh(v) & 1-f(v)
           \end{pmatrix}
\eneq

Then $e(u,v)$ is a projection. There is $\dt_0>0$ (independent of
unitaries $u, v$ and $A$)  such that if $\|[u,v]\|<\dt_0,$ the
spectrum of positive element $e(u,v)$ has a gap at $1/2.$ The bott
element of $u$ and $v$ is an element in $K_0(A)$ as defined by
in \cite{EL1}  and \cite{EL2} is
\beq\label{bot2}
{\rm bott}_1(u,v)=[\chi_{1/2, \infty}(e(u,v))]-[\begin{pmatrix} 1 & 0\\
0 & 0
\end{pmatrix}].
 \eneq

Note that $\chi_{1/2, \infty}$ is a continuous function on ${\rm
sp}(e(u,v)).$ Suppose that \\${\rm sp}(e(u,v))\subset (-\infty, a]
\cup [1-a,\infty)$ for some $0<a<1/2.$ Then $\chi_{1/2, \infty}$
can be replaced by any other positive continuous function $F$ such
that $F(t)=0$ if $t\le a$ and $F(t)=1$ if $t\ge 1/2.$

}

\end{NN}

\begin{NN}
{\rm

Let $X$ be a finite CW complex and let $A$ be a unital separable
stably finite \CA. Recall that an element $\af\in
Hom_{\Lambda}(\underline{K}(C(X)), \underline{K}(A))$ is said to
be positive, if
$$\af(K_0(C(X))_+\setminus \{0\})\subset
K_0(A)_+\setminus \{0\}.$$ These positive elements may be denoted
by $KL(C(X), A)_+.$

}
\end{NN}

\begin{NN}\label{Y-X}
{\rm Let $X$ be a connected finite CW complex and let $\xi_X\in X$
be a base point. Put $Y_X=X\setminus \{\xi_X\}.$ }
\end{NN}

\begin{NN}
{\rm Let $X$ be a compact metric space with metric $\di(-,-).$ In
what follows, we will use
$$
\di((x,t),(y,s))=\sqrt{\di(x,y)^2+|t-s|^2}
$$
for $(x,t), (y,s)\in X\times S^1.$ }
\end{NN}

\begin{NN}
{\rm  Let $A$ and $B$ be two \CA s and $\phi: A\to B$ be a
contractive positive linear map. Let ${\mathcal G}\subset A$ be a
subset and let $\sigma>0.$ We say that $\phi$ is $\sigma$-${\mathcal
G}$-injective if
$$
\|\phi(a)\|\ge \sigma \|a\|
$$
for all $a\in {\mathcal G}.$

}
\end{NN}

\begin{NN}\label{Dspec}
{\rm Let $X$ be a compact metric space and let $A$ be a unital
\CA. Suppose that $\phi: C(X)\to A$ is a \hm. There is a compact
subset $F\subset X$ such that
$$
{\rm ker}\phi=\{f\in C(X): f|_F=0\}.
$$
Thus there is a monomorphism $h_1: C(F)\to A$ such that
$h=h_1\circ \pi,$ where $\pi: C(X)\to C(F)$ is the quotient map.
We say the spectrum of $\phi$ is $F.$ }
\end{NN}

\begin{NN}
{\rm Let $C=A\otimes B,$ where $A$ and $B$ are unital \CA s. Let
$D$ be another unital \CA. Suppose that $\phi: C\to D$ is a \hm.
Throughout this work, by $\phi|_A: A\to D$ we mean the \hm\,
defined by $\phi|_A(a)=\phi(a\otimes 1)$ for all $a\in A.$ }

\end{NN}

\begin{NN}
{\rm

Let $L_1, L_2: A\to B$ be two maps, let $\ep>0$  and ${\mathcal F}\subset A$ be a subset.
We write
$$
L_1\approx_{\ep} L_2\,\,\,\text{on}\,\,\,{\mathcal F},
$$
if
$$
\|L_1(f)-L_2(f)\|<\ep\rforal f\in {\mathcal F}.
$$
}
\end{NN}

\begin{NN}\label{sigma}
{\rm Let $X$ be a compact metric space, $\dt>0$ and ${\mathcal
F}\subset C(X)$ be a finite subset. Define $\sigma=\sigma_{X,\dt,
{\mathcal F}}$ to be the largest positive number satisfying the
following:
$$
|f(x)-f(x')|<\dt\rforal f\in {\mathcal F},
$$
provided ${\rm dist}(x, x')<\sigma.$

}
\end{NN}

\begin{NN}

{\rm Let $\tau$ be a  state on $C(X).$  Denote by
 $\mu_{\tau}$  the probability Borel measure induced by  $\tau.$

}

\end{NN}

\begin{NN}
{\rm Let $X$ be a metric  space and let $x\in X.$ Suppose that
$r>0.$ We will use $O(x, r)$ for
$$
\{y\in X: \di(x,y)<r\}.
$$
}
\end{NN}

\begin{NN} {\rm (see 1.2 of \cite{LP})}
{\rm Let $X$ be a compact metric space, let $A$ be a  unital \CA\,
and let $L: C(X)\to A$ be a  \morp. Suppose $\ep>0$ and ${\mathcal
F}$ is a finite subset of $C(X).$  Denote by $\Sigma_{\ep}(L,
{\mathcal F})$ the closure of the subset of those points
$\lambda\in X$ for which there is a non-zero hereditary \SCA\, $B$
of $A$ satisfying
$$
\|(f(\lambda)-L(f))b\|<\ep\andeqn \|b(f(\lambda)-L(f))\|<\ep
$$
for all $f\in {\mathcal F}$ and $b\in B$ with $\|b\|\le 1.$ Note that
if $\ep<\sigma,$ then $\Sigma_{\ep}(L, {\mathcal F})\subset
\Sigma_{\sigma}(L, {\mathcal F}).$ }

\end{NN}

\begin{lem}\label{inj}
Let $X$ be a compact metric space, let $\ep>0,$  $\sigma>0$ and
let ${\mathcal F}\subset C(X)$ be a finite subset of the unit
ball. There is $\dt>0$ and a finite subset ${\mathcal G}\subset
C(X)$ satisfying the following: For any unital \CA\, $A,$ any
unital $\dt$-${\mathcal G}$-multiplicative \morp\, $L: C(X)\to A,$
if $L$ is also $1/2$-${\mathcal G}$-injective, then
$\Sigma_{\ep}(L, {\mathcal F})$ is $\sigma$-dense in $X.$
\end{lem}

\begin{proof}
Choose a $\sigma$-dense finite subset $\{x_1,x_2,...,x_m\}$ in $X.$
Let $\eta=\sigma_{X, \ep/3, {\mathcal F}}.$  In particular,
\beq\label{inj-0}
|f(x)-f(x')|<\ep/3\rforal f\in {\mathcal F}
\eneq
whenever $\di(x,x')<\eta.$ Choose non-negative functions
$g_1,g_2,...,g_m, g_1', g_2',...,g_m'\in C(X)$ such that
$g_i(x)\le 1,$ $g_i(x)=1$ if $x\in O(x_i, \eta/4)$ and $g_i(x)=0$
if $x\not\in O(x_i, \eta/2)$ and $g_i'(x)\le 1,$ $g_i'(x)=1$ if
$x\in O(x_i, \eta/8)$ and $g_i'(x)=0$ if $x\in O(x_i, \eta/4),$
$i=1,2,...,m.$ Choose $\dt<\ep/4$ and
$$
{\mathcal G}={\mathcal F}\cup\{g_1,g_2,..., g_m, g_1',g_2',...,g_m'\}.
$$
Now suppose that $L: C(X)\to A$ is $\dt$-${\mathcal G}$-multiplicative and
$1/2$-${\mathcal G}$-injective.
Let $B_i=\overline{L(g_i')AL(g_i')},$ $i=1,2,...,m.$ Since
$L$ is $1/2$-${\mathcal G}$-injective, $B_i\not=\{0\},$ $i=1,2,...,m.$
Let $b=L(g_i')cL(g_i')$ for some $c\in A$ such that $\|b\|\le 1.$
Then, by (\ref{inj-0}),
\beq\label{inj-1}\nonumber
\|(f(x_i)-L(f))b\| &\le
&\|(f(x_1)-L(f))L(g_i)L(g_i')cL(g_i')\|+\dt\\\nonumber
&\le &\| L((f(x_i)-f)g_i)b\|+\dt+\dt\\
&<& \ep/ 3+\dt+\dt<5\ep/6.
\eneq
It follows that
$$
\|(f(x_i)-L(f))b\|<\ep
$$
for all $b\in B_i,$ $i=1,2,...,m.$ Similarly,
$$
\|b(f(x_i)-L(f))\|<\ep
$$
for $i=1,2,...,m.$ Therefore $x_i\in \Sigma_{\ep}(L, {\mathcal F}),$
$i=1,2,...,m.$ It follows that $\Sigma_{\ep}(L,{\mathcal F})$ is
$\sigma$-dense in $X.$

\end{proof}

\begin{lem}\label{inj2}
Let $X$ be a compact metric space and let $d>0.$ There is a finite
subset ${\mathcal G}\subset C(X)$ satisfying the following:

For any unital \CA\, $A$ and a unital \hm\, $h: C(X)\to A$ which
is $1/2$-${\mathcal G}$-injective, then spectrum $F$ of $h$ is
$d$-dense in $X.$

\end{lem}

\begin{proof}
Let $\{x_1,x_2,...,x_n\}$ be a $d/2$-dense subset of $X.$ Let
$f_i\in C(X)$ such that $0\le f_i\le 1,$ $f_i(x)=1$ if $x\in
O(x_i, d/4)$ and $f_i(x)=0$ if $x\not\in O(x_i, d/2).$ Put
$$
{\mathcal F}=\{f_1,f_2,...,f_n\}.
$$

 Let
$\eta=\sigma_{X, 1/4, {\mathcal F}}.$  Let $g_i, g_i'\in C(X)$ be as
in the proof of \ref{inj}, $i=2,...,n.$ Suppose that $h: C(X)\to
A$ is a unital \hm\, which is $1/2$-${\mathcal G}$-injective. Let $F$
be the spectrum of $h.$

If $F$ were not $d$-dense, then there is $i$ such that
$O(x_i,d/2)\cap F=\emptyset.$ Then $f_i(x)=0$ for all $x\in F.$ It
follows that $h(f_i)=0.$ However, as in the proof of \ref{inj}, we
have
$$
\|b\|=\|f(x_i)b\|=\|(f(x_i)-h(f))b\|<1/2
$$
for any $b\in B_1$ with $\|b\|= 1.$ A contradiction.

\end{proof}

\begin{lem}\label{triv}
{\rm (1)}\,\,\, Let $X$ be a connected finite CW complex  and let
$A$ be a unital anti-liminal   \CA. Then there exists a unital
monomorphism $\phi_0: C(X)\to A$ and a unital monomorphism
$\phi_1: C([0,1])\to A$ and a unital monomorphism $\psi: C(X)\to
C([0,1])$ such that
\beq\label{triv1}
\phi_0=\phi_1\circ \psi.
\eneq

{\rm (2)}\,\,\, If $X$ is a finite CW complex with $k$ connected
components and $A$ is a unital \CA\, with $k$ mutually orthogonal
projections $p_1,p_2,...,p_k$ such that each $p_iAp_i$ is
anti-liminal, then there is unital monomorphism $\psi_0: C(X)\to
A$ such that
$$
[\phi_0]=[\psi_{00}]\,\,\,\text{in}\,\,\, KK(C(X), A),
$$
where $\psi_{00}$ is a point-evaluation on $k$ points with one
point on each component.

{\rm (3)}\,\,\,If $X$ is a compact metric space and if  $A$ is a
unital simple \CA\, with real rank zero, then there is a a unital
monomorphism $\phi_0: C(X)\to A$ and a unital monomorphism
$\phi_1: C(\Omega)\to A$ and a unital monomorphism $\psi: C(X)\to
C(\Omega)$ such that
\beq\label{triv2}
\phi_0=\phi_1\circ \psi,
\eneq
where $\Omega$ is the Cantor set.

\end{lem}

\begin{proof}
It is clear that part (2) of the lemma follows from  part (1). So
we may assume that $X$ is connected.  It is well known that there
is a surjective continuous map $s: [0,1]\to X.$ This, in turn,
gives a unital monomorphism $\psi: C(X)\to C([0,1]).$ It follows
from \cite{AS} that there is a positive element $a\in A$ such that
$sp(a)=[0,1].$ Define $\phi_1: C([0,1])\to A$ by $\phi_1(f)=f(a)$
for all $f\in C([0,1]).$ Define $\phi_0=\phi_1\circ \psi.$

To see part (3), we note that there is unital monomorphism
$\psi_0: C(\Omega)\to A$ and a unital monomorphism $\psi: C(X)\to
C(\Omega).$ Then define $\phi_1=\psi_0\circ \psi.$

\end{proof}

\begin{NN}

{\rm Two projections $p$ and $q$ in a \CA\, $A$ are siad to be
equivalent if there exists $w\in A$ such that $w^*w=p$ and
$ww^*=q.$ }

\end{NN}

\begin{NN}
{\rm Let $A$ be a stably finite \CA. We will use $T(A)$ for the
tracial state space of $A.$ Denote by $Aff(T(A))$ the real
continuous affine functions on $T(A).$ Let $\tau\in T(A),$ we will
also use $\tau$ for $\tau\otimes Tr$ on $A\otimes M_n,$ where $Tr$
is the standard trace on $M_n.$

There is a (positive) \hm\, $\rho_A: K_0(A)\to Aff(T(A))$ defined
by $\rho([p])=\tau(p)$ for projection $p\in M_n(A).$ Denote by
${\rm ker}\rho_A$ (or ${\rm ker}\rho$)  the kernel of the map
$\rho_A.$

}
\end{NN}

\begin{NN}
{\rm A unital \CA\, $A$ is purely infinite and simple, if $A\not=\C,$ and for any
$a\not=0,$ there are $x, y\in A$ such that $xay=1.$
We refer the reader to \cite{Cu}, \cite{Ro2},\cite{Ro}  \cite{Ph1} and \cite{P} , for example,
for some  related results about purely infinite simple \CA s
used in this paper.

A unital separable simple \CA\, $A$ with tracial rank zero has
real rank zero, stable rank one and weakly unperforated $K_0(A).$
When $A$ has tracial rank zero, we write $TR(A)=0.$ We refer the
reader to \cite{Lntrd}, \cite{Lntr0}, \cite{Lnan}, \cite{Lnduke},
\cite{Br}, \cite{W1} and \cite{W2} for more information. }
\end{NN}

Most of this section will not be used until Chapter 2.
\vspace{0.2in}



\section{The Basic Homotopy Lemma for dim$(X)\le 1$}

{\it In what follows, we use ${\bf B}$ for the class of unital \CA
s which are simple \CA s with real rank zero and stable rank one,
or purely infinite simple \CA s. }

\begin{NN}\label{factd1}
{\rm When $X$ is a connected finite CW complex of covering
dimension 1, $X\times S^1$ is of dimension 2. Then the followings
are true.

1)\,$K_0(C(X))=\Z,$  $K_0(C(X\times S^1))=K_0(C(X))\oplus
\bt_1(K_1(C(X)))=\Z\oplus \bt_1(K_1(C(X)))$ and $K_1(C(X\times
S^1))=K_1(C(X))\oplus \bt_0(K_0(C(X)))=K_1(C(X))\oplus \Z.$ In
particular, $K_i(C(X\times S^1))$ is torsion free and
$\text{ker}\rho_{C(X\times S^1)}=\bt_1(K_1(C(X))).$


2) Suppose that $h: C(X)\to A$ is a unital \hm\, and $u\in A$ is a
unitary. Suppose that $\|uh(a)-h(a)u\|<\dt$ for all $a\in {\mathcal
G},$ where ${\mathcal G}\subset C(X)$ is a finite subset and $\dt>0.$
Suppose that $\text{bott}_1(h,u)$ is well defined and
$\text{bott}_1(h,u)=\{0\}$ and $[u]=\{0\}.$ Let ${\mathcal P}\subset
\underline{K}(C(X\times S^1).$ Then, with sufficiently large
${\mathcal G}$ and sufficiently small $\dt.$
$$
[\psi]|_{\mathcal P}=[\phi_0]|_{\mathcal P}
\,\,\,\text{in}\,\,\,KL(C(X\times S^1),A),
$$
where $\psi(f\otimes g)=h(f)g(u)$ for all $f\in C(X)$ and $g\in
C(S^1)$ and $\phi_0(a)=a(\xi)\cdot 1_A$ for all $a\in C(X\times
S^1),$ where $\xi\in X\times S^1$ is a point.

3)\,\,\,Let $Y\subset X\times S^1$ be a compact subset.  It is
known and easy to see that $s_{*0}$ maps ${\rm ker}\rho_{C(X\times
S^1)}$ onto ${\rm ker}\rho_{C(Y)}$ (see Lemma 2.2 of \cite{GL1}).
In particular, if $h: C(X\times S^1)\to B$ is a \hm\, such that
$(h_{*0})|_{{\rm ker}\rho_{C(X\times S^1)}}=0$ then
$(h_1)_{*0}|_{{\rm ker}\rho_{C(Y)}}=0,$ where $h_1: C(Y)\to B$ is
a \hm\, defined by $h=h_1\circ s.$

 (4)\,\,\, If $Y$ is a compact subset of $X\times S^1,$
then $K_i(C(Y))$ is also torsion free. Let $s: C(X\times S^1)\to
C(Y)$ is the quotient map. Let $\phi: C(Y)\to B$ be a unital \hm\,
for some unital \CA\, $B.$ Then $[\phi]\in {\mathcal N}$ (see 0.4 of
\cite{Lnalm}) if and only if
$$
\phi_{*1}=0\andeqn \phi_{*0}(s_*(\bt_1(K_1(C(X)))))=0.
$$
 }
\end{NN}

Recall that a unital \CA\, $A$ is said to be $K_1$-simple, if
$K_1(A)=U(A)/U_0(A)$ and if $[p]=[q]$ in $K_0(A)$ for two
projections $p$ and $q$ in $A,$ then there exists $w\in A$ such
that $w^*w=p$ and $ww^*=q.$

\begin{lem}\label{1Kemb}
Let $X$ be a finite CW complex with torsion free $K_1(C(X))$ and
let $B$ be a unital $K_1$-simple \CA\, of real rank zero.
 Suppose that $\lambda: K_1(C(X))\to K_1(B)$ is a \hm.
Then, for any $\dt>0,$ any finite subset ${\mathcal G}\subset C(X),$
 any finite subset ${\mathcal P}_1\subset K_1(C(X))$ and
any finite subset ${\mathcal P}_0\subset {\rm ker}\rho_{C(X)},$
 there exists
$\dt$-${\mathcal G}$-multiplicative \morp\, $L: C(X)\to B$ such that
$$
[L]|_{{\mathcal P}_1}=\lambda|_{{\mathcal P}_1}\andeqn [L]|_{{\mathcal
P}_0}=0.
$$

If, in addition, $X$ has dimension 1, then there is a unital
monomorphism $h: C(X)\to B$ such that
$$
h_{1*}=0.
$$
\end{lem}

\begin{proof}
If $B$ is a finite dimensional \CA\, we define $L(f)=f(\xi)1_B$
for all $f\in C(X),$ where $\xi\in X$ is a fixed point. We now
assume that $B$ is non-elementary. It is also clear that one may
reduce the general case to the case that $X$ is connected.

 Write $K_1(C(X))=\Z^k$ for some
integer $k>0.$ Let $C$ be a unital simple A$\T$-algebra of real
rank zero with $K_1(C)=\Z^k $ and ${\rm ker}\rho_C=\{0\}.$ For the
convenience, one may assume that $C=\overline{\cup_{n=1}^{\infty}
C_n},$ where
$$
C_n=\oplus_{i=1}^k M_{r(i)}(C(S^1)), n=1,2....
$$
The construction of such A$\T$-algebra is standard.

It follows from  a result of \cite{Li} that there is a unital \hm\,
$h_0: C(X)\to C$ such that
\beq\label{Kemb-1}
(h_0)_{*1}={\rm id}_{\Z^k}\andeqn ((h_0)_{*0})|_{{\rm
ker}\rho_{C(X)}}=\{0\}.
\eneq
(see also \cite{LnK}). For any finite subset ${\mathcal G}\subset
C(X),$ we may assume that $h_0({\mathcal G})\subset C_n$ for $n\ge
1.$ Nuclearity  of \CA s involved implies that there is a
$\dt$-${\mathcal G}$-multiplicative \morp\, $L': C(X)\to C_n$ such
that
$$
\|h_0(f)-L'(f)\|<\dt\rforal {\mathcal G}.
$$
We may also assume (by choosing a larger $n$) that
$\overline{{\mathcal Q}_{\dt,{\mathcal G}}}\supset {\mathcal P}_0\cup {\mathcal
P}_1$ and that
$$
[L']|_{{\mathcal P}_1}={\rm id}_{\Z^k}|_{{\mathcal P}_1}\andeqn
[L']|_{{\mathcal P}_0}=\{0\}.
$$
There are non-zero mutually orthogonal projections $e_{j, i}\in B$
($j=1,2,...,r(i),$ $i=1,2,..,k$) such that for fixed $k,$
$\{e_{j,k}: i=1,2,...,r(i)\}$ is a set of non-zero mutually
orthogonal projections.

Let  $g_i,i=1,2,...,k$ be the standard generator $\Z^k$ and let
$z_i=\lambda(g_i),$ $i=1,2,...,k.$ There is a unitary $v_i\in
e_{1,i}Be_{1,i}$ such that $[v_i+(1-e_{1,i})]=z_i,$ $i=1,2,....$
Put $p=\sum_{i=1}^k\sum_{j=1}^{r(i)}e_{i,j}.$  We then obtain a
unital  monomorphism $h_1: C_n\to pBp$ such that
$$
(h_1)_{*1}=\lambda
$$
as a \hm\, from $\Z^k$ to $K_1(B).$ Define $h_{00}: C(X)\to
(1-p)A(1-p)$ by $h_{00}(f)=f(\xi)(1-p)$ for all $f\in C(X),$ where
$\xi\in X$ is a fixed point in $X.$  Define $L=(h_1\circ L'\circ
h_0)\oplus h_{00}.$ where $\xi\in X$ is a fixed point in $X.$  It
is clear that so defined $L$ meets the requirements.

By \ref{triv}, there is a unital monomorphism $h_{00}': C(X)\to
(1-p)B(1-p)$ such that $[h_{00}']=[h_{00}]$ in $KK(C(X), B).$ The
last statement follows from the fact that $C(X)$ is
semi-projective (by 5.1 of \cite{Lo}) and replacing $h_1\circ
L'\circ h_0$ by a \hm\, and replacing $h_{00}$ by $h_{00}'.$

\end{proof}

\begin{lem}\label{d1f}
Let $X$ be a connected finite CW complex with dimension 1 (with a fixed metric), let
$A\in {\bf B}$
 and let $h: C(X)\to A$ be a monomorphism. Suppose that
there is a unitary $u\in A$ such that
\beq\label{d1f1}
h(a)u=uh(a)\tforal \, a\in C(X),\,\,\,[u]=0\,\,\,\text{in}\,\,\,
K_1(A) \andeqn {\rm{bott}_1}(h,u)=0.
\eneq
Suppose also  that $\psi: C(X\times S^1)\to A$ defined by
$\psi(f\otimes g)=h(f)g(u)$ for $f\in C(X)$ and $g\in C(S^1)$ is a
monomorphism. Then, for any $\ep>0$ and any finite subset ${\mathcal
F}\subset C(X),$
 there exists a rectifiable continuous path of unitaries
$\{u_t: t\in [0,1]\}$ of $A$ such that
\beq\label{d1f2}
u_0=u,\,\,\,u_1=1_A\andeqn \|[h(a), u_t]\|<\ep
\eneq
for all $a\in {\mathcal F}$ and all $t\in [0,1].$ Moreover,
\beq\label{d1f3}
{\rm{Length}}(\{u_t\})\le \pi+\ep\pi
\eneq
\end{lem}

\begin{proof}
Let $\ep>0$ and ${\mathcal F}\subset C(X)$ be a finite subset. We may
assume that $1_{C(X)}\in {\mathcal F}.$ Let
$$
{\mathcal F}'=\{f\otimes a: f\in {\mathcal F}, a=1_{C(S^1)},
\text{or}\,\,\,a=z\}.
$$
Let $\sigma_1=\sigma(X\times S^1, {\mathcal F}', \ep/32)$ (see
\ref{sigma}).
Let $\dt_1>0,$  ${\mathcal G}\subset C(X\times S^1)$ and  ${\mathcal
P}\subset \underline{K}(C(X\times S^1))$ be required by Theorem
1.12 of \cite{Lnalm} associated with $X\times S^1,$ ${\mathcal F}',$
$\ep/16$ (in place of $\ep$). Here we also assume that any two
$\dt_1$-${\mathcal G}$-multiplicative \morp s $L_1, L_2$ from
$C(X\times S^1),$
\beq\label{d1f5}
[L_1]|_{\mathcal P}=[L_2]|_{\mathcal P},\,\,\,\text{if}\,\,\,
\eneq
\beq\label{d1f6}
L_1\approx_{\dt_1} L_2 \,\,\,\,\text{on}\,\,\,{\mathcal G}
\eneq

Without loss of generality, we may assume that ${\mathcal G}={\mathcal
G}_1\otimes S,$ where ${\mathcal G}_1$ is in the unit ball of $C(X)$
and $S=\{1_{C(S^1)}, z\}.$ We may also assume that $\dt_1<\ep/2$
 and that ${\mathcal F}\subset {\mathcal G}_1.$

Let $\eta_1=(1/2)\sigma_{X\times S^1,\min\{\ep/32, \dt_1/32\},
{\mathcal G}}.$ So if ${\rm dist}(x, x')<\eta_1$ and if ${\rm
dist}(t,t')<\eta_1,$ then
\beq\label{d1f4}
|f(x\times t)-f(x'\times t')|<\min\{\ep/32, \dt_1/32\}\,\rforal\,
f\in {\mathcal G}.
\eneq
Let $\{x_1,x_2,...,x_m\}$ be $\sigma_1/4$-dense in $X$ and let
$\{t_1,t_2,...,t_l\}$ divide the unit circle into $l$ arcs with
the  length with $2\pi/l<\eta_1.$  We may assume that
$\eta_1<\sigma_1/4.$ We also assume that $t_1=1.$ In particular, $
\{x_i\times t_j: i=1,2,...,m, j=1,2,...,l\} $ is
$\sigma_1/2$-dense.

Let $g_{i,j}$ be nonnegative functions in $C(X\times S^1)$ with
$0\le g_{i,j}\le 1$ such that $g_{i,j}(\xi)=1$ if ${\rm dist}(\xi,
x_i\times t_j)<\eta_1/4s$ and $g_{i,j}(\xi)=0$ if ${\rm dist}(\xi,
x_i\times t_j)\ge \eta_1/2s.$ Since $A$ is a unital simple \CA\,
with real rank zero and $\psi$ is injective, there are non-zero,
mutually orthogonal and mutually equivalent projections
$E_{i,j},E_{i,j}'\in \overline{\psi(g_{i,j})A\psi(g_{i,j}))}$ for
each $i$ and $j.$ Since $A$ is simple and has real rank zero, then
(by repeated application of (1) of Lemma 3.5.6 of \cite{Lnbk}, for
example) there is a non-zero projection $E_0\le E_{1,1}$ such that
\beq\label{d1f6+}
[E_0]\le [E_{i,j}],\,\,\,i=1,2,...,m\andeqn j=1,2,...,l.
\eneq
Put
\beq\label{d1f7}
\Phi(f)=\sum_{i,j} f(x_i\times t_j)(E_{i,j}+E_{i,j}')+{\tilde
\phi}(f)\,\rforal\, f\in C(X\times S^1),
\eneq
where
\beq\label{dlf8-} &&{\tilde
\phi}(f)=(1-\sum_{i,j}(E_{i,j}+E_{i,j}'))\psi(f)(1-\sum_{i,j}(E_{i,j}+E_{i,j}'))\andeqn \\
\nonumber
 &&\hspace{-0.6in}\Phi_1(f)=f(x_1\times
t_1)(E_{1,1}'+E_{1,1}-E_0)+\sum_{(i,j)\not=(1,1)} f(x_i\times
t_j)(E_{i,j}+E_{i,j}')+\\\label{d1f8}
&& +{\tilde \phi}(f)
\eneq
for $f\in C(X\times S^1).$ Thus, by the choice of $\sigma_1,$
\beq\label{d1f8+}
\|\psi(f)-\Phi(f)\|<\ep/32\rforal f\in {\mathcal G}.
\eneq
By the choice of $\eta_1,$  one sees that  ${\tilde \phi}$ is
$\dt_1$-${\mathcal G}$-multiplicative. So $\Phi_1$ is $\dt_1$-${\mathcal
G}$-multiplicative. Moreover, $\Sigma_{\dt_1}(\Phi_1, {\mathcal F})$
is $\sigma_1/2$-dense,
 since $\{(x_i\times t_j): i,j\}$ is $\sigma_1/2$-dense. There are
two non-zero mutually orthogonal projections $e_1, e_2\in E_0AE_0$
such that $e_1+e_2=E_{0}.$ Since both $e_1Ae_1$ and $e_2Ae_2$ are
simple \CA s in ${\bf B},$
 and $X$ is a connected finite CW complex, it is easy to
construct unital monomorphisms $h_1: C(X)\to e_1Ae_1$ and $h_2:
C(X)\to e_2Ae_2$ such that
\beq\label{d1f9}
(h_1)_{*1}=h_{*1}\andeqn (h_2)_{*1}+h_{*1}=\{0\}
\eneq
(by \ref{1Kemb}). It follows from Cor 1.14 of \cite{Lnalm} that
there are mutually orthogonal projections
$\{p_1,p_2,...,p_K\}\subset E_{0}AE_{0}$ and points
 $\{y_1,y_2,...,y_K\}\subset X$ such that $\sum_{i=1}^Kp_i=E_0$ and
 \beq\label{d1f10}
 \|(h_1(f)+h_2(f))-\sum_{k=1}^Kf(y_k)p_k\|<\ep/32\rforal f\in {\mathcal G}_1.
 \eneq

Define $\psi_1(f\otimes g)=h_1(f)g(1)e_1\andeqn \psi_2(f\otimes
g)=h_2(f)g(1)e_2$ for all $f\in C(X)$ and $g\in C(S^1).$ Define
$\psi_0: C(X\times S^1)\to E_0AE_0$ by
\beq\label{d1f11}
\psi_0(f\otimes g)=\sum_{k=1}^Kf(y_k)g(1)p_k
\eneq
for all $f\in C(X)$ and $g\in C(S^1).$ By replacing $\ep/32$ by
$\ep/16,$ in (\ref{d1f10}), we may assume that
$\{y_1,y_2,...,y_K\}\subset \{x_1,x_2,...,x_{K'}\}$ with $K'\le
m.$ Keep this in mind, since $p_k\le e_2,$
  by (\ref{d1f6+}), we obtain a unitary $w_1\in A$ such
that
\beq\label{d1f12}
{\rm ad}\, w_1\circ (\psi_0(f)+\Phi_1(f))= \Phi(f)
\eneq
for all $f\in C(X\times S^1).$ It follows that (by (\ref{d1f10})
\beq\label{d1f12+}
{\rm ad}\, w_1\circ( (\psi_1\oplus \psi_2)+\Phi_1)\approx_{\ep/16}
\Phi\,\,\,\text{on}\,\,\, {\mathcal G}
\eneq
and (by \ref{d1f8+}))
\beq\label{d1f12++}
{\rm ad}\, w_1\circ( (\psi_1\oplus
\psi_2)+\Phi_1)\approx_{3\ep/32} \psi\,\,\,\text{on}\,\,\, {\mathcal
G}.
\eneq
It follows from (\ref{d1f9}) and the condition that
$\text{bott}_1(h,u)=0$ and $[u]=0$  that
\beq\label{d1f13}
[{\rm ad}\, w_1\circ (\psi_2\oplus \Phi_1)]|_{\mathcal
P}=[\Phi_{00}]|_{\mathcal P},
\eneq
where $\Phi_{00}=f(\xi_X\times 1)(1-e_1)$ for all $f\in C(X\times
S^1).$  We also check that ${\rm ad}\, w_1\circ (\psi_2\oplus
\Phi_1)$ is $\dt_1$-${\mathcal G}$-multiplicative and
$\Sigma_{\dt_1}(w_1\circ (\psi_2\oplus \Phi_1),{\mathcal F}')$ is
$\sigma_1/2$-dense.  It follows from Theorem 1.12 of
\cite{Lnalm} that there is a \hm\, $\Phi_0: C(X\times S^1)\to
(1-e_1)A(1-e_1)$ with finite dimensional range such that
\beq\label{d1f14}
{\rm ad}\,w_1\circ (\psi_2\oplus \Phi_1)\approx_{\ep/16}
\Phi_0\,\,\,\,\text{on}\,\,\, {\mathcal F}'.
\eneq
In the finite dimensional commutative \SCA\, $\Phi_0(C(X\times
S^1),$ there is a continuous rectifiable path $\{U_t: t\in
[0,1]\}$ of $(1-e)A(1-e)$ such that
\beq\label{d1f15}
U_0=\Phi_0(1\otimes z), U_1=1-w_1^*e_1w_1\andeqn \Phi_0(f\otimes
1)U_t=U_t\Phi_0(f\otimes 1)
\eneq
for all $t\in [0,1]$ and $f\in C(X).$ Moreover
\beq\label{d1f16}
{\rm{Length}}(\{U_t\})\le \pi.
\eneq
Define $v_t=w_1^*e_1w_1\oplus U_t$ for $t\in [0,1].$ Then, (note
that $\psi_1(1\otimes z)=e_1$) by (\ref{d1f14}) and
(\ref{d1f12++}),
\beq\label{d1f17}\nonumber
\|v_0-u\| &=& \|w_1^*e_1w_1\oplus \Phi_0(1\otimes
z)-u\|\\\nonumber &\le & \|w_1^*e_1w_1\oplus [\Phi_0(1\otimes
z)-w_1^*(\psi_2\oplus \Phi_1)(1\otimes z)w_1]\|\\\nonumber &+&
\|w_1^*e_1w_1\oplus w_1^*(\psi_2 \oplus \Phi_1)(1\otimes
z)w_1-u\|\\\nonumber
&<&\ep/16+\|{\rm ad}\, w_1^*((\psi_1\oplus \psi_2)(1\otimes z)
\oplus \Phi_1(1\otimes z)w_1-\psi(1\otimes z)\|\\
&<&\ep/16+3\ep/32=\ep/8.
\eneq
We estimate that, by (\ref{d1f12++}) and (\ref{d1f14}) and
(\ref{d1f15}),
\beq\label{d1f18}
\|[h(f),v_t]\|&<&6\ep/32+\|[{\rm ad}\, w_1\circ ((\psi_1+\psi_2)+\Phi_1)(f\otimes 1),v_t]\|\\
&<&3\ep/16+\ep/8+\|[{\rm ad}\, w_1\circ (\psi_1+\Phi_0)(f\otimes 1),
v_t]\|\\
&=&5\ep/16
\eneq
for all $f\in {\mathcal F}.$ Combing this with (\ref{d1f17}), we
obtain a continuous rectifiable path of unitaries $\{u_t:t\in
[0,1]\}$ of $A$ such that
\beq\label{d1f19}
u_0=u,\,\,\,u_1=1_A\andeqn \|[h(a),u_t]\|<\ep\rforal\, a\in {\mathcal
F}
\eneq
and all $t\in [0,1].$ Moreover,
\beq\label{d1f20}
{\rm{Length}}(\{u_t\})\le \pi+\ep\pi.
\eneq

\end{proof}

\begin{lem}\label{d1hf}
Let $X_1$ be a connected finite CW complex with dimension 1. Let
$\ep>0$ and ${\mathcal F}_1\subset C(X_1)$ be a finite subset. There
is $\sigma>0$ satisfying the following: Suppose that $A\in {\bf
B}$ is an infinite dimensional unital \CA, suppose that  $h_1:
C(X_1)\to A$ is a \hm\, and suppose that  $u\in A$ is a unitary
with $[u]=\{0\}$ in $K_1(A)$ such that
\beq\label{d1hf1}
h_1(a)u=uh_1(a)\tforal\, a\in C(X_1)\andeqn {\rm{bott}_1}(
h_1,u)=\{0\},
\eneq
and suppose that $X$ is a  subset of $X_1$ which is a finite CW
complex and the spectrum  of  $h_1$ is $Y\subset X$ which is
$\sigma$-dense in $X.$ Then there exists a unital monomorphism
${\bar h}: C(X\times S^1)\to A$ and a rectifiable continuous path
of unitaries $\{u_t: t\in [0,1]\}$ such that
\beq\label{d1hf2}
u_0=u,\,\,\, u_1={\bar h}(1\otimes z),\,\,\,\|[h_1(a), u_t]\|<\ep
\andeqn
\eneq
\beq\label{dihf3}
\|{\bar h}(s(a)\otimes 1)-h_1(a)\|<\ep\rforal \, a\in {\mathcal F}_1,
\eneq
where $s: C(X_1)\to C(X)$ is defined by $s(f)=f|_X.$   Moreover
\beq\label{d1hf4}
{\rm{Length}}(\{u_t\})\le \pi+\ep\pi
\eneq

\end{lem}

\begin{proof}
There is a \hm\, $h: C(X)\to A$ such that $h=h_1\circ s,$ where
$s: C(X_1)\to C(X)$ is the quotient map. Put ${\mathcal F}=s({\mathcal
F})$ and
$$
{\mathcal F}'=\{f\otimes g: f\in {\mathcal F}\andeqn
g=1,\,\,\,\text{or}\,\,\, g=z\}.
$$

Define $\psi: C(X\times S^1)\to A$ by $\psi(f\otimes g)=h(f)g(u)$
for $f\in C(X)$ and $g\in C(S^1).$ By (\ref{d1hf1}),
$[\psi]|_{{\rm ker}\rho_{C(X\times S^1)}}=\{0\}.$

Large part of the argument is the same as that  used in \ref{d1f}.
So we will keep the notation there and keep all the proof from the
beginning of that proof to the equation (\ref{d1f4}). Moreover, we
will refer to Theorem 2.5 of \cite{Lnalm} instead of Lemma 1.12 of
\cite{Lnalm} and $\dt_1,$ ${\mathcal G}$ and ${\mathcal P}$ are now stated
in Theorem 2.5 of \cite{Lnalm} (for $X\times S^1$). Note that we
still assume that ${\mathcal F}\subset {\mathcal G}_1.$ Note that we have
$\overline{{\mathcal Q}_{\dt_1, {\mathcal G}} }\subset {\mathcal P}.$
  It should be noted (see 1.7 of
\cite{Lnalm}) that
$$\
 \sigma_{X_1\times S^1, {\mathcal F}_1, \ep/32}\le \sigma_{X\times S^1, {\mathcal F}, \ep/32}
 =\sigma_1.
$$
Put $\sigma=\sigma_{X_1\times S^1, {\mathcal F}_1, \ep/32}/4.$ Let $Y$
be $\sigma$-dense in $X.$  We  assume that
$\{x_1,x_2,....,x_m\}\subset Y$ is $\sigma$ -dense in $X.$

Without loss of generality, we may still assume, for each $i,$
that there is $j$ such that
 $\psi(g_{i,j})\not=0.$
 We may replace $\{t_1,t_2,...,t_l\}$ by
 $\{t_{i(1)},t_{i(2)},...t_{i(k(i))}\}$
 with $t_{i(1)}=1.$
 We assume that $\{x_i\times t_{i,j}:
 j=1,2,...,k(i),i=1,2,...,m\}$ is $\sigma$ -dense in $Y.$
We now replace $g_{i,j}$ by $g_{i,1},$ $E_{i,j}$ by $E_{i,1}$ and
$E_{i,j}'$ by $E_{i,1}'.$  Accordingly,  (\ref{d1f6+}),
(\ref{d1f7}),(\ref{dlf8-}),(\ref{d1f8}) and (\ref{d1f8+}) hold
where we replace all $j$ by $1,$ $t_1$ by $t_{1(1)}$ and $t_j$ by
$t_{i(1)}=1,$ respectively.

Now, instead of taking two projections in $E_0AE_0,$ we  choose
$e_1=E_0.$
 We also introduce another  \hm\, $F_3$ as follows.
By \ref{triv}, there are  monomorphisms $F_1: C(X\times S^1)\to
C(D)$ and $F_2: C(D)\to e_1Ae_1,$ where $D$ is a finite union of
closed intervals. Put $F_3=F_2\circ F_1.$ It follows from Cor.
1.14 of \cite{Lnalm} that there is a unital \hm\, $\psi_0:
C(X\times S^1)\to e_1Ae_1$ with finite dimensional range such that
\beq\label{d1hf5}
\|F_3(f)-\psi_0(f)\|<\min\{\ep/32,\dt_1/32\}
\eneq
for all $f\in {\mathcal G}.$ Without loss of generality, by changing
$\ep/32$ to $2\ep/32=\ep/16,$ we may assume that
\beq\label{d1hf6}
\psi_0(f\otimes g)=\sum_{i,j}f(x_i)g(t_j)p_{i,j}
\eneq
for all $f\in C(X)$ and $g\in C(S^1)$ and $\{p_{i,j}\}$ is a set
of finitely many mutually orthogonal projections in $e_1Ae_1$ with
$\sum_{i,j}p_{i,j}=e_1.$  Define another \hm\, $\psi_{00}:
C(X\times S^1)\to e_1Ae_1$ such that
\beq\label{d1hf7}
\psi_{00}(f\otimes g)=\sum_{i=1}^{m} f(x_i)g(1)q_i
\eneq
for all $f\in C(X)$ and $g\in C(S^1),$ where $q_i=\sum_j p_{i,j},$
$i=1,2,...,m.$ Working in the  finite dimensional commutative
\SCA\, $\psi_0(C(X\otimes S^1)),$ it is easy to obtain a
continuous rectifiable path of unitaries $\{v_t: t\in [0,1]\}$ of
$e_1Ae_1$ such that
\beq\label{d1hf8}
v_0=\psi_{00}(1\otimes z),\,\,\, v_1=\psi_0(1\otimes z)\andeqn
v_t\psi_{00}(f)=\psi_{00}(f)v_t
\eneq
for all $f\in C(X\times S^1)$ and $t\in [0,1].$ Moreover,
\beq\label{d1hf9}
{\rm{Length}}(\{v_t\})\le \pi.
\eneq

Since $q_i\le e_1=E_0,$ the same argument in the proof of
\ref{d1f} provides a unitary $w_1\in A$ such that
\beq\label{d1hf10}
{\rm ad}\, w_1\circ (\psi_{00})+\Phi_1)(f)=\Phi(f)
\eneq
for all $f\in C(X\times S^1).$ It follows, as in the proof of
\ref{d1f}, that
\beq\label{d1hf11}
{\rm ad}\, w_1\circ (\psi_{00}+\Phi_1)(f)\approx_{\ep/16}
\psi(f)\,\,\,\text{on}\,\,\, f\in {\mathcal G}
\eneq

Put $V_t=w_1^*[v_t\oplus \Phi_1(1\otimes z)]w_1$ for $t\in [0,1].$
Then by (\ref{d1hf8}) and (\ref{d1hf11}),
\beq\label{d1hf12}
\|[V_t, \psi(f\otimes 1)]\|<\ep/8
\eneq
for all $f\in {\mathcal F}.$ We have, by (\ref{d1hf11}),
\beq\label{d1hf13}
\|V_0-u\| = \|w_1^* (v_0\oplus \Phi_1(1\otimes
z))w_1-\psi(1\otimes z)\|<\ep/8
\eneq
Note  that
$$
[{\rm ad} w_1\circ( \psi_{00})+\Phi_1]|_{\mathcal P}=[\psi]|_{\mathcal P}.
$$
Since $\psi_{00}$ has finite dimensional range and $[\psi]|_{{\rm
ker}\rho_{C(X)}}=\{0\},$ by applying Theorem 2.5 of \cite{Lnalm},
one obtains a \hm\, $\Phi_0: C(X\times S^1)\to (1-E_0)A(1-E_0)$
such that
\beq\label{d1hf14}
\Phi_1\approx_{\ep/16} \Phi_0 \,\,\,{\rm on}\,\,\,{\mathcal F}'.
\eneq
Now define ${\bar h}: C(X\times S^1)\to A$ by
 \beq\label{d1hf15}
 {\bar h}(f)={\rm ad}\,w_1\circ (F_3\oplus \Phi_0)(f)\,\,\,\rforal\, f\in C(X\times S^1).
 \eneq
 By (\ref{d1hf14}) and (\ref{d1hf11}), we have
 \beq\label{d1hf15+}
 \|{\bar h}(f\otimes 1)-h(f)\|<\ep/8\rforal f\in {\mathcal F}.
 \eneq

We estimate that, by (\ref{d1hf5}) and (\ref{d1hf14}),
\beq\label{d1hf16}
\hspace{-0.3in}\|V_1-{\bar h}(1\otimes
z)\|&=&\|w_1^*(\psi_0(1\otimes z)\oplus
\Phi_1(1\otimes z))w_1-{\bar h}(1\otimes z)\|\\
\hspace{-0.3in}&< &\ep/32+\ep/16
\eneq
Note that $u=\psi(1\otimes z).$  By connecting $V_1$ with ${\bar
h}(1\otimes z)$ and $V_0$ to $u$  appropriately, we obtain a
continuous rectifiable path of unitaries $\{u_t: t\in [0,1]\}$
such that
\beq\label{d1hf17}
u_0=u,\,\,\, u_1={\bar h}(1\otimes z)\andeqn \|[u_t, h(f)]\|<\ep
\eneq
for all $f\in {\mathcal F}$ and $t\in [0,1].$ Moreover,
\beq
\text{Length}(\{u_t\})\le \pi +\ep\pi
\eneq
Since $F_3$ is a monomorphism, we conclude that ${\bar h}$ is also
a monomorphism.

\end{proof}

\begin{rem}\label{remhp}
{\rm In the statement of \ref{d1hf}, define $h_2: C(X_1)\to A$ by
$h_2(f)={\bar h}(s(f)\otimes 1)$ for all $f\in C(X_1).$  Note
that, if $\ep$ is small enough and ${\mathcal F}$ is large enough,
then $\text{bott}_1(h_2, u_1)=\text{bott}_1(h,u)=0$ (see
\ref{Khp}). Since $X$ has dimension $1,$
$s_{*1}(K_1(C(X_1)))=K_1(C(X)).$ This implies that
$\text{bott}_1(h_2',u_1)=0,$ where $h_2': C(X)\to A$ by
$h_2'(f)={\bar h}(f\otimes 1)$ for $f\in C(X).$

It should be noted that in the statement of \ref{d1f} and
\ref{d1hf}, the length of $\{u_t\}$ can be controlled by
$\pi+\ep.$ }
\end{rem}

\begin{lem}\label{Ld1}
Let $X$ be a connected finite CW complex with dimension 1.
Then,
for any $\ep>0$ and any finite subset ${\mathcal F}\subset C(X),$
there is $\sigma>0$ satisfying the following:
Suppose that $A\in {\bf B}$ is not finite dimensional
and that $h: C(X)\to A$ is a unital \hm\, whose spectrum is $\sigma$-dense in $X$
and suppose that there is a unitary $u\in A$ with $[u]=\{0\}$ in
$K_1(A)$ such that
\beq\label{2tt1}
h(a)u=uh(a)\,\tforal\, a\in C(X) \andeqn {\rm{bott}_1}(h,u)=0.
\eneq
Then,
there exists a rectifiable continuous path of unitaries
$\{u_t: t\in [0,1]\}$ of $A$ such that
\beq\label{2tt2}
u_0=u,\,\,\,u_1=1_A\andeqn \|[h(a), u_t]\|<\ep
\eneq
for all $a\in {\mathcal F}$ and all $t\in [0,1].$ Moreover,
\beq\label{2tt3}
{\rm{Length}}(\{u_t\})\le 2\pi+\ep\pi.
\eneq
\end{lem}

\begin{proof}
  Let $\ep>0$ and ${\mathcal F}\subset C(X)$ be a
finite subset. We may assume that ${\mathcal F}$ is in the unit ball
of $C(X).$ Let ${\mathcal P}\subset K_1(C(X))$ be a finite subset
containing a set of generators. Let $\ep_1>0$ and ${\mathcal
F}_1\subset C(X)$ be a finite subset so that  $\|[h'(a),
w]\|<\ep_1$ for all $a\in {\mathcal F}_1$ implies that
$\text{bott}_1(h,w)$ is well defined for any unital \hm \, $h':
C(X)\to A$ and any unitary $w\in A.$ We may assume that
$\ep_1<\ep$ and ${\mathcal F}\subset {\mathcal F}_1.$

Let $\sigma$ be in \ref{d1hf} corresponding to $\ep_1/4$ and
${\mathcal F}_1.$ There is a subset $Y\subset X$ which is a finite CW
complex so that the spectrum of $h$ is $\sigma$-dense in $Y.$ Then
it is clear that the lemma follows from \ref{d1f} and \ref{d1hf}
(see also \ref{remhp}).

\end{proof}

\begin{thm}\label{dmyi}
Let $X$ be a finite CW complex with dimension 1. Then, for any
$\ep>0$ and any finite subset ${\mathcal F}\subset C(X),$ there exists
$\dt>0,$  a finite subset ${\mathcal G}\subset C(X)$ and $\sigma>0$ satisfying the
following:

Suppose that $A\in {\bf B},$ suppose that $h: C(X)\to A$ is a
unital \hm\, whose spectrum is $\sigma$-dense in $X$ and suppose that there is a unitary $u\in A$
such that
\beq\label{T2tt1}
\|[h(a), u]\|<\dt\,\tforal\, a\in {\mathcal G},\,\,\, {\rm
bott}_0(h,u)=0\andeqn {\rm{bott}_1}(h,u)=0.
\eneq

Then there exists a rectifiable continuous path of unitaries
$\{u_t: t\in [0,1]\}$ of $A$ such that
\beq\label{T2tt2}
u_0=u,\,\,\,u_1=1_A\andeqn \|[h(a), u_t]\|<\ep
\eneq
for all $a\in {\mathcal F}$ and all $t\in [0,1].$ Moreover,
\beq\label{T2tt3}
{\rm{Length}}(\{u_t\})\le 2\pi+\ep.
\eneq
\end{thm}

\begin{proof}
We assume that $A$ is not finite dimensional. So $A$ is not
elementary. The case that $A$ is finite dimensional will be dealt
next in \ref{Dm1finite}. By considering each connected component
of $X,$  since we have assume that ${\rm bott}_0(h,u)=0,$
one easily reduces the general case to the case that $X$
is connected (see also the beginning of the proof of \ref{MT1}).
So for the rest of this proof we assume that $X$ is connected and
$[u]=0$ in $K_1(A).$  Fix
$\ep>0$ and a finite subset ${\mathcal F}\subset C(X).$ We may assume
that $1\in {\mathcal F}.$ Let ${\mathcal F}'=\{f\otimes g: f\in {\mathcal
F}\andeqn g=1,\text{or}\,\,\,g=z\}.$ Let $\eta_1>0$ (in place of
$\eta$) be in \ref{Ld1} corresponding to $\ep/2$ and ${\mathcal F}.$
Put $\eta=\min\{\eta_1/2, \ep/2\}.$


Let $\dt_1>0,$ ${\mathcal G}_1\subset C(X\times S^1)$ and ${\mathcal
P}\subset \underline{K}(C(X\times S^1))$ be a finite subset
required by Theorem 2.5 of \cite{Lnalm} corresponding to $\eta/2$
and ${\mathcal F}'.$ We may assume that ${\mathcal P}$ contains a set of
generators of $\bt_1(K_1(C(X))) $ and $\dt_1<\eta.$ We may assume
that $\dt_1$ is sufficiently small and ${\mathcal G}_2$ is sufficiently
large so that for any two $\dt_1$-${\mathcal G}_2$-multiplicative \morp
s $L_i: C(X\times S^1)\to B$ (for any unital \CA\,), $[L_i]|_{\mathcal
P}$ is well-defined and
\beq\label{dny1-1}
[L_1]|_{\mathcal P}=[L_2]|_{\mathcal P},
\eneq
provided that
$$
\|L_1(f)-L_2(f)\|<\dt_1\rforal f\in {\mathcal G}_1.
$$
Moreover, by choosing even smaller $\dt_1,$ we may assume that
${\mathcal G}_1={\mathcal G}'\otimes S,$ where ${\mathcal G}'\subset C(X)$ and
$S=\{1_{C(S^1)}, z\}.$

Let $\dt>0$ and ${\mathcal G}\subset C(X)$ (in place of ${\mathcal F}_1$)
be a finite subset required in \ref{appn} for $\dt_1$ (in place of
$\ep$) and ${\mathcal G}'$ (in place of ${\mathcal F}$).

Now suppose that $h$ and $u$ satisfy the conditions in the theorem
for the above $\dt$ and ${\mathcal G}.$ Let $\phi: C(X\times S^1)\to
A$ be defined by $\phi(f\otimes g)=h(f)g(u)$ for all $f\in C(X)$
and $g\in C(S^1).$ The condition that $\text{bott}_1(h,u)=0,$
implies that
\beq\label{dmy1-2}
[\psi]|_{\mathcal P}=\af({\mathcal P})
\eneq
for some $\af\in {\mathcal N}k$ (see 2.1 of \cite{Lnalm}). By  \ref{appn} and Theorem
2.5 of \cite{Lnalm}, there is a unital \hm\, $H: C(X\times S^1)\to
A$ such that
\beq\label{dym1-3}
\|H(f\otimes g)-\psi(f\otimes g)\|<\eta_1/2\rforal f\in {\mathcal
G}_1\andeqn g=1,\,\,\,\text{or}\,\,\, g=z.
\eneq
Put $h_1(f)=H(f\otimes 1)$ (for $f\in C(X)$ and $v=H(1\otimes z).$
Then
$$
[h_1,v]=0,\,\,\, \text{bott}_1(h_1,v)=0\andeqn [v]=0
\,\,\,\text{in}\,\,\,K_1(A).
$$
Thus \ref{Ld1} applies to $h_1$ and $v$ (with $\ep/2$). By
(\ref{dym1-3}),
\beq\label{dym1-4}
h_1\approx_{\eta_1/2} h\,\,\,\text{on}\,\,\,{\mathcal G}_1.
\eneq
The lemma follows.

\end{proof}


We actually prove the following:

\begin{cor}\label{Cdmyi}
Let $X_1$ be a finite CW complex with dimension 1. Then, for any
$\ep>0$ and any finite subset ${\mathcal F}\subset C(X),$ there
exists $\dt>0,$  a finite subset ${\mathcal G}\subset C(X)$ and
$\sigma>0$ satisfying the following:

Suppose $X$ is a finite CW complex which is also a compact subset
of a finite CW complex $X\subset X_1,$ suppose that $A\in {\bf
B},$ suppose that $h: C(X)\to A$ is a unital \hm\, whose spectrum
is $\sigma$-dense in $X$ and suppose that there is a unitary $u\in
A$ such that
\beq\label{C2tt1}
\|[h(a), u]\|<\dt\,\tforal\, a\in {\mathcal G},\,\,\, {{\rm
bott}_0}(h,u)=0\andeqn {\rm{bott}_1}(h,u)=0.
\eneq

Then there exists a rectifiable continuous path of unitaries
$\{u_t: t\in [0,1]\}$ of $A$ such that
\beq\label{C2tt2}
u_0=u,\,\,\,u_1=1_A\andeqn \|[h(a), u_t]\|<\ep
\eneq
for all $a\in {\mathcal F}$ and all $t\in [0,1].$ Moreover,
\beq\label{C2tt3}
{\rm{Length}}(\{u_t\})\le 2\pi+\ep.
\eneq
\end{cor}

\begin{cor}\label{Acor}
Let $X$ be a finite CW complex with dimension 1. Then, for any
$\ep>0$ and any finite subset ${\mathcal F}\subset C(X),$ there exists
$\dt>0$ and   a finite subset ${\mathcal G}\subset C(X)$  satisfying the
following:

Suppose that $A\in {\bf B},$ suppose that $h: C(X)\to A$ is a
unital monomorphism  and suppose that there is a unitary $u\in A$
such that
\beq\label{CT2tt1}
\|[h(a), u]\|<\dt\,\tforal\, a\in {\mathcal G},\,\,\, {{\rm
bott}_0}(h,u)=0\andeqn {\rm{bott}_1}(h,u)=0.
\eneq

Then there exists a rectifiable continuous path of unitaries
$\{u_t: t\in [0,1]\}$ of $A$ such that
\beq\label{CT2tt2}
u_0=u,\,\,\,u_1=1_A\andeqn \|[h(a), u_t]\|<\ep
\eneq
for all $a\in {\mathcal F}$ and all $t\in [0,1].$ Moreover,
\beq\label{CT2tt3}
{\rm{Length}}(\{u_t\})\le 2\pi+\ep.
\eneq

\end{cor}

The following lemma, in particular, deals with the case that \CA s
are finite dimensional. The other reason to include the following
is that, in the case that $K_1(A)=0,$ the bound for the length can
be made  shortest possible.

\begin{thm}\label{Dm1finite}
Let $X$ be a finite CW complex of dimension 1. Then, for any
$\ep>0$ and any finite subset ${\mathcal F}\subset C(X),$ there exists
$\dt>0$ and a finite subset ${\mathcal G}\subset C(X)$ satisfying the
following:

Let $A\in {\bf B}$ with $K_1(A)=\{0\},$  let $h: C(X)\to A$ be a
unital \hm\, and let $u\in A$ be a unitary  such that
$$
\|[h(g),u]\|<\dt\tforal f\in {\mathcal G}\andeqn {\rm{bott}_1}(h,
u)=0.
$$
Then there exists a continuous rectifiable path of unitaries
$\{u_t: t\in [0,1]\}$ such that
$$
u_0=u,\,\,\, u_1=1_A\andeqn \|[h(f), u_t]\|<\ep\tforal f\in
{\mathcal F}.
$$
Moreover,
$$
{\rm{Length}}(\{u_t\})\le \pi +\ep.
$$
\end{thm}

\begin{proof}
Let $\ep>0$ and ${\mathcal F}\subset C(X)$ be a finite subset. Let
${\mathcal P}_0\subset K_0(C(X\times S^1))$ be a finite subset
which contains a set of generators of ${\rm ker}\rho_{C(X\times
S^1)}.$ There exists $\dt_0>0$ and a finite subset ${\mathcal
G}_0\subset C(X\times S^1)$ such that for any two unital
$\dt_0$-${\mathcal G}_0$-multiplicative \morp s $L_1, L_2:
C(X\times S^1)\to A$
 with
 $$
 \|L_1(f)-L_2(f)\|<\dt_0\rforal f\in {\mathcal G}_0,
 $$
one has
$$
[L_1]|_{{\mathcal P}_0}=[L_2]|_{{\mathcal P}_0}.
$$
Without loss of generality, we may assume that $ {\mathcal G}_0={\mathcal
G}_1\otimes S,$ where ${\mathcal G}_1\subset C(X)$ is a finite subset
and $S=\{1_{C(S^1)}, z\}.$ Put $\ep_1=\min\{\ep/4, \dt_0/2\}$ and
${\mathcal F}_1={\mathcal F}\cup {\mathcal G}_1.$

Let $\dt_1>0,$ ${\mathcal G}_2\subset C(X)$ and
 ${\mathcal P}_1\subset \underline{K}(C(X\times S^1))$ be finite subsets  required in 2.5 of \cite{Lnalm} corresponding to $\ep_1$ and
${\mathcal F}_1.$ Let $\af\in KK(C(X\times S^1), A)$ such that
$$
\af|_{{\rm ker}\rho_{C(X\times S^1)} }=\{0\}\andeqn
\af([1_{C(X)}])=[1_A].
$$
(Note that ${\rm ker}\rho_{C(X\times S^1)}=\bt_1(K_1(C(X))).$)
Without loss of generality, we may assume that $ {\mathcal G}_2={\mathcal
G}\otimes S,$ where ${\mathcal G}\subset C(X)$ is a finite subset. Put
${\mathcal P}={\mathcal P}_1\cup{\mathcal P}_0.$ We may also assume that
${\mathcal G}_1\subset {\mathcal G}_2$ and $\dt_1<\dt_0/2.$

 Define $\psi: C(X\times S^1)\to A$ by $\psi(f\otimes
g)=h(f)g(u)$ for all $f\in C(X)$ and $g\in C(S^1).$ Then,
to simplify the notation, by applying
\ref{appn}, we may assume that  $\psi$ is a
$\dt_1$-${\mathcal G}_1$-multiplicative \morp. The condition that
$K_1(A)=\{0\}$ and $\text{bott}_1(h,u)=0$ implies that
$$
[\psi]|_{{\mathcal P}}=\af|_{{\mathcal P}}.
$$
It follows from 2.5 of \cite{Lnalm} that there exists a unital
\hm\, $H_0: C(X\times S^1)\to A$ such that
\beq\label{DM1f-1}
\|h(f)g(u)-H_0(f\otimes g)\|<\ep_1\rforal f\in {\mathcal F}_1
\eneq
and $g\in S.$ Let $Y$ be the spectrum of $H_0.$ Then $Y$ is a
compact subset of $X\times S^1.$  By (3) and (4) of \ref{factd1},
the choice of $\dt_0$ and ${\mathcal G}_0$ and the condition that
$K_1(A)=\{0\},$ we have
\beq\label{DM1f-2}
[H_0]\in {\mathcal N}
\eneq

It follows from 1.14 of \cite{Lnalm} that there is a unital \hm\,
$H: C(Y)\to A$ with finite finite dimensional range such that
\beq\label{DM1f-3}
\|H_0(f\otimes g)-H(f\otimes g)\|<\ep/4\rforal f\in {\mathcal F}
\eneq
and $g=1$ or $g=z.$ In the finite dimensional commutative \SCA\,\\
$H(C(X\times S^1)),$ we find a continuous path of unitaries
$\{v_t: t\in [0,1]\}$ such that
\beq\label{DM1f-4}
v_0=H(1\otimes z),\,\,\, v_1=1\andeqn \text{Length}(\{v_t\})\le
\pi.
\eneq
The lemma then follows easily.

\end{proof}


\begin{cor}\label{bbeek}
Let $\ep>0.$ There is $\dt>0$ satisfying the following:

For any two unitaries $u$ and $v$ in a unital \CA\, $A\in {\bf
B}$ with  $K_1(A)=\{0\}$ and if
$$
\|[u,\,v]\|<\dt \andeqn \text{bott}_1(u,v)=0,
$$
then there exists a continuous path of unitaries $\{u_t: t\in
[0,1]\}$ of $A$ such that
$$
u_0=u,\,\,\, u_1=1\andeqn \|[u_t,\,v]\|<\ep.
$$
Moreover,
$$
{\rm{Length}}(\{u_t\})\le \pi +\ep.
$$
\end{cor}

\begin{thm}\label{dmyiT}
Let $X$ be a compact metric space  with dimension no more than 1.
Then, for any $\ep>0$ and any finite subset ${\mathcal F}\subset
C(X),$ there exists $\dt>0,$ $\sigma>0,$ a finite subset
${\mathcal G}\subset C(X)$ and a finite subset ${\mathcal
P}_0\subset K_0(C(X))$ and ${\mathcal P}_1\subset K_1(C(X))$
satisfying the following:

Suppose that $A$ is a unital  simple \CA\, in {\bf B},
  suppose that $h: C(X)\to A$ is a
unital \hm\, whose spectrum is $\sigma$-dense in $X$ and suppose
that there is a unitary $u\in A$  such that
\beq\label{TT2tt1}
\hspace{-0.2in}\|[h(a), u]\|<\dt\,\tforal\, a\in {\mathcal
G},\,\,\, {\rm{bott}_0}(h,u)|_{{\mathcal P}_0}=0\andeqn
{\rm{bott}_1}(h,u)|_{{\mathcal P}_1}=0.
\eneq

Then there exists a rectifiable continuous path of unitaries
$\{u_t: t\in [0,1]\}$ of $A$ such that
\beq\label{TT2tt2}
u_0=u,\,\,\,u_1=1_A\andeqn \|[h(a), u_t]\|<\ep
\eneq
for all $a\in {\mathcal F}$ and all $t\in [0,1].$ Moreover,
\beq\label{TT2tt3}
{\rm{Length}}(\{u_t\})\le 2\pi+\ep\pi
\eneq
\end{thm}

\begin{proof}
There is  (by \cite{Mar}) a sequence of one-dimensional finite CW
complex $X_n$ such that $C(X)=\lim_{n\to\infty} (C(X_n), \phi_n).$
For any $\ep>0$ and any finite subset ${\mathcal F}\subset C(X),$
there is an integer $N$ and a finite subset ${\mathcal F}_1\subset
C(X_N)$ such that for any $f\in {\mathcal F},$ there exists $g_f\in
{\mathcal F}_1$ such that
\beq\label{DMT-1}
\|f-\phi_N(g_f)\|<\ep/4
\eneq
Replacing $h$ by $h\circ \phi_N,$ we see the theorem is reduced to
the case that $X$ is a compact subset of a  finite CW complex.

Now we assume that $X$ is a compact subset of a finite CW complex.
Then, there is a sequence of finite CW complexes $X_n\supset X$ such that
$\cap_{n=1}^{\infty}X_n=X.$ In particular,
$$
\lim_{n\to\infty}{\rm dist}(X, X_n)=0.
$$
Then the theorem follows from  \ref{Cdmyi}.

\end{proof}

\begin{rem}

{\rm When \CA\, $A$ is  in ${\bf B},$ both Theorem \ref{dmyi} and
Theorem \ref{dmyiT} improve the Basic Homotopy Lemma of
\cite{BEEK}.

(1) The Basic Homotopy Lemma in \cite{BEEK}  deals with the case
that $X=S^1.$ Theorem \ref{dmyi} allows  more general spaces $X$
which is a general one-dimensional finite CW complex and may not
be embedded into the $\R^2.$ Furthermore, Theorem \ref{dmyiT}
allows any compact metric space with dimension no more than 1. It
should be noted that spaces such as Hawaii ear ring which is one
dimensional but is not a compact subset of any one-dimensional
finite CW complex.

\vspace{0.1in}

(2) The constant $\dt$ in Theorem \ref{Dm1finite} only depends on
$X,$ $\ep$ and ${\mathcal F}.$ It does not depend on the spectrum
of $h$ which is a compact subset of $X.$ In other words, $\dt$ is
independent of the choices of compact subsets of $X.$ This case is
a significant improvement since in the Basic Homotopy Lemma of
\cite{BEEK} the spectrum of the unitary $u$ is assumed to be
$\dt$-dense in $S^1$ with a possible gap. It should also be noted
that if we do not care about the independence of $\dt$ on the
subsets of $X$ the proofs of this section could  be further
simplified.

\vspace{0.1in}

 (3) In both Theorem \ref{dmyiT}  and Theorem
\ref{dmyi}, the length of the path $\{u_t: t\in [0,1]\}$ is
reduced to $2\pi+\ep.$ This is a less than half of $5\pi+1$
required by the Basic Homotopy Lemma of \cite{BEEK}. The Lemma
\ref{Dm1finite} shows that, at least for the case that
$K_1(A)=\{0\}$ such as $A$ is a simple AF-algebra, the length of
the path can be further reduced to $\pi+\ep$ which  is the best
possible  bound.
However, at this point, we do not known that the bound $2\pi+\ep$
can be further reduced in general.

 (5) As indicated in \cite{BEEK}, the path to the proof of  Basic Homotopy
 Lemma in \cite{BEEK} is long. Here we present a  shorter route. The idea
 to establish something  like \ref{d1f} and then \ref{d1hf} is
 taken from our earlier result in \cite{LS} (see also Lemma 2.3 of \cite{Ln0ki}). However, the execution
 of this idea in this section relies heavily on the results in
 \cite{Lnalm}. When the dimension of $X$ becomes more than one,
 things become much more complicated and some of
 them may not have been foreseen as we will see in the next few
 sections.

}
\end{rem}

\chapter{The Basic Homotopy Lemma for higher dimensional spaces}

\setcounter{section}{3}

\section{$K$-theory and traces}

\begin{lem}\label{1LK}
Let $C$ be a unital amenable \CA\,  and let $A$ be a unital \CA.
Let $\phi: C\otimes C(S^1)\to A$ be a unital \hm\, and let
$u=\phi(1\otimes z).$

Suppose that $\psi: C\otimes C( S^1)\to A$ is another \hm\, and
$v=\psi(1\otimes z)$ such that
$$
[\phi|_{C}]=[\psi|_{C}]\,\,\,\text{in}\,\,\, KL(C, A),\andeqn
{\rm{Bott}}(\phi|_C, u)= {\rm{Bott}}(\psi|_C, v).
$$

Then
$$
[\phi]=[\psi]\,\,\,\text{in}\,\,\,KL(C\otimes C(S^1),A).
$$

In particular,
\beq\label{1lk1}
{\rm{Bott}}(\phi|_C, u)=0
\eneq
if and only if
\beq\label{1lk2}
[\psi_0]=[\phi]\,\,\,\text{in}\,\,\, KL(C\otimes C(S^1), A)
\eneq
where $\psi_0(a\otimes f(z))=\phi(a)f(1_A)$ for all $f\in C(S^1).$

Moreover, we have the following:

{\rm (1) } Suppose that $C=C(X)$ for some connected finite CW
complex and suppose that $\text{ker}\rho_{C(X)}=\{0\}$ and
$K_1(C(X))$ has no torsion, if $[u]=0$ in $K_1(A)$ and if
\beq\label{1lk0}
{\rm{bott}_1}(\phi|_{C(X)}, u)=0,
\eneq
then
\beq\label{1lk2+1}
[\psi_0]=[\phi]\,\,\,\text{in}\,\,\, KL(C(X\times S^1), A).
\eneq

{\rm (2)} In the case that $K_0(C)$ has no torsion, or $K_0(A)$ is
divisible,  and at the same time $K_1(C)$ has no torsion, or
$K_1(A)$ is divisible, if
\beq\label{1lk1+}
{\rm{Bott}_i}(\phi, u)=0,\,\,\, i=0,1,
\eneq
then
\beq\label{1lk2+2}
[\psi_0]=[\phi]\,\,\,\text{in}\,\,\, KL(C(X\times S^1), A)
\eneq

{\rm (3)} In the case that  $K_0(A)$ is divisible and
$K_1(A)=\{0\}$ and if
\beq\label{1lk1++}
{\rm{bott}_1}(\phi|_{C(X)}, u)=0
\eneq
then
\beq\label{1lk2+3}
[\psi_0]=[\phi]\,\,\,\text{in}\,\,\, KL(C(X\times S^1), A).
\eneq

{\rm (4)}  When $K_1(A)=\{0\}$ and $K_0(C(X))$ has no torsion, and
if
\beq\label{1lk1+++}
{\rm{bott}_1}(\phi|_{C(X)}, u)=0
\eneq
then
\beq\label{1lk2+4}
[\psi_0]=[\phi]\,\,\,\text{in}\,\,\, KL(C(X\times S^1), A).
\eneq

\end{lem}

\begin{proof}
Let
$$
\bt^{(0)}: K_0(C)\to K_1(C\otimes C(S^1)), \bt^{(1)}: K_1(C)\to
K_0(C\otimes C(S^1)),
$$
$$
\bt_k^{(0)}: K_0(C, \Z/k\Z)\to K_1(C\otimes C(S^1), \Z/k\Z)\andeqn
$$
$$
\bt_k^{(1)} : K_1(C, \Z/k\Z)\to K_0(C\otimes C(S^1)), \Z/k\Z)
$$
be defined by the Bott map, $k=2,3,...$ (as in \ref{Dbot2}). Put
$B=C\otimes C(S^1).$ We have the following two commutative
diagrams.
$$
{\small \put(-160,0){$K_0(B)$} \put(0,0){$K_0(B,{\Z }/k\Z )$}
\put(180,0){$K_1(B)$} \put(-85,-40){$K_0(A)$}
\put(0,-40){$K_0(A,{\Z}/k{\Z})$} \put(105,-40){$K_1(A)$} \put(-85,
-70){$K_0(A)$} \put(0,-70){$K_1(A, {\Z}/k{\Z})$}
\put(105,-70){$K_1(A)$} \put(-160,-110){$K_0(B)$}
\put(0,-110){$K_1(B,{\Z}/k{\Z})$} \put(180,-110){$K_1(B)$}
\put(-120, 2){\vector(1,0){95}} \put(70,1){\vector(1,0){95}}
\put(-123,-3){\vector(1,-1){30}} \put(30,-3){\vector(0,-1){25}}
\put(180,-2){\vector(-1,-1){30}} \put(-45,-38){\vector(1,0){35}}
\put(70,-38){\vector(1,0){25}} \put(-147, -90){\vector(0,1){85}}
\put(-75,-60){\vector(0,1){15}} \put(115, -45){\vector(0,-1){15}}
\put(190,-7){\vector(0,-1){85}} \put(-7,-68){\vector(-1,0){35}}
\put(95,-68){\vector(-1,0){25}} \put(-123,-102){\vector(1,1){30}}
\put(175, -105){\vector(-1,1){30}} \put(30,-104){\vector(0,1){30}}
\put(-5, -108){\vector(-1,0){100}}
\put(170,-108){\vector(-1,0){95}} \put(-112,-12){$[\phi]$}
\put(15, -15){$[\phi] $} \put(145, -12){$[\phi]$}
\put(-135,-92){$[\phi]$} \put(15, -88){$[\phi]$} \put(160,
-88){$[\phi]$} }
$$
and
$$
{\small \put(-160,0){$K_0(B)$} \put(0,0){$K_0(B,{\Z }/k\Z )$}
\put(180,0){$K_1(B)$} \put(-85,-40){$K_0(A)$}
\put(0,-40){$K_0(A,{\Z}/k{\Z})$} \put(105,-40){$K_1(A)$} \put(-85,
-70){$K_0(A)$} \put(0,-70){$K_1(A, {\Z}/k{\Z})$}
\put(105,-70){$K_1(A)$} \put(-160,-110){$K_0(B)$}
\put(0,-110){$K_1(B,{\Z}/k{\Z})$} \put(180,-110){$K_1(B)$}
\put(-120, 2){\vector(1,0){95}} \put(70,1){\vector(1,0){95}}
\put(-123,-3){\vector(1,-1){30}} \put(30,-3){\vector(0,-1){25}}
\put(180,-2){\vector(-1,-1){30}} \put(-45,-38){\vector(1,0){35}}
\put(70,-38){\vector(1,0){25}} \put(-147, -90){\vector(0,1){85}}
\put(-75,-60){\vector(0,1){15}} \put(115, -45){\vector(0,-1){15}}
\put(190,-7){\vector(0,-1){85}} \put(-7,-68){\vector(-1,0){35}}
\put(95,-68){\vector(-1,0){25}} \put(-123,-102){\vector(1,1){30}}
\put(175, -105){\vector(-1,1){30}} \put(30,-104){\vector(0,1){30}}
\put(-5, -108){\vector(-1,0){100}}
\put(170,-108){\vector(-1,0){95}} \put(-112,-12){$[\psi]$}
\put(15, -15){$[\psi] $} \put(145, -12){$[\psi]$}
\put(-135,-92){$[\psi]$} \put(15, -88){$[\psi]$} \put(160,
-88){$[\psi]$} }
$$
Put $\phi'=\phi|_C$ and $\psi'=\psi|_C.$ Then
$$
K_i(C\otimes C(S^1)=K_i(C)\oplus \bt^{(i-1)}(K_{i-1}(C)\andeqn
$$
$$
K_i(C\otimes C(S^1), \Z/k\Z)=K_i(C\otimes C(S^1),\Z/k\Z)\oplus
\bt_k^{(i-1)}(K_{i-1}(C\otimes C(S^1),\Z/k\Z)),
$$
$i=0,1$ and $k=2,3,....$

 This gives
the following two commutative diagrams:
$$
{\small \put(-170,0){$\bt^{(1)}(K_1(C))$}
\put(-10,0){$\bt_k^{(1)}(K_1(C,{\Z}/k{\Z}))$}
\put(180,0){$\bt^{(0)}(K_0(C))$} \put(-85,-40){$K_0(A)$}
\put(0,-40){$K_0(A,{\Z}/k{\Z})$} \put(105,-40){$K_1(A)$} \put(-85,
-70){$K_0(A)$} \put(0,-70){$K_1(A, {\Z}/k{\Z})$}
\put(105,-70){$K_1(A)$} \put(-170,-110){$\bt^{(1)}(K_1(C))$}
\put(-10,-110){$\bt_k^{(0)}(K_0(C,{\Z}/k{\Z}))$}
\put(180,-110){$\bt^{(0)}(K_0(C))$} \put(-110,
2){\vector(1,0){95}} \put(75,1){\vector(1,0){95}}
\put(-122,-4){\vector(1,-1){30}} \put(30,-3){\vector(0,-1){25}}
\put(178,-3){\vector(-1,-1){30}} \put(-45,-38){\vector(1,0){35}}
\put(70,-38){\vector(1,0){25}} \put(-147, -90){\vector(0,1){85}}
\put(-75,-60){\vector(0,1){15}} \put(115, -45){\vector(0,-1){15}}
\put(190,-7){\vector(0,-1){85}} \put(-7,-68){\vector(-1,0){35}}
\put(95,-68){\vector(-1,0){25}} \put(-123,-100){\vector(1,1){30}}
\put(175, -105){\vector(-1,1){30}} \put(30,-104){\vector(0,1){30}}
\put(-12, -108){\vector(-1,0){98}}
\put(170,-108){\vector(-1,0){95}} \put(-112,-12){$[\phi]$}
\put(15, -15){$[\phi] $} \put(145, -12){$[\phi]$}
\put(-135,-92){$[\phi]$} \put(15, -88){$[\phi]$} \put(160,
-88){$[\phi]$} }
$$
and
$$
{\small \put(-170,0){$\bt^{(1)}(K_1(C))$}
\put(-10,0){$\bt_k^{(1)}(K_1(C,{\Z}/k{\Z}))$}
\put(180,0){$\bt^{(0)}(K_0(C))$} \put(-85,-40){$K_0(A)$}
\put(0,-40){$K_0(A,{\Z}/k{\Z})$} \put(105,-40){$K_1(A)$} \put(-85,
-70){$K_0(A)$} \put(0,-70){$K_1(A, {\Z}/k{\Z})$}
\put(105,-70){$K_1(A)$} \put(-170,-110){$\bt^{(1)}(K_1(C))$}
\put(-10,-110){$\bt_k^{(0)}(K_0(C,{\Z}/k{\Z}))$}
\put(180,-110){$\bt^{(0)}(K_0(C))$} \put(-110,
2){\vector(1,0){95}} \put(75,1){\vector(1,0){95}}
\put(-122,-4){\vector(1,-1){30}} \put(30,-3){\vector(0,-1){25}}
\put(178,-3){\vector(-1,-1){30}} \put(-45,-38){\vector(1,0){35}}
\put(70,-38){\vector(1,0){25}} \put(-147, -90){\vector(0,1){85}}
\put(-75,-60){\vector(0,1){15}} \put(115, -45){\vector(0,-1){15}}
\put(190,-7){\vector(0,-1){85}} \put(-7,-68){\vector(-1,0){35}}
\put(95,-68){\vector(-1,0){25}} \put(-123,-100){\vector(1,1){30}}
\put(175, -105){\vector(-1,1){30}} \put(30,-104){\vector(0,1){30}}
\put(-12, -108){\vector(-1,0){98}}
\put(170,-108){\vector(-1,0){95}} \put(-112,-12){$[\psi]$}
\put(15, -15){$[\psi] $} \put(145, -12){$[\psi]$}
\put(-135,-92){$[\psi]$} \put(15, -88){$[\psi]$} \put(160,
-88){$[\psi]$} }.
$$
By the assumption that
$$
\text{Bott}(\phi', u)=\text{Bott}(\psi', v),
$$
all the corresponding maps in the above two diagrams from outside
six terms into the inside six terms are the same. Then, by
combining the assumption that
$$
[\phi']=[\psi'],
$$
we conclude that
$$
[\phi]=[\psi].
$$
This proves the first part of the statement.

Note that
$$
[\psi_0|_C]=[\phi']\,\,\,\text{in}\,\,\, KL(C, A)\andeqn
$$
$$
\text{Bott}(\psi_0|_{C}, 1_A)=0.
$$
Thus, by what we have proved,
$$
[\psi_0]=[\phi]\,\,\,\text{in}\,\,\,KL(B,A)
$$
if and only if $\text{Bott}(\phi', u)=0.$

Now consider the special cases.

(1)  When $\text{ker}\rho_{C(X)}=\{0\},$ $K_0(C(X))=\Z.$ Thus if
$[u]=0$ in $K_1(A),$
\beq\label{1lkn1}
[\phi]\circ \bt^{(0)}=\{0\}
\eneq
The condition that
 $\text{bott}_1(\phi|_{C(X)}, u)=0$ implies that
 \beq\label{1lkn2}
 [\phi]\circ \bt^{(1)}=\{0\}.
 \eneq
 Moreover, $\text{ker}\rho_{C(X)}=\{0\}$ implies that $K_0(C(X))$ is torsion
free. Thus by chasing the upper half of the second diagram above,
using the exactness, we see that
\beq\label{1lkn3}
[\phi]|_{\bt_k^{(1)}(K_1(C(X),\Z/k\Z))}=\{0\},\,\,\, k=2,3,....
\eneq
Similarly, since $K_1(C(X))$ has no torsion, we have
\beq\label{1lkn4}
[\phi]|_{\bt_k^{(0)}(K_0(C(X), \Z/k\Z))}=\{0\},\,\,\,k=2,3,....
\eneq
Therefore (\ref{1lk2+1}) holds in this case.

The proof of the rest cases are similar. Note that in all these
cases,\\ $[\phi]|_{\bt^{(i)}(K_i(C(X))}=\{0\},$ $i=0,1.$ All
additional conditions, by chasing the second diagram, imply the
the maps in the middle are zero, i.e.,
\beq\label{1lkn5}
[\phi]|_{\bt_k^{(i)}(K_i(C(X),\Z/k\Z))}=\{0\},\,\,\,i=0,1,k=2,3,....
\eneq

\end{proof}

\begin{lem}\label{shk}
Let  $X$ be a connected finite CW complex and let $A$ be a unital
separable simple \CA\, with stable rank one, real rank zero and
weakly unperforated $K_0(A).$ Suppose that $\phi: C(X)\to A$ is a
unital \hm. Then, for any non-zero projection $p\in A,$ there
exists a unital monomorphisms  $h_1, h_2: C(X)\to pAp$ such that
\beq\label{shk1}
[h_1|_{C_0(Y_X)}]=[\phi_{C_0(Y_X)}] \andeqn
[h_2|_{C_0(Y_X)}]+[\phi|_{C_0(Y_X)}]=\{0\}.
\eneq

\end{lem}

\begin{proof}
Recall $Y_X=X\setminus \xi_X$ for a point $\xi_X\in X$ (see
\ref{Y-X}).

We have
\beq\label{shk2}
K_0(C(X))=\Z\oplus K_0(C_0(Y_X))\andeqn K_1(C(X))=K_1(C_0(Y_X)).
\eneq
Moreover, we may write
\beq\label{shk3}
K_0(C(X),\Z/k\Z)&=&\Z/k\Z\oplus K_0(C_0(Y_X),
\Z/k\Z)\andeqn\\\label{shk4} K_1(C(X),
\Z/k\Z)&=&K_1(C_0(Y_X),\Z/k\Z).
\eneq

Let $q\le p$ be a nonzero projection such that $p-q\not=0.$
Define an element\\
 $\gamma_1\in Hom_{\Lambda}(\underline{K}(C(X)),
\underline{K}(A))$ by $\gamma_1([1_{C(X)}])=[q]$ and
$\gamma_1(\overline{[1_{C(X)}]})=\overline{[q]},$ where
$\overline{[1_{C(X)}]}$ is the image of $[1_{C(X)}]$ under the map
from $K_0(C(X))\to K_1(C(X),\Z/k\Z)$ in $\Z/k\Z\subset
K_0(C(X))/kK_0(C(X))$ and $\overline{[q]}$ is the image of $[q]$
under the map $K_0(A)\to K_0(A, \Z/k\Z)$ in $K_0(A)/kK_0(A),$ and
define $\gamma_1(x)=[\phi](x)$ for $x\in K_0(C_0(Y_X)),$ for $x\in
K_1(C(X))=K_1(C_0(Y_X)),$ for $x\in K_0(C_0(Y_X), \Z/k\Z)$ and for
$x\in K_1(C(X),\Z/k\Z)=K_1(C_0(Y_X),\Z/k\Z),$ $k=2,3,....$ It
follows from (\ref{shk3}) and (\ref{shk4}) that a straightforward
computation shows that $\gamma_1\in
Hom_{\Lambda}(\underline{K}(C(X)), \underline{K}(A)).$ Moreover,
$\gamma_1$ is an positive element since $[\phi]$ is. Similarly,
define $\gamma_2([1_{C(X)}])=\gamma_1([1_{C(X)}]),$
$\gamma_2(\overline{[1_{C(X)}]})=\gamma_1(\overline{[1_{C(X)}]})=\overline{[q]}$
and define $\gamma_2(x)=-\gamma_1(x)$ for $x\in K_0(C_0(Y_X)),$
for $x\in K_1(C(X))=K_1(C_0(Y_X)),$ for $x\in K_0(C_0(Y_X),
\Z/k\Z)$ and for $x\in K_1(C(X),\Z/k\Z)=K_1(C_0(Y_X),\Z/k\Z),$
$k=2,3,....$ The same argument shows that $\gamma_2\in KL(C(X),
A)_+.$ It now follows from Theorem 4.7 of \cite{LnK} that there
are \hm\, $h'_1, h'_2: C(X)\to qAq$ such that
\beq\label{shk6}
[h_1'|_{C_0(Y_X)}]=[\phi|_{C_0(Y_X)}] \andeqn
[h_2'|_{C_0(Y_X)}]+[\phi|_{C_0(Y_X)}]=\{0\}.
\eneq
By applying \ref{triv}, we obtain a unital monomorphism $h_{00}: C(X)\to (p-q)A(p-q)$
which factors through $C([0,1]).$
In particular,
$$
[h_{00}|_{C_0(Y_X)}]=0\,\,\,\text{in}\,\,\,KK(C(X),A).
$$
Now define $h_1=h_{00}+h_1'$ and $h_2=h_{00}+h_2''.$

\end{proof}

\begin{lem}\label{mmm}
Let $X$ be a compact metric space and let $A$ be a unital
separable simple \CA\, with real rank zero, stable rank one and
weakly unperforated $K_0(A).$ Suppose that $\Lambda:C(X)_{s.a}\to
Aff(T(A))$ is a unital positive linear map. Then, for any
$\gamma>0$ and any finite subset ${\mathcal G}\subset C(X),$ there
exists a unital \hm\, $\phi: C(X)\to A$ with finite dimensional
range such that
\beq\label{mmm1}
|\tau\circ \phi(f)-\Lambda(f)(\tau)|<\gamma\,\tforal \, f\in
{\mathcal G}\andeqn\, \tforal\,\tau\in T(A).
\eneq
\end{lem}

The proof is contained in the proof  of Theorem 3.6 of
\cite{HLX3}. We present here for the convenience of the reader.

\begin{proof}
Let $\gamma>0$ and ${\mathcal G}\subset C(X)$ be given. Without loss
of generality, we may assume that ${\mathcal G}\subset C(X)_{s.a}.$
 Choose
$\dt=\sigma_{X. \gamma/3, {\mathcal G}}.$
 Suppose
that $\{y_1,y_2,...,y_n\}$ is $\dt/2$-dense in $X.$ Let
$$
Y_j=\{x\in X: {\rm dist}(y,y_j)< \dt \}, j=1,2,...,n.
$$

Let $\{g_1,\cdots,g_n\}$ be a partition of unity, i.e., a subset
of nonnegative functions in $C(X)$ satisfying the following
conditions:

(i) $g_i(x)=0$ for all $x\not\in Y_i,i=1,2\cdots,n$ and

(ii) $\sum_{i=1}^ng_i(x)=1$ for all $x\in X$.

For each $i$, define $\hat{g_{i}}\in\Aff(T(A))$  by

$$
\hat{g_i}(\tau)=\tau(\Lambda(g_i)).
$$
If $\Lambda(g_i)\not=0,$ then
$$
\inf \{\hat{g_{i}}(\tau): \tau\in T(A)\}>0,
$$
since $A$ is simple and $T(A)$ is compact. Put
$$
\gamma_1=\inf\{\hat{g_{i}}(\tau):\tau\in T(A),g_i\in {\mathcal G},
\Lambda(g_i)\not=0\}.
$$

Let $\gamma_2=\min\{\gamma,\gamma_1\}$. If $\hat{g_{i}}=0,$ we
choose $q_i=0.$ Since $A$ has real rank zero and the range of the
mapping $\rho_A$ is dense in $\Aff(T(A)),$  if
$\hat{g_{i}}\not=0,$ there exists a projection $q_{i}\in A$ such
that
$$
\|\rho_A(q_{i})-(\hat{g_{i}}-{\frac {\gamma_2}{6n}})\|<{\frac
{\gamma_2}{6n}}.
$$
So for all $\tau\in T(A),$ if $\hat{g_i}\not=0,$
$$
\tau(\Lambda(g_{i}))-{\frac{\gamma_2}{3n}}<\tau(q_{i})
<\tau(\Lambda(g_{i})).
$$
 Since $\Lambda\not=0,$ this implies for each $i$,
$$
1-{\frac {\gamma_2} 3}=\sum_{i=1}^n\tau(\Lambda(g_{i}))-{\frac
{\gamma_2}3} <\sum_{i=1}^n\tau(q_{i}) <
\sum_{i=1}^n\tau(\Lambda(g_{i}))=1.
$$
Let $Q=\oplus_{i=1}^n q_i$ be a projection in $M_n(A).$ It follows
that $\tau(Q) < \tau(1_A)$ for all $\tau\in T(A).$ Since
$TR(M_n(A))=0,$ $M_n(A)$ has the Fundamental Comparison Property
(see \cite{Lntrd}). So from  $\tau(Q)<\tau(1_A)$ for all $\tau\in
T(A),$ one obtains mutually orthogonal projections
$q_1',q_2',...,q_n'\in A$ such that $[q_i']=[q_i],$ $i=1,2,...,n.$
Put $p_i=q_i',$ $i=1,2,...,n-1$ and $p_n=1-\sum_{i=1}^{n-1}q_i'.$

For any $\tau\in T(A)$,
$$
\tau(\Lambda(g_n)) =1-\sum_{i=1}^{n-1}\tau(\Lambda(g_i))\le
1-\sum_{i=1}^{n-1}\tau(q_i)=\tau(p_n)
$$
and
$$
\tau(p_n)=1-\sum_{i=1}^{n-1}\tau(q_i)<1-(\sum_{i=1}^{n-1}\tau(\Lambda(g_i))-{\frac{\gamma_2}3})<\tau
(\Lambda(g_n))+{\frac {\gamma_2} 3}.
$$
So, for all $\tau\in T(A)$,
$$
|\tau(p_i)-\tau(\Lambda(g_i))|\le {\frac{\gamma_2}{3n}},\,\,\,
i=1,\cdots,n-1,
$$
and
$$
|\tau(p_n)-\tau(\Lambda(g_n)|<{\frac{\gamma_2}3}.
$$

Define homomorphism $\phi$ by $\phi(f)=\sum_{i=1}^nf(y_{i})p_i$
for $f\in C(X).$
 Then $\psi$ is a unital homomorphism from $C(X)$ into $A$ with finite
dimensional range. Now, we have,
for any $\tau\in T(A)$, that

\begin{eqnarray*}
|\tau(\Lambda(f))-\tau(\phi(f))|
&=&|\tau(\Lambda(f))-\sum_{k=1}^nf(y_i)\tau(p_i)|\cr &\le
&|\tau(\Lambda(f))-\sum_{i=1}^nf(y_i)\tau(\Lambda(g_i))|+
|\sum_{i=1}^nf(y_i)(\tau(p_i)-\tau(\Lambda(g_i)))|\cr
&\le&\sum_{i=1}^n\int_{Y_i}|f(x)-f(y_i)|g_i(x)d\mu_{\tau\circ\Lambda}+(n-1){\frac
{\gamma_2} {3n}}+{\frac {\gamma_2} 3}\cr
&<&{\frac{\gamma}{3}}+{\frac {2\gamma_2}3}\le\gamma.\cr
\end{eqnarray*}

\end{proof}

\begin{cor}\label{cmm}
Let $X$ be a compact metric space and let $A$ be a unital
separable simple \CA\, with real rank zero, stable rank one and
weakly unperforated $K_0(A).$ Suppose that $L: C(X)\to A$ is a
unital positive linear map. Then, for any $\gamma>0$ and any
finite subset ${\mathcal G}\subset C(X),$ there exists a unital \hm\,
$\phi: C(X)\to A$ such that
\beq\label{cmm1}
|\tau\circ \phi(f)-\tau\circ L(f)|<\gamma\,\rforal \, f\in {\mathcal
G}\andeqn\, \rforal\,\tau\in T(A).
\eneq
\end{cor}

\vspace{0.2in}

\section{Some finite dimensional approximations}

\begin{lem}\label{L1dig}
Let $X$ be a compact metric space, let $\ep>0,$ $\gamma>0$ and let
${\mathcal F}\subset C(X)$ be a finite subset. There exists $\dt>0$
and there exists a finite subset ${\mathcal G}\subset C(X)$ satisfy
the following: Suppose that $A$ is a unital \CA\, of real rank
zero, $\tau: A\to \C$ is a state on $A$ and $\phi: C(X)\to A$ is a
unital $\dt$-${\mathcal G}$-multiplicative \morp. Then, there is a
projection $p\in A,$ a unital $\ep$-${\mathcal F}$-multiplicative
\morp\, $\psi: C(X)\to (1-p)A(1-p)$ and a unital \hm\, $h: C(X)\to
pAp$ with finite dimensional range such that
$$
\|\phi(f)-(\psi(f)\oplus
h(f))\|<\ep\,\,\,\text{for\,\,\,all}\,\,\, f\in {\mathcal F}
$$
and
$$
\tau(1-p)<\gamma.
$$
\end{lem}

\begin{proof}
Fix $\ep_0>0$ and $\gamma_0>0.$ Fix a finite subset ${\mathcal
F}_0\subset C(X).$ Suppose that the lemma is false. There exists a
sequence of unital \CA s $\{A_n\}$ of real rank zero, a sequence
of unital $\dt_n$-${\mathcal G}_n$-multiplicative \morp s $\phi_n:
C(X)\to A_n$ and a sequence of states $\tau_n$ of $A_n,$ where
$\sum_{n=1}^{\infty}\dt_n<\infty$ and $\cup_{n=1}^{\infty} {\mathcal
G}_n$ is dense in $C(X),$ satisfying the following:
\beq\label{Ldige1}
\inf\{\sup\{\|\phi_n(f)-[(1-p_n)\phi_n(f)(1-p_n)+h_n(f)]\|: f\in
{\mathcal F}_0\}\}\ge \ep_0,
\eneq
where the infimum is taken among all possible projections $p_n\in
A_n$ with
$$
\tau_n(1-p_n)<\gamma_0\andeqn
$$
all possible unital \hm s $h_n: C(X)\to p_nAp_n,$ $n=1,2,....$

Define $\Phi: C(X)\to l^{\infty}(\{A_n\})$ by
$\Phi(f)=\{\phi_n(f)\}$ and define $q: l^{\infty}(\{B_n\})\to
q_{\infty}(\{B_n\})=l^{\infty}(\{B_n\})/c_0(\{B_n\}).$ Then
$q\circ \Phi: C(X)\to q^{\infty}(\{B_n\})$ is a \hm. There is a
compact subset $Y\subset X$ and monomorphism $\bar\Psi: C(Y)\to
q_{\infty}(\{B_n\})$ such that
$$
\bar \Psi\circ \pi=q\circ \Phi,
$$
where $\pi: C(X)\to C(Y)$ is the quotient map. Define
$$
t_n(\{a_n\})=\tau_n(a_n)
$$
for $\{a_n\}\in l^{\infty}(\{B_n\}).$ Let $\tau$ be a weak limit
of $\{\tau_n\}.$ By passing to a subsequence, if necessary, we may
assume that
$$
\lim_{n\to\infty} t_n(\{a_n\})=\tau(\{a_n\})
$$
for all $\{a_n\}\in \Phi(C(X)).$ Moreover, since for each
$\{a_n\}\in c_0(\{B_n\}),$ $\lim_{n\to\infty} t_n(a_n)=0,$ we may
view $\tau$ as a state of $q_{\infty}(\{B_n\}).$ Note that
$q_{\infty}(\{B_n\})$ has real rank zero. It follows from Lemma
2.11 of \cite{Lnalm} that there is a projection ${\bar p}\in
q_{\infty}(\{B_n\})$ satisfying the following:
\beq\label{l1dig-1}
\|{\bar \Psi}(f)-[(1-{\bar p}){\bar \Psi}(f)(1-{\bar
p})+\sum_{i=1}^mf(x_i){\bar p}_i]\|<\ep_0/3\andeqn
\eneq
\beq\label{l1dig-2}
\|(1-{\bar p})\Psi(f)-\Psi(f)(1-{\bar p})\|<\ep_0/3
\eneq
for all $f\in {\mathcal F}_0,$ where $\{{\bar p}_1, {\bar
p}_2,...,{\bar p}_m\}$ is a set of mutually orthogonal projections
and $\{x_1, x_2,...,x_m\}$ is a finite subset of $X,$ and where
$$
\tau(1-{\bar p})<\gamma_0/3.
$$
It follows easily that there is a projection $P\in
l^{\infty}(\{B_n\}))$ and are mutually orthogonal projections
$\{P_1, P_2,...,P_m\}\subset l^{\infty}(\{B_n\})$ such that
$$
q(P)={\bar p}\andeqn q(P_i)={\bar p}_i,\,\,\,i=1,2,...,m.
$$
There are projections $p^{(n)}, p^{(n)}_i\in A_n,$ $i=1,2,...,m,$
and $n=1,2,...$ such that
$$
P=\{p^{(n)}\}\andeqn P_i=\{p^{(n)}_i\},\,\,\,i=1,2,...,m.
$$
It follows that, for all sufficiently large $n,$
$$
\|\phi_n(f)-[(1-p^{(n)})\phi_n(f)(1-p^{(n)})-\sum_{i=1}^m
f(x_i)p^{(n)}_i]\|<\ep_0/2
$$
for all $f\in {\mathcal F}_0$ and
$$
\tau_n(1-p^{(n)})<\gamma_0/2.
$$
This contradicts with (\ref{Ldige1}).

\end{proof}

\begin{cor}\label{CL1dig}
Let $X$ be a compact metric space, let $\ep>0,$ $\gamma>0$ and let
${\mathcal F}\subset C(X)$ be a finite subset. There exists
$\dt>0$ and there exists a finite subset ${\mathcal G}\subset
C(X)$ satisfy the following: Suppose that $A$ is a unital \CA\, of
real rank zero, $\tau: A\to \C$ is a state on $A$ and $\phi:
C(X)\to A$ is a unital $\dt$-${\mathcal G}$-multiplicative \morp\,
and $u\in A$ is an element such that \\ $\|u^*u-1\|<\dt,$
$\|uu^*-1\|<\dt$ and
$$
\|[\phi(g), u]\|<\dt\rforal g\in {\mathcal G}.
$$

 Then, there is a
projection $p\in A,$
and a unital \hm\, $h: C(X)\to
pAp$ with finite dimensional range such that
\beq\label{l1dig-0}
\|\phi(f)g(u)-((1-p)\phi(f)g(u)(1-p)\oplus
h(f\otimes g))\|<\ep\andeqn
\eneq
\beq\label{l1dig-00}
\|[(1-p), \phi(f)]\|<\dt
\eneq
for all $f\in {\mathcal F}$ and $g=1_{C(X)}$ or $g=z,$
where
\beq\label{cl1dig-1}
h(f\otimes g)=\sum_{i=1}^n \sum_{j=1}^{m(i)}f(y_i)g(t_{i,j}) p_{i,j} \rforal f\in C(X), g\in C(S^1),
\eneq
where $y_i\in X$ and $t_{i,j}\in S^1$ and $\{p_{i,j}: i,j\}$ is a set of mutually orthogonal non-zero projections
with $\sum_{i=1}\sum_{j=1}^{m(i)}p_{i,j}=p$
and
$$
\tau(1-p)<\gamma.
$$

Moreover, for any $\eta>0,$ we may assume that
\beq\label{cL1dig-2}
\sum_{j=1}^{m(i)}p_{i,j}\subset \overline{\phi(f_i\otimes 1)A\phi(f_i\otimes 1)},
\eneq
where $f_i\in C(X_1)$ such that $0\le f_i\le 1,$ $f_i(y)=1$ if $y\in O(y_i, \eta/2)\subset X_1$ and $f_i(y)=0$ if
$y\not\in O(y_i, \eta)\subset X_1,$ $i=1,2,...,n.$
\end{cor}

\begin{proof}
We make a few modification of the proof of \ref{L1dig}.
In the proof of \ref{L1dig},
define a linear map $L_n: C(X\otimes S^1)\to B_n$ by $L_n(f\otimes g)=\phi_n(f)g(u_n),$
where $\phi_n: C(X)\to B_n$ is a sequence of $\dt_n$-${\mathcal G}_n$-multipolicative \morp s
and $u_n\in B_n$ with
$$
\lim_{n\to\infty}\|u_n^*u_n-1_{B_n}\|=0=\lim_{n\to\infty}\|u_nu_n^*-1_{B_n}\|=0\andeqn
$$
$$
\lim_{n\to\infty}\|[\phi_n(f),u_n]\|=0\rforal f\in C(X).
$$

Redefine $\Phi: C(X\times S^1)\to l^{\infty}(\{B_n\})$ by
$\Psi(f\otimes g)=\{L_n(f\otimes g)\}$ for $f\in C(X)$ and $g\in
C(S^1).$ Then ${\bar \Psi}=q\circ \Psi: C(X\times S^1)\to
q_{\infty}(\{B_n\})$ is a unital \hm, where $q:
l^{\infty}(\{B_n\})\to q_{\infty}(\{B_n\})$ is the quotient map.

In this case,  in  (\ref{l1dig-1}), we may replace (\ref{l1dig-1})
by the following:

\beq\label{cl1dig-2}
\|{\bar \Psi}(f)-[(1-{\bar p}){\bar \Psi}(f)(1-{\bar
p})+\sum_{i=1}^n\sum_{j=1}^{m(i)}f(y_i)g(t_{i,j}){\bar
p_{i,j}}]\|<\ep_0/3
\eneq
for all $f\in {\mathcal F}\andeqn g\in S,$ where $S=\{1_{C(S^1)},z\}.$
As in proof of \ref{L1dig}, we conclude that (\ref{l1dig-0}) and
(\ref{l1dig-00}) hold.

For any $\eta>0,$ from the the proof of 2.11 of \cite{Lnalm}, we may assume
that
\beq\label{cl1dig-3}
{\bar p_{i,j}}\in B_{G_i},\,\,\,j=1,2,...,m(i), i=1,2,...,n
\eneq
where $G_i=O(y_i,\eta/2)$  and $B_{G_i}$ is the hereditary \SCA\, generated by
${\bar \Psi}(f_i\otimes 1),$ $i=1,2,...,n.$

Let $B_n$ be as in the proof of \ref{L1dig}. Let $p_{i,j}^{(n)}\in
B_n$ be projections such that $q\circ (\{p_{i,j}^{(n)}\})={\bar
p_{i,j}},$ $j=1,2,...,m(i)$ and $i=1,2,...,n.$ Then
\beq\label{c1dig-4}
\lim_{n\to\infty}\|\phi_n(f_i\otimes 1)p_{i,j}^{(n)}\phi_n(f_i\otimes 1)-p_{i,j}^{(n)}\|=0.
\eneq
It follows from 2.5.3 of \cite{Lnbk} that, for each sufficiently large $n,$  there is a projection
$q_{i,j}^{(n)}\in \overline{\phi_n(f_i\otimes 1) A\phi_n(f_i\otimes 1)}$ such that
\beq\label{c1dig-5}
\lim_{n\to\infty}\|p_{i,j}^{(n)}-q_{i,j}^{(n)}\|=0.
\eneq
It follows that we may replace $p_{i,j}^{(n)}$ by $q_{i,j}^{(n)}$ and
then the corollary follows.

\end{proof}

\begin{lem}\label{L2dig}
Let $X$ be a compact metric space, let $\ep>0,$ $\gamma>0$ and let
${\mathcal F}\subset C(X)$ be a finite subset. There exists $\dt>0$
and there exists a finite subset ${\mathcal G}\subset C(X)$ satisfy
the following: Suppose that $A$ is a unital \CA\, of tracial rank
zero and $\phi: C(X)\to A$ is a unital $\dt$-${\mathcal
G}$-multiplicative \morp. Then, there is a projection $p\in A,$ a
unital $\ep$-${\mathcal F}$-multiplicative \morp\, $\psi: C(X)\to
(1-p)A(1-p)$ and a unital \hm\, $h: C(X)\to pAp$ with finite
dimensional range such that
$$
\|\phi(f)-(\psi(f)\oplus
h(f))\|<\ep\,\,\,\text{for\,\,\,all}\,\,\, f\in {\mathcal F}
$$
and
$$
\tau(1-p)<\gamma\,\,\,\text{for\,\,\,all}\,\,\, \tau\in T(A).
$$
In particular, $\psi$ can be chosen to be
$\psi(f)=(1-p)\phi(f)(1-p)$ for all $f\in C(X).$

\end{lem}

\begin{proof}

Fix $\ep>0,$ $\gamma>0$ and a finite subset ${\mathcal F}\subset
C(X).$

Let $\dt_0>0$ and ${\mathcal G}_0\subset C(X)$ be finite subset
required by \ref{L1dig} corresponding to $\ep/2$ and $\gamma/2.$

Choose $\dt=\dt_0/2$ and ${\mathcal G}={\mathcal G}_0$ and let $\phi:
C(X)\to A$ be a \hm\, which satisfies the conditions of the lemma.

Since $TR(A)=0,$ there exists a sequence of finite dimensional
\SCA s $B_n$ with $e_n=1_{B_n}$ and a sequence of \morp s $\phi_n:
A\to B_n$ such that

(1) $\lim_{n\to\infty}\|e_na-ae_n\|= 0$ for all $a\in A,$

(2) $\lim_{n\to\infty}\|\phi_n(a)-e_nae_n\|=0$ for all $a\in A$
and $\phi_n(1)=e_n;$

(3) $\lim_{n\to\infty}\|\phi_n(ab)-\phi_n(a)\phi_n(b)\|=0$ for all
$a, b\in A$ and

(4) $\tau(1-e_n)\to 0$ uniformly on $T(A).$

 We write $B_n=\oplus_{i=1}^{r(n)}D(i,n),$ where
each $D(i,n)$ is a simple finite dimensional \CA, a full matrix
algebra. Denote by $\Phi(i,n): A\to D(i,n)$ the map which is the
composition of the projection map from $B_n$ onto $D(i,n)$ with
$\phi_n.$ Denote by $\tau(i,n)$ the standard normalized trace on
$D(i,n).$
Put $\phi_{(i,n)}=\Phi(i,n)\circ \phi.$
From (1), (2), (3), (4) above, by applying \ref{L1dig} to each
$\phi_{(i,n)}, $ for each $i$ ( and $\tau(i,n)$) and all
sufficiently large $n,$ we have

\begin{eqnarray}\label{eMP1}\nonumber
\|\phi_{(i,n)}(f)-[(1_{D(i,n)}-p_{i,n})\phi_{(i,n)}(f\otimes
z)(1_{D(i,n)}-p_{i,n})+
\sum_{j=1}^{m(i)}f(x_{i,j})p_{j,i,n}]\|<\ep/4\\
\end{eqnarray}
for all $f\in {\mathcal F},$  where $x_j\in X,$ $ \{p_{j,i, n}: j,i\}$
are mutually orthogonal projections in $D(i,n),$
$\sum_{j=1}^{m(i)} p_{j,i,n}=p_{i,n}$ with
\beq\label{f77}
\tau(i,n)(e_n-p_{i,n})<\gamma/3
\eneq
and $p_{i,n}=\sum_{j=1}^m p_{j,i,n},$ $i=1,2,...,r(n).$

Put $q_n=\sum_{i=1}^{r(n)}p_{i,n}.$ For any $\tau\in T(A),$
$\tau|_{B_n}$ has the form
$$
\tau|_{B_n}=\sum_{i=1}^{r(n)} \af_{i,n} \tau(i,n),
$$
where $\af_{i,n}\ge 0$ and $\sum_{i=1}^{r(n)}\af_{i,n}<1.$ Thus,
by (\ref{f77}),
$$
\tau(e_n-q_n)<\sum_{i=1}^{r(n)}\af_{i,n}(\gamma/3)=\gamma/3
$$
Therefore
\beq\label{L2dige2}
\tau(1-p_n)<\gamma/3+\gamma/3<\gamma\,\,\,\text{for\,\,\,all}\,\,\,\tau\in
T(A).
\eneq

Define $h_n: C(X)\to B_n$ by
$$
h_n(f)=\sum_{i=1}^{r(n)} \sum_{j=1}^{m(i)}f(x_{i,j})p_{j,i,n}
$$
for all $f\in C(X)$ and define $\psi_n: C(X)\to (1-q_n)A(1-q_n)$
by
$$
\psi_n(f)=(1-q_n)\psi(f)(1-q_n)
$$
for $f\in C(X).$ From (\ref{eMP1}) and (\ref{L2dige2}), the lemma
follows by choosing $h=h_n$ and $\psi=\psi_n$ for sufficiently
large $n$ and with $\dt<\ep/2$ and ${\mathcal G}\supset {\mathcal F}.$

\end{proof}

\begin{lem}\label{CL2dig}
Let $X$ be a compact metric space, let $\ep>0,$ $\gamma>0$ and let
${\mathcal F}\subset C(X)$ be a finite subset. There exists $\dt>0$
and there exists a finite subset ${\mathcal G}\subset C(X)$ satisfy
the following: Suppose that $A$ is a unital \CA\, of tracial rank
zero and $\phi: C(X)\to A$ is a unital \hm\, and $u\in A$ is a unitary
such that
$$
\|[\phi(g),\,u]\|<\dt\rforal g\in {\mathcal G}.
$$
Then, for any $\eta>0$ and any finite subset ${\mathcal F}_1\subset C(X),$ there is a projection $p\in A$  and a unital \hm\, $h: C(X)\to pAp$ with finite
dimensional range such that
\beq\label{cl2dig-01}
\|\phi(f)g(u)-((1-p)\phi(f)g(u) (1-p)\oplus h(f\otimes g))\|<\ep
\eneq
for all $f\in {\mathcal F}, g\in S,$
\beq\label{cl2dig-02}
\|\phi(f)-(1-p)\phi(f)(1-p)\oplus h(f\otimes 1)\|<\eta\rforal f\in {\mathcal F}_1,
\eneq
\beq\label{cl2dig-03}
\|[(1-p),\,\phi(f)]\|<\eta\rforal f\in {\mathcal F}_1
\eneq
and
$$
\tau(1-p)<\gamma\,\,\,\text{for\,\,\,all}\,\,\, \tau\in T(A),
$$
where $S=\{1_{C(S^1)}, z\}.$
\end{lem}

\begin{proof}
We will use the proof of \ref{L2dig} and \ref{CL1dig} We proceed
the proof of \ref{L2dig}. By applying \ref{CL1dig}, in stead of
(\ref{eMP1}) we can have
\beq\label{C2L-1}
\hspace{-0.4in}&&\|\Phi(i,n)(\phi(f)g(u))\\
\hspace{-0.4in}&&\hspace{-0.4in}-[(1_{D(i,n)}-p_{i,n})\Phi(i,n)(\phi(f)g(u))(1_{D(i,n)}-p_{i,n})+H_{i,n}(f\otimes
g)]\|<\ep/4
\eneq
for all $f\in {\mathcal F}$ and $g\in S,$
where
\beq\label{C2L-2}
H_{i,n}(f\otimes g)=\sum_{j=1}^{L(i)}
\sum_{k=1}^{m(j)}f(x_{i,j})g(t_{i,j,k})p_{i,j,k}
\eneq
for all $f\in C(X)\andeqn g\in C(S^1),$ and where $\{p_{i,j,k}:
i,j,k\}$ are mutually orthogonal projections with
$\sum_{j,k}p_{i,j,k}=p_{i,n}.$ For any $\eta>0,$ choose
$\sigma=\sigma_{X, \ep/8, {\mathcal F}_1}.$

By applying \ref{CL1dig}, we may assume that
\beq\label{C2L-3}
p_{i,j,k}\in \overline{\Phi(i,n)\circ \phi(f_{i,j})A\Phi(i,n)\circ \phi(f_{i,j})},
\eneq
where $f_{i,j}\in C(X)$ such that $0\le f_{ij}\le 1,$ $f_{i,j}(x)=1$ if $x\in O(x_{i,j}, \sigma/2)$ and
$f_{i,j}(x)=0$ if $x\not\in O(x_{i,j},\sigma).$
By the choice of $\sigma,$ it follows, since $\phi$ is a \hm, for all sufficiently large $n,$
$$
\|\Phi(i,n)\circ \phi(f)-[(1_{D(i,n)}-p_{i,n})\Phi(i,n)\circ
\phi(f)(1_{D(i,n)}-p_{i,n})+H_{i,n}(f\otimes 1)]\|<\eta
$$
and
$$
\|[\Phi(i,n)\circ \phi(f),\Phi(i,n)(u)]\|<\eta
$$
for all $f\in {\mathcal F}_1.$
The lemma then follows.

\end{proof}

\begin{cor}\label{pcldig}
Let $X$ be a compact metric space, let $\ep>0$ and let ${\mathcal
F}\subset C(X)$ be a finite subset. There exists $\dt>0$ and there
exists a finite subset ${\mathcal G}\subset C(X)$ satisfy the
following: Suppose that $A$ is a unital \CA\, of real rank zero
and $\phi: C(X)\to A$ is a unital \hm\, and $u\in A$ is a unitary
such that
$$
\|[\phi(g),\,u]\|<\dt\rforal g\in {\mathcal G}.
$$
Then, for any $\eta>0$ and any finite subset ${\mathcal F}_1\subset
C(X),$ there is a non-zero projection $p\in A$  and a unital \hm\,
$h: C(X)\to pAp$ with finite dimensional range such that
\beq\label{pcldig-01}
\|\phi(f)g(u)-((1-p)\phi(f)g(u) (1-p)\oplus h(f\otimes
g))\|<\ep
\eneq
for all $f\in {\mathcal F}, g\in S,$
\beq\label{pcldig-02}
\|\phi(f)-(1-p)\phi(f)(1-p)\oplus h(f\otimes 1)\|<\eta\rforal f\in
{\mathcal F}_1,
\eneq
\beq\label{pcldig-03}
\|[(1-p),\,\phi(f)]\|<\eta\rforal f\in {\mathcal F}_1.
\eneq

\end{cor}

\begin{proof}
This follows from \ref{CL1dig} as in the proof of \ref{CL2dig}.
But we do not need to consider the trace.
\end{proof}

\begin{rem}\label{Rdig}
{\rm It is essentially important in both \ref{pcldig} and
\ref{CL2dig} that $\eta$ and ${\mathcal F}_1$ are independent of the
choice of $\dt$ and ${\mathcal G}$ so that $\eta$ can be made
arbitrarily small and ${\mathcal F}_1$ can be made arbitrarily large
even after $\dt$ and ${\mathcal G}$ have  been determined.

In \ref{pcldig}, $h(f)=\sum_{j=1}^mf(\xi_j)p_j$ for all $f\in
C(X\times S^1),$ where $x_j\in X\times S^1$ and $p_1,p_2,...,p_m$
is a set of mutually orthogonal non-zero projections. In
particular, we may assume that $h(f)=f(x_1)p_1,$ by replacing
$(1-p)\phi(1-p)$ by $(1-p)\phi(1-p)+\sum_{j=2}^mf(x_j)p_j.$ But
the latter term is close to $(1-p_1)\phi(1-p_1)$ within $2\ep$ in
(\ref{pcldig-01})  and within $2\eta$ in (\ref{pcldig-02}).
Moreover, (\ref{pcldig-03}) holds if one replaces $p$ by $p_1$ and
$\eta$ by $2\eta.$ Therefore, the lemma holds, if $h$ has the form
$h(f)=f(x)p$ for a single point. Moreover, it works for any
projection $q\le p.$

}\end{rem}

We need the following statement of a result of S. Zhang (1.3 of
\cite{Zhr}).

\begin{lem}{\rm (Zhang's Riesz Interpolation)}\label{zhr1}
Let $A$ be a unital \CA\, with real rank zero. If
$p_1,p_2,...,p_n$ are mutually orthogonal projections in $A$ and
$q$ is another projection such that $q$ is equivalent to a
projection in $(p_1+p_2+\cdots p_n)A(p_1+p_2+\cdots p_n),$ then
there are projections $q'$ and $q_1,q_2,...,q_n$ in $A$ such that
$q$ and $q'$ are equivalent and $q'=q_1+q_2+\cdots +q_n,$ $q_i\le
p_i,$ $i=1,2,...,n.$
\end{lem}

\begin{proof}
This follows from 1.3 of \cite{Zhr} immediately.
There is a partial isometry $v\in A$ such that
$$
v^*v=q\andeqn vv^*\le p_1+p_2\cdots +p_n.
$$
Put $e=vv^*.$ By 1.3 of \cite{Zhr}, one obtains mutually
orthogonal projections $e_1, e_2,...,e_n\in A$ such that
$$
e=e_1+e_2\cdots +e_n
$$
and there are partial isometries $v_i\in A$ such that
$$
v_i^*v_i=e_i\andeqn v_iv_i^*\le p_i,\,\,\,i=1,2,...,n.
$$
Now define $q_i=v_iv_i^*,$ $i=1,2,...,n$ and $q'=\sum_{i=1}^n q_i.$

\end{proof}
\begin{lem}\label{zhr2}
Let $A$ be a unital \CA\, of real rank zero and let $p_1, p_2,...,p_n$ be a set of mutually
orthogonal projections with $p=\sum_{i=1}^np_i.$
Suppose that $e_1,e_2,...,e_m$ are mutually orthogonal projections such that
$e=\sum_{j=1}^me_j$ and  $e$ is equivalent to $p.$ Then
there exists a unital commutative finite dimensional \SCA\, $B$ of $pAp$ which contains
$p_1,p_2,...,p_n$ and $q_1,q_2,...,q_m$ such that
$q_j$ is equivalent to $e_j,$ $j=1,2,...,m.$
\end{lem}

\begin{proof}
There is a projection $e_j'\in pAp$ such that $e_i'$ is equivalent to $e_j,$ $i=1,2,...,m.$
By Zhang's Riesz interpolation theorem (\ref{zhr1}), there are projections
$q_{1,i}\le p_i,$ $i=1,2,...,n$ such that
$q_1$ is equivalent to $e_1,$ where
$q_1=\sum_{i=1}^nq_{1,i}.$ Note that $q_{1,i}$ commutes with each $p_i,$ $i=1,2,...,n.$
We then repeat this argument to $e_j',$ $j=2,3,...,m$ and $p_i-q_{1,i},$ $i=1,2,....,n.$
The lemma then follows.

\end{proof}

\begin{lem}\label{ndig}
Let $X$ be a compact metric space, let $\ep>0,$ $\gamma>0$ and
${\mathcal F}\subset C(X)$ be a finite subset. There exists $\dt>0$
and there exists a finite subset ${\mathcal G}\subset C(X)$ satisfying
the following.

Suppose that $A$ is a unital separable simple \CA\, with tracial
rank zero and $\psi_1, \psi_2: C(X)\to A$ are two $\dt$-${\mathcal
G}$-multiplicative \morp s.
 Then
 \beq\label{ndig1}
\| \psi_i(f)-(\phi_i(f)\oplus h_i(f))\|<\ep\,\rforal f\in {\mathcal
F},
\eneq
 where $\phi_i: C(X)\to (1-p_i)A(1-p_i)$ is a unital $\ep$-${\mathcal F}$-multiplicative \morp\, and
 $h_i: C(X)\to p_iAp_i$ is a unital \hm\, with finite dimensional range ($i=1,2$) for
 two unitarily equivalent
 projections
 $p_1, p_2\in A$ for which
 \beq\label{ndig2}
 \tau(1-p_i)<\sigma\,\rforal \, \tau\in T(A).
 \eneq

\end{lem}

\begin{proof}
By applying \ref{L2dig} to $\psi_1$ and $\psi_2,$ we obtain that
$\phi_i': C(X)\to (1-p_i)A(1-p_i)$ is a unital $\ep$-${\mathcal
F}$-multiplicative \morp\, and
 $h_i': C(X)\to p_iAp_i$ is a unital \hm\, with finite dimensional range ($i=1,2$)
 such that
 $$
\| \psi_i(f)-(\phi_i'(f)\oplus h_i'(f))\|<\ep\,\rforal f\in {\mathcal
F},\,\,\, i=1,2, \andeqn
$$
$$
 \tau(1-p_i)<\sigma/2\,\rforal \, \tau\in T(A).
 $$
Write
$$
h_i'(f)=\sum_{k=1}^{m(i)} f(x_k^{(i)})p^{(i)}_k\rforal\, f\in
C(X),
$$
where $\{p^{(i)}_k: k=1,2,...,m(i)\}$ is a set of mutually
orthogonal projections and $x_k^{(i)}\in X,$ $k=1,2,...,m(i),$
$i=1,2.$ Since $A$ is a simple \CA\, with tracial rank zero, there
is a projection $q^{(i)}\le p^{(i)}$ such that
$$
[q^{(1)}]=[q^{(2)}]\andeqn \tau(p^{(i)}-q^{(i)})<\sigma/3
$$
for all $\tau\in T(A).$ By Zhang's Riesz interpolation
(\ref{zhr1}), there are $e_k^{(i)}\le p_k^{(i)}$ such that
\beq\label{ndig-0-1}
[\sum_{k=1}^{m(i)}q^{(i)}e_k^{(i)}]=[q^{(i)}],\,\,\,i=1,2.
\eneq
Let $p_i=\sum_{k=1}^{m(i)}e_k^{(i)},$ put
$$
\phi_i(f)=\phi_i'(f)+\sum_{k=1}^{m(i)}f(x_k^{(i)})(p^{(i)}_k-e_k^{(i)})\andeqn
$$
 $$
 h_i(f)=\sum_{k=1}^{m(i)}f(x_k^{(i)})e_k^{(i)}
 $$
 for all $f\in C(X),$ $i=1,2.$
 It is clear that so defined $\phi_i$ and $h_i$ satisfy the requirements.

\end{proof}



\section{The Basic Homotopy Lemma --- full spectrum}

\begin{thm}\label{1T}
Let $X$ be a connected finite CW complex (with a fixed metric), let $A$ be a unital
separable simple \CA\, with tracial rank zero and let $h: C(X)\to
A$ be a unital monomorphism. Suppose that $u\in A$ is a unitary
such that
\beq\label{1t1}
h(a)u=uh(a)\tforal \, a\in C(X)\,\,\text{and}\,\,\,
{\rm{Bott}}(h,u)=0.
\eneq
Suppose also that the \hm\,\,$\psi: C(X\times S^1)\to A$ induced
by $\psi(a\otimes z)=h(a)u$ is injective. Then, for any $\ep>0$
and any finite subset ${\mathcal F}\subset C(X),$ there is a
continuous rectifiable path of unitaries $\{u_t:t\in [0,1]\}$ of
$A$ such that
\beq\label{1t2}
\|[h(a),u_t]\|<\ep\,\,\,\text{for all}\,\,\, t\in [0,1],\,\, a\in
{\mathcal F}\,\,\,\text{and}\,\,\, u_0=u,\,\,u_1=1_A.
\eneq
Moreover,
\beq\label{1t3}
{\rm{Length}}(\{u_t\})\le \pi+\ep\pi.
\eneq
\end{thm}

\begin{proof}
Let $\ep>0$ and ${\mathcal F}\subset C(X)$ be a finite subset.
Without loss of generality, we may assume that $1_{C(X)}\in
{\mathcal F}.$ Let ${\mathcal F}_1=\{a\otimes b: a\in {\mathcal
F}, b=u,b=1\}.$ Let $\eta=\sigma_{X\times S^1, \ep/16, {\mathcal
F}_1}$ (be as in \ref{sigma}).
 Let $\{x_1,x_2,...,x_m\}$ be $\eta/2$-dense in
$X\times S^1$ and let $s\ge 1$ be an integer such that $O_i\cap
O_j=\emptyset,$ if $i\not=j,$ where
$$
O_i=\{x\in X\times S^1: {\rm dist}(x,
x_i)<\eta/2s\},\,i=1,2,...,m.
$$
Since $A$ is simple, there exists $\sigma>0$ such that
\beq\label{1t4}
\mu_{\tau\circ\psi}(O_i)\ge 2\sigma\eta,\,\,\,i=1,2,...,m
\eneq
for all $\tau\in T(A).$ Let $1/2>\gamma>0,$ $\dt>0,$ ${\mathcal G}\subset
C(X\times S^1)$ be a finite subset of $C(X\times S^1)$ and ${\mathcal
P}\subset \underline{K}(C(X\times S^1))$ be a finite subset
required by Theorem 4.6 of \cite{Lncd} associated with $\ep/16,$ $\sigma/2$
and ${\mathcal F}_1.$ We may assume that $\dt$ and ${\mathcal G}$ are so
chosen that $\overline{{\mathcal Q}_{\dt, {\mathcal G}}}\supset {\mathcal P}.$
In other words, for any ${\mathcal G}$-$\dt$-multiplicative \morp\,
$L: C(X\times S^1)\to B$ (any unital  \CA\, $B$), $[L]|_{{\mathcal
P}}$ is well-defined.

Choose $x\in X$ and let $\xi=x\times 1.$ Let $Y_X=X\setminus
\{x\}.$ Choose $\eta_1>0$ such that
$|f(\zeta)-f(\zeta')|<\min\{\dt/2, \ep/16\}$ if ${\rm
dist}(\zeta,\zeta')<\eta_1$ for all $f\in {\mathcal G}.$ Let $g\in
C(X\times S^1)$ be a nonnegative continuous function with $0\le
g\le 1$ such that $g(\zeta)=1$ if ${\rm dist}(\zeta,\xi)<\eta_1/2$
and $g(\zeta)=0$ if ${\rm dist}(\zeta, \xi)\ge \eta_1.$ Since $A$
is simple and $\psi$ is injective,
${\overline{\psi(g)A\psi(g)}}\not=\{0\}.$ Since $TR(A)=0,$ $A$ has
real rank zero. Choose two non-zero mutually orthogonal and
mutually equivalent  projections $e_1$ and $e_2$ in the simple
\CA\, ${\overline{\psi(g)A\psi(g)}}\,\text{(}\not=\{0\}\text{)}$
for which
\beq\label{1t5}
\tau(e_1+e_2)<\min\{\gamma/4, \sigma\eta/4\} \,\rforal\, \tau\in
T(A).
\eneq

It follows from \ref{shk} that there are unital \hm s $h_1:
C(X)\to e_1Ae_1$ and $h_2: C(X)\to e_2Ae_2$ such that
\beq\label{1t7}
[h_1|_{C_0(Y_X)}]=[h|_{C_0(Y_X)}]\andeqn
[h_2|_{C_0(Y_X)}]+[h_{C_0(Y_X)}]=\{0\}.
\eneq
Define $\psi_1: C(X\times S^1)\to e_1Ae_1$ by $\psi_1(a\otimes
f)=h_1(a) f(1)e_1$ for $a\in C(X)$ and $f\in C(S^1)$ and define
$\psi_2: C(X\times S^1)\to e_2Ae_2$ by  $\psi_2(a\otimes f)=h_2(a)
f(1)e_2$ for $a\in C(X)$ and $f\in C(S^1).$ Put $e=e_1+e_2.$
Define $\Psi_0: C(X\times S^1)\to (1-e_1)A(1-e_1)$ by
$\Psi_0(a\otimes f)=\psi_2(a\otimes f)\oplus (1-e)h(a)g(u)(1-e)$
for all $a\in C(X)$ and $f\in C(S^1)$ and define $\Psi: C(X)\to A$
by
\beq\label{1t8}
\Psi(b)=\psi_1(a\otimes f)\oplus \Psi_0(b)\,\rforal\, b\in
C(X\times S^1).
\eneq
Given the choice of $\eta_1,$ both $\Psi_0$ and $\Psi$ are
$\dt$-${\mathcal G}$-multiplicative.

Since $\tau(e_1)<\min\{\gamma/4, \sigma\cdot \eta/4\},$ we also
have that
\beq\label{1t9-1}
\mu_{\tau\circ \Psi}(O_i)\ge 7\sigma\eta/4,\,\, \,i=1,2,...,m.
\eneq
By \ref{cmm}, there is a unital \hm\, $\phi_0: C(X\times S^1)\to
(1-e_1)A(1-e_1)$ with finite dimensional range such that
\beq\label{1t-m}
|\tau\circ \phi_0(b)-\tau\circ \Psi_0(b)|<\gamma/2\,\rforal b\in
{\mathcal G}\andeqn
\eneq
\beq\label{1t-m2}
\mu_{\tau\circ \phi_0}(O_i)\ge \sigma\eta,\,i=1,2,...,m
\eneq
for all $\tau\in T(A).$

Note that
\beq\label{1t-n1}
\psi\approx_{\dt/2} f(\xi)(e_1+e_2)\oplus (1-e) \psi(1-e)\,\,\,
{\rm on}\,\,\, {\mathcal G}.
\eneq
Thus, by applying \ref{1LK} and by (\ref{1t7}), since
$\overline{{\mathcal P}_{\dt,{\mathcal G}}}\supset {\mathcal P},$
\beq\label{1t-mk}
[\Psi_0]|_{\mathcal P}=[\phi_0]|_{\mathcal P}.
\eneq
Note that $\tau(1-e_1)\ge 3/4$ and  ${\gamma/2\over{\tau(1-e_1)}}<\gamma$ for all $\tau\in T(A).$
It then follows from Theorem 4.6 of \cite{Lncd} and (\ref{1t-m}), (\ref{1t-m2}) and
(\ref{1t-mk}) that there exists a
unitary $w_1\in (1-e_1)A(1-e_1)$ such that
\beq\label{1t-m3}
{\rm ad}\, w_1\circ \phi_0\approx_{\ep/16}
\Psi_0\,\,\,\text{on}\,\,\, {\mathcal F}_1.
\eneq

On the other hand, we have
\beq\label{1t9}
\mu_{\tau\circ \Psi}(O_i)\ge 7\sigma\eta/4,\,\rforal\, \tau\in
T(A), \,i=1,2,...,m.
\eneq
Moreover,
\beq\label{1t10}
|\tau\circ \Psi(b)-\tau\circ \psi(b)|<\gamma\,\rforal\, b\in
{\mathcal G}\andeqn  \rforal \, \tau\in T(A).
\eneq
Furthermore, by \ref{1LK} and by (\ref{1t7}),
\beq\label{1t11}
[\Psi]|_{\mathcal P}=[\psi]|_{\mathcal P}.
\eneq
It follows from Theorem 4.6 of \cite{Lncd} that there exists a
uniatry $w_2\in A$ such that
\beq\label{1t12}
{\rm ad}\, w_2\circ \Psi\approx_{\ep/16}
\psi\,\,\,\,\text{on}\,\,\, {\mathcal F}_1
\eneq
Let $\omega=\phi_0(1\otimes z).$ Put $A_0=\phi_0(C(X\times S^1)).$
Since $\phi_0$ has finite dimensional range, $A_0$ is a finite
dimensional commutative \CA. Thus there is a continuous
rectifiable path of unitaries $\{\omega_t: t\in [0,1]\}$ in
$A_0\subset (1-e_1)A(1-e_1)$ such that
\beq\label{1t13}
\omega_0=\omega,\,\,\, \omega_1=1-e_1\andeqn \|[\phi_0(a\otimes
1),\omega_t]\|=0
\eneq
for all $a\in C(X)$ and for all $t\in [0,1],$ and
\beq\label{1t13+}
\rm{Length}(\{\omega_t\})\le \pi.
\eneq
Define
\beq\label{1t14}
U_t=w_2^*(e_1+w_1^*(\omega_t)w_1)w_2\rforal \, t\in [0,1].
\eneq
Clearly $U_1=1_A.$ Then, by (\ref{1t-m3}) and (\ref{1t12}),
\beq\nonumber
\|U_0-u\|&=&\|{\rm ad}\, w_2\circ  (\psi_1(1\otimes z)\oplus {\rm
ad}\, w_1\circ \phi_0(1\otimes z))-u\|\\\nonumber &<&\ep/16+\|{\rm
ad}\, w_2\circ (\psi_1(1\otimes z)\oplus \Psi_0(1\otimes
z))-u\|\\\label{1t15+} &<&\ep/16+\ep/16+\|\psi(1\otimes
z)-u\|=\ep/8.
\eneq
We also have, by (\ref{1t-m3}), (\ref{1t12}) (\ref{1t15+}) and
(\ref{1t13}),
\beq\label{1t15}
\|[h(a),U_t]\|=\|[\psi(a\otimes 1), U_t]\|<2(\ep/8)+\ep/16=5\ep/16
\eneq
for all $a\in {\mathcal F}.$ Moreover,
\beq\label{1t16}
{\rm{Length}}(\{U_t\})\le \pi.
\eneq
By (\ref{1t16}) and (\ref{1t15}), we obtain a continuous
rectifiable path of unitaries $\{u_t:t\in [0,1]\}$ of $A$ such
that
\beq\label{1t17}
u=u_0,\,\,\, u_1=1_A\andeqn \|[h(a), u_t]\|<\ep\,\rforal \, t\in
[0,1].
\eneq
Furthermore,
\beq\label{1t18}
{\rm{Length}}(\{u_t\})\le \pi+\ep\pi.
\eneq

\end{proof}

\begin{lem}\label{Ldig}
Let $X$ be a compact metric space, $\ep>0$ and ${\mathcal F}\subset
C(X)$ be a finite subset. Let $L\ge 1$ be an integer and let
$0<\eta\le \sigma_{X, {\mathcal F}, \ep/8}.$
 Then, for any integer
$s>0,$
any finite $\eta/2$-dense subset $\{x_1,x_2,...,x_m\}$ of $X$ for
which $O_i\cap O_j=\emptyset,\,\,{\rm if}\,\,\, i\not=j, $ where
$$
O_i=\{\xi\in X: {\rm dist}(\xi, x_i)<\eta/2s\}
$$
and any  $1/2s>\sigma>0,$ there exist a finite subset ${\mathcal
G}\subset C(X\times S^1),$ $\dt>0$ and an integer $l>0$ with
$8\pi/l<\ep$ satisfying the following:

For any unital separable simple  \CA\, $A$ with tracial rank
zero
and any $\dt$-${\mathcal G}$-multiplicative \morp\, $\phi: C(X\times
S^1)\to A,$ if $\mu_{\tau\circ \phi}(O_i\times S^1)>\sigma\cdot
\eta$ for all $\tau\in T(A)$ and  for all $i,$ then there are mutually orthogonal projections
$p_{i,1},p_{i,2},...,p_{m, l}$ in $A$ such that
$$
\|\phi(f\otimes g(z))-((1-p)\phi(f\otimes
g(z))(1-p)+\sum_{i=1}^m\sum_{j=1}^{J(i)}f(x_i)g(z_{i,j})p_{i,j}\|<\ep\rforal
f\in {\mathcal F}
$$
$$
 \andeqn \|(1-p)\phi(f\otimes g)-\phi(f\otimes g)(1-p)\|<\ep,
$$
where $g=1_{C(S^1)}$ or $g(z)=z,$ $p=\sum_{i=1}^mp_i,$
$p_i=\sum_{j=1}^k p_{i,j}$ and $z_{i,j}$ are points on the unit
circle, $1\le J(i)\le l,$
$$
\tau(p_k)>{3\sigma\over{4}}\cdot \eta\,\,\tforal \tau\in
T(A),\,k=1,2,...,m.
$$
\end{lem}

\begin{proof}
There are $f_i\in C(X)$ such that $0\le f\le 1,$ $f(x)=0$ if
$x\not\in O_i$ and
\beq\label{Ldigne1}
\tau\circ \phi(f_i)>\mu_{\phi}(O_i)-\sigma\cdot \eta/8
\eneq
for all $\tau\in T(A),$ $i=1,2,...,m.$ Put ${\mathcal F}_1={\mathcal
F}\cup\{f_1,f_2,...,f_m\},$ $S=\{1_{C(S^1)}, z\},$ and ${\mathcal
F}'={\mathcal F}_1\otimes S.$

Put $\gamma=\sigma\eta/16.$ Let $\dt>0$ and ${\mathcal G}\subset
C(X\times S^1)$ be a finite subset which are required by
\ref{L2dig} for $\min\{\ep/2,\gamma/2\}$ (in place of $\ep$),
$\gamma$ and ${\mathcal F}'$ (in place of ${\mathcal F}$).

Suppose that $\phi$ satisfies the assumption of the lemma for
$\dt$ and ${\mathcal G}$ above.

 Applying \ref{L2dig}, we obtain a projection $p,$ a \hm\,
 $H: C(X\times S^1)\to pAp$ with finite dimensional range and an
 $\ep/2$-${\mathcal G}$-multiplicative \morp\, $\psi_0: C(X\times S^1)\to (1-p)A(1-p)$ such that
 \beq\label{Ldigne2}
 \|\phi(f)-[\psi_0(f)+H(f)]\|<\min\{\ep/2, \gamma/2\}
 \eneq
 for all $f\in {\mathcal F}'$ and
 \beq\label{Ldigne3}
 \tau(1-p)<\gamma\rforal \, \tau\in T(A).
 \eneq
 Clearly we may assume that $\psi_0(f)=(1-p)\phi(f)(1-p)$ for all
 $f\in C(X).$
 Define $h(f)=H(f\otimes 1)$ for $f\in C(X).$ Then we may write that
 $$
 h(f)=\sum_{j=1}^kf(y_j)e_j\rforal f\in C(X),
 $$
 where $y_j\in X$ and $\{e_1,e_2,...,e_k\}$ is a set of mutually
 orthogonal projections.
 In particular,
\beq\label{Ldigne4}
h(f_i) =\sum_{y_j\in O_j} f(y_j)e_j,\,\,\,i=1,2,...,m.
\eneq
Put $p_i'=\sum_{y_j\in O_j}e_j.$
 It follows from (\ref{Ldigne1}),
(\ref{Ldigne2}), (\ref{Ldigne3}) and (\ref{Ldigne4})  that
\beq\label{Ldigne5}
\tau(p_i')>{3\over{4}}\sigma\cdot \eta\rforal \, \tau\in
T(A),\,\,\,i=1,2,...,m.
\eneq
Note that, if $d(x, x')<\eta,$
$$
|f(x)-f(x')|<\ep/8
$$
for all $f\in {\mathcal F}.$ It is then easy to see (by the choice of
$l$) that we may assume that
$$
H(f\otimes g)=\sum_{i=1}^m \sum_{j=1}^{J(i)}f(x_i)g(z_{i,j})
p_{i,j}
$$
for all $f\in C(X)$ and $g\in C(S^1)$ as required (with $J(i)\le
l$).

\end{proof}

\begin{rem}\label{R2dig}
In Theorem \ref{Ldig}, the integer $l$ depends only on $\ep$ and
${\mathcal F}.$ In fact, it is easy to see that $l$ can be chosen so
that $\pi M/l<\ep/2,$ where $M=\max \{\|f\|: f\in {\mathcal F}\}.$
\end{rem}

\begin{lem}\label{fullsp}
Let $X$ be a compact metric space, $\ep>0$ and ${\mathcal F}\subset
C(X)$ be a finite subset. Let $l$ be a positive integer for which
$256\pi M/l<\ep,$ where $M=\max\{1, \max\{\|f\|: f\in {\mathcal F}\}
\}.$
Let $\eta=\sigma_{X,{\mathcal F},\ep/32}$ be as \ref{sigma}. Then,  for
any finite $\eta/2$-dense subset $\{x_1,x_2,...,x_m\}$ of $X$ for
which $O_i\cap O_j=\emptyset,$ where
$$
O_i=\{x\in X: {\rm dist}(x, x_i)<\eta/2s\}
$$
for some integer $s\ge 1$ and for any $\sigma_1>$ for which
$\sigma_1<1/2s,$ and for any $\dt_0>0$ and any finite subset
${\mathcal G}_0\subset C(X\otimes S^1),$
 there exist
a finite subset ${\mathcal G}\subset C(X),$ $\dt>0$
satisfying the following:

Suppose that $A$ is a unital separable simple \CA\, with tracial
rank zero, $h: C(X)\to A$ is a unital monomorphism and $u\in A$ is
a unitary such that
\beq\label{f1}
\|[h(a), u]\|<\dt\,\,\,for\,\,all \, \,a\in {\mathcal G}\andeqn
\mu_{\tau\circ h}(O_i)\ge \sigma_1\eta\,\,for\,\,all \,\, \tau\in
T(A).
\eneq

Then there is  a $\dt_0$-${\mathcal G}_0$-multiplicative \morp\,
$\phi: C(X)\otimes C(S^1)\to A$ and a rectifiable continuous path
$\{u_t:t\in [0,1]\}$ such that
\beq\label{f2}
u_0=u,\,\, \,\|[\phi(a\otimes 1), u_t]\|<\ep\,\rforal \, a\in
{\mathcal F},
\eneq
\beq\label{f3}
\|\phi(a\otimes 1)-h(a)\|<\ep,\,\,\, \|\phi(a\otimes
z)-h(a)u\|<\ep\tforal \, a\in {\mathcal F},
\eneq
where $z\in C(S^1)$ is the standard unitary generator of $C(S^1),$
and
\beq\label{f4}
\mu_{\tau\circ \phi}(O(x_i\times
t_j))>{\sigma_1\over{2l}}\eta,\,\,\, i=1,2,...,m, j=1,2,...,l
\eneq
for all $\tau\in T(A),$ where $t_1,t_2, ...,t_l$ are $l$ points on
the unit circle which divide $S^1$ into $l$ arcs evenly and where
$$
O(x_i\times t_j)=\{x\times t\in X\times S^1: {\rm dist}(x,
x_i)<\eta/2s \andeqn {\rm dist}(t,t_j)<\pi/4sl\}
$$
for all $\tau\in T(A)$
 {\rm (}so that $O(x_i\times t_j)\cap O(x_{i'}\times t_{j'})=\emptyset$
if $(i,j)\not=(i', j')${\rm )}.

Moreover,
\beq\label{f5}
{\rm Lenghth}(\{u_t\})\le \pi+\ep\pi.
\eneq
\end{lem}

\begin{proof}
Fix $\ep$ and ${\mathcal F}.$ Without loss of generality, we may
assume that ${\mathcal F}$ is in the unit ball of $C(X).$ We may
assume that $\ep<1/5$ and $1_{C(X)}\in {\mathcal F}.$ Let $\eta>0$
such that $|f(x)-f(x')|<\ep/32$ for all $f\in {\mathcal F}$ if ${\rm
dist}(x,x')<\eta.$ Let $\{x_1,x_2,...,x_m\}$ be an $\eta/2$-dense
set. Let $s>0$ such that $O_i\cap O_j=\emptyset$ if $i\not=j,$
where
\beq\label{f6-}
O_i=O(x_i, \eta/2s)=\{x\in X: {\rm dist}(x_i,
x)<\eta/2s\},\,\,\,i=1,2,....,m.
\eneq
Suppose also $1/2s>\sigma_1>0$ and $\mu_{\tau\circ
\phi}(O_i)>\sigma_1\eta$ for all $\tau\in T(A).$  Let
$\dt_0$ and ${\mathcal G}_0$ be given. Without loss of generality (by
using smaller $\dt_0$), we may assume that ${\mathcal G}_0={\mathcal
F}_0\otimes S,$ where ${\mathcal G}_0\subset C(X)$ is a finite subset
and $S=\{1_{C(S^1)}, z\}.$ Put $\ep_1=\min\{\ep, \dt_0\}$ and
${\mathcal F}_1={\mathcal F}_0\cup{\mathcal F}.$ We may assume that
${\mathcal F}_1$ is a subset of the unit ball.

Let $\dt_1>0$ and ${\mathcal G}_1\subset C(X)$ be a finite subset (and
$8\pi /l<\ep/32$)
as required by \ref{Ldig} corresponding to $\ep_1/32,$ ${\mathcal
F}_1,$ $s$ and $\sigma_1$ above (instead of $\ep,$ ${\mathcal F},$ $s$
and $\sigma$). We may assume that $\dt_1<\ep/32$ and $1_{C(X)}\in
{\mathcal G}_1.$ We may further assume that
\beq\label{add1}
\dt_1/2\le \dt_0\andeqn {\mathcal F}_0\subset {\mathcal G}_1.
\eneq

Let $B=C(X).$ Let ${\mathcal G}_2$ a finite subset (in place of ${\mathcal
F}_1$) and $\dt_2>0$ ( in place of $\dt$) be required by Lemma
\ref{appn} associated with $\dt_1/2$ and ${\mathcal G}_1$ (in place  of
$\ep$ and ${\mathcal F}_0$).

Suppose that $A$ is a unital separable simple \CA\, with tracial
rank zero, suppose that $h: C(X)\to A$ is a unital monomorphism,
suppose that $u\in A$ is a unitary such that
\beq\label{f6}
\|[h(a), u]\|<\dt_2\, \rforal \, a\in {\mathcal G}_2\andeqn
\mu_{\tau\circ h}(O_i)>\sigma_1\eta
\eneq
for all $\tau\in T(A),$ $i=1,2,...,m.$

We may assume that $\dt_2<\dt_1$ and ${\mathcal G}_1\subset
{\mathcal G}_2.$ It follows from \ref{appn} that there exists a
$\dt_1$-${\mathcal G}_1$-multiplicative \morp\, $\psi: C(X)\otimes
C(S^1)\to A$ such that
\beq\label{f7-}
\|\psi(a)-h(a)\|<\dt_1/2\andeqn \|\psi(a\otimes z)-h(a)u\|<\dt_1/2
\eneq
for all $a\in {\mathcal G}_1.$

By applying \ref{Ldig}, we obtain a projection $p\in A$ and
mutually orthogonal projections $\{p_{i,j}: j=1,2,...,J(i),
i=1,2,...,m\}$  with $\sum_{i,j}p_{i,j}=p$ such that $J(i)\le l,$
\beq\label{nf8}
\hspace{-0.3in}\|\psi(f\otimes g)-[(1-p)\psi(f\otimes
g)(1-p)+\sum_{i=1}^m\sum_{j=1}^{J(i)}f(x_i)g(z_{i,j})p_{i,j}]\|<\ep_1/32
\eneq
for all  $ f\in {\mathcal F}_1$ and for $g\in S$ and
\beq\label{nf9}
\tau(p_i)>3\sigma_1\cdot \eta/4\,\rforal\, \tau\in T(A),
\eneq
where $p_i=\sum_{j=1}^{J(i)}p_{i,j},$ $i=1,2,...,m.$

Since $\ep<1/5,$ by (\ref{nf8}), (see for example Lemma 2.5.8 of
\cite{Lnbk}) there are unitaries $w_1\in (1-p)A(1-p),$ such that
\beq\label{nf10}
\|\psi(1\otimes z)-w_1\oplus \sum_{i=1}^{m} \sum_{j=1}^{J(i)}
z_{i,j}p_{i,j}\| <\ep/32+\ep/4
\eneq
for all $f\in {\mathcal F}_1.$ Fix $i.$ For each $j,$ since
$p_{i}Ap_{i}$ is a unital simple \CA\, with real rank zero and
stable rank one and $K_0(A)$ is weakly unperforated,
it follows from \ref{zhr2} that there is a finite dimensional commutative  unital
\SCA\, $B_i\subset p_iAp_i$ which contains projections
$p_{i,j},$ $j=1,2,...,J(i),$ and mutually orthogonal projections
$q_{i,j}\in B_i,$ $j=1,2,...,l,$ such that
$\sum_{j=1}^l q_{i,j}=p_i$ and
\beq\label{nf12}
\tau(q_{i,j})>{\tau(p_{i})\over{l}}-{\sigma_1\eta\over{16\cdot
l^2}} \,\rforal \, \tau\in T(A),
\eneq
$i=1,2,...,m.$ Define $B=\oplus_{i=1}^m B_i$ and define a unital
\hm\,\\ $\Psi_0: C(X\times S^1)\to B$ such that
\beq\label{f12+}
\Psi_0(f\otimes g)=\sum_{i=1}^m f(x_i)\sum_{j=1}^l g(t_{j})q_{i,j}\rforal f\in C(X)\andeqn g\in C(S^1).
\eneq


There is a continuous path of unitaries $\{v_t: t\in [0,1]\}$ in $B$ such that
\beq\label{f11}
v_0=w_1\oplus \sum_{i=1}^m\sum_{j=1}^l t_{i,j}q_{i,j},\,\,\,
v_1=\Psi_0(1\otimes z)\andeqn {\rm{Length}}(\{v_t\})\le \pi.
\eneq
 Define
$U_t=w_1\oplus v_t$ for $t\in [0,1]$ and define
\beq\label{nf18}
\phi(f\otimes g(z))=(1-p)\psi(f\otimes z)(1-p)+
\sum_{i=1}^m\sum_{j=1}^lf(x_i)g(t_j)q_{i,j}
\eneq
for all $f\in C(X)$ and $g\in C(S^1).$ By (\ref{nf8}), we know
that $\phi$ is $\dt_0$-${\mathcal G}_1$-multiplicative.

We compute that
\beq\label{f13}
U_0=w_1\oplus v_0,
U_1=w_1\oplus v_1,
\eneq
\beq\label{f14+}
\|U_0-u\|&\le &\|U_0-\psi(1\otimes z)\|+\|\psi(1\otimes z)-u\|\\
&<&\ep/32+\ep/4 +\dt_1<\ep/16+\ep/4
\eneq
(by (\ref{f7-}) and \ref{nf10}). Also, by  (\ref{f7-}) and
(\ref{nf8}),
\beq\label{f15}
\|\phi(a\otimes 1)-h(a)\|<\ep/2\andeqn \|\phi(a\otimes
z)-h(a)U_1\|<\ep/2
\eneq
for all $a\in {\mathcal F},$ and by (\ref{nf10}),
\beq\label{f16}
\|[U_t, \phi(f\otimes 1)]\|<\ep/32+\ep/4+\ep/32\rforal t\in [0,1].
\eneq
Moreover, by (\ref{f17}),
\beq\label{f16++}
\text{Length}(\{U_t\})\le \pi.
\eneq
Furthermore,  by (\ref{nf12}) and (\ref{nf9})
 we have
\beq\label{f17}
\tau(q_{i,s})>3\sigma_1\cdot\eta/4l-\sigma_1\cdot \eta/16l\ge
\sigma_1\cdot\eta/2l
\eneq
for all $\tau\in T(A)$
and
\beq\label{f17+}
\mu_{\tau\circ \phi}(O(x_i\times t_j))\ge {\sigma_1\over{2l}}\eta,
\eneq
$i=1,2,...,m$ and $j=1,2,...,l.$ By (\ref{f14+}), we may write
\beq\label{f18}
u(U_0)^*=\exp(ib)\,\,\,\,{\rm for\,\,\, some\,\,\,} b\in A_{s.a}
\,\,\,{\rm with}\,\,\,\|b\|\le \ep\pi
\eneq
 We define a continuous
rectifiable path of unitaries $\{u_t: t\in [0,1]\}\subset A$ by
\beq\label{f19}
u_t=exp(i(1-2t)b)U_0\,\,\,{\rm for}\,\,\, t\in [0,1/2]\andeqn
u_t=U_{2t-1}.
\eneq
So $u_0=u$ and $u_1=U_1.$ Moreover, by (\ref{f16}), (\ref{f14+})
and (\ref{f18}),
\beq\label{f20}
\|[\phi(f\otimes 1), u_t]\|<\ep\,\rforal\, t\in [0,1]
\eneq
By (\ref{f18}),
\beq\label{f21}
{\text{Length}}(\{u_t: t\in [0,1/2]\})\le \ep\pi.
\eneq
Finally, by (\ref{f16++}) and (\ref{f21}), we obtain that
\beq
{\rm Length}(\{u_t\})\le \pi+\ep\pi.
\eneq

\end{proof}

\section{The Basic Homotopy Lemma --- finite CW compleces}

\begin{df}\label{mdn}
{\rm Let $X$ be a compact metric space, let ${\varDelta}: (0,1)\to
(0,1)$ be an increasing map and let $\mu$ be a Borel probability
measure. We say $\mu$ is ${\varDelta}$- distributed, if for any
$\eta\in (0,1),$
$$
\mu(O(x,\eta))\ge \Delta(\eta)\eta
$$ for all $x\in X.$

Recall that a Borel measure on $X$ is said to be {\it strictly
positive} if for any non-empty open subset $O\subset X,$
$\mu(O)>0.$ If $\mu$ is a strictly positive Borel probability
measure, then there is always a ${\varDelta}: (0,1)\to (0,1)$ such
that $\mu$ is ${\varDelta}$-distributed. To see this, fix $\eta\in
(0,1).$ Note that there are $x_1,x_2,...,x_m\in X$ such that
$$
\cup_{i=1}^m O(x_i, \eta/2)\supset X.
$$
Let
$$
\varDelta(\eta)={\min\{\mu(O(x,\eta/2)): i=1,2,...,m\}\over{\eta}}
$$
Then, for any $x\in X,$ there is $i$ such that $x\in
O(x_i,\eta/2).$ Thus
$$
O(x, \eta)\supset O(x_i, \eta/2).
$$
Therefore
$$
\mu(O(x,\eta))\ge \varDelta(\eta)\cdot \eta\rforal x\in X.
$$
Thus $\mu$ is $\varDelta$-distributed. }

\end{df}

\begin{prop}\label{Sdis}
Let $A$ be a unital simple \CA\, with tracial state space $T(A).$
Let $\phi: C(X)\to A$ be a monomorphism. Then there is an
increasing map $\varDelta: (0,1)\to (0,1)$ such that
$\mu_{\tau\circ \phi}$ is $\varDelta$-distributed for all $\tau\in
T(A).$
\end{prop}

\begin{proof}

Fix $\eta\in (0,1).$ For each $x\in X,$  define $f_{x,\eta}\in
C(X)_+$ with $0\le f_{x,\eta}\le 1$ such that $f_{x,\eta}(x)=1$ if
$x\in O(x, \eta/4)$ and $f_{x, \eta}(x)=0$ if $x\not\in O(x,
\eta/2).$ Since $A$ is simple, there is $d(x)>0$ such that
$$
\tau(\phi(f_{x, \eta}))\ge d(x)/\eta\rforal \tau\in T(A).
$$
Thus
$$
\mu_{\tau\circ \phi}(O(x, \eta/2))\ge d(x)\cdot \eta
$$
for all $x\in X.$ Then, as shown above, one sees  that
$\mu_{\tau\circ \phi}$ is $\varDelta$-distributed for
$$
\varDelta(\eta)=\min\{d(x_i):i=1,2,...,m\}
$$
for $\eta\in (0,1)$ and for all $\tau\in T(A).$
\end{proof}

\begin{df}\label{dist}
{\rm Let $\ep>0$ and let ${\mathcal F}\subset C(X)$ be a finite
subset. Let $1>\Delta>0$ be a positive number. We say that a Borel
probability measure is ($\ep,$ ${\mathcal F}$, ${\mathcal D},$
$\Delta$)-distributed, if
there exists  an $\eta/2$-dense finite subset ${\mathcal
D}=\{x_1,x_2,...,x_m\}\subset X,$ where $0<\eta<\min\{
{\ep\over{2+\max\{\|f\|: f\in {\mathcal F}\}}}\},
\sigma_{X,\ep/32, {\mathcal F}},$ and an integer $s \ge1$
satisfies the following:

(1) $O_i\cap O_j=\emptyset,$  if $i\not=j,$

(2) $\mu(O_i)\ge \eta\cdot \varDelta,\,\,i=1,2,...,m,$ where
$$
O_i=\{x\in X: {\rm dist}(x, x_i)<\eta/2s\},\,i=1,2,...,m.
$$
Here we assume that
$$
\eta/s\le \min\{\di(x_i,x_j),i\not=j\}.
$$

If there is an increasing map $\varDelta: (0,1)\to (0,1)$ such
that $\mu$ is $\varDelta$-distributed, then, for any $\ep>0$ and
any finite subset ${\mathcal F}\subset C(X),$ there exists $\eta>0,$ a
finite subset ${\mathcal D}$ which is $\eta/2$ -dense, an integer
$s\ge 1$
 such that $\mu$ is ($\ep$, ${\mathcal F}$, ${\mathcal D},$
$\varDelta(\eta/2s)$)-distributed.
 }
\end{df}

\begin{thm}\label{MT1}
Let $X$ be a  finite CW-complex with a fixed metric. For any
$\ep>0,$ any finite subset ${\mathcal F}\subset C(X),$
$\varDelta>0,$ $\eta>0$ and $s\ge 1,$ there exists $\dt>0$ and a
finite subset ${\mathcal G}\subset C(X)$ satisfying the following:

Suppose that $A$ is a unital separable simple \CA\, with tracial
rank zero, $h: C(X)\to A$ is a unital monomorphism and $u\in A$ is
a unitary
 such that $\mu_{\tau\circ h}$ is ($\ep$, ${\mathcal F}$,${\mathcal D},$ $\varDelta$)-distributed
 for some finite ($\eta/2$-dense) subset ${\mathcal D}$ and for all $\tau\in T(A),$
\beq\label{mt1}
\|[h(a), u]\|<\dt\rforal\, a\in {\mathcal G}\andeqn
{\rm{Bott}}(h,u)=0.
\eneq
Then, there exists a continuous rectifiable path of unitaries
$\{u_t:t\in [0,1]\}$ of $A$ such that
\beq\label{mt2}
u_0=u,\,\,\,u_1=1_A\andeqn \|[h(a),u_t]\|<\ep\tforal\, a\in
{\mathcal F} \andeqn t\in [0,1].
\eneq
Moreover,
\beq\label{mt3}
{\rm{Length}}(\{u_t\})\le 2\pi+\ep\pi.
\eneq
\end{thm}

\begin{proof}
We first show that we may reduce the general case to the case that
$X$ is a connected finite CW complex.
 We may assume that
$X=\cup_{i=1}^k X_i$ which is a disjoint union of connected finite
CW complex $X_i.$ Let $e_i$ be the function in $C(X)$ which is 1
in $X_i$ and zero elsewhere. Suppose that we have shown the
theorem holds for any connected finite CW complex.

Let $\dt_i$ be required for $X_i,$ $i=1,2,...,k.$  Put
$\dt=\min\{\dt_i/2: i=1,2,...,k\}.$ Put  $P_i=h(e_i).$ Then with a
sufficiently large ${\mathcal G},$ we may assume that $\|[P_i,
u]\|<\dt.$ There is, for each $i,$ a unitary $u_i\in P_iAP_i$ such
that
$$
\|P_iuP_i-u_i\|<2\dt\andeqn \|u-\sum_{i=1}^nu_i\|<2\dt.
$$
Note that the condition $Bott(h, u)=0$ implies that
$$
Bott(h_i, u_i)=0,\,\,\,i=1,2,...,k,
$$
provided that $\dt$ is small enough. These lines of argument lead
to the reduction of the general case to the case that $X$ is a
connected finite CW complex.

For the rest of the proof, we will assume that $X$ is a connected
finite CW complex.

 Let $\ep>0$ and ${\mathcal F}\subset C(X)$ be a
finite subset. We may assume that $1\in {\mathcal F}.$ Let
$$
{\mathcal F}'=\{f\otimes a: f\in {\mathcal F}\andeqn
a=z,\,\,\,\text{or}\,\,\,a=1\}\subset C(X\times S^1).
$$
Let $\eta$ be positive so that
$$\eta<\min\{{\ep\over{1+\max\{\|f\|: f\in {\mathcal
F}\}}}, \sigma_{X,\ep/32, {\mathcal F}}\}.
$$
It follows that if $g\in {\mathcal F}',$ then
$$
|g(z)-g(z')|<\ep/32\rforal g\in {\mathcal F}'
$$
if ${\rm dist}(z,z')<\eta$ ($z,z'\in X\times S^1$). Let
$s\ge 1$ be given and $\varDelta>0.$

To simplify notation, without loss of generality, we may assume
that ${\mathcal F}$ is in the unit ball of $C(X).$  Choose $l$ to
be an integer such that $256\pi/l<\min\{\eta/4s, \ep/16\}.$ Let
$$\sigma=\varDelta/8l.$$

Let $\gamma>0,$ $\dt_1>0,$ ${\mathcal G}_1\subset C(X\times S^1)$ be a
finite subset and ${\mathcal P}\subset \underline{K}(C(X\times S^1))$
be a finite subset which are required by Theorem 4.6 of
\cite{Lncd} associated with $\ep/16,$ ${\mathcal F}'$ and $\sigma$
above. Without loss of generality, we may assume (by choosing
smaller $\dt_1$) that ${\mathcal G}_1={\mathcal F}_1\otimes S,$ where
${\mathcal F}_1\subset C(X)$ is a finite subset of $C(X)$ and $S=\{z,
1_{C(S^1)}\}.$ We may also assume that ${\mathcal F}_1\supset {\mathcal
F}$ and  $\dt_1<\ep/32.$

Let ${\mathcal P}_{C(X)}\subset \underline{K}(C(X))$ be associated with
$\text{Bott}$ map as defined in \ref{Dbot2}. We may assume that
${\boldsymbol{ \bt}}({\mathcal P}_{C(X)})\subset {\mathcal P}.$

Let $\dt>0$ and ${\mathcal G}\subset C(X)$ be a finite subset required
by Lemma (\ref{fullsp}) associated with $\dt_1$ (in place of
$\ep$), ${\mathcal F}_1$ (in place of ${\mathcal F}$), $\eta,$
$\sigma_1=\varDelta$ and with $\dt_0=\dt_1$ and ${\mathcal G}_0={\mathcal
G}_1.$

Let $h$ be a \hm\,  and $u\in A$ be a unitary in the theorem
associated with the above $\dt,$ ${\mathcal G}.$  Since $\tau\circ h$
is ($\ep$, ${\mathcal F},$ ${\mathcal D},$ $\varDelta$)-distributed, there
exists an $\eta/2$-dense subset ${\mathcal D}=\{x_1,x_2,...,x_m\}$ of
$X$ such that $O_i\cap O_j=\emptyset,$ if $i\not=j,$ where
$$
O_i=\{x\in X: {\rm dist}(x, x_i)<\eta/2s\},\,i=1,2,...,m
$$
and such that
\beq\label{mt1+}
\mu_{\tau\circ h}(O_i)\ge \eta\varDelta,\,\,\, i=1,2,...,m
\eneq
for all $\tau\in T(A).$

By Lemma \ref{fullsp}, there is a $\dt_1$-${\mathcal
G}_1$-multiplicative \morp\, $\phi: C(X\times S^1)\to A$ and a
rectifiable continuous path of unitaries $\{w_t: t\in [0,1]\}$ of
$A$ such that
\beq\label{mt2++}
w_0=u,\,\,\, \|[\phi(a),w_t]\|<\dt_1\,\rforal\, a\in {\mathcal G}_1,
\eneq
\beq\label{mt3--}
\text{Length}(\{w_t\})\le \pi+\ep/4\pi,
\eneq
\beq\label{mt3+}
\|h(a)-\phi(a\otimes 1)\|<\dt_1,\,\,\, \|\phi(a\otimes
z)-h(a)w_1\|<\dt_1 \rforal a\in {\mathcal F}_1\,\andeqn
\eneq
\beq\label{mt4}
\mu_{\tau\circ \phi}(O(x_i\times
t_j))>{\varDelta\eta\over{2l}}\,\rforal\, \tau\in T(A),
\eneq
where
$$
O(x_i\times t_j)=\{(x\times t)\in X\times S^1: {\rm dist}(x,
x_i)<\eta/2s\andeqn {\rm dist}(t,t_j)<\pi/4l\},
$$
$i=1,2,...,m, j=1,2,...,l$ and $\{t_1,t_2,...,t_l\}$ divides the
unit circle into $l$ arcs with the same length. We note that, by
(\ref{mt3+}),
\beq\label{mt41}
\|w_1-\phi(1\otimes z)\|<\dt_1.
\eneq

Put $O_{i,j}=O(x_i\times t_j),$ $i=1,2,...,m$ and $j=1,2,...,l.$
Then $O_{i,j}\cap O_{i',j'}=\emptyset$ if $(i,j)\not=(i',j').$ We
also know that $\{x_i\times t_j: i=1,2,...,m,j=1,2,...,l\}$ is
$\eta/2$-dense in $X\times S^1.$ Moreover,
\beq\label{mt5}
\mu_{\tau\circ \phi}(O_{i,j})>{4\sigma}\eta, i=1,2,...,m\andeqn
j=1,2,...,l
\eneq

Let $e_1, e_2 \in A$ be two non-zero  mutually orthogonal
projections such that
\beq\label{mt6}
\tau(e_i)<\min\{\gamma/4, {\sigma\over{4}}\eta\}\,\rforal\,
\tau\in T(A),\,\,\,i=1,2.
\eneq
 By \ref{shk}, there is a unital \hm\,
$h_1': C(X)\to e_1Ae_1$ such that
\beq\label{mt7}
[h_1'|_{C_0(Y_X)}]=[h|_{C_0(Y_X)}]\,\,\,\text{in}\,\,\,KK(C(X),
A).
\eneq
By \ref{triv},  there is a unital monomorphism $\psi_1': C(X\times
S^1)\to e_2Ae_2$ such that
\beq\label{mt7+}
[\psi_1'|_{C_0(Y_X)}]=\{0\}\,\,\,\text{in}\,\,\, KK(C(X), A),
\eneq
since it factors through $C([0,1]).$

Define $\psi_1:C(X\times S^1)\to (e_1+e_2)A(e_1+e_2)$ by
$\psi_1(a\otimes g)=h_1'(a)g(1)e_1\oplus \psi_1'(a\otimes g)$ for
all $a\in C(X)$ and $g\in C(S^1).$

By \ref{cmm}, there is a unital \hm\, $\psi_2: C(X\times S^1)\to
(1-e)A(1-e)$ with finite dimensional range such that
\beq\label{mt8}
|\tau\circ \phi(a)-\tau\circ \psi_2(a)|<\gamma/2 \andeqn
\mu_{\tau\circ \psi_2}(O_{i,j})>{3\sigma}\eta \,\rforal\, \tau\in
T(A)
\eneq
and for all $a\in {\mathcal G}_1.$

Define $\psi: C(X\times S^1)\to A$ by $\psi(a)=\psi_1(a)\oplus
\psi_2(a)$ for all $a\in C(X).$ It follows that (see  \ref{Dbot2})
\beq\label{mt9}
[\psi]|_{\mathcal P}=[\phi]|_{\mathcal P}.
\eneq

By (\ref{mt8}) and (\ref{mt6}),
\beq\label{mt10}
|\tau\circ \psi(a)-\tau\circ \phi(a)|<\gamma\,\rforal\, \tau\in
T(A)
\eneq
and for all $a\in {\mathcal G}_1,$ and
\beq\label{mt11}
\mu_{\tau\circ \psi}(O_{i,j})>{\sigma}\eta \,\rforal\, \tau\in
T(A),
\eneq
$i=1,2,...,m$ and $j=1,2,...,l.$

By Theorem 4.6 of \cite{Lncd} and the choices of $\dt_1$ and
${\mathcal G}_1,$ we obtain a unitary $Z\in A$ such that
\beq\label{mt12}
{\rm ad}\, Z\circ \psi\approx_{\ep/16} \phi\,\,\,\text{on}\,\,\,
{\mathcal F}'.
\eneq
Note that ${\rm ad}\,Z\circ \psi$ is a unital monomorphism. Put
$H_1={\rm ad}\, Z\circ \psi$ and $h_2: C(X)\to A$ by
$h_2(f)=H_1(f\otimes 1)$ for all $f\in C(X).$  So $h_2$ is a monomorphism. By (\ref{mt9}) and
${\boldsymbol{ \bt}}({\mathcal P}_{C(X)})\subset {\mathcal P},$ we conclude that
\beq\label{MT-1-}
\text{Bott}(h_2(f), H_1(1\otimes z))=0.
\eneq

 We now apply \ref{1T} to the monomorphism  $h_2$ (for $f\in C(X)$) and
 the unitary $u'=H_1(1\otimes z).$ We obtain a continuous path of unitaries $\{v_t:
t\in [0,1]\}\subset A$ such that
\beq\label{MTshort1}
v_0=H_1(1\otimes z),\,\,\, v_1=1_A,\,\,\, \|[\psi(f\otimes 1),
v_t]\|<\ep/16
\eneq
for all $t\in [0,1]$ and for all $f\in {\mathcal F}.$ Moreover,
\beq\label{MTshort2}
\text{Length}(\{v_t\})\le \pi+\ep\pi/8.
\eneq
By (\ref{mt12}) and (\ref{mt3+}),
\beq\label{MTshort3}
\|[h(a), v_t]\|<\ep/8\rforal a\in {\mathcal F}\andeqn \rforal t\in
[0,1].
\eneq
Then by (\ref{mt41})
\beq\label{MTshort4}
\|v_1-w_1\|<\ep/8.
\eneq
 By connecting the path $\{w_t\}$ and $\{v_t\}$
appropriately, we see that the theorem follows.

\end{proof}

\begin{rem}
{\rm When $X,$ $\ep$ and ${\mathcal F}$ are given, $\eta$ is also
determined. Given $\eta,$ one can find an integer $s\ge 1$ that
works for at least one $\eta/2$-dense subset
$\{x_1,x_2,...,x_m\}.$ The statement of the Theorem \ref{MT1} says
that $\dt$ does not depend on ${\mathcal D}$ as long as ${\mathcal
D}$ is $\eta/2$-dense and $\eta/s<\min\{\di(x_i,x_j): i\not=j\}.$

Such ${\mathcal D}$ is called $\eta/2$-dense and $s$-separated. If we
fix one such $s(\eta)$ which depends on the choice of $\eta,$ and
since $\eta$ is determined by $X,$ $\ep$ and ${\mathcal F},$ we may
rewrite Theorem \ref{MT1} as follows:}
\end{rem}

\begin{thm}\label{MTC}
Let $X$ be a  finite CW-complex with a fixed metric. For any
$\ep>0,$ any finite subset ${\mathcal F}\subset C(X),$
$\varDelta>0,$  there exists $\dt>0$ and a finite subset ${\mathcal
G}\subset C(X)$ satisfying the following:

Suppose that $A$ is a unital separable simple \CA\, with tracial
rank zero, $h: C(X)\to A$ is a unital monomorphism and $u\in A$ is
a unitary
 such that $\mu_{\tau\circ h}$ is ($\ep$, ${\mathcal F}$,${\mathcal D},$ $\varDelta$)-distributed
 for some finite  subset ${\mathcal D}$ (which is $\eta/2$-dense and
 $s(\eta)$-separated) and for all $\tau\in T(A),$
\beq\label{cmt1}
\|[h(a), u]\|<\dt\rforal\, a\in {\mathcal G}\andeqn
{\rm{Bott}}(h,u)=0.
\eneq
Then, there exists a continuous rectifiable path of unitaries
$\{u_t:t\in [0,1]\}$ of $A$ such that
\beq\label{cmt2}
u_0=u,\,\,\,u_1=1_A\andeqn \|[h(a),u_t]\|<\ep\tforal\, a\in
{\mathcal F}\andeqn t\in [0,1].
\eneq
Moreover,
\beq\label{cmt3}
{\rm{Length}}(\{u_t\})\le 2\pi+\ep\pi.
\eneq
\end{thm}

%

\begin{cor}\label{CT0}
Let $X$ be a connected finite CW-complex with a fixed metric and
let $\varDelta: (0,1)\to (0,1)$ be an increasing map. For any
$\ep>0,$ any finite subset ${\mathcal F}\subset C(X),$ a unital
separable simple \CA\, $A$ with tracial rank zero and a unital
monomorphism $h: C(X)\to A$ for which $\mu_{\tau\circ h}$ is
$\varDelta$-distributed for all $\tau\in T(A),$  there exists
$\dt>0$ and a finite subset ${\mathcal G}\subset C(X)$ satisfying the
following: If  $u\in A$ is a unitary  such that
\beq\label{ct0-1} \|[h(a),
u]\|<\dt\rforal\, a\in {\mathcal G}\andeqn {\rm{Bott}}(h,u)=0.
\eneq
Then, there exists a continuous rectifiable path of unitaries
$\{u_t:t\in [0,1]\}$ of $A$ such that
\beq\label{ct0-2}
u_0=u,\,\,\,u_1=1_A \andeqn \|[h(a),u_t]\|<\ep\tforal\, a\in
{\mathcal F}, t\in [0,1].
\eneq
Moreover,
\beq\label{ct0-3}
{\rm{Length}}(\{u_t\})\le 2\pi+\ep\pi.
\eneq
\end{cor}

\vspace{0.2in}

Of course we can omit the mention of $\varDelta,$ if we allows
$\dt$ depends on $h.$

\begin{cor}\label{CT1}
Let $X$ be a  finite CW-complex with a fixed metric. For any
$\ep>0,$ any finite subset ${\mathcal F}\subset C(X),$ a unital
separable simple \CA\, $A$ with tracial rank zero and a unital
monomorphism $h: C(X)\to A,$ there exists $\dt>0$ and a finite
subset ${\mathcal G}\subset C(X)$ satisfying the following: If  $u\in
A$ is a unitary  such that
\beq\nonumber
\|[h(a), u]\|<\dt\rforal\, a\in {\mathcal G}\andeqn
{\rm{Bott}}(h,u)=0.
\eneq
Then, there exists a continuous rectifiable path of unitaries
$\{u_t:t\in [0,1]\}$ of $A$ such that
\beq\nonumber
u_0=u,\,\,\,u_1=1_A\andeqn \|[h(a),u_t]\|<\ep\rforal\, a\in
{\mathcal F}\andeqn t\in [0,1].
\eneq
Moreover,
\beq\nonumber
{\rm{Length}}(\{u_t\})\le 2\pi+\ep.
\eneq
\end{cor}

\vspace{0.2in}

\section{The Basic Homotopy Lemma --- compact metric spaces}

\begin{thm}\label{MT2}
Let $X$ be a compact metric space which is a compact subset of a
finite CW complex. For any $\ep>0,$ any finite subset ${\mathcal
F}\subset C(X)$ and any map ${\varDelta}: (0,1)\to (0,1),$ there
exists $\dt>0,$ a finite subset ${\mathcal G}\subset C(X)$  and a
finite subset ${\mathcal P}\subset \underline{K}(C(X))$ satisfying the
following:

Suppose that $A$ is a unital separable simple \CA\, with tracial
rank zero, $h: C(X)\to A$ is a unital monomorphism and $u\in A$ is
a unitary
 such that $\mu_{\tau\circ h}$ is $\varDelta$-distributed for all $\tau\in T(A),$
\beq\label{2mt1}
\|[h(a), u]\|<\dt\tforal\, a\in {\mathcal G} \andeqn
{\rm{Bott}}(h,u)|_{\mathcal P}=0.
\eneq
Then, there exists a continuous rectifiable path of unitaries
$\{u_t:t\in [0,1]\}$ of $A$ such that
\beq\label{2mt2}
u_0=u,\,\,\,u_1=1_A\andeqn \|[h(a),u_t]\|<\ep\tforal\, a\in
{\mathcal F}\andeqn t\in [0,1].
\eneq
Moreover,
\beq\label{2mt3}
{\rm{Length}}(\{u_t\})\le 2\pi+\ep\pi
\eneq
\end{thm}

\begin{proof}
Without of loss of generality, we may assume that there is a
finite CW complex $Y'$ such that $X\subset Y'$ is a compact subset
(with a fixed metric).

 Let $\ep>0$ and ${\mathcal F}\subset C(X)$ be a
finite subset. We may assume that $1\in {\mathcal F}.$  Let
$S=\{1_{C(S^1)}, z\}.$ To simplify notation, without loss of
generality, we may assume that ${\mathcal F}$ is a subset of the unit
ball of $C(X).$

Let
$$
{\mathcal F}'=\{f\otimes a: f\in {\mathcal F}\andeqn a\in S\}\subset
C(X\times S^1).
$$
Let $\eta_1>0$ so that
$$\eta_1<\min\{{\ep\over{64(1+\max\{\|f\|: f\in {\mathcal
F}\})}}, \sigma_{X,\ep/64, {\mathcal F}}\}.
$$
It follows that if $g\in {\mathcal F}',$ then
$$
|g(z)-g(z')|<\ep/64\rforal g\in {\mathcal F}'
$$
if ${\rm dist}(z,z')<\eta$ ($z,z'\in X\times S^1$).

By the Tietz Extension Theorem, one has $g_f\in C(Y')$ such that
$(g_f)|_{X}=f$ and $\|g_f\|=\|f\|$ for all $f\in {\mathcal F}.$
Put ${\mathcal G}_{00}'=\{g_f\in C(Y'): f\in {\mathcal F}\}.$ Let
$\eta_0=\sigma_{Y', \ep/64,{\mathcal G}_{00}'}.$ It is easy to
construct a finite CW complex $Y\subset Y'$ such that $X\subset Y$
and $X$ is $\min\{\eta_0/2, \eta_1/4\}$-dense in $Y.$ Put
$${\mathcal G}_0=\{g\in C(Y): g|_Y=g_f\,\,\,\text{for\,\,\,some}\,\,\, f\in {\mathcal F}\}.$$
We estimate that
\beq\label{sigma-1}
\sigma_{Y, \ep/32,{\mathcal G}_0}\ge \sigma_{X, \ep/64,{\mathcal F}}.
\eneq

Let $\eta=\eta_1/2.$ Put
$$
{\mathcal G}_0'=\{g\otimes  f\in C(Y\times S^1): f\in {\mathcal G}_0
\andeqn g\in S\}.
$$

It follows that if $g\in {\mathcal G}_0',$ then
$$
|g(z)-g(z')|<\ep/32\rforal g\in {\mathcal G}_0'
$$
if ${\rm dist}(z,z')<\eta$ ($z,z'\in Y\times S^1$).

Suppose that $X$ contains an $\eta/2$-dense subset which is
$s(\eta)$-separated. Let $\varDelta: (0,1)\to (0,1)$ be an
increasing map and $\varDelta=\varDelta(\eta/2s).$

 Choose an integer $l$ such that $256\pi/l<\min\{\eta/4s, \ep/4\}.$ Let
$$\sigma=\varDelta/8l.$$

Let $\gamma>0,$ $\dt_1>0,$ ${\mathcal G}_1\subset C(Y\times S^1)$ be a
finite subset and ${\mathcal P}_1\subset \underline{K}(C(Y\times
S^1))$ be a finite subset which are required by Theorem 4.6 of
\cite{Lncd} associated with $\ep/4,$ ${\mathcal G}_0'$ and $\sigma$
above. Without loss of generality, we may assume (by choosing
smaller $\dt_1$) that ${\mathcal G}_1={\mathcal F}_1\otimes S,$ where
${\mathcal F}_1\subset C(Y)$ is a finite subset of $C(Y).$ We may also
assume that ${\mathcal F}_1\supset {\mathcal F}$ and  $\dt_1<\ep/4.$

Let $\dt>0$ and ${\mathcal G}\subset C(X)$ be a finite subset required
by Lemma (\ref{fullsp}) associated with $\dt_1$ (in place of
$\ep$), ${\mathcal F}_1$ (in place of ${\mathcal F}$), $\eta,$
$\sigma_1=\varDelta$ and with $\dt_1$  (in place of $\dt_0$)  and
${\mathcal G}_1$ (in place of ${\mathcal G}_0$).

Let ${\mathcal P}=[\theta]({{\mathcal P}_1}),$ where $\theta:
C(Y)\to C(X)$ is the quotient map. Let $h$ be a \hm\, and $u\in A$
be a unitary in the theorem associated with the above $\dt,$
${\mathcal G}$ and ${\mathcal P}.$ Moreover, we assume that
$\tau\circ h$ is ($\ep$, ${\mathcal F},$ ${\mathcal D},$
$\varDelta$)-distributed, so there exists an $\eta/2$-dense subset
${\mathcal D}=\{x_1,x_2,...,x_m\}$ of $X$ such that $O_i\cap
O_j=\emptyset,$ if $i\not=j,$ where
$$
O_i=\{x\in X: {\rm dist}(x, x_i)<\eta/2s\},\,\,\,i=1,2,...,m
$$
for $s=s(\eta)>0$ given above and such that
\beq\label{2mt-1+}
\mu_{\tau\circ h}(O_i)\ge \eta\varDelta,\,\,\, i=1,2,...,m
\eneq
for all $\tau\in T(A).$ Note here we assume that
$\eta/s<\min\{\di(x_i,x_j): i\not=j\}.$

Note also that $\{x_1,x_2,...,x_m\}$ is $\eta_1/2$-dense in $Y.$

By Lemma \ref{fullsp}, there is a ${\mathcal
G}_1$-$\dt_1$-multiplicative \morp\, $\phi: C(X\times S^1)\to A$
and a rectifiable continuous path of unitaries $\{w_t: t\in
[0,1]\}$ of $A$ such that
\beq\label{2mt2++}
w_0=u,\,\,\, \|[\phi(a),w_t]\|<\dt_1\,\rforal\, a\in {\mathcal G}_1,
\eneq
\beq\label{2mt3--}
\text{Length}(\{w_t\})\le \pi+\ep\cdot \pi/4,
\eneq
\beq\label{2mt3+}
\|h(a)-\phi(a\otimes 1)\|<\dt_1,\,\,\, \|\phi(a\otimes
z)-h(a)w_1\|<\dt_1 \rforal a\in {\mathcal F}_1\,\andeqn
\eneq
\beq\label{2mt4}
\mu_{\tau\circ \phi}(O(x_i\times
t_j))>{\varDelta\eta\over{2l}}\,\rforal\, \tau\in T(A),
\eneq
where
$$
O(x_i\times t_j)=\{(x\times t)\in X\times S^1: {\rm dist}(x,
x_i)<\eta/2s\andeqn {\rm dist}(t,t_j)<\pi/4sl\},
$$
$i=1,2,...,m, j=1,2,...,l$ and $\{t_1,t_2,...,t_l\}$ divides the
unit circle into $l$ arcs with the same length.

Put
$$
{\tilde O_{i,j}}=\{( y\times t)\in Y\times S^1:
\di(y,x_i)<\eta/2s\andeqn \di(t,t_j)<\pi/4sl\},
$$
$i=1,2,...,m$ and $j=1,2,...,l.$ Then ${\tilde  O}_{i,j}\cap
{\tilde  O}_{i',j'}=\emptyset$ if $(i,j)\not=(i',j').$ We also
know that $\{x_i\times t_j: i=1,2,...,m,j=1,2,...,l\}$ is
$\eta/2$-dense in $X\times S^1.$ Define $\Phi: C(Y\times S^1)\to
A$ by $\Phi(f)=\phi(\theta(f))$ for all $f\in C(Y\times S^1).$
Moreover,
\beq\label{2mt5}
\mu_{\tau\circ \Phi}({\tilde O}_{i,j})>{4\sigma}\eta,
i=1,2,...,m\andeqn j=1,2,...,l.
\eneq

Write $Y=\sqcup_{k=1}^RY_k,$ where each $Y_k$ is a connected
finite CW complex. Fix a  point $\xi_k\in  Y_k$ for each $k.$ Let
$E_k$ be the characteristic function for $Y_k,$ $k=1,2,...,R.$ For
the convenience, one may assume that $[E_k]\in {\mathcal P}_1.$ Denote
by $E_k'=E_k\times 1_{C(S^1)},$ $k=1,2,...,R.$

Put $e_k=h(\theta(E_k)),$ $k=1,2,...,R.$ Choose nonzero projection
$d_k\le e_k$ in $A,$ $k=1,2,...,R$ such that
\beq\label{2mt6}
\tau(d_k)<\min\{{\gamma\over{4R}}, {\sigma\cdot
\eta\over{4R}}\}\,\rforal\, \tau\in T(A),
\eneq
$k=1,2,...,R.$ We may assume that $e_k-d_k\not=0,$ $k=1,2,...,R.$

By \ref{cmm},  there is a unital \hm\, $\psi_{2,k}: C(Y_i\times
S^1)\to (e_k-d_k)A(e_k-d_k)$ with finite dimensional range such
that
\beq\label{2mt8}
|\tau\circ (\phi(E_k'aE_k'))-\tau\circ
\psi_{2,k}(a)|<\min\{\gamma/2R, {\sigma\cdot \eta\over{2R}}\}
\rforal \tau\in T(A)
\eneq
and for all $a\in {\mathcal G}_1,$ $k=1,2,...,R.$ Set $e=\sum_{i=1}^Rd_i.$ Then, by
(\ref{2mt6}),
\beq\label{2mt6+}
\tau(e)<\min\{{\gamma\over{4}},{\sigma\cdot \eta\over{4}}\}\rforal
\tau\in T(A).
\eneq

Define $\psi_2: C(Y\times S^1)\to (1-e)A(1-e)$ by
$\psi_2(f)=\sum_{k=1}^R\psi_{2,k}(E_k'f)$ for $f\in C(Y\times
S^1).$ Then we have
\beq\label{2mt8+}
|\tau\circ (\phi\circ\theta(a))-\tau\circ
\psi_{2}(a)|<\min\{\gamma/2, {\sigma\cdot \eta\over{2}}\}\andeqn
\mu_{\tau\circ \psi_2}({\tilde O}_{i,j})>{3\sigma}\eta
\eneq
 for all $\tau\in  T(A)$
 and for all $a\in {\mathcal G}_1.$

By \ref{shk}, there is a unital monomorphism  $h_{1,k}: C(Y)\to
d_kAd_k$ such that
\beq\label{2mt7}
[h_{1,k}|_{C_0(Y_k\setminus \{\xi_k\})}]=[(h\circ
\theta)|_{C_0(Y_k\setminus
\{\xi_k\})}]\,\,\,\text{in}\,\,\,KL(C(Y_k), A),
\eneq
$k=1,2,....,R.$ Define $h_1(f)=\oplus_{k=1}^Rh_{1,k}(E_kf)$ for
$f\in C(Y).$

Define $\psi_1:C(Y\times S^1)\to eAe$ by $\psi_1(a\otimes
g)=h_1(a)g(1)e$ for all $a\in C(Y)$ and $g\in C(S^1).$

Define $\psi: C(Y\times S^1)\to A$ by $\psi(a)=\psi_1(a)\oplus
\psi_2(a)$ for all $a\in C(Y).$ We have $h_1(E_k)=e_k=h\circ
\theta(E_k),$ $k=1,2,...,R.$ It follows from this, (\ref{2mt7})
and \ref{1LK} that
\beq\label{2mt9}
[\psi]|_{{\mathcal P}_1}=[\phi\circ \theta]|_{{\mathcal P}_1}
\eneq

By (\ref{2mt8+}) and (\ref{2mt6+}),
\beq\label{2mt10}
|\tau\circ \psi(a)-\tau\circ \phi\circ
\theta(a)|<\gamma\,\rforal\, \tau\in T(A)
\eneq
and for all $a\in {\mathcal G}_1,$ and
\beq\label{2mt11}
\mu_{\tau\circ \psi}({\tilde O_{i,j}})>{\sigma}\eta \,\rforal\,
\tau\in T(A),
\eneq
$i=1,2,...,m$ and $j=1,2,...,l.$

By Theorem 4.6 of \cite{Lncd} and the choices of $\dt_1$ and
${\mathcal G}_1,$ we obtain a unitary $Z\in A$ such that
\beq\label{2mt12}
{\rm ad}\, Z\circ \psi\approx_{\ep/4} \phi\circ
\theta\,\,\,\text{on}\,\,\, {\mathcal G}_0'.
\eneq

Since $\psi_2$ has finite dimensional range, from the definition
of $\psi_1$ and $\psi,$ we obtain a rectifiable continuous path of
unitaries $\{U_t: t\in [0,1]\}$ of $A$ with
\beq\label{2addleng1}
\text{Length}(\{U_t\})\le \pi+\ep\pi/2
\eneq
such that
\beq\label{2mt13}
U_0=\psi(1\otimes z),\,\,\, U_1=1_A\andeqn \|[\psi(a\otimes 1),
U_t]\|=0
\eneq
for all $a\in {\mathcal F}$ and for all $t\in [0,1].$

Put $W_t=Z^*U_tZ$ for $t\in [0,1].$ Then $W_1=1_A.$ We also have
that (by (\ref{2mt3+}) and (\ref{mt12}))
\beq\label{2mt14}
\|W_0-u\|<\dt_1+\|W_0-\phi(1\otimes z)\|<\ep/4+\ep/4=\ep/2.
\eneq
Moreover,
\beq\label{2mt15}
\|[W_t, \phi(a\otimes 1)\|<\ep/2\,\rforal\, a\in {\mathcal F}.
\eneq
By  (\ref{2mt3+}), (\ref{2mt12}) and using the fact (\ref{2mt14})
and (\ref{2addleng1}), and  using the path $\{w_t\}$ and the path
$\{W_t:t\in [0,1]\},$ we obtain a rectifiable continuous path of
unitaries $\{u_t: t\in [0,1]\}$ of $A$ such that
\beq\label{2mt16}
u_0=u,\,\,\, u_1=1_A,\,\,\,\|[h(a), u_t]\|<\ep
\eneq
for all $a\in {\mathcal F}$ and for all $t\in [0,1]$ and
\beq\label{2mt17}
\text{Length}(\{u_t\})\le 2\pi+\ep\pi.
\eneq

\end{proof}

\begin{rem}\label{Rmt3}
{\rm One can absorb the proof of Theorem \ref{MT1} into that of
Theorem \ref{MT2}. However, it does not seem necessarily helpful
to mix the issue of dependence of $\dt$ on a local measure
distribution ($X, {\mathcal F}, {\mathcal D}, \varDelta$) together
with other constants such as $\eta$ and $s$ with the issue of
reduction of a compact subset $X$ of a finite CW complex to the
special situation that $X$ is a finite CW complex. }
\end{rem}

\begin{thm}\label{MT3}
Let $X$ be a compact metric space. For any $\ep>0,$ any finite
subset ${\mathcal F}\subset C(X)$ and any map ${\varDelta}: (0,1)\to
(0,1),$ there exists $\dt>0,$ a finite subset ${\mathcal G}\subset
C(X)$  and a finite subset ${\mathcal P}\subset \underline{K}(C(X))$
satisfying the following:

Suppose that $A$ is a unital separable simple \CA\, with tracial
rank zero, $h: C(X)\to A$ is a unital monomorphism and $u\in A$ is
a unitary
 such that $\mu_{\tau\circ h}$ is $\varDelta$-distributed for all $\tau\in T(A),$
\beq\label{3mt1}
\|[h(a), u]\|<\dt\rforal\, a\in {\mathcal G} \andeqn
{\rm{Bott}}(h,u)|_{\mathcal P}=0.
\eneq
Then, there exists a continuous rectifiable path of unitaries
$\{u_t:t\in [0,1]\}$ of $A$ such that
\beq\label{3mt2}
u_0=u,\,\,\,u_1=1_A\andeqn \|[h(a),u_t]\|<\ep\tforal\, a\in
{\mathcal F}\andeqn t\in [0,1].
\eneq
Moreover,
\beq\label{3mt3}
{\rm{Length}}(\{u_t\})\le 2\pi+\ep\pi.
\eneq
\end{thm}

\begin{proof}
It is well known (see, for example, \cite{Mar}) that there exists a
sequence of finite CW complex $Y_n$ such that
$C(X)=\lim_{n\to\infty}(C(Y_n), \lambda_n).$ By replacing $Y_n$ by
one of its compact subset $X_n,$ we may write
$C(X)=\lim_{n\to\infty}(C(X_n), \lambda_n),$ where each
$\lambda_n$ is a monomorphism form $C(X_n)$ into $C(X).$ Thus
there is a surjective homeomorphism $\af_n: X\to X_n$ such that
$\lambda_n(f)=f\circ \af_n$ for all $f\in C(X_n).$

Now fix an $\ep>0$ and a finite subset ${\mathcal F}\subset C(X).$ Let
$\eta_1=\sigma_{X, \ep/64,{\mathcal F}}.$

There is an integer $N\ge 1$ and a finite subset ${\mathcal
F}_1\subset C(X_N)$ such that, for each $f\in {\mathcal F},$ there
exists $g_f\in C(X_N)$ such that
\beq\label{3mt-1}
\|f-\lambda_N(g_f)\|<\ep/4.
\eneq

Suppose that $\varDelta: (0,1)\to (0,1)$ is an increasing map. For
each $\eta>0$ and $y\in X_N,$ there is $x(y)\in X$ and
$r(y,\eta)>0$ such that
\beq\label{3mt-2}
\af_N(O(x(y), r(y)))\subset O(y,\eta).
\eneq
Since
$$
\cup_{x\in {X_N}}O(x,\eta/2)\supset X_N,
$$
there are $y_1,y_2,...,y_m\in X_N$ such that
$$
\cup_{k=1}^m O(y_i, \eta/2)\supset X_N.
$$
Now for any $y\in X_N,$ $y\in O(y_k,\eta/2)$ for some $k\in
\{1,2,...,m\}.$ Then
$$
O(y,\eta)\supset O(y_k,\eta/2).
$$
Put
$$
{\bar \eta}=\min\{r(y_k,\eta): k=1,2,...,m\}>0
$$
Define $\varDelta_1: (0,1)\to (0,1)$ by
$$
\varDelta_1(\eta)=\varDelta({\bar \eta})\rforal \eta\in(0,1).
$$

Let $\dt>0,$  ${\mathcal G}_1\subset C(X_N)$ be a finite subset and
${\mathcal P}_1\subset \underline{K}(C(X_N))$ be a finite subset
required by \ref{MT2} associated with $\ep/4$ (in place of $\ep$),
$\lambda_N({\mathcal F}_1)$ (in place of ${\mathcal F}$) and $\varDelta_1:
(0,1)\to (0,1).$

Put ${\mathcal P}=[(\lambda_N)]({\mathcal P}_1)$ in $\underline{K}(C(X)).$

Suppose that $h: C(X)\to A$ is a unital monomorphism and $u\in A$
is a unitary satisfying the conditions in the theorem associated
with the above $\dt,$ ${\mathcal G},$ ${\mathcal P}.$

Define $h_1: C(Y)\to A$ by $h_1\circ \lambda_N=h.$ Then
$$
\|[h_1(g), u]\|<\dt\rforal g\in {\mathcal G}_1,
$$
$$
\text{Bott}(h_1,u)|_{{\mathcal P}_1}=0.
$$
Moreover one checks that $\mu_{\tau\circ h_1}$ is $\varDelta_1$-distributed. By
the choices of $\dt,$ ${\mathcal G}_1$ and ${\mathcal P}_1,$ by applying
\ref{MT2}, one obtains a continuous path of unitaries $\{u_t: t\in
[0,1]\}\subset A$ with
$$
\text{length}(\{u_t\})\le 2\pi+\ep\pi
$$
such that
\beq\label{3mt-4}
u_0=u,\,\,\, u_1=1\andeqn \|[h_1(f), u_t]\|<\ep/4
\eneq
for all $f\in {\mathcal F}_1.$ It follows from (\ref{3mt-4}) and
(\ref{3mt-1}) that
$$
\|[h(f),u_t]\|<\ep\rforal f\in {\mathcal F}.
$$

\end{proof}

\begin{cor}\label{CT2}
Let $X$ be a  compact metric. For any $\ep>0,$ any finite subset
${\mathcal F}\subset C(X),$ a unital separable simple \CA\, $A$ with
tracial rank zero and a unital monomorphism $h: C(X)\to A,$ there exists
$\dt>0,$ a finite subset ${\mathcal G}\subset C(X)$ and a finite
subset ${\mathcal P}\subset \underline{K}(C(X))$ satisfying the
following: If $u\in A$ is a unitary such that
\beq\label{c2cmt1} \|[h(a),
u]\|<\dt\rforal\, a\in {\mathcal G}\andeqn \rm{Bott}(h,u)|_{\mathcal
P}=0.
\eneq
Then, there exists a continuous rectifiable path of unitaries
$\{u_t:t\in [0,1]\}$ of $A$ such that
\beq\label{c2cmt2}
\|[h(a),u_t]\|<\ep\rforal\, a\in {\mathcal F}, u_0=u\andeqn u_1=1_A.
\eneq
Moreover,
\beq\label{c2cmt3}
{\rm{Length}}(\{u_t\})\le 2\pi+\ep.
\eneq
\end{cor}

\begin{cor}\label{CT3}
Let $X$ be a compact  metric space. For any $\ep>0,$ any finite
subset ${\mathcal F}\subset C(X),$ a unital simple separable simple
\CA\, $A$ with tracial rank zero and a unital monomorphism $h: C(X)\to
A.$

If there
 is a unitary  $u\in A$ such that
\beq\label{3cmt1}
h(a)u=uh(a)\rforal\, a\in C(X)\andeqn {\rm{Bott}}(h,u)=0.
\eneq
Then, there exists a continuous rectifiable path of unitaries
$\{u_t:t\in [0,1]\}$ of $A$ such that
\beq\label{3cmt2}
\|[h(a),u_t]\|<\ep\rforal\, a\in {\mathcal F}, u_0=u\andeqn u_1=1_A.
\eneq
Moreover,
\beq\label{3cmt3}
{\rm{Length}}(\{u_t\})\le 2\pi+\ep.
\eneq
\end{cor}

\section{The constant $\dt$ and an obstruction behind the measure distribution}

The reader undoubtedly notices that the constant $\dt$ in Theorem
\ref{dmyi}  and in \ref{dmyiT} is universal. But the constant
$\dt$ in Theorem \ref{MTC} and in Theorem \ref{MT3} is not quite
universal, it depends on the measure distribution $\mu_{\tau\circ
h}.$ The introduction of the concept of measure distribution is
not just for the convenience of the proof. In fact it is
essential.  In this section we will demonstrate that when
$\text{dim}(X)\ge 2,$ the constant can not be universal for any
such space $X$ whenever $A$ is a separable non-elementary simple
\CA\, with real rank zero and stable rank one. There is a
topological obstruction for the choice of $\dt$ hidden behind the
measure distribution.

\begin{lem}\label{ideal}
Let $A$ be a unital separable simple \CA\, with stable rank one.
Suppose that $\{d_n\}$ and $\{e_n\}$ are two sequences of non-zero
projections such that there are $k(n)$ mutually orthogonal and
mutually equivalent projections $q(n,1),q(n,2),$ $...,q(n,k(n))$
in $e_nAe_n$ each of which is equivalent to $d_n,$ $n=1,2,....$

Let $d=\{d_n\}\in l^{\infty}(A)$ and let $I_d$ be the (closed
two-sided) ideal of $l^{\infty}(A)$ generated by $d.$ Then
$e=\{e_n\}\not\in I_d$ if $lim_{k\to\infty}k(n)=\infty.$ Moreover,
$\pi(e)\not\in \pi(I_d), $ where $\pi: l^{\infty}(A)\to
q_{\infty}(A)$ is the quotient map.

\end{lem}

\begin{proof}
Let $\tau$ be a quasitrace of $A.$ Define
$t_n(a)=\tau(a)/\tau(d_n)$ for all $a\in A_{s.a}.$ Suppose that
$\pi(e)\in \pi(I_d).$ There there are $x_1,x_2,...,x_m\in
l^{\infty}(A)$ such that
\beq\label{ideal1}
a=\sum_{i=1}^mx_i^*dx_i-e\in c_0(A).
\eneq
Suppose that $a=\{a_n\},$ $x_i=\{x(i,n)\}$ and $L=\max\{\|x_i\|:
i=1,2,...,m\}.$ Then
\beq\label{ideal2}
\lim_{n\to\infty}\|a_n\|=0.
\eneq
Thus, for all sufficiently large $n,$
\beq\label{ideal3}
Lm+1\ge t_n(e_n)\ge k(n).
\eneq
This is impossible.

\end{proof}

\begin{lem}\label{kv}
Let $X$ be a connected finite CW complex with dimension at least 2
and let $Y\subset X$ be a compact subset which is homeomorphic to
a $k$-cell with $k\ge 2.$ Let $A$ be a unital \CA\, and let $h_1:
C(Y)\to A$ be a unital \hm. Define $h=h_1\circ \pi,$ where $\pi:
C(X)\to C(Y)$ is the quotient map. Then, for any sufficiently
small $\dt,$ any sufficiently large finite subset ${\mathcal G}\subset
C(X)$ and any unitary $u\in A$ with $[u]=\{0\},$
\beq\label{kv1}
{\rm{Bott}}(h,u)=\{0\},
\eneq
provided that $\|[h(a),u]\|<\dt$ for all $a\in {\mathcal G}.$

\end{lem}

\begin{proof}
Let $\{p_i\in M_{n(i)}(C(X)): i=1,2,...,k\}$ be $k$  projections
which generate $K_0(C(X)).$ Assume that $p_i$ has rank $r(i)\le
n(i),$ $i=1,2,...,k.$ Since $Y$ is contractive, there is a unitary
$w_i\in M_{n(i)}(C(Y))$ such that
\beq\label{kv2}
w_i^*\pi(p_i)w_i=1_{r(i)},\,\,\,i=1,2,...,k,
\eneq
where $1_{r(i)}$ is a diagonal matrix in $M_{n(i)}$ with $r(i)$
1's and $n(i)-r(i)$ zeros. There is $a_i\in M_{n(i)}( C(X))$ such
that $\pi(a_i)=w_i,$ $i=1,2,...,k.$

Put
$$
\overline{u}^{n(i)}={\text{diag}}(\overbrace{u,u,...,u}^{n(i)}).
$$
Then, for sufficiently small $\dt$ and sufficiently large ${\mathcal
G},$
\beq\label{kv4-}
h_1(w_i^*)(h(p_i)\overline{u}^{n(i)})h_1(w_i)&=& h(a_ip_i)\overline{u}^{n(i)}h(a_i)\\
&\approx_{1/16}& 1_{r(i)}\overline{u}^{n(i)}.
\eneq
It implies that
\beq\label{kv3}
[h(p_i)\overline{u}^{n(i)}]=\{0\}\,\,\,\,
\text{(in}\,\,\,K_1(A)\,\text{)}.
\eneq
This implies that
$$
\text{bott}_0(h, u)=0.
$$

The fact that
$$
\text{bott}_1(h,u)=0
$$
follows from the fact that $K_1(C(Y))=\{0\}.$

Since $K_i(C(Y\times S^1))$ is torsion free ($i=0,1$), we conclude
that
$$
\text{Bott}(h,u)=0.
$$

%
\end{proof}

A continuous path with length $L$ may have infinite Lipschitz constant.
 Note in the following lemma, $\{v_t: t\in [0,1]\}$ not only satisfies
 the Lipschitz condition but also its  length does not increase.

\begin{lem}\label{eqexp}
Let $L>0$ be a positive number, let $C$ be a unital \CA\, and let
$M>0.$
Then, for any $\ep>0,$ there is $\dt>0$ satisfying the following:

Let $A$ be a unital \CA, let ${\mathcal G}\subset C$ be a self-adjoint
subset with $\|f\|\le M$ for all $f\in {\mathcal G},$ $\{u_t: t\in
[0,1]\}$ be a continuous path of unitaries and let $h: C\to A$ be
a monomorphism such that
$$
u_1=1,\,\,\, \|[h(g),u_t]\|<\dt\tforal \, g\in {\mathcal G}
\andeqn \,\tforal \, t\in [0,1]
$$
and
\beq\label{eqex01}
{\rm Length}(\{u_t\})\le L.
\eneq
 Then there is a path of unitaries
$\{v_t: t\in [0,1]\}$ such that
$$
v_0=u_0,\,\,\, v_1=1,
$$
$$
\|h(g),v_t]\|<\ep\tforal g\in {\mathcal G}\andeqn t\in [0,1],
$$
$$
{\rm Length}(\{v_t\})\le L.
$$
Moreover,
$$
\|v_t-v_{t'}\|\le L|t-t'|
$$
for all $t, t'\in [0,1].$
\end{lem}

\begin{proof}
Let $0=t_0<t_1<\cdots <t_m=1$ be a partition of $[0,1]$ such that
\beq\label{eqexp02}
0<\text{Length}(\{u_t: t\in [t_{i-1}, t_i]\}\le \pi/2,
\eneq
$i=1,2,...,m.$

 Put $w_i=u_{t_{i-1}}^*u_{t_i},$ $i=1,2,...,m.$
 Then there is a rectifiable continuous path from $w_i$ to $1$
 with length no more than $\pi/2.$ It follows that
$$
sp(w_i)\subset  S=\{e^{i\pi t}: t\in [-1/2, 1/2]\}.
$$
 Then
 $$
 F(e^{i\pi t})=t
 $$
 defines a continuous function from $S$ to $ [-1/2, 1/2].$
 Then
$$
w_i=\exp(\sqrt{-1}\,\pi F(w_i)),\,\,\,i=1,2,...,m.
$$
Function $F$ does not depend on $w_i$ nor $A.$

Since ${\mathcal G}$ is self-adjoint and $\|f\|\le M$ for all $f\in
{\mathcal G},$ there is a positive number $\dt$ depending on $\ep,$
$M$ and $F$ only, such that, if
$$
\|[h(g), u_t]\|<\dt\rforal g\in {\mathcal G}\andeqn t\in [0,1],
$$
then
\beq\label{eqexp04}
\|[h(g), \exp(\sqrt{-1}\,\pi t\,F(w_i))]\|<\ep
\eneq
for all $t\in [0,1]$ (cf. 2.5.11 of \cite{Lnbk} and 2.6.10 of
\cite{Lnbk}). Let
$$
l_i=\text{Length}(\{u_t:t\in [t_{i-1}, t_i]\}),\,\,\, i=1,2,...,m.
$$
Then
$$
\|\pi F(w_i)\|=l_i,\,\,\, i=1,2,...,m.
$$
 Now, define
$$
v_0=u_0,\,\,\, v_t=u_{i-1}(exp(\sqrt{-1}\pi
{t-s_{i-1}\over{s_i-s_{i-1}}}F(w_i))) \rforal \, t\in
[s_{i-1},s_i],
$$
where $s_0=0,$ $s_i=\sum_{j=1}^il_j/L,$ $i=1,2,...,m.$ Clearly, if
$t,t'\in [s_{i-1},s_i],$ then
\beq\label{eqexp041}\nonumber
\|v_t-v_{t'}\|&=&\|\exp(\sqrt{-1}\pi
{t-s_{i-1}\over{s_i-s_{i-1}}}F(w_i)) -\exp(\sqrt{-1}\pi
{t'-s_{i-1}\over{s_i-s_{i-1}}}F(w_i)) \|\\
&\le & \|\pi F(w_i)\|{|t-t'|\over{l_i/L}}=L|t-t'|.
\eneq

One then computes that
\beq\label{eqexp05}
\|v_t-v_{t'}\|\le L |t-t'|\rforal t,t'\in [0,1]
\eneq

\end{proof}

The following follows from Zhang's approximately halving
projection lemma.

\begin{lem}\label{half}
Let $A$ be a non-elementary simple \CA\, of real rank zero and let
$p\in A$ be a non-zero projection. Then, for any integer $n\ge 1,$
there exist mutually orthogonal projections $p_1,p_2,...,p_{n+1}$
such that $p_1$ is equivalent to $p_i,$ $i=1,2,...,n,$  $p_{n+1}$
is equivalent to a projection in $p_1Ap_1$
and
$$
p=\sum_{i=1}^{n+1}p_i.
$$
\end{lem}

\begin{proof}
Fix $n.$
  There is $m\ge 1$ such that $2^m\ge 2n^2.$
Write
$$
2^m=ns+r,\,\,\, r<s,
$$
where $s\ge 1$ is a positive integer and $r$ is a nonnegative
integer. It follows that $s\ge n.$

By Theorem 1.1 of \cite{Zhhalf}, there are mutually orthogonal
projections $q_1,q_2,...,q_{2^m+1}$ such that $q_1$ is equivalent
to $q_k,$ $k=1,2,..., 2^m,$ and $q_{2^m+1}$ is equivalent to a
projection in $q_1Aq_1.$
Put
$$
p_i=\sum_{k=(i-1)s+1}^{is} q_k,,\,\,\,i=1,2,...,n\andeqn
$$
$$
p_{n+1}=q_{2^m+1}+\sum_{k=ns+1}^{2ns+r}q_k.
$$
One checks that so defined $p_1,p_2,...,p_{n+1}$ meet the
requirements.

\end{proof}

\begin{thm}\label{CEX}
Let $X$ be a connected finite CW complex with ${\rm dim}(X)\ge 2$
and let $A$ be a unital separable simple \CA\, with real rank zero
and stable rank one.

Then there exists a sequence of unital monomorphisms $h_n: C(X)\to
A$ and a sequence of unitaries $\{u_n\}\subset A$ satisfying the
following:
\beq\label{cex1}
\lim_{n\to\infty}\|[u_n, h_n(a)]\|=0\andeqn {\rm{Bott}}(h_n,
u_n)=\{0\}, n=1,2,...,
\eneq
for all $a\in C(X),$ and, for any $L>0,$ there exists $\ep_0$ and
a finite subset ${\mathcal F}\subset C(X)$ such that

\beq\label{cex2}
\sup\{ \sup\{ \|[u_n(t), h_n(a)]\|: t\in [0,1]\}: a\in {\mathcal
F}\} \ge \ep_0
\eneq
for any sequence of continuous rectifiable path $\{u_n(t): t\in
[0,1]\}$ of unitaries of $A$ for which
\beq\label{cex3}
u_n(0)=u_n, \,\,\, u_n(1)=1_A \andeqn \rm{Length}(\{u_t\}) \le L.
\eneq

\end{thm}

\begin{proof}
For each integer $n,$ there are mutually orthogonal projections
$e_n,$ $p(n,1),p(n,2),...,p(n,n^3)$ in $A$ such that $p(n,1)$ is
equivalent to every $p(n,i)$,$i=2,3,...,n^3,$ and
$$
[e_n]\le [p(n,1)]\andeqn e_n+\sum_{i=1}^{n^3}p(n,i)=1_A.
$$
 by \ref{half}. Since $\text{dim}X\ge 2,$ there is a
subset $X_0\subset X$ which is homeomorphic to a (closed) $k$-cell
with dimension of $X_0$ is at least $2.$ There is a compact subset
$Y\subset X_0$ such that $Y$ is homeomorphic to $S^1.$ Without
loss of generality, we may write $Y=S^1.$ By \ref{triv}, we define
a monomorphism $h_n^{(0)}: C([0,1])\to e_nAe_n$ and let $h_{00}:
C(X)\to C([0,1])$ be a monomorphism. Put $\psi_n=h_n^{(0)}\circ
h_{00}.$ Let $\{e_{i,j}: 1\le i,j\le n^3\}$ be a system of matrix
unit for $M_{n^3}.$ By mapping $e_{1,1}$ to $p(n,1),$ we obtain
 a unital monomorphism $j_n: M_{n^3}\to (1-e_n)A(1-e_n).$
Put
$$
E(n,i)=\sum_{j=n(i-1)+1}^{ni} p(n,j),\,\,\,i=1,2,...,n^2.
$$
$j_n$ induces a unital \hm\, $j_{n,i}: M_n\to E(n,i)AE(n,i),$
$i=1,2,...,n^2.$ We will again use $j_{n,i}$ for the extension
from $M_2(M_n)$ to $M_2(E(n,i)AE(n,i)),$ $i=1,2,...,n^2.$
 As in
\cite{EL1}, \cite{EL2} and \cite{Lo2} there are two sequences of unitaries $\{u_n'\}$ and
$\{v_n'\}$ in $M_{n}$ such that
\beq\label{cex5}
\lim_{n\to\infty}\|u_n'v_n'-v_n'u_n'\|=0 \andeqn
\text{rank}\{[e(u_n',v_n')]\}=n-1.
\eneq
Denote $p_n'=\begin{pmatrix} 1 &  0\\
                                              0 & 0
                                              \end{pmatrix}$ and
                                              $q_n'=e(u_n', v_n')$  in $M_2(M_n).$
There are partial isometries $w_n'\in M_2(M_n)$ such that
$p_n'-(w_n')^*q_n'w_n'$ is a rank one projection in $M_2(M_n).$

Define $p_n=j_{n,i}(p_n'),$ $q_n=j_{n,i}(q_n')$ and
$w_n=j_{n,i}(w_n'),$ $i=1,2,...,n^2$ and $n=1,2,....$ Define
$\phi_{n,i}': C(S^1)\to p(n,j)Ap(n,j)$ by
$\phi_{n,i}(g)=g(j_{n,i}(v_n'))$ for $g\in C(S^1),$ $n=1,2,...,$
and define $\phi_{n,i}=\phi_{n,i}'\circ \pi',$ where $\pi':
C(X)\to C(S^1)$ is the quotient map which first maps $C(X)$ onto
$C(X_0)$ and then maps $C(X_0)$ onto $C(S^1).$ Define $h_n:
C(X)\to A$ by
\beq\label{cex6}
h_n(f)=\psi_n(f)\oplus \sum_{i=1}^{n^2}\phi_{n,i}(f)
\eneq
for $f\in C(X).$ Note that $h_n$ are monomorphisms. Define
\beq\label{cex7}
u_n=e_n\oplus \sum_{i=1}^{n^2} j_{n,i}(u_n'),\,\,\, n=1,2,....
\eneq
Since each $j_{n,i}(u_n)$ is in a finite dimensional \SCA,
\beq\label{cex8}
[u_n]=\{0\} \,\,\,\text{in}\,\,\, K_1(A).
\eneq
We also have
\beq\label{cex9}
\lim_{n\to\infty}\|[h_n(f),u_n]\| =0\,\rforal\, f\in C(X).
\eneq
Since $\psi_n$ factors through $C([0,1])$ and $\phi_{n,i}$ factors
through $C(X_0),$ by  (\ref{cex8}) and  \ref{kv},
\beq\label{cex10}
\text{Bott}(h_n, u_n)=\{0\},\,\,\, n=1,2,....
\eneq

Denote $B=l^{\infty}(A)$ and $Q=q_{\infty}(A).$ Let $H: C(X)\to B$
be defined by $\{h_n\}.$ Let $I_d$ be the (closed two-sided) ideal
of $B$ generated by
$$
d=\{ E(n,1)\}.
$$
Put
\beq\nonumber
p&=&\{\text{diag}(0, j_{n,1}(p_n'), j_{n,2}(p_n'),...,
j_{n,n^2}(p_n')\},\\\nonumber
 q&=&\{\text{diag}(0, j_{n,1}(q_n'),
j_{n,2}(q_n'),...,j_{n,n^2}(q_n')\}\,\,\,\andeqn\\
w&=&\{\text{diag}(0,
j_{n,1}(w_n'),j_{n,2}(w_n'),...,j_{n,n^2}(w_n')\}.
\eneq
Put
\beq\label{cex11}
e=p-w^*qw.
\eneq
Note that $[E(n,1)]=n[p_n-w_n^*q_nw_n].$  Hence
$$
 [e]=n[E(n,1)].
$$
 It follows from \ref{ideal} that
\beq\label{cex12}
e\not\in M_2(I_d)
\eneq

Now suppose that the lemma fails. Then, by passing to a
subsequence if necessary, we obtain  a constant $L>0$ and a
sequence of continuous rectifiable paths $\{z_t^{(n)}: t\in
[0,1]\}$ of $A$ such that
\beq\label{cex13}
z_0^{(n)}=u_n,\,\,\, z_1^{(n)}=1_A \andeqn
\lim_{n\to\infty}\|[z_t^{(n)}, h_n(f)]\|=0.
\eneq
Moreover, for each $n,$
\beq\label{cex14}
\text{Length}(\{z_t^{(n)}\})\le L
\eneq
It follows from \ref{eqexp} that there exists, for each $n,$ a
path of unitaries $\{v_t(n): t\in [0,1]\}$ in $A$ such that
$$
v_0(n)=z_0^{(n)}, \lim_{n\to\infty}\|[h_n(g),v_t(n)]\|=0
$$
for all $g\in C(X)$ and
\beq\label{cex15-1}
\|v_t(n)-v_{t'}(n)\|\le L|t-t'|
\eneq
for all $t,t'\in [0,1]$ and all $n.$ In particular, $\{v_t(n)\}$
is equi-continuous on $[0,1].$

Let $\pi: B\to Q$ be the quotient map. Put $Z_t=\pi(\{v_t\})$ for
$t\in [0,1].$ By (\ref{cex15-1}), $\{Z_t:t\in [0,1]\}$ is a
continuous path of unitaries in $Q.$ Moreover
\beq\label{cex15}
Z_0=\pi(\{u_n\}),\,\,\, Z_1=1_B\andeqn  Z_t(\pi\circ
H(f))=(\pi\circ H(f))Z_t
\eneq
for all $f\in C(X)$ and $t\in [0,1].$ Let ${\bar \pi}_I: Q\to
Q/\pi(I_d).$ Then ${\bar \pi}_I\circ \pi\circ H$ induces a \hm\,
$\Psi: C(S^1)\to Q/\pi(I_d)$ since $\{e_n\}\in I_d.$ We also have
\beq\label{cex16}
{\bar \pi}_I(Z_0)={\bar \pi}(\{u_n\}),\,\,\, {\bar
\pi}_I(Z_1)=1\andeqn  {\bar \pi}_I(Z_t)\Psi(g)=\Psi(g){\bar
\pi}_I(Z_t\})
\eneq
for all $g\in C(S^1).$ However, since $e\not\in M_2(I_d),$ we
compute that
\beq\label{cex17}
\text{Bott}(\Psi, {\bar \pi}_I(\{u_n\})\not=\{0\}.
\eneq
This contradicts (\ref{cex16}).

\end{proof}

\begin{NN}\label{Lip-path}

{\rm As we know that a continuous rectifiable path with length $L$
may fail to be Lipschitz, for the future use, applying
\ref{eqexp}, Theorem \ref{MT3} may be written as the following. We
also remark that Theorem \ref{dmyiT} may also be written
similarly.

}

\end{NN}
\begin{thm}\label{MT3+Lip}
Let $X$ be a compact metric space. For any $\ep>0,$ any finite
subset ${\mathcal F}\subset C(X)$ and any (increasing) map
${\varDelta}: (0,1)\to (0,1),$ there exists $\dt>0,$ a finite
subset ${\mathcal G}\subset C(X)$  and a finite subset ${\mathcal
P}\subset \underline{K}(C(X))$ satisfying the following:

Suppose that $A$ is a unital separable simple \CA\, with tracial
rank zero, $h: C(X)\to A$ is a unital monomorphism and $u\in A$ is
a unitary
 such that $\mu_{\tau\circ h}$ is $\varDelta$-distributed for all $\tau\in T(A),$
\beq\label{Lip3mt1}
\|[h(a), u]\|<\dt\tforal\, a\in {\mathcal G} \andeqn
{\rm{Bott}}(h,u)|_{\mathcal P}=0.
\eneq
Then, there exists a continuous rectifiable path of unitaries
$\{u_t:t\in [0,1]\}$ of $A$ such that
\beq\label{Lip3mt2}
u_0=u,\,\,\,u_1=1_A\andeqn \|[h(a),u_t]\|<\ep\rforal\, a\in
{\mathcal F}\andeqn t\in [0,1].
\eneq
Moreover,
\beq\label{Lip3mt3}
&&\|u_t-u_{t'}\| \le  (2\pi+\ep)|t-t'|\tforal t, t'\in [0,1]\\
&&\hspace{-0.2in}\andeqn
{\rm{Length}}(\{u_t\}) \le  2\pi+\ep.
\eneq
\end{thm}

\vspace{0.2in}

\begin{rem}\label{Rdelta}
{\rm Let $X$ be a connected finite CW complex, let $A$ be a unital
simple \CA\, and let $h: C(X)\to A$ be a unital monomorphism.
Suppose that $u\in A$ is a unitary such that $[h(f),u]=0$ for all
$f\in C(X).$ Denote by $\phi: C(X\times S^1)\to A$ the \hm\,
induced by $h$ and $u.$ Suppose there is a rectifiable continuous
path of unitaries $\{u_t: t\in [0,1]\}$ of $A$ such that $\|[h(a),
u]\|<\ep<1/2$ on a set of generators. We know that it is necessary
that  $\text{Bott}(h, u)=0.$ Theorem \ref{1T} states that under
the assumption that $\phi$ is injective, $\text{Bott}(h,u)=0$ is
also sufficient to have such path of unitaries provided that $A$
has tracial rank zero.

As in the proof of \ref{CEX}, if there is a sequence of unitaries
$u_n\in A$  such that
$$
\lim_{n\to\infty}\|[h(a), u_n]\|=0\rforal f\in C(X),
$$
then we obtain a \hm\, $\Psi: C(X\times S^1)\to q_{\infty}(A).$
The existence of the desired paths of unitaries  implies (at
least) that
$$
\Psi_{*0}|_{\bt_1(K_1(C(X)))}=0.
$$
However, $q_{\infty}(A)$ is not simple. Suppose that $I\subset
q_{\infty}(A)$ is a closed two-sided ideal and ${\bar \pi}:
q_{\infty}(A)\to q_{\infty}(A)/I$ is the quotient map. Even if
$\Psi$ is a monomorphism, ${\rm ker} {\bar \pi}\circ \Psi$ may not
be zero. The trouble starts here. Suppose that $F$ is a compact
subset of $X$ such that
$$
{\rm ker}\Psi\supset \{ f\in C(X\times S^1): f|_{F\times S^1}=0\}.
$$
We obtain a \hm\, $\Psi_1: C(F\times S^1)\to q_{\infty}(A)/I.$ The
existence of the paths  also requires that
$$
(\Psi_1)_{*0}(\bt_1(K_1(C(F))))=0.
$$
This can not be guaranteed by $\Psi_{*0}|_{\bt_1(K_1(C(X)))}=0$
when ${\rm dim}X\ge 2$ as we see in the proof. The relevant ideals
are those  associated  with diminishing tracial states. This
topological obstruction can be revealed this way. The
requirement to have a fixed measure distribution as in \ref{MTC} guarantees that
${\rm ker}{\bar \pi}\circ \Psi=\{0\}$ for those ideals. Therefore
the condition of vanishing  $\text{Bott}$ maps is sufficient  at least
for the case that $A$ is assumed to have tracial rank zero.

When ${\rm dim}X\le 1,$ $s_{*1}(K_1(C(X)))=K_1(C(F)),$ where $s:
C(X)\to C(F)$ is the quotient map. Thus ideals with diminishing
tracial states do not produce anything surprising. Therefore, when
${\rm dim}X\le 1,$ as one sees in the section 3, no measure
distribution are considered.

In Section 11,  we will consider the case that $A$ is assumed to
be purely infinite and simple. In that case, $A$ has no tarcial
states. So there is nothing about measure distribution. In turns
out, there is no hidden obstruction as we will see in the section
after the next.

}

\end{rem}

\chapter{Purely infinite simple $C^*$-algebras}

\setcounter{section}{9}

\section{Purely infinite simple $C^*$-alegbras}

\begin{lem}{\rm (Lemma 1.5 of \cite{LP})}\label{LP}
Let $X$ be a compact metric space, let ${\mathcal F}$ be a finite
subset of the unit ball of $C(X).$ For any $\ep>0$ and any
$\sigma>0,$  there exists $\eta>0$ such that for any unital \CA\,
$A$ of real rank zero and any contractive unital $*$-preserving
linear map $\psi: C(X)\to A,$ there is a finite  subset  $\{x_1,
x_2,...,x_n\}\subset \Sigma_{\eta}(\psi, {\mathcal F})$ which is
$\sigma$-dense in $\Sigma_{\eta}(\psi, {\mathcal F})$ and mutually
orthogonal projections $p_1,p_2,...,p_n$ in $A$ satisfying the
following:
$$
\|\psi(f)-\sum_{i=1}^n f(x_i)p_i-p\psi(f)p\|<\ep
$$
for all $f\in {\mathcal F}.$
\end{lem}

\begin{proof}
The statement is slightly different from that of Lemma 1.5 of \cite{LP}. One can deduce this from
Lemma 1.5 of \cite{LP}. But it also follows from the  proof of Lemma 1.5 of \cite{LP}.

\end{proof}

\begin{lem}\label{Jsp}
Let $X$ be a compact metric space and let ${\mathcal F}\subset C(X)$ be a finite subset.
For any $\ep>0$ and any $\sigma>0,$  there exists $\dt>0$ and a finite subset ${\mathcal G}\subset C(X)$ satisfying the following:

For any unital \CA\, $B$ and any unital monomorphism
$h: C(X)\to B,$ if there exists a unitary $u\in B$ such that
$$
\|[h(g), u]\|<\dt\rforal g\in {\mathcal G},
$$
then there is a finite $\sigma$-dense subset
$\{x_1,x_2,...,x_m\}\subset X$ and $s_1,s_2,...,s_m$ (may not be
distinct) such that $(x_i, s_i)\in\Sigma_{\ep}(\psi, {\mathcal
F}\otimes S),$ $i=1,2,...,m,$ where $S=\{1_{C(S^1)}, z\}$ and where
$\psi(g\otimes f)=h(g)f(u)$ for all $g\in C(X)$ and $f\in C(S^1).$

\end{lem}

\begin{thm}\label{NEW}
Let $X$ be a compact metric space and let $\ep>0$ and let ${\mathcal
F}\subset C(X)$ be a finite subset. There is $\dt>0$
and a finite subset ${\mathcal G}\subset C(X)$ and a finite subset
${\mathcal P}\subset \underline{K}(C(X))$ satisfying the following:

Suppose that $A$ is a  unital  purely infinite simple
\CA\,  and $\psi, \phi: C(X)\to A$ are two  unital $\dt$-${\mathcal G}$-multiplicative \morp s
 for which
\beq\label{old1}
[\psi]|_{\mathcal P}=[\phi]|_{\mathcal P},
\eneq
and $\psi$ and $\phi$ are both  $1/2$-${\mathcal G}$-injective,
then  exists  a unitary $u\in A$ such that
\beq\label{old2}
\|{\rm ad}\, u\circ \psi(f)-\phi(f)\|<\ep\,\rforal\, f\in {\mathcal
F}.
\eneq
 \end{thm}

\begin{proof}
If $A$ is also assumed to be amenable, this theorem follows from
the combination of a theorem of Dadarlat (Theorem 1.7 of
\cite{D1}) and Cor 7.6 of \cite{Lnbdf}.  In general case, we
proceed as follows.

Let $\dt_1>0$ and ${\mathcal G}_1$ and ${\mathcal P}$ be required as in
Theorem 3.1 of \cite{GL2} for  $\ep/16$ and ${\mathcal F}.$ Note
(as stated before Theorem 3.1 of \cite{GL2}) that Theorem 3.1 of
\cite{GL2} applies to the case that $A$ is a unital purely
infinite simple \CA.
We may assume that ${\mathcal F}\subset {\mathcal G}_1.$

By choosing smaller $\dt_1$ and larger ${\mathcal G}_1,$ we may assume that
$$
[L_1]|_{\mathcal P}=[L_2]|_{\mathcal P}
$$
if both $L_1$ and $L_2$ are $\dt_1$-${\mathcal G}_1$-multiplicative \morp\, from $C(X)$ to $A,$ and if
$$
L_1\approx_{\dt_1} L_2 \,\,\,\text{on}\,\,\,{\mathcal G}_1.
$$

Let $\eta_1=\min\{\ep/16, \dt_1/4\}$  and $\eta_2<\sigma_{X, \eta_1, {\mathcal G}_1}.$
Therefore
$|f(x)-f(y)|<\eta_1$ for all
$f\in {\mathcal G}_1,$ whenever ${\rm dist}(x,y)<\eta_2.$
Put  $\eta_3=\min\{\eta_1, \eta_2\}.$

Let $\eta>0$ be required in \ref{LP} corresponding to $\eta_3$ (in place of $\ep$), $\eta_2/2$
(in place of $\sigma$)
and ${\mathcal G}_1$ (in place of ${\mathcal F}$). We may assume that $\eta<\eta_3/2.$

Let  $\dt>0$ and ${\mathcal G}\subset C(X)$ be a finite subset  required by \ref{inj} associated with $\eta$
(in place of $\ep$) and ${\mathcal G}_1$  (in place of ${\mathcal F}$).  We may assume that ${\mathcal G}_1\subset {\mathcal G}$ and
$\dt<\eta.$

Now suppose that $\psi,\phi: C(X)\to A$ are two unital $\dt$-${\mathcal G}$-multiplicative \morp s
which are both $1/2$-${\mathcal G}$-injective such that
\beq\label{old3-}
[\psi]|_{\mathcal P}=[\phi]|_{\mathcal P}.
\eneq

By Lemma 1.5 of \cite{LP} and \ref{inj},
\beq\label{old3}
\hspace{-0.2in}\|\phi(f)- \sum_{j=1}^N
f(\xi_j)p_j+\phi_1(f)\|<\eta_3\andeqn \|\psi(f)-\sum_{i=1}^{N'}
f(\xi_i')p_i'+\psi_1(f)\|<\eta_3
\eneq
 for all $ f\in {\mathcal G}_1,$
where $\{p_1,p_2,...,p_N\}$ and $\{p_1,p_2,...,p_{N'}\}$ are  two sets  of non-zero mutually
orthogonal projections, $P=\sum_{j=1}^N p_j,$  $P'=\sum_{i=1}^{N'}p_i',$ $\{\xi_1,\xi_2,...,
\xi_N\}$ and $\{\xi_1',\xi_2',...,\xi_{N'}'\}$ are two $\eta_2$-dense subsets  in $X$ and $\psi_1: C(X)\to (1-P)A(1-P)$
and $\psi_1': C(X)\to (1-P')A(1-P')$ are two  unital $\dt$-${\mathcal G}$-multiplicative \morp s.

By the choice of $\eta_2,$ by replacing $\eta_3$ by $\eta_1$ in (\ref{old3}),
we may assume that $N=N'$ and $\xi_i=\xi_i',$ $i=1,2,...,N.$
Since $A$ is purely infinite, we may write $p_i=e_i+e_i',$ where $e_i,$ $e_i'$ are two mutually
orthogonal projections
such that $[e_i]=[p_i'],$ $i=1,2,...,N.$ Put $\phi_1'(f)=\sum_{i=1}^Nf(\xi_i)e_i'+\phi_1(f)$ for all
$f\in C(X).$ There is a unitary $U\in A$ such that
$$
U^*p_i'U=e_i,\,\,\,i=1,2,...,N.
$$
Replacing $\phi_1$ by $\phi_1'$ and $\psi$ by ${\rm ad}\, U\circ \psi,$ in (\ref{old3}), we may assume
that $p_i=p_i'$ and $P=P'.$  In other words, we may assume that
\beq
\label{old4} \hspace{-0.2in}\|\phi(f)- \sum_{j=1}^N
f(\xi_j)p_j+\phi_1(f)\|<\eta_1\andeqn \|\psi(f)-\sum_{i=1}^{N}
f(\xi_i)p_i+\psi_1(f)\|<\eta_1
\eneq
 for all $ f\in {\mathcal G}_1, $

Define $h_{0}: C(X)\to PAP$ by $h_0(f)=\sum_{i=1}^{N} f(\xi_i)p_i$ for all $f\in C(X).$
By the choice of $\dt_1$ and ${\mathcal G}_1$ and by (\ref{old3-}),
 $\phi_1$ and $\psi_1$ are $\dt_1/2$-${\mathcal G}_1$-multiplicative
and
\beq\label{old5}
[\phi_1]|_{\mathcal P}=[\psi_1]|_{\mathcal P}.
\eneq
By Theorem 3.1 of \cite{GL2}, there is an integer $l>0,$ a
unitary\\  $w\in M_{l+1}((1-P)A(1-P))$ and a unital \hm\,\\
$h_{00}: C(X)\to  M_{l}((1-P)A(1-P))$ with finite dimensional
range such that

\beq\label{old7}
\text{diag}(\psi_1, h_{00})\approx_{\ep/16} {\text{ad}}\, w\circ
\text{diag}(\phi_1, h_{00})\,\,\,\,\text{on}\,\,\, {\mathcal F}.
\eneq

By the choice of $\eta_2$ and by replacing $\ep/16$ by $\ep/8$ in
(\ref{old7}), we may assume that
\beq\label{old8}
h_{00}(f)=\sum_{j=1}^N f(\xi_j) q_j\,\rforal\, f\in C(X)
\eneq
where $\{q_1, q_2,...,q_N\}$ are mutually orthogonal
projections in $M_l((1-P)A(1-P)).$
Since $A$ is purely infinite and simple, we may write that  $p_i=q_i'+q_i'',$ where
$q_i$ and $q_i'$ are two mutually orthogonal projections such that
$[q_i]=[q_i'],$ $1=2,...,N.$ There is a unitary $w_1\in  M_{l+1}(A)$ such that
\beq\label{old9}
w_1^*q_iw=q_i'\le p_i,\,\,\,\,\,i=1,2,...,N.
\eneq
Thus, by (\ref{old7}), we obtain a unitary $u\in A$ such that
\beq\label{old10}
{\rm ad}\, u\circ (h_0+\phi_1)\approx_{\ep/4}h_0+\psi_1\,\,\,{\rm on}\,\,\, {\mathcal F}.
\eneq
By (\ref{old4}), we have
$$
{\rm ad}\, u\circ \psi\approx_{\ep} \phi\,\,\,{\rm on}\,\,\, {\mathcal F}.
$$
\end{proof}

\begin{cor}\label{presemi}
Let $X$ be a compact metric space and let ${\mathcal F}\subset C(X)$ be a finite subset. For any $\ep>0,$ there are
$\dt>0,$ a finite subset ${\mathcal G}\subset C(X)$ and a finite subset ${\mathcal P}\subset\underline{K}(C(X))$ satisfying the following:

For any unital purely infinite simple \CA\, $A$ and any unital $\dt$-${\mathcal G}$-multiplicative
\morp\, $\psi: C(X)\to A$ which is $1/2$-${\mathcal G}$-injective such that
$$
[\psi]|_{\mathcal P}=[h]|_{\mathcal P},
$$
for some \hm\, $h: C(X)\to A,$ there is a unital monomorphism $\phi: C(X)\to A$ such that
$$
\|\psi(f)-\phi(f)\|<\ep\rforal f\in {\mathcal F}.
$$
\end{cor}

\begin{proof}
For sufficiently small $\dt$ and large ${\mathcal G},$ we may assume that $[h(1_{C(X)})]=[1_A].$
To make $h$ injective, let $e\in A$ be a projection with $[e]=[h(1_{C(X)})]$ such that $1-e\not=0.$
So we may assume that $h$ maps $C(X)$ unitally to $eAe.$
But then $[1-e]=0.$ There is a unital embedding $\phi_0: {\mathcal O}_2\to eAe.$ By applying \ref{triv},
we obtain a unital embedding $h_0: C(X)\to {\mathcal O}_2.$ Put $h_{00}=\phi_0\circ h_0.$
Then $h_1=h+h_{00}.$ Then $h_1$ is a unital monomorphism.
We then apply \ref{NEW} to obtain $\phi.$
\end{proof}

The following is a combination of a result of M. R\o rdam \cite{Ro}, a result of \cite{ER} and
the classification theorem of Kirchberg-Phillips .

\begin{lem}\label{rod}
Let $A$ be a unital separable amenable purely infinite simple \CA\, satisfying the Universal Coefficient Theorem
and $B$ be any unital purely infinite simple \CA.
Then for any $\af\in KL(A, B),$ there is a monomorphism $\psi: A\to B$ such that
$$
\af=[\psi].
$$
Moreover, if $\af([1_A])=[1_B],$ $\psi$ can be made into a unital \hm.
\end{lem}

\begin{proof}
Let ${\bar {\mathcal C}}$ be the class of those unital separable purely infinite simple \CA s for which the lemma holds.
Let $(G_0', g_0', G_1')$ be a triple, where $G_0'$ and $G_1'$ are finitely generated abelian groups
and $g_0'\in G_0'$ is an element. The proof of 5.6 of \cite{ER} shows that the class ${\mathcal C}$ contains
one unital separable amenable purely infinite simple \CA\, $A$ satisfying the UCT such that
$(K_0(A), [1_A], K_1(A))=(G_0', g_0', G_1').$
Now let $(G_0, g_0, G_1)$ be a triple for arbitrary countable abelian groups. Then it is an inductive limit of those triples with finitely generated abelian groups.
By 8.2 of \cite{Ro} (see also 5.9 of \cite{Ro} and 5.4 of \cite{ER}),
${\bar C}$ contains a unital separable amenable purely infinite simple \CA\, $A$ satisfying the UCT such that
$(K_0(A), [1_A], K_1(A))=(G_0, g_0, G_1).$  By the classification of those purely infinite simple \CA s
of Kirchberg and Phillips (see \cite{P}, for example), there is only one such unital amenable separable purely
infinite simple \CA\, satisfying the UCT.
\end{proof}

The following also follows from a result of E. Kirchberg (in a coming paper).

\begin{thm}\label{ext}
Let $B$ be a unital separable amenable \CA\, which satisfies the UCT and let $A$ be a unital purely infinite
simple \CA.  Then, for any $\kappa\in KL(B,A),$ there is a
monomorphism $h: C(X)\to A$ such that
$$
[h]=\kappa.
$$
If $A$ is unital and $\kappa([1_{B}])=[1_A],$ then $h$ can be chosen to be unital.
\end{thm}

\begin{proof}
There is a unital amenable separable purely infinite simple \CA\, $C$ satisfying the UCT such that
$$
(K_0(A), [1_A], K_1(A))=(K_0(B), [1_{B}], K_1(B)).
$$
Thus, by 7.2 of \cite{RS}, there is an invertible element $x\in KK(B, C).$
Thus, by Theorem 6.7 of \cite{Lnbdf}, there is a unital \hm\, $h_1: B\to  C$ such that
$[h_1]={\bar x},$ where ${\bar x}\in KL(B, C)$ is the image of $x.$ Since $x$ is invertible, we obtain an isomorphism
$$
KK(B, A)\cong KK(C, A).
$$
Let ${\tilde \kappa}\in KK(B,A)$ such that the image of ${\tilde \kappa}$ in $KL(B,A)$ is $\kappa.$
Let ${\tilde \kappa}'\in KK(C, A)$ be the element which is the image of ${\tilde \kappa}$ under the above isomorphism.
It follows that there is ${\tilde \af}\in KK(C,C)$ such that
\beq\label{ext-1}
[h_1]\times {\tilde \af} \times {\tilde \kappa}'={\tilde \kappa},
\eneq
where $\times$ is the Kasparov product.
Let $\af$ be the image of ${\tilde \af}$ in $KL(C,C)$ and let $\kappa'$ be the image of ${\tilde \kappa}'$ in $KL(C,A),$
respectively.
By \ref{rod}, there are monomorphisms $h_2: C\to A$ and \hm\, $h_3: C\to C$ such that
\beq\label{ext-2}
[h_2]=\kappa'\andeqn [h_3]=\af.
\eneq

Define $h=h_2\circ h_3\circ h_1.$ Then
$$
[h]=\kappa\,\,\,\text{in}\,\,\, KL(B, A).
$$

To make sure that $h'$ is a monomorphism,
consider two non-zero projections $p_1, p_2\in A$ such that
$p_1\le p_2,$ $p_2-p_1\not=0$ and
$$
[p_1]=[p_2]=\kappa([1_{B}]).
$$
By applying what we have shown to $p_1Ap_1,$ we obtain a \hm\,
$h':B\to p_1Ap_1$ such that $[h']=\kappa.$

There is a unital embedding  $h_0$ from ${\mathcal O}_2$ into
$(p_2-p_1)A(p_2-p_1)$ (note that \\$[p_2-p_1]=0$ in $K_0(A)$). It
follows from 2.3 of \cite{KP1} that there is a unital monomorphism
$h_0': B\to {\mathcal O}_2.$  Let $h_{00}=h_0\circ h_0'.$ Then $$
[h_{00}]=0\,\,\,\text{in}\,\,\, KK(B, A).
$$
Define $h=h_{00}+h'.$ Then $h$ is a monomorphism. Moreover,
if $\kappa([1_{B}])=[1_A],$ we
can choose $p_2=1_A.$

\end{proof}

\begin{thm}\label{semi}
Let $X$ be a finite CW complex and let $\ep>0$ and let ${\mathcal
F}\subset C(X)$ be a finite subset. There exists $\dt>0$
and a finite subset ${\mathcal G}\subset C(X)$ satisfying
the following.

For any  purely infinite simple \CA\, $A$ and any unital
$\dt$-${\mathcal G}$ -multiplicative \morp\, $\psi: C(X)\to A$ which
is $1/2$-${\mathcal G}$-injective, there is a unital monomorphism $h: C(X)\to A$
such that
\beq\label{semi1}
\|\psi(f)-h(f)\|<\ep\,\rforal\, f\in {\mathcal F}.
\eneq
\end{thm}

\begin{proof}
If we also assume that $A\otimes{\mathcal O}_{\infty}\cong A,$ for
example, $A$ is amenable, this follows from Cor 7.6 of
\cite{Lnbdf}. The general case follows from \ref{presemi}.
Since $K_i(C(X))$ is finitely generated ($i=0,1$),
it is easy to see that, with sufficiently small $\dt$ and
sufficiently large ${\mathcal G},$ $[\psi]$ gives an element in
$KL(C(X), A).$ This certainly follows from \ref{apK1}. The next thing is to show that, for any $\kappa\in
KL(C(X),A),$ there is a unital monomorphism $h: C(X)\to A$ such that
$[h]=\kappa.$ This follows from Theorem \ref{ext}. Then
\ref{presemi} applies.
\end{proof}

\begin{lem}\label{pi1}
Let $X$ be a compact metric space, let $\ep>0,$ let $\sigma>0,$
and let ${\mathcal F}\subset C(X)$ be a finite subset. Suppose that
$A$ is a unital purely infinite simple \CA\, and suppose that $h:
C(X)\to A$ is a unital \hm\, with the spectrum $F\subset X.$ Then
there are nonzero  mutually orthogonal projections
$e_1,e_2,...,e_n$ and a $\sigma$-dense subset
$\{x_1,x_2,...,x_n\}\subset F$ and a unital \hm\, $h_1: C(X)\to
(1-e)A(1-e)$ {\rm (}where $e=\sum_{i=1}^ne_i${\rm )}  with the spectrum $F$ such that
$$
\|h(f)-(h_1(f)+\sum_{i=1}^nf(x_i)e_i)\|<\ep\rforal f\in {\mathcal F}.
$$

Moreover, we may choose $[e_i]=0$ in $K_0(A),$ $i=1,2,..,n$ and
$[h_1]=[h]$ in $KL(C(X),A).$
\end{lem}

\begin{proof}
To simplify the notation, without loss of generality, we may assume that $X=F.$
Fix $\ep>0$ and a finite subset ${\mathcal F}\subset C(X).$ Without
loss of generality, we may assume that ${\mathcal F}$ is in the unit
ball of $C(X).$ Let $\dt>0$ and ${\mathcal G}\subset C(X)$ be a finite
subset required by \ref{presemi} associated with
$\ep/4$ and ${\mathcal F}.$ Let $\dt_1=\min\{\dt, \ep/4\}.$ Let
$\eta=\sigma_{X, \dt_1, {\mathcal G}}.$

By choosing smaller $\dt$ and larger ${\mathcal G},$ we may also assume that
$$
[L_1]|_{\mathcal P}=[L_2]|_{\mathcal P}
$$
for any two unital $\dt$-${\mathcal G}$-multiplicative \morp s from $C(X)$ to $A,$ provided that
$$
L_1\approx_{\dt} L_2\,\,\,{\rm on}\,\,\,{\mathcal G}.
$$

Let $\{x_1,x_2,...,x_n\}$ be a
$\sigma$-dense subset of $X=F$ such that
\beq\label{pi1-01}
\max\{|g(x_i)|: i=1,2,...,n\}\ge (3/4)\|g\|\rforal g\in {\mathcal G}.
\eneq


It follows from \ref{LP} that there are non-zero mutually orthogonal projections such that
\beq\label{pi1-1}
\|h(f)-[(1-p)h(f)(1-p)+\sum_{i=1}^nf(x_i)e_i]\|<\dt_1\andeqn
\eneq
\beq\label{pi1-2}
\|[1-p, h(f)]\|<\dt_1
\eneq
for all $f\in {\mathcal G}.$ Put $e=\sum_{i=1}^n e_i.$ Define $\psi:
C(X)\to (1-e)A(1-e)$ by $\psi(f)=(1-p)h(f)(1-p)$ for $f\in C(X).$
Then $\psi$ is $\dt_1$-${\mathcal G}$-multiplicative, by
(\ref{pi1-2}).
Write $e_i=e_i'+e_i'',$ where $e_i'$ and $e_i''$ are nonzero mutually orthogonal projections
such that
\beq\label{pi1-2k}
[e_i']=0 \,\,\,{\rm in}\,\,\, K_0(A),\,\,\,
i=1,2,...,n.
\eneq
This is possible since $A$ is purely infinite and simple.
By replacing $(1-p)h(f)(1-p)$ by  $(1-p)h(f)(1-p)+\sum_{i=1}^nf(x_i)e_i''$ and
replacing $e_i$ by $e_i',$ we may  assume, by (\ref{pi1-01}),  that $\psi$ is $1/2$-${\mathcal G}$-injective
and $[e_i]=0$ in $K_0(A),$ $i=1,2,...,n.$
Define $h_0: C(X)\to  eAe$ by $h_0(f)=\sum_{i=1}^n f(x_i)e_i$ for all $f\in C(X).$
Since $[e_i]=0$ in $K_0(A),$ we have that
\beq\label{pi1-K}
[h_0]=0\,\,\,\,\text{in}\,\,\, KL(C(X), A).
\eneq
Thus, by (\ref{pi1-1}), (\ref{pi1-2k}) and by (\ref{pi1-K}),
\beq\label{pi1-K2}
[\psi]|_{\mathcal P}=[h]|_{\mathcal P}.
\eneq

It follows from
\ref{presemi}
 that
there is a unital monomorphism $h_1: C(X)\to (1-e)A(1-e)$ such that
\beq\label{pi1-3}
\|h_1(f)-\psi(f)\|<\ep/4\rforal f\in {\mathcal F}.
\eneq
Combining (\ref{pi1-1}) and (\ref{pi1-3}), we obtain that
$$
\|h(f)-[h_1(f)+\sum_{i=1}^nf(x_i)e_i]\|<\ep\rforal f\in {\mathcal F}.
$$

This proves the first part of the lemma. Note that $h_1$ was
chosen to have $[h_1]=[h]$ as in the proof of \ref{presemi}. Thus
the above proof also implies the last part of the lemma.

\end{proof}

\begin{cor}\label{pi2}
Let $X$ be a compact metric space, let $\ep>0,$ let $\sigma>0,$
and let ${\mathcal F}\subset C(X)$ be a finite subset. Suppose that
$A$ is a unital purely infinite simple \CA\, and suppose that $h:
C(X)\to A$ is a unital \hm\, with the spectrum $F\subset X.$ Then, for any
finite subset $\{x_1,x_2,...,x_n\}\subset F,$ any $z_1,z_2,...,z_n\in K_0(A),$
there are nonzero  mutually orthogonal projections
$e_1,e_2,...,e_n$
with $[e_i]=z_i,$ $i=1,2,...,n,$  and a unital \hm\, $h_1: C(X)\to
(1-e)A(1-e)$ (where $e=\sum_{i=1}^ne_i$)  with the spectrum $F$ such that
$$
\|h(f)-(h_1(f)+\sum_{i=1}^nf(x_i)e_i)\|<\ep\rforal f\in {\mathcal F}.
$$
\end{cor}

\begin{proof}
The only thing that needs to be added to \ref{pi1} is that there is a non-zero
projection $e_i'\le e_i$  (in $A$ ) such that $[e_i']=z_i,$ $i=1,2,...,n.$
\end{proof}
\vspace{0.2in}

\section{Basic Homotopy Lemma in purely infinite simple $C^*$-algebras}

Now we establish the Basic Homotopy Lemma in purely infinite
simple \CA s.

\begin{lem}\label{pihf}
Let $X$ be a compact metric space  and let $A$ be a unital
purely infinite simple \CA\,and $h: C(X)\to A$ be a
monomorphism. Suppose that there is a unitary $u\in A$
 such that
\beq\label{pihf01}
h(a)u=uh(a)\tforal \, a\in C(X)\andeqn {\rm{Bott}}(h,u)=0.
\eneq
Suppose also  that $\psi: C(X\times S^1)\to A$ defined by
$\psi(f\otimes g)=h(f)g(u)$ for $f\in C(X)$ and $g\in C(S^1)$ is a
monomorphism. Then, for any $\ep>0$ and any finite subset ${\mathcal
F}\subset C(X),$
 there exists a rectifiable continuous path of unitaries
$\{u_t: t\in [0,1]\}$ of $A$ such that
\beq\label{pihf02}
u_0=u,\,\,\,u_1=1_A\andeqn \|[h(a), u_t]\|<\ep
\eneq
for all $a\in {\mathcal F}$ and all $t\in [0,1].$ Moreover,
\beq\label{pihf03}
{\rm{Length}}(\{u_t\})\le \pi+\ep\pi.
\eneq
\end{lem}

\begin{proof}

Let $\ep>0$ and ${\mathcal F}\subset C(X)$ be a finite subset. Without
loss of generality, we may assume that $1_{C(X)}\in {\mathcal F}.$ Let
${\mathcal F}_1=\{a\otimes b: a\in {\mathcal F}, b=u,b=1\}.$ Let $\eta=\sigma_{X\times S^1, \ep/16, {\mathcal F}}.$
Let  $\dt>0,$ ${\mathcal G}\subset C(X\times S^1)$ be a
finite subset of $C(X\times S^1)$ and ${\mathcal P}\subset
\underline{K}(C(X\times S^1))$ be a finite subset required by
\ref{NEW} associated with $\ep/16$ and ${\mathcal F}_1.$

We may assume that $\dt$ and ${\mathcal G}$ are so chosen that,
 for
any two $\dt_1$-${\mathcal G}$-multiplicative \morp\, $L_1, L_2:
C(X\times S^1)\to B$ (any unital  \CA\, $B$),
$$
[L_1]|_{{\mathcal P}}=[L_2]|_{\mathcal P},
$$
provided that
$$
L_1\approx_{\dt}L_2\,\,\,\text{on}\,\,\,{\mathcal G}.
$$

With even smaller $\dt$ we may also assume that ${\mathcal G}={\mathcal
F}_2\otimes {\mathcal G}_1,$ where ${\mathcal F}_2\subset C(X)$ and ${\mathcal
G}_1\subset C(S^1)$ are finite subsets and $1_{C(X)}\subset {\mathcal
F}_2$ and $1_{C(S^1)}\subset {\mathcal G}_1.$

 Choose $x\in X$ and let
$\xi=x\times 1.$

It follows from \ref{pi2} that there is a non-zero projection
$e\in A$ with $[e]=0$ in $K_0(A)$ and a unital monomorphism
$\psi': C(X\times S^1)\to (1-e)A(1-e)$ such that
\beq\label{d1ff-1}
\|\psi(f)- (\psi'(f)+f(\xi)e)\|<\dt_1/2\rforal f\in {\mathcal G}.
\eneq
Write $e=e_1+e_2,$ where $e_1$ and $e_2$ are two non-zero mutually orthogonal projections
with $[e_1]=[1_A]$ and $[e_2]=-[1_A]$ in $K_0(A).$

It follows from \ref{ext} that there are unital monomorphisms
$h_1: C(X)\to e_1Ae_1$ and $h_2: C(X)\to e_2Ae_2$ such that
\beq\label{d1ff-2}
[h_1]=[h]\andeqn [h_2]=-[h]\,\,\,\,\text{in}\,\,\, KL(C(X), A).
\eneq

Define $\psi_1: C(X\times S^1)\to e_1Ae_1$ by $\psi_1(f\otimes g)=h_1(f)g(1)e_1$ for all $f\in C(X)$ and
$g\in C(S^1)$ and defined
$\psi_2: C(X\times S^1)\to e_2Ae_2$ by $\psi_2(f\otimes g)=h_2(f)g(1)e_2$ for all $f\in C(X).$
Define $\Psi_0: C(X\times S^1)\to (1-e_1)A(1-e_1)$ by
$\Psi_0=\psi_2+\psi'$ and define $\Psi: C(X\times S^1)\to A$ by $\Psi=\psi_1+\Psi_0.$
Note that
\beq\label{d1ff-3}
\|h(f)-(f(\xi)(e_1+e_2)+\psi'(f\otimes 1))\|<\dt_1/2\rforal f\in
C(X).
\eneq
Define $\phi_0: C(X\times S^1)\to (1-e_1)A(1-e_1)$ by $\phi_0(f)=\sum_{j=1}^mf(x_i)q_i$ for all $f\in C(X\times S^1),$
where $\{x_1,x_2,...,x_m\}$ are in $X$ and $\{q_1, q_2,...,q_m\}$ is a set of nonzero mutually orthogonal
projections such that $[q_i]=0$ in $K_0(A)$ and $\sum_{i=1}^mq_i=(1-e_1).$
This is possible since $[1-e_1]=0$ in $K_0(A).$
Then
\beq\label{d1ff-4}
[\phi_0]=0\,\,\,\text{in}\,\,\, KL(C(X), A).
\eneq
Moreover, by choosing more points and large $m,$ we may assume
that $\phi_0$ is $1/2$-${\mathcal G}$-injective.

Note that $\text{Bott}(h,u)=0.$ Thus , by \ref{Dbot2} and by
(\ref{d1ff-1}), (\ref{d1ff-2}) and (\ref{d1ff-4}),
\beq\label{d1ff-5}
[\Psi_0]|_{P}=[\phi_0]|_{\mathcal P}.
\eneq

It follows from \ref{NEW} that there exists a unitary $w_1\in
(1-e_1)A(1-e_1)$ such that
\beq\label{n1t-m3}
{\rm ad}\, w_1\circ \phi_0\approx_{\ep/16}
\Psi_0\,\,\,\text{on}\,\,\, {\mathcal F}_1.
\eneq

By \ref{Dbot2} and by  (\ref{d1ff-1}), (\ref{d1ff-2}) and (\ref{d1ff-4})
\beq\label{n1t11}
[\Psi]|_{\mathcal P}=[\psi]|_{\mathcal P}.
\eneq
It follows from \ref{NEW} that there exists a uniatry $w_2\in A$
such that
\beq\label{n1t12}
{\rm ad}\, w_2\circ \Psi\approx_{\ep/16}
\psi\,\,\,\,\text{on}\,\,\, {\mathcal F}_1.
\eneq
Let $v=\phi_0(1\otimes z).$ In the finite dimensional commutative
\SCA\, \\$\phi_0(C(X\times S^1), $ it is easy to find a continuous
rectifiable path of unitaries \\$\{v_t: t\in [0,1]\}$ in
$(1-e_1)A(1-e_1)$ such that
\beq\label{n1t13}
v_0=v,\,\,\, v_1=1-e_1\andeqn \|[\phi_0(a\otimes 1),v_t]\|=0
\eneq
for all $t\in [0,1]$ and
\beq\label{n1t13+}
\text{Length}(\{v_t\})\le \pi.
\eneq
Define
\beq\label{n1t14}
U_t=w_2^*(e_1+w_1^*(v_t)w_1)w_2\rforal \, t\in [0,1].
\eneq
Clearly $U_1=1_A.$ Then, by (\ref{n1t-m3}) and (\ref{n1t12}),
\beq\nonumber
\|U_0-u\|&=&\|{\rm ad}\, w_2\circ  (\psi_1(1\otimes z)\oplus {\rm
ad}\, w_1\circ \phi_0(1\otimes z))-u\|\\\nonumber &<&\ep/16+\|{\rm
ad}\, w_2\circ (\psi_1(1\otimes z)\oplus \Psi_0(1\otimes
z))-u\|\\\label{n1t15-} &<&\ep/16+\ep/16+\|\psi(1\otimes
z)-u\|=\ep/8.
\eneq
We also have, by (\ref{n1t-m3}), (\ref{n1t12}) and (\ref{n1t13}),
\beq\label{n1t15}
\|[h(a),U_t]\|=2(\ep/16)+\|[{\rm ad}\, w_2\circ \Psi(a\otimes 1),
U_t]\|<\ep/8+\ep/8=\ep/4
\eneq
for all $a\in {\mathcal F}.$ Moreover,
\beq\label{n1t16}
\text{Length}(\{U_t\})\le \pi.
\eneq
By (\ref{n1t16}), (\ref{n1t15-}) and  (\ref{n1t15}), we obtain a
continuous rectifiable path of unitaries $\{u_t:t\in [0,1]\}$ of
$A$ such that
\beq\label{n1t17}
u=u_0,\,\,\, u_1=1_A\andeqn \|[h(a), u_t]\|<\ep\,\rforal \, t\in
[0,1].
\eneq
Furthermore,
\beq\label{n1t18}
\text{Length}(\{u_t\})\le \pi+\ep\pi.
\eneq

\end{proof}

\begin{lem}\label{pihf2}
Let $X$ be a compact metric space and let ${\mathcal F}\subset C(X)$ be a finite subset.
For any $\ep>0,$ there exists $\dt>0$ and a finite subset ${\mathcal G}\subset C(X)$
satisfying
the following:
  Suppose that $A$ is a unital purely infinite simple \CA,
 that $h: C(X)\to A$ is  a  unital monomorphism and that $u\in
A$ is a unitary  such that
\beq\label{pd1hf1}
\|[h(a), u]\|<\dt \tforal\, a\in {\mathcal G}.
\eneq
Then there exists a
unital monomorphism $H: C(X\times S^1)\to A$ and a rectifiable
continuous path of unitaries $\{u_t: t\in [0,1]\}$ such that
\beq\label{pd1hf2}
u_0=u,\,\,\, u_1=H(1\otimes z),
\eneq
\beq\label{pdihf3}
\|H(a\otimes 1)-h(a)\|<\ep \andeqn \|[h(a),
u_t]\|<\ep\,\,\,\tforal t\in [0,1]
\eneq
and for all $a\in {\mathcal F}.$  Moreover
\beq\label{pd1hf4}
{\rm{Length}}(\{u_t\})\le \pi+\ep.
\eneq
Furthermore,  there is a finite subset ${\mathcal P}\subset \underline{K}(C(X)),$
if in addition to {\rm (\ref{pd1hf1})},
$$
{\rm{Bott}}(h,u)|_{\mathcal P}=0,
$$
we may also require that
\beq\label{pdahf4+}
 [H]|_{{\boldsymbol{ \bt}}(\underline{K}(C(X)))}=0.
\eneq

\end{lem}
\begin{proof}
Fix $\ep>0$ and a finite subset ${\mathcal F}\subset C(X).$  To
simplify the notation, without loss of generality, we may assume
that ${\mathcal F}$ is in the unit ball of $C(X).$ Let $\dt_1>0,$
${\mathcal G}_1\subset C(X\times S^1)$ and ${\mathcal P}_1\subset \underline{K}(C(X\times S^1))$ be a finite subset which are
required in \ref{presemi} associated with $\ep/32$ and ${\mathcal F}\otimes S,$
where $S=\{1_{C(S^1)}, z\}.$
Without loss of generality, we may assume that ${\mathcal G}_1={\mathcal
F}_1\otimes S,$ where ${\mathcal F}_1\subset C(X)$ is a finite subset.

We may assume that $\dt_1<\min\{1/2,\ep/32\}$ and ${\mathcal F}\subset
{\mathcal F}_1.$
Put
$${\mathcal P}_2=\{({\rm id}-\hat{{\boldsymbol{ \bt}}})(x), \hat{{\boldsymbol{ \bt}}}(x), x: x\in {\mathcal P}_1\}.$$
By choosing possibly even smaller $\dt_1$ and larger ${\mathcal F}_1,$ we may assume that
$$
[L_1]|_{{\mathcal P}_2}=[L_2]|_{{\mathcal P}_2}
$$
for any pair of $\dt_1$-${\mathcal F}_1\otimes S$-multiplicative \morp s from $C(X\otimes S^1)$ to
any unital \CA\, $B,$ provided that
$$
L_1\approx_{\dt_1} L_2\,\,\,\text{on}\,\,\, {\mathcal F}_1\otimes S.
$$

There is $\sigma>0$ such that for any $\sigma$-dense subset
 $\{x_1,x_2,...,x_n\}\subset X$  and a finite subset $\{t_1,t_2,...,t_m\}$ such that
\beq\label{pp1}
&&\max\{\|f(x_i)\|: i=1,2,...,n\}\ge (3/4)\|f\|
\eneq
for all $ f\in {\mathcal F}_1$ and
\beq\label{pp1+}
&&\max\{\|g(x_i\times t_j)\|: 1\le i\le n\andeqn 1\le j\le m\} \ge (3/4)\|g\|
\eneq
for all $ g\in {\mathcal G}_1.$

Let $\dt_2>0$ and ${\mathcal G}_2\subset C(X\times S^1)$ be a finite
subset required in \ref{inj} for $\dt_1/4$ (in place of $\ep$),
$\sigma/2$ (in place of $\sigma$) and ${\mathcal G}_1$ (in place of
${\mathcal F}$).  Without loss of generality, we may assume that
$\dt_2<\dt_1/4$ and ${\mathcal G}_2={\mathcal F}_2\otimes S$ with ${\mathcal
F}_2\subset C(X)$ being a finite subset containing ${\mathcal F}_1.$

Let $\eta>0$ be required in \ref{LP} for $\dt_2/2$ (in place of $\ep$),
$\sigma/2$ and ${\mathcal F}_2.$ We may assume that $\eta<\dt_2/2.$

Let $\dt>0$ and ${\mathcal G}\subset C(X)$ (in place of ${\mathcal F}_1$)
be a finite subset required by \ref{appn} for $\dt_2/2$ and ${\mathcal
F}_2.$ We may assume ${\mathcal F}_2\subset {\mathcal G}.$ We may further
assume that $\dt$ is smaller than that in \ref{Jsp} and ${\mathcal G}$
is larger than that in \ref{Jsp} for $\eta$ (in place of $\ep$)
and $\sigma/2$ (in place of $\sigma$) and ${\mathcal F}_2$ (in place
of ${\mathcal F}$).

Let ${\mathcal P}\subset \underline{K}(C(X))$ be a finite subset such
that ${\boldsymbol{ \bt}}({\mathcal P})\supset \hat{{\boldsymbol{ \bt}}}({\mathcal P}_1).$ Suppose that
$h$ and $u$ satisfy the conditions in the theorem for above $\dt,$
${\mathcal G}$ and ${\mathcal P}.$ It follows from \ref{appn} that there
is a $\dt_2/2$-${\mathcal G}_2$-multiplicative \morp\, $\psi:
C(X\times S^1)\to A$ such that
\beq\label{pp1-1}
\|\psi(f\otimes g)-h(f)g(u)\|<\dt_2/2 \rforal f\in {\mathcal F}_2\andeqn g\in S.
\eneq

By applying \ref{LP} as well as \ref{Jsp}, we have
\beq\label{pp2}
\|\psi(f\otimes g)-(\Psi_1(f)+\sum_{i=1}^n f(x_i)g(s_i)e_i)\|<\dt_2/2
\eneq
for all $f\in {\mathcal F}_2$ and $g\in S,$ where $\Psi_1: C(X\times
S^1)\to (1-e)A(1-e)$ is a unital $\dt_2$-${\mathcal
G}_2$-multiplicative \morp, $e_1,e_2,...,e_n$ are non-zero
mutually orthogonal projections, $e=\sum_{i=1}^n e_i,$
$[e_i]=0$
in $K_0(A)$ ($i=1,2,...,n$),
$\{x_1,x_2,...,x_n\}$ is $\sigma$-dense, and $s_1, s_2,...,s_n$
are points in $S^1$ (not necessarily distinct).
Define $\Psi_1': C(X\otimes S^1)\to eAe$ by $\Psi_1'(f)=\sum_{i=1}^nf(x_i\times s_i)e_i$
for all $f\in C(X\times S^1).$ It follows that
\beq\label{pp2KK}
[\Psi_1']=0\,\,\,\text{in}\,\,\, KL(C(X), A).
\eneq

There are $m$
non-zero mutually orthogonal projections $e_{i,1},
e_{i,2},...,e_{i,m}$ in $e_iAe_i$
such that
$[e_{i,j}]=0$ in $K_0(A),$ $j=1,2,...,m,$
for each $i.$ Define
$\Psi_0: C(X\times S^1)\to eAe$ by $\Psi_0(f\otimes
g)=\sum_{i=1}^n (\sum_{j=1}^mf(x_i)g(t_j)e_{i,j})$ for all $f\in
C(X)$ and $g\in C(S^1).$

There is a continuous path of unitaries $\{w_t: t\in [0,1]\}$ in the
finite dimensional commutative
\SCA\, $\Psi_0(C(X\times S^1))$ such that
\beq\label{pp3}
{\rm{Length}}(\{w_t\})\le \pi,\,\,\, w_0=\sum_{i=1}^ns_ie_i,\,\,\,
w_1=\Psi_0(1\otimes z).
\eneq
Moreover,
\beq\label{pp4}
w_t\Psi_0(f\otimes 1)=\Psi_0(f\otimes 1)w_t\rforal t\in [0,1]
\eneq
and for all $f\in C(X).$
Note also
\beq\label{pp4KK}
[\Psi_0]=0\,\,\,\text{in}\,\,\,KL(C(X\times S^1), A).
\eneq
By (\ref{pp1}), $\Psi_0$ is $(3/4)$-${\mathcal G}_1$-injective. So
$\Psi_1+\Psi_0$ is $\dt_2$-${\mathcal G}_1$-multiplicative and
$(3/4)$-${\mathcal G}_1$-injective. Let $\Phi_0: C(X\times S^1)\to A$
by $\Phi_0(f\otimes g)=h(f)g(1)$ for all $f\in C(X)$ and $g\in
C(S^1).$ The condition that $\text{Bott}(h,u)|_{{\mathcal P}}=\{0\}$
(see \ref{1LK}) together with the choice of $\dt_1$ and ${\mathcal
F}_1$ implies that
\beq\label{ppKK}
[\Phi_0]|_{{\mathcal P}_2}=[\psi]|_{{\mathcal P}_2}.
\eneq
By (\ref{pp2}), (\ref{pp2KK}), (\ref{pp4KK}) and (\ref{ppKK}) and the choice of $\dt_1$ and ${\mathcal F}_1,$ we have
\beq\label{pp5LKK}
[\Phi_0]|_{{\mathcal P}_2}=[\Psi_1+\Psi_0]|_{{\mathcal P}_2}.
\eneq

By applying \ref{presemi}, we obtain
a unital monomorphism $H: C(X\times S^1)\to A$ such that
\beq\label{pp5}
\Psi_1+\Psi_0\approx_{\ep/32} H\,\,\,\text{on}\,\,\,{\mathcal F}.
\eneq
Furthermore, from the proof of \ref{presemi}, if
$\text{Bott}(h,u)|_{{\mathcal P}}=0$ is assumed, we can also assume
that
\beq\label{pp5+}
[H]|_{{\boldsymbol{ \bt}} (\underline{K}(C(X)))}=0.
\eneq
There is a unitary $w_0\in (1-e)A(1-e)$ such that
\beq\label{pp5++}
\|\Psi_1(1\otimes z)-w_0\|<2\dt_2<\dt_1<\ep/32.
\eneq

Define  $U_t=w_0+w_t$ for $ t\in [0,1].$ Then $\{U_t: t\in [0,1]\}$ is a
continuous path of unitaries in $A$ such that
\beq\label{pp6}
U_0=w_0+\sum_{i=1}^n s_ie_i,\,\,\, U_1=H(1\otimes
z)\andeqn\\\label{pp6+}
 \|[U_t,H(f\otimes
1)]\|<2\dt_2+\ep/32<\ep/16 \rforal t\in [0,1]
\eneq
and for all $ f\in C(X).$

Moreover,
\beq\label{pp7}
\text{Length}(\{U_t\})\le \pi.
\eneq

Note also we have, by (\ref{pp2}) and (\ref{pp5}), and by (\ref{pp5++})
\beq\label{pp8}
H(f\otimes 1)\approx_{\ep/16} h(f)\,\,\,\text{on}\,\,\,{\mathcal
F}\andeqn
\eneq
 \beq\label{pp9}
 \|u-U_0\|<\ep/4.
 \eneq
 We also have, by (\ref{pp8}) and (\ref{pp6+}),
 \beq\label{pp10}
 \|[h(f), U_t]\|<\ep/8
 \eneq
 The lemma follows if we connect $u$ with $U_0$ properly.

\end{proof}

\begin{thm}\label{Tpi}
Let $X$ be a compact metric space. For any $\ep>0$ and any
finite subset ${\mathcal F}\subset C(X),$ there exists $\dt>0,$  a finite subset
${\mathcal G}\subset C(X)$ and a finite subset ${\mathcal P}\subset \underline{K}(C(X))$ satisfying the following:

Suppose that $A$ is a unital  purely infinite simple \CA\, and
suppose that $h: C(X)\to A$ is a unital monomorphism. Suppose that
there is a unitary $u\in A$ such that
\beq\label{tpi1}
\|[h(a), u]\|<\dt\,\tforal\, a\in {\mathcal G}\andeqn
{\rm{Bott}}(h,u)|_{\mathcal P}=\{0\}.
\eneq
Then there is a continuous rectifiable path of unitaries $\{u_t:
t\in [0,1]\}$ of $A$ such that
\beq\label{tpi2}
u_0=u,\,\,\,u_1=1_A\andeqn \|[h(a), u_t]\|<\ep\,\rforal \, a\in
{\mathcal F}\andeqn t\in [0,1].
\eneq
Moreover,
\beq\label{tpi3}
&&\|u_t-u_{t'}\|\le (2\pi+\ep)|t-t'|\tforal t,t'\in [0,1]\andeqn\\
&&{\rm{Length}}(\{u_t\})\le 2\pi+\ep.
\eneq

\end{thm}

\begin{proof}

Fix $\ep>0$ and finite subset ${\mathcal F}\subset C(X).$ Let $\dt>0,$
${\mathcal G}\subset C(X)$ be a finite subset and ${\mathcal P}\subset
\underline{K}(C(X))$ be a finite subset required in \ref{pihf2}
for $\ep/4$ and ${\mathcal F}.$

Suppose that $h$ and $u$ satisfy the conditions in the theorem. By
\ref{pihf2}, there is a unital monomorphism $H: C(X\times S^1)\to
A$ and there is a continuous rectifiable path of unitaries $\{v_t:
t\in [0,1]\}$ in $ A$ satisfying the following:
\beq\label{Tpi1-1}
v_0=u,\,\,\,v_1=H(1\otimes z),
\eneq
\beq\label{Tpi1-2}
\|H(f\otimes 1)-h(f)\|<\ep/4,\,\,\, \|[h(f),v_t]\|<\ep/4
\eneq
for all $t\in [0,1]$ and for all $ f\in {\mathcal F},$
\beq\label{Tpi1-3}
\text{Length}(\{v_t\})\le \pi+\ep/4\andeqn
[H]|_{{\boldsymbol{ \bt}}(\underline{K}(C(X)))}=\{0\}.
\eneq
Put $h_1(f)=H(f\otimes 1)$ for all $f\in C(X).$ By (\ref{Tpi1-3}),
\beq\label{Tpi1-4}
\text{Bott}(h_1, v_1)=0.
\eneq
It follows from \ref{pihf} that there is a continuous rectifiable
path $\{w_t: t\in [0,1]\}$ in $A$ such that
\beq\label{Tpi1-5}
w_0=v_1,\,\,\,w_1=1_A, \|[h_1(f), w_t]\|<\ep/4
\eneq
for all $ t\in [0,1]$ and $f\in {\mathcal F},$ and
\beq\label{Tpi1-6}
\text{Length}(\{w_t\})\le \pi+\ep\pi/4.
\eneq

 Now define $u_t=v_{2t}$ for $t\in [0,1/2]$ and $u_t=w_{2t-1}$
for $t\in [1/2,1].$ Then,
$$
u_0=v_0=u,\,\,\, u_1=w_1=1_A.
$$
By (\ref{Tpi1-2}) and  (\ref{Tpi1-5}),
$$
\|[h(f), u_t]\|<\ep\rforal t\in [0,1]\andeqn \rforal f\in {\mathcal
F}.
$$
Moreover, by (\ref{Tpi1-3}) and (\ref{Tpi1-6})
$$
\text{length}(\{u_t\})\le 2\pi+\ep.
$$

To get (\ref{tpi3}), we apply \ref{eqexp} and what we have proved
above. Note that we may have to choose a different $\dt.$
\end{proof}

\begin{cor}\label{Cpi}
Let $X$ be a finite CW complex and let ${\mathcal F}\subset C(X)$ be a
finite subset. For any $\ep>0,$  there exists $\dt>0,$ a finite
subset ${\mathcal G}\subset C(X)$ satisfying the following:

Suppose that $A$ is a unital  purely infinite simple \CA\, and
suppose that $h: C(X)\to A$ is a unital monomorphism. Suppose that
there is a unitary $u\in A$ such that
\beq\label{ctpi1}
\|[h(a), u]\|<\dt\,\rforal\, a\in {\mathcal G}\andeqn
\rm{Bott}(h,u)=\{0\}.
\eneq
Then there is a continuous rectifiable path of unitaries $\{u_t:
t\in [0,1]\}$ of $A$ such that
\beq\label{ctpi2}
u_0=u,\,\,\,u_1=1_A\andeqn \|[h(a), u_t]\|<\ep\,\rforal \, a\in
{\mathcal F}\andeqn t\in [0,1].
\eneq
Moreover,
\beq\label{ctpi3}
\rm{Length}(\{u_t\})\le 2\pi+\ep\pi.
\eneq

\end{cor}

\begin{cor}\label{Cpi2}
Let $X$ be a compact metric space and let $A$ be  a unital purely
infinite simple \CA.
 Suppose that $h: C(X)\to A$ is a unital
monomorphism and that  there is a unitary $u\in A$ such that
\beq\label{cctpi1}
uh(f)=h(f)u \tforal f\in C(X)\andeqn {\rm{Bott}}(h,u)=0.
\eneq
Then there is a continuous rectifiable path of unitaries $\{u_t:
t\in [0,1]\}$ of $A$ such that
\beq\label{cctpi2}
u_0=u,\,\,\,u_1=1_A\andeqn \|[h(a), u_t]\|<\ep\,\tforal \, a\in
{\mathcal F}\andeqn t\in [0,1].
\eneq
Moreover,
\beq\label{cctpi3}
{\rm{Length}}(\{u_t\})\le 2\pi+\ep\pi.
\eneq

\end{cor}

\begin{rem}\label{reminj}

{\rm As we have seen in the case that $A$ is a unital simple \CA\,
of real rank zero and stable rank one, one essential difference of
the original Basic Homotopy Lemma from the the general form is
that the constant $\dt$ could no longer be universal whenever the
dimension of $X$ is at least two. In that case, a measure
distribution becomes part of the statement. This additional factor
disappears when $A$ is purely infinite. (There is no measure/trace
in this case!) However, there is a second difference. In both
purely infinite and finite cases, when the dimension of $X$ is at
least two, we need to assume that $h$ is a monomorphism, and when
the dimension of $X$ is no more than one, $h$ is only assumed to
be a \hm\,  for the case that $K_1(A)=\{0\}.$

It is much easier to reveal the  topological obstruction in this
case. Just consider a monomorphism $\phi: C(S^1\times S^1)\to A$
with $\text{bott}_1(\phi)\not=0$ and $\text{bott}_0(\phi)=0.$ In
other words,
$$
\text{bott}_1(u,v)\not=0\andeqn [v]=0,
$$
where $u=\phi(z\otimes 1)$ and $v=\phi(1\otimes z).$ This is
possible if $A$ is a unital purely infinite simple \CA\, or $A$ is
a unital separable simple \CA\, with tracial rank zero so that
${\rm ker}\rho_A\not=\{0\}.$ Let $D$ be the closed unit disk and
view $S^1$ as a compact subset of $D.$ Define $h: C(D)\to A$ by
$h(f)=\phi(f|_{S^1}(z\otimes 1))$ for all $f\in C(D).$ Since $D$
is contractive,
$$
\text{Bott}(h, v)=0.
$$
But one can not find any continuous path of unitaries $\{v_t: t\in [0,1]\}$ of $A$ such that
$$
v_0=v,\,\,\, v_1=1_A\andeqn \|[h(f),v_t]\|<\ep
$$
for all $t\in [0,1]$ and for some small $\ep,$ where $f(z)=z$ for
all $z\in D,$  since $\text{bott}_1(u,v)\not=0$ and $h(f)=u.$
Counterexample for the case that $A$ is a unital separable simple
\CA\, with tracial rank zero for which ${\rm ker}\rho_A=\{0\}$ can
also be made but not for commuting pairs. One can  begin with a
sequence of unitaries $\{u_n\}$ and $\{v_n\}$ with $v_n\in U_0(A)$
such that $\lim_{n\to\infty}\|[u_n, v_n]\|=0$ and
$\text{bott}_1(u_n, v_n)\not=0.$ Then define $h_n: C(D)\to A$ by
$h_n(f)=f(u_n)$ for $f\in C(D).$ Again $\text{Bott}(h_n, v_n)=0$
and $\lim_{n\to\infty}\|[h_n(f), v_n]\|=0$ for all $f\in C(D).$
This will give a counterexample. In summary, the condition that
$h$ is monomorphism can not be replaced by the condition that $h$
is a \hm\, in the Basic Homotopy Lemma whenever $\text{dim} X\ge
2.$

}
\end{rem}


\chapter{Approximate homotopy}

\setcounter{section}{11}

\section{Homotopy length}
In the next two sections, we will study when two maps from $C(X)$
to $A,$ a unital purely infinite simple \CA, or a separable simple
\CA\, of tracial rank zero, are approximately homotopic. When two
maps are approximately homotopic, we will give an estimate for the
length of the homotopy. In this section we will discuss the notion
of the  length of a homotopy.

\begin{lem}\label{uo}
Let $A\in {\bf B}$ be a unital simple \CA\, and let $X$ be a
compact metric space. Let $\ep>0$ and ${\mathcal F}\subset C(X)$ be a
finite subset. Suppose that $h_0, h_1: C(X)\to A$ are two unital
\hm s such that there is unitary $u\in A$ such that
$$
 {\rm ad}\, u\circ
 h_1(f)\approx_{\ep/2}h_0(f)\,\,\,\text{on}\,\,\,
{\mathcal F}.
$$
Then, there is a unitary $v\in U_0(A)$ such that
$$
{\rm ad}\, v\circ h_1\approx_{\ep} h_0\,\,\,\,\text{on}\,\,\,
{\mathcal F}.
$$
\end{lem}

\begin{proof}
 Fix $\ep$ and ${\mathcal F}\subset C(X).$
Let $\eta=\sigma_{X, \ep/2,{\mathcal F}}.$
 For a point $\xi\in X,$ let $g\in C(X)_+$ with support lies in
$O(\xi, \eta)=\{x\in X: {\rm dist}(x, \xi)<\eta\}.$ We may choose
such $\xi$ in the spectrum of $h_0$ so that $h_0(g)\not=0.$ Then
there is a non-zero projection $e\in \overline{h_0(g)Ah_0(g)}.$

Since $A$ is assumed to be purely infinite or has real rank zero
and stable rank one, there is a unitary $w\in eAe$ so that
$[(1-e)+w]=[u^*]$ in $K_1(A).$ Moreover, $u((1-e)+w)\in U_0(A).$
Put $v=u((1-e)+w).$ Note that, by the choice of $\eta,$
$$
\|((1-e)+w)^*h_0(f)-h_0(f)((1-e)+w)\|<\ep/2
$$
for all $f\in {\mathcal F}.$ Thus
$$
{\rm ad}\, v\circ h_1\approx_{\ep} h_0\,\,\,{\rm on}\,\,\,{\mathcal
F}.
$$

\end{proof}

\begin{df}\label{hp1def}
{\rm Let $A$ and $B$ be unital \CA s. Let $H: A\to C([0,1],B)$ be
a \morp. One has the following notion of the length of the
homotopy. For any partition $P\,: 0=t_0<t_1<\cdots<t_m=1$ and
$f\in A,$ put
$$
L_H(P,f)=\sum_{i=1}^m \|\pi_{t_i}\circ H(f)-\pi_{t_{i-1}}\circ
H(f)\|\andeqn
$$
$$
L_H(f)=\sup_{P}L_H(P, f).
$$
 Define
$$
\text{Length}(\{\pi_t\circ H\})=\sup_{f\in A, \|f\|\le 1} L_H(f).
$$

}
\end{df}

\vspace{0.1in}

\begin{thm}\label{Thp}
Let $A$ be a separable unital simple \CA\, of tracial rank zero
 and let $X$ be a
compact metric space. Suppose that $h_0, h_1: C(X)\to A$ are two
unital monomorphisms such that
\beq\label{thp1}
[h_0]=[h_1]\,\,\,\text{in}\,\,\,KL(C(X),A)\andeqn
\eneq
\beq\label{thp2}
\tau\circ h_1=\tau\circ h_2\,\rforal\, \tau\in T(A)
\eneq
Then, for any $\ep>0$ and any finite subset ${\mathcal F}\subset
C(X),$ there is a unital  \hm\, $\Phi: C(X)\to C([0,1], A)$ such
that $\pi_t\circ H$ is a unital monomorphism for all $t\in [0,1],$
\beq\label{thp3}
h_0\approx_{\ep} \pi_0\circ \Phi\,\,\,\text{on}\,\,\,{\mathcal
F}\andeqn \pi_1\circ \Phi=h_1.
\eneq
Moreover,
\beq\label{thp5}
\rm{Length}(\{\pi_t\circ \Phi\})\le 2\pi.
\eneq

\end{thm}

\begin{proof}
It follows from Theorem 3.3 of \cite{Lncd} that there exists a
sequence of unitaries $w\in A$ such that
$$
\lim_{n\to\infty}{\rm ad}\, w_n\circ h_1(f)= h_0(f)\,\rforal\,
f\in C(X).
$$
It follows from \ref{uo} that there is $z\in U_0(A)$ such that
$$
{\rm ad}\, z\circ h_1\approx_{\ep/2} h_0\,\,\,{\rm on}\,\,\, {\mathcal
F}.
$$
It follows from   \cite{LnSF}  that, for any $\ep>0,$  there is a
continuous path of unitaries $\{z_t:t\in [0,1]\}$ in $A$ such that
$$
z_0=1,\,\,\, \|z_1-z\|<\ep/2M,
$$
where $M=\max\{\|f\|: f\in {\mathcal F}\},$ and
$$
\text{Length}(\{z_t\})\le \pi
$$

 Define $H: C(X)\to C([0,1],A)$ by
$$
\pi_t\circ H(f)={\rm ad}\, z_t\circ h_1.
$$
Then
$$
\text{Length}(\{\pi_t\circ H\})\le 2\pi.
$$

\end{proof}



\begin{NN}

{\rm Let $\xi_1, \xi_2\in X.$ Suppose that there is a continuous
rectifiable path $\gamma: [0,1]\to X$ such that $\gamma(0)=\xi_1$
and $\gamma(1)=\xi_2.$ Let $A$ be a unital \CA. Define $\psi_1,
\psi_2: C(X)\to A$ by $\psi_i(f)=f(\xi_i)\cdot 1_A.$ These two \hm
s are homotopic. In fact one can define $H: C(X)\to C([0,1],A)$ by
$\pi_t\circ H(f)=f(\gamma(t))$ for all $f\in C(X).$ However, with
definition \ref{hp1def}, $\text{Length}(\{\pi_t\circ H\})=\infty$
even in the case that $X=[0,1],$ $\xi_1=0$ and $\xi_2=1.$ This is
because the presence of continuous functions  of unbounded
variation.  Nevertheless $\{\pi_t\circ H\}$ should be regarded
good homotopy. This leads us to consider a different notion of the
length. } \end{NN}

\begin{df}\label{lDef}
{\rm Let $X$ be a path connected compact metric space and fix a
base point $\xi_X.$ Denote by $P_t(x, y)$ a continuous path in $X$
which starts at $x$ and ends at $y.$ By {\it universal homotopy
length}, we mean the following constant:
\beq\label{ldef1}
L(X,\xi_X)=\sup_{y\in X}\inf\{\text{Length}(\{P_t(\xi_X,y)\}):
P_t(\xi_X, y)\}.
\eneq
Define
$$
\underline{L}_p(X)=\inf_{\xi_X\in X}L(X, \xi_X)\andeqn {\bar
L}_p(X)=\sup_{\xi_X\in X}L(X, \xi_X).
$$

The {\it total length} of $X,$ denote by $L(X),$ is the following
constant:
\beq\label{tol}
L(X)=\sup_{\xi, \zeta\in X}\inf\{\text{Length}(\{P_t(\xi,
\zeta)\}: P_t(\xi,\zeta)\}.
\eneq

It is clear that
\beq\label{to2}
\underline{L}_p(X)\le L(X, \xi_X)\le {\bar L}_p(X)\le L(X)\le 2L(X,\xi_X).
\eneq
If $D$ is the closed unit disk, then
\beq\label{ldef1+}
\underline{L}_p(D)=L(D, \{0\})=1\andeqn {\bar L}_p(D)=L(D)=2.
\eneq
If $S^1$ is the unit circle, then
\beq\label{ldef2}
\underline{L}_p(S^1,1)=L(S^1, 1)={\bar L}_p(S^1)=L(S^1)=\pi.
\eneq

There are plenty of examples that $L(X, \xi_X)=\infty$ for any
$\xi_X\in X.$ For example, let $X$ be the closure of the image of
map $f(t)=t\sin(1/t)$ for $t\in (0,1].$

}

\end{df}

\begin{df}\label{Ldef}
{\rm Let $X$ be a path connected compact metric space. We fix a
metric $d(-,-): X\times X\to \R_+.$ A function $f\in C(X)$ is sad
to be Lipschitz if
$$
\sup_{x, x'\in X, x\not=x'} {|f(x)-f(x')|\over{d(x, x')}}<\infty.
$$

For Lipschitz function $f,$ define
\beq\label{lipdef1}
{\bar D}_f=\sup_{x\in X} \lim_{\dt\to 0}\sup_{x',  x''\in O(x,
\dt)} {|f(x')-f(x'')|\over{d(x',x'')}}.
\eneq
Note that
$$
{\bar D}_f\le \sup_{x, x'\in X, x\not=x'} {|f(x)-f(x')|\over{d(x,
x')}}<\infty.
$$
}

\end{df}

\begin{df}

{\rm Let $A$ be a unital \CA\, and let $H: C(X)\to C([0,1],A)$ be
a homotopy path from $\pi_0\circ H$ to $\pi_1\circ H.$

By the homotopy length of $\{\pi_t\circ H\},$ denote by
$$
{\overline{ \text{Length}}}(\{\pi_t\circ H\}),
$$
we mean
$$
\sup\{L_H(f): f\in C(X), {\bar D_f}\le 1\}.
$$
}
\end{df}

\begin{prop}\label{Lnprop1}
Let $X$ be a connected compact metric space. For any unital \CA\,
$A$ and unital \hm\, $H: C(X)\to C([0,1], A),$
\beq\label{Lnprop2}
\overline{\rm{Length}}(\{\pi_t\circ H\})\le
\underline{L}_p(X)\rm{Length}(\{\pi_t\circ H\})
\eneq
for any $\xi_X\in X.$
\end{prop}

\begin{proof}
Note that if $\text{Length}(\{\pi_t\circ H\})=0,$ then
$\overline{\text{Length}}(\{\pi_t\circ H\})=0.$ So we may assume
that $\text{Length}(\{\pi_t\circ H\})\not=0.$ In this case, it is
clear that we may assume that $\underline{L}_p(X)<\infty.$ It
follows that $L(X)<\infty$ (by (\ref{to2})).

Fix $\xi_X\in X.$ We first assume that
$$
\overline{\text{Length}}(\{\pi_t\circ H\})<\infty.
$$

Fix $\ep>0.$ There is a Lipschitz function $f\in C(X)$ with ${\bar
D}_f\le 1$ such that
\beq\label{Lnprope1}
L_H(f)>\overline{\text{Length}}(\{\pi_t\circ H\})-\ep.
\eneq
Fix $x\in X.$ Let $\eta>0.$ There is a continuous rectifiable path
$Q_t(x, \xi_X)$ which connects $x$ with $\xi_X$ such that
\beq\label{Lnprope2}
\text{Length}(\{Q_t(x, \xi_X)\})< \inf\{ \text{Length}(\{P_t(x,
\xi_X)\}): P_t(x, \xi_X)\}+\eta.
\eneq

Denote by $Q^*$ the image of the path $\{Q_t: t\in [0,1]\}.$ It is
a compact subset of $X.$ For each $y\in Q^*,$ there is a $\dt_y>0$
such that
\beq\label{Lnprope2-}
|f(y')-f(y'')|\le (1+\eta)d(y,y')
\eneq
for any $y', y''\in O(y, \dt_y)$ (Note that ${\bar D}_f\le 1$).
Since $\cup_{y\in Q^*} O(y, \dt_y/2)\supset Q^*,$ one easily
obtains a partition $\{ t_0=0<t_1<t_2<\cdots <t_n=1\}$ such that,
if $x\in [t_{i-1}, t_i],$ then
\beq\label{Lnprope3}
|f(x)-f(Q_{t_i})|\le (1+\eta) d(x, Q_{t_i})\andeqn
\eneq
\beq\label{Lnaprope3+}
|f(Q_{t_i})-f(Q_{t_{i-1}})|\le  (1+\eta) d(Q_{t_i-1}, Q_{t_i}),
\eneq
$i=1,2,...,n.$

Put $g=f(x)-f(\xi_X).$ Then
\beq\label{Lnprope4}\nonumber
|g(x)|&=&|f(x)-f(\xi_X)| \le
\sum_{i=1}^n|f(Q_{t_i})-f(Q_{t_{i-1}})|\\\nonumber
&\le &
(1+\eta)\sum_{i=1}^n d(Q_{t_i}, Q_{t_{i-1}})\\\nonumber
&\le & (1+\eta)\text{Length}(\{Q_t(x, \xi)\})\\
&\le & (1+\eta)L(X, \xi_X)+\eta+\eta^2
\eneq
for all $\eta>0.$ It follows that $\|g\|\le L(X, \xi_X).$ Note
that
\beq\label{Lnprope5}
\pi_t\circ H(g)=\pi_t\circ H(f)-f(\xi_X)\cdot 1_A.
\eneq
It follows from (\ref{Lnprope5}) and (\ref{Lnprope1}) that
\beq\label{Lnprope6}
L_H(g)=L_H(f)>\overline{\text{Length}}(\{\pi_t\circ H\})-\ep.
\eneq
Therefore
\beq\label{Lnprope7}
L(X,\xi_X)\text{Length}(\{\pi_t\circ H\}) \ge L_H(g)
>\overline{\text{Length}}(\{\pi_t\circ H\})-\ep\\
\eneq
for all $\ep>0.$  Let $\ep\to 0.$ We then obtain
\beq\label{Lnprope8}
\overline{\text{Length}}(\{\pi_t\circ H\})\le L(X,
\xi_X)\text{Length}(\{\pi_t\circ H\})
\eneq
for any $\xi_X\in X.$ Thus (\ref{Lnprop2}) holds.

If $\overline{\text{Length}}(\{\pi_t\circ H\})=\infty,$ let $L>0.$
We replace (\ref{Lnprope1} by
\beq\label{Lnprope9}
L_H(f)>L.
\eneq
Then, instead of (\ref{Lnprope6}), we obtain that
\beq\label{Lnprope10}
L_H(g)=L_H(f)>L.
\eneq
Therefore
\beq\label{Lnprope11}
L(X, \xi_X)\text{Length}(\{\pi_t\circ H\}) \ge L_H(g)>L.
\eneq
It follows that
$$
\text{Length}(\{\pi_t\circ H\})=\infty.
$$

\end{proof}

 Fix a base point
$\xi_X$ and a point $x\in X.$ Let $\gamma: [0,1]\to X$ be a
continuous path such that $\gamma(0)=x$ and $\gamma(1)=\xi_X.$ Let
$A$ be a unital \CA. Define $h_0(f)=f(x)1_A$ and
$h_1(f)=f(\xi_X)1_A$ for all $f\in C(X).$ One obtains a homotopy
path $H: C(X)\to C([0,1],A)$ by $\pi_t\circ H(f)=f(\gamma(t))1_A$
for $f\in C(X).$  Then one has the following easy fact:

\begin{prop}\label{lnprop1}
\beq\label{Lng}
{\overline{\rm{Length}}}(\{\pi_t\circ H\})\le
{\rm{Length}}(\{\gamma(t): t\in [0,1]\}).
\eneq

\end{prop}

\begin{proof}
Let $\ep>0.$ There is $f\in C(X)$ with ${\bar L}_f\le 1$ such that
\beq\label{lnprope1}
L_H(f)>{\overline{\text{Length}}}(\{\pi_t\circ H\})-\ep/2.
\eneq
There is a partition $P:$
$$
0=t_0<t_1<\cdots <t_n=1
$$
such that
\beq\label{lnprope2}
\sum_{j=1}^n| f(\gamma(t_j))-f(\gamma(t_{j-1}))|>L_H(f)-\ep/2.
\eneq
Let $\eta>0.$ For each $t\in [0,1],$ there is $\dt_t>0$ such that
\beq\label{lnprope3}
|f(\gamma(t'))-f(\gamma(t''))|\le
(1+\eta)d(\gamma(t'),\gamma(t''))
\eneq
if $\gamma(t'),\gamma(t'')\in O(\gamma(t), \dt_t).$ From this, one
obtains a finer partition
$$
0=s_0<s_1<\cdots <s_m=1
$$
such that
\beq\label{lnprope4}
|f(\gamma(s_i))-f(\gamma(s_{i-1}))| &\le& (1+\eta)d(\gamma(s_i),
\gamma(s_{i-1})),
\eneq
$i=1,2,...,n$ and
\beq\label{lnprope4+}
 \sum_{i=1}^m
|f(\gamma(s_i))-f(\gamma(s_{i-1}))|&>&L_H(f)-\ep/2.
\eneq
It follows that
\beq\nonumber
L_H(f)-\ep/2&<&\sum_{i=1}^m (1+\eta)d(\gamma(s_i), \gamma(s_{i-1})\\
&\le & (1+\eta)\text{Length}(\{\gamma(t): t\in [0,1]\})
\eneq
for any $\eta>0.$ Let $\eta\to 0,$ one  obtains
\beq\label{lnprope5}
L_H(f)-\ep/2 \le \text{Length}(\{\gamma(t): t\in [0,1]\})
\eneq
Combining (\ref{lnprope1}), one has
\beq\label{lnprope6}
\overline{\text{Length}}(\{\pi_t\circ H\})-\ep \le
\text{Length}(\{\gamma(t): t\in [0,1]\}).
\eneq
Let $\ep\to 0.$ One conclude that
\beq\label{lnprope7}
\overline{\text{Length}}(\{\pi_t\circ H\}) \le
\text{Length}(\{\gamma(t): t\in [0,1]\}).
\eneq

\end{proof}

\begin{rem}
{\rm  If $X$ is a path connected compact subset of $\C^n,$ with
$$\di(\xi,\zeta)=\max_{1\le i\le n}|x_i-y_i|,$$ where
$\xi=(x_1,x_2,...,x_n)$ and $y=(y_1,y_2,...,y_n).$ One easily
shows that (\ref{Lng}) becomes an equality.
 When  $X=\gamma$ and $\gamma$
is smooth, then it is also easy to check that (\ref{Lng}) becomes
equality.}
\end{rem}

\begin{lem}\label{LengthL}
Let $X$ be a path connected compact metric space and let $\xi_X\in
X$ be a point.

{\rm (1)}\,\,\, Then $L(X\times S^1,\xi_X\times 1)\le \sqrt{L(X,
\xi)^2+\pi^2}.$

{\rm (2)} Suppose that $A$ is a unital \CA\, and $h: C(X)\to A$ is
unital \hm\, with finite dimensional range. Let $\psi: C(X)\to A$
be defined by $\psi(f)=f(\xi_X)\cdot 1_A.$ Then, for any $\ep>0,$
and finite subset ${\mathcal F}\subset C(X),$  there exist two \hm s
$H_1, H_2: C(X)\to C([0,1], A)$ such that
\beq\label{LenL1-}
\pi_0\circ H_1=h, \,\,\pi_1\circ H_1=\psi\rforal f\in C(X),\\
\overline{\rm{Length}}(\{\pi_t\circ H\})\le L(X, \xi_X)+\ep
\eneq
and
\beq\label{lenL1--}
\pi_0\circ H_2=h, \,\,\,\pi_1\circ H_2\approx_{\ep} \psi\rforal
f\in {\mathcal F}\andeqn\\
\overline{\rm{Length}}(\{\pi_t\circ H_2\})\le L(X, \xi_X).
\eneq

{\rm (3)} Suppose that $h_1, h_2: C(X)\to A$ are two unital \hm s
so that $h_1(A), h_2(A)\subset B\subset A,$ where $B$ is a unital
commutative finite dimensional \SCA. Then, for any $\ep>0,$  there
is a unital \hm\, $H: C(X)\to C([0,1],B)$ such that
\beq\label{lenL1-01}
\pi_t\circ H=h_1,\,\,\, \pi_1\circ H=h_2\andeqn\\
\overline{\rm{Length}}(\{\pi_t\circ H\})\le L(X)+\ep.
\eneq

 {\rm (4)}\,\,\, Suppose
that $h: C(X\times S^1)\to A$ is a unital \hm\, with finite
dimensional range, where $A$ is a unital \CA. Then, for any
$\ep>0,$ there exists a \hm\, $H: C(X\times S^1)\to C([0,1], A)$
such that
\beq\label{LenL1}
\pi_0\circ H=h, \,\,\pi_1\circ H(f)=f(\xi_X\times 1)1_A\rforal
f\in C(X\times S^1),
\eneq
\beq\label{LenL2}
{\overline{\rm{Length}}}(\{\pi_t\circ H|_{C(X)\otimes 1}\})\le
L(X, \xi_X)+\ep,
\eneq
\beq\label{LenL3}
{\overline{\rm{Length}}}(\{\pi_t\circ H|_{1\otimes C(S^1)}\})\le
\pi\andeqn
\eneq
\beq\label{LenL4}
[\pi_t\circ H(1\otimes z), \pi_{t'}\circ H(f\otimes 1)]=0
\eneq
for any $f\in C(X)$ and $t, t'\in [0,1].$

{\rm (5)} If $h_1, h_2: C(X\times S^1)\to B$ are two unital \hm s, where $B$ is a unital commutative
finite dimensional \CA,
then
for any
$\ep>0,$ there exists a \hm\, $H: C(X\times S^1)\to C([0,1], B)$
such that
\beq\label{LL-51}
\pi_0\circ H=h_1, \,\,\pi_1\circ H=h_2
\eneq
\beq\label{LL-52}
{\overline{\rm{Length}}}(\{\pi_t\circ H|_{C(X)\otimes 1}\})\le
L(X)+\ep,
\eneq
\beq\label{LL-53}
{\overline{\rm{Length}}}(\{\pi_t\circ H|_{1\otimes C(S^1)}\})\le
\pi.
\eneq

\end{lem}

\begin{proof}

We may assume that $L(X, \xi_X)<\infty.$ Then $L(X)<\infty.$

(1) is obvious.

For (2), let
$$
h(f)=\sum_{i=1}^n f(x_i)e_i\rforal f\in C(X).
$$
where $x_1,x_2,...,x_n\in X$ and $\{e_i: i=1,2,...,n\}$ is a set
of mutually orthogonal projections such that $\sum_{i=1}^n
e_i=1_A.$ For any $\ep>0,$ define
$$
\eta=\sigma_{X, \ep, {\mathcal F}}\andeqn \ep_1=\min\{\ep, \eta\}.
$$
For each $i,$ there is a rectifiable continuous path $\gamma_i:
[0,1]\to X$ such that $\gamma_i(0)=x_i$ and $\gamma_i(1)=\xi_X$
such that
$$
\text{Length}(\{\gamma_i\})< L(X,\xi_X)+\ep_1,\,\,\,i=1,2,...,n,
$$

Define $H_1: C(X\times S^1)\to C([0,1],A)$ by
$$
\pi_t\circ H_1(f)=\sum_{i=1}^n f(\gamma_i(t))e_{i}
$$
for all $f\in C(X\times S^1)$ and $t\in [0,1].$

Fix a  partition ${\mathcal P}: 0=t_0<t_1<\cdots <t_k=1.$ For any
$f\in C(X)$ with ${\bar L}_f\le 1,$ by applying \ref{lnprop1},
\beq\label{LL-1}\nonumber
\sum_{j=1}^k\|\pi_{t_j}\circ H_1(f)-\pi_{t_{j-1}}\circ H_1(f)\|
&=& \sum_{j=1}^k
(\sum_{i=1}^n |f(\gamma_i(t_j))-f(\gamma_i(t_{j-1}))|e_i)\\\nonumber
&=& \sum_{i=1}^n
\sum_{j=1}^k|f(\gamma_i(t_j))-f(\gamma_i(t_{j-1}))|e_i)\\\nonumber
&\le
&\max_{1\le i\le n}\{\sum_{j=1}^k|f(\gamma_i(t_j))-f(\gamma_i(t_{j-1}))|\}\\
&\le &\ L(X,\xi_X)+\ep_1.
\eneq
It follows that
\beq\label{LL-2}
\overline{\text{Length}}(\{\pi_t\circ H_1\})\le L(X, \xi_X)+\ep.
\eneq

 By re-parameterizing the
paths $\gamma_i,$ we may assume that, for some $a\in (0,1),$
\beq\nonumber
\text{Length}(\{\gamma_i(t): t\in [0, a]\})\le L(X,\xi_X)
\andeqn\\
\text{Length}(\{\gamma_i(t): t\in [a,1]\})<\eta,
\eneq
$i=1,2,...,n.$ Define $H_2: C(X)\to C([0,1],A)$ by $\pi_t\circ
H_2=\pi_{t\over{a}}\circ H_1.$  We have that
$$
\di(\xi_X, \gamma_i(a))<\eta,\,\,\,i=1,2,...,n.
$$
It follows that, for $f\in {\mathcal F},$
\beq\label{LL-4}\nonumber
\|\pi_1\circ H_2(f)-\psi(f)\|&=&\|\pi_a\circ H_1(f)-f(\xi_X)\cdot 1_A\|\\
&=&\|\sum_{i=1}^nf(\gamma_i(a))e_i-\sum_{i=1}^n f(\xi_X)
e_i\|<\ep.
\eneq
Moreover,
\beq\label{LL-5}
\overline{\text{Length}}(\{\pi_t\circ H_2\})\le L(X, \xi_X).
\eneq

For (3), suppose that $B$ is the commutative finite dimensional
\SCA\, generated by $\{e_1,e_2,...,e_m\},$  where $e_1,e_2,...,e_m$ are non-zero mutually orthogonal
projections.
We may write
$$
h_i(f)=\sum_{k=1}^m f(x_{i,k})e_k\rforal f\in C(X),
$$
where $\{x_{i,k}: 1\le k\le m\}$ is a set of (not
necessarily distinct) points in $X,$ $i=1,2.$
 For any $\ep>0,$ for each $k,$ there
exists a continuous path $\{\gamma_k: t\in [0,1]\}$
 in $X$ such that
 $$
 \gamma_k(0)=x_{1,k},\,\,\,\gamma_k(1)=x_{2,k}\andeqn
 $$
 $$
 \text{Length}(\{\gamma_k(t)\})\le L(X)+\ep.
 $$
Define $H: C(X)\to C([0,1],B)$ by
$$
\pi_t\circ H(f)=\sum_{k=1}^m f(\gamma_k(t))e_k\rforal f\in C(X).
$$
Similar to the estimate (\ref{LL-1}), we obtain that
$$
\overline{\text{Length}}(\{\pi_t\circ H\})\le L(X)+\ep.
$$

 For (4),
we write
$$
h(f)=\sum_{i=1}^n\sum_{j=1}^m f(\xi_i\times t_j)e_{i,j}
$$
for all $f\in C(X\times S^1),$ where $x_1,x_2,...,x_n\in X$ and
$t_1,t_2,...,t_m\in S^1$, and $\{e_{i,j}, i=1,2,...,n,
j=1,2,...,m\}$ is a set of mutually orthogonal projections such
that $\sum_{i,j}e_{i,j}=1_A.$ For each $i,$ there is a continuous
path $\gamma_i: [0,1/2]\to X$ such that $\gamma_i(0)=x_i$ and
$\gamma_i(1)=\xi_X.$ Define $h_t: C(X\times S^1)\to A$ by
$$
h_t(f)=\sum_{i=1}^n \sum_{j=1}^m f(\gamma_i(t)\times t_j)e_{i,j}
$$
for all $f\in C(X\times S^1)$ and $t\in [0,1/2].$ If $L(X,
\xi_X)<\infty,$ then, for any $\ep>0,$ we may also assume that
$$
\text{Length}(\{\gamma_i\})\le L(X,\xi_X)+\ep,\,\,\,i=1,2,...,n.
$$

There is for each $j,$ a rectifiable continuous path $\lambda_j:
[1/2,1]\to S^1$ such that
$$
\lambda_j(0)=t_j,\,\,\,\,\lambda_j(1)=1\andeqn
\text{Length}(\{\lambda_j\})\le \pi
$$
$j=1,2,...,m.$ For $t\in [1/2, 1],$ define
$$
h_t(f)=\sum_{i=1}^n\sum_{j=1}^m f(\xi_X\times \lambda(t))e_{i,j}
$$
for all $f\in C(X\times S^1).$

Note that
\beq\label{LenL01}
h_t(f)h_{t'}(f)=h_{t'}(f)h_t(f)\,\rforal f\in C(X\times S^1)
\eneq
and any $t, t'\in [0,1].$ We then define $H: C(X\times S^1)\to
C([0,1],A)$ by
$$
H(f)(t)=h_t(f)\,\rforal\, f\in C(X) \andeqn t\in [0,1]
$$
By the proof of (2), it is easy to check  that (\ref{LenL1}),
(\ref{LenL2}),(\ref{LenL3}) and (\ref{LenL4}) hold.

Part (5) follows from the combination of the proof of (3) and (4).
\end{proof}

\vspace{0.2in}

\section{Approximate homotopy  for \hm s}

\begin{thm}\label{GHT1}
Let $X$ be a path connected compact metric space with the base
point $\xi_X$ and let $A$ be a unital separable simple \CA\, which
has tracial rank zero, or is purely infinite. Suppose that $h:
C(X)\to A$ is a unital monomorphism such that
\beq\label{lt1}
[h|_{C_0(Y_X)}]=\{0\}\,\,\,\,\text{in}\,\,\,KL(C(X), A).
\eneq
Then, for any $\ep>0$ and any compact subset ${\mathcal F}\subset
C(X),$ there is \hm\, $H: C(X) \to  C([0,1],A)$ such that
\beq\label{lt2}
\pi_0\circ H\approx_{\ep} h\,\,\,\,\text{on}\,\,\, {\mathcal
F}\andeqn \pi_1\circ H=\psi,
\eneq
where $\psi(f)=f(\xi_X)\cdot 1_A.$ Moreover,
\beq\label{lt3}
{\overline{\rm{Length}}}(\{\pi_t\circ H\})\le L(X,\xi_X).
\eneq

\end{thm}

\begin{proof}
Note that $[\psi|_{C_0(Y_X)}]=0.$ Thus $[\psi]=[h]$ in $KL(C(X),A).$
By Theorem 3.8 of \cite{HLX3} (in the case that $TR(A)=0$), or by
Theorem 1.7 of \cite{D1} (in the case that $A$ is purely infinite),
there is a \hm\, $h_0: C(X)\to A$ with finite dimensional range
such that
$$
h_1\approx_{\ep/2} h_0\,\,\,{\rm on}\,\,\, {\mathcal F}.
$$
 Since $h_0$ has finite dimensional range, by
\ref{LengthL}, there is  a \hm\, $H: C(X)\to C([0,1],A)$ such that
$$
\pi_0\circ H\approx_{\ep/2}h_0\,\,\,\text{on}\,\,\,{\mathcal
F}\,\andeqn \pi_1\circ H=\psi.
$$
Moreover, $${\overline{\text{Length}}}(\{\pi_t\circ H\})\le
L(X,\xi_X).$$

\end{proof}

\begin{prop}\label{GHT0}
Let $A$ and $B$ be two unital \CA s and let $h_1, h_2: A\to B$ be
two unital \hm s. Suppose that, for any $\ep>0$ and any finite
subset ${\mathcal F}\subset A,$  there is a unital \hm\, $H: A\to
C([0,1], B)$ such that
$$
\pi_0\circ H=h_1\andeqn \pi_1\circ H\approx_{\ep}
h_2\,\,\,\text{on}\,\,\,{\mathcal F}.
$$
Then
$$
[h_1]=[h_2]\,\,\,\text{in}\,\,\,KL(A,B).
$$
\end{prop}

\begin{proof}
This follows straightforward  from the definition and see \ref{Khp}.

\end{proof}

Proposition \ref{GHT0} states an obvious fact that if $h_1$ and
$h_2$ are homotopic then $[h_1]=[h_2]$ in $KL(A,B).$ We will show
that, when $A=C(X)$ and $B$ is a unital purely infinite simple
\CA\, or $B$ is a unital separable simple \CA\, of tracial rank
zero, at least for monomorphisms, $[h]$ is a complete
approximately homotopy invariant.

\begin{thm}\label{GHT2}
Let $X$ be a compact metric space and let $A$ be a unital purely
infinite simple \CA. Suppose that $h_1, h_2: C(X)\to A$ are two
unital monomorphisms such that
\beq\label{ghp2-1}
[h_1]=[h_2]\,\,\,\text{in}\,\,\,KL(C(X),A).
\eneq
Then, for any $\ep>0$ and any finite subset ${\mathcal F}\subset
C(X),$  there is a unital \hm\, $H: C(X)\to C([0,1],A)$ such that
\beq\label{ghp2-2}
\pi_0\circ H=h_1\andeqn \pi_1\circ H\approx_{\ep}
h_2\,\,\,\text{on}\,\,\,{\mathcal F},
\eneq
Moreover, each $\pi_t\circ H$ is a monomorphism {\rm (}for $t\in
[0,1]${\rm )} and
\beq\label{ghp2-3}
\rm{Length}(\{\pi_t\circ H\})\le 2\pi+\ep \andeqn
\eneq
If, in addition, $X$ is path connected,
then
\beq\label{ghp2-3+}
\overline{\rm{Length}}(\{\pi_t\circ H\})\le 2\pi
\underline{L}_p(X)+\ep.
\eneq

Furthermore, the converse also holds.

\end{thm}

\begin{proof}
It follows from a theorem of Dadarlat \cite{D1} (see \ref{NEW})
that $h_1$ and $h_2$ are approximately unitarily equivalent, by
applying \ref{uo}, we obtain a unitary $u\in A$ with $[u]=0$ in
$K_1(A)$ such that
$$
{\rm ad}\, u\circ h_1\approx_{\ep} h_2.
$$
By a theorem of N. C. Phillips (\cite{Ph1}), there is a continuous
unitaries $\{u_t: t\in [0,1]\}$ (see \ref{GHT1}, or see Lemma
4.4.1 of \cite{Lnbk} for the exact statement) such that
$$
u_0=1,\,\,\,u_1=u\andeqn \text{Length}(\{u_t\})\le \pi+\ep/2.
$$
Define $H: C(X)\to C([0,1],A)$ by $H(f)(t)={\rm ad}\, u_t\circ
h_1(f)$ for $f\in C(X).$  Then,
$$
\pi_0\circ H=h_1\andeqn \pi_1\circ H={\rm ad}\, u\circ h_1.
$$
Moreover, we compute  that
\beq\label{GHT2-1}
\text{Length}(\{\pi_t\circ H\})\le 2\pi+\ep .
\eneq

If $X$ is path connected and $L(X,\xi_X)<\infty,$ we may choose
$\{u_t: t\in [0,1]\}$ so that
$$
\text{Length}(\{u_t\})\le \pi+{\ep\over{4(1+L(X))}}.
$$
Then (\ref{GHT2-1}) becomes
$$
\text{Length}(\{\pi_t\circ H\})\le 2\pi+{\ep\over{2(1+L(X))}}.
$$
By \ref{Lnprop1},
$$
\overline{\text{Lnegth}}(\{\pi_t\circ H\})\le 2\pi
\underline{L}_p(X)+\ep.
$$

\end{proof}

\begin{rem}\label{remconn}
{\rm  The assumption that both $h_1$ and $h_2$ are monomorphism is
important. Suppose that $X$ has at least two path connected
components, say $X_1$ and $X_2.$ Fix two points, $\xi_i\in X_i,$
$i=1,2.$ Suppose that $A$ is a unital purely infinite simple \CA\,
with $[1_A]=0$ in $K_0(A).$ Define $h_i: C(X)\to A$ by
$h_i(f)=f(\xi_i)\cdot 1_A$ for all $f\in C(X),$ $i=1,2.$ Then
$[h_1]=[h_2]=0$ in $KL(C(X),A).$ It is clear that they are not
approximately  homotopic.  However, when $X$ is path connected, we
have the following approximate homotopty result. }

\end{rem}

\begin{thm}\label{GHP3}
Let $X$ be a path connected compact metric space and let $A$ be a
unital purely infinite simple \CA. Suppose that $h_1, h_2: C(X)\to
A$ are two unital \hm s such that
$$
[h_1]=[h_2]\,\,\,\text{in}\,\,\, KL(C(X),A).
$$
Then, for any $\ep>0$ and any finite subset ${\mathcal F}\subset
C(X),$ there exist two  unital \hm s $H_1,\, H_2: C(X)\to
C([0,1],A)$ such that
$$
\pi_0\circ H_1\approx_{\ep/3}h_1,\,\,\, \pi_1\circ
H_1\approx_{\ep/3} \pi_0\circ H_2 \andeqn \pi_1\circ
H_2\approx_{\ep/3} h_2\,\,\,\text{on}\,\,\,{\mathcal F}.
$$
Moreover,
$$
\overline{\rm{Length}}(\{\pi_t\circ H_1\})\le
L(X,\xi_1)(1+2\pi)+\ep/2 \andeqn
$$
$$
\overline{\rm{Length}}(\{\pi_t\circ H_2\})\le L(X,\xi_2)+\ep/2,
$$
where $\xi_i$ is (any) point in the spectrum of $h_i,$ $i=1,2.$

{\rm (Note $L(X, \xi_i)\le {\bar L}_p(X)$)}

\end{thm}

\begin{proof}
Let $\ep>0$ and ${\mathcal F}\subset C(X)$ be finite subset. Let
$\dt>0$ and ${\mathcal G}\subset C(X)$ be a finite subset required in
\ref{NEW} for $\ep/3$ and ${\mathcal F}.$

Suppose that $\xi_i$ is a point in the spectrum of $h_i,$ $i=1,2.$
It follows from virtue of \ref{pi1} and \ref{pi2} that we may
write that
\beq\label{GHT3-1}
&&\|h_1(f)-(h_1'(f)+f(\xi_1)p_1)\|<\ep/3\andeqn \\
&&\|h_2(f)-(h_2'(f)+f(\xi_2)p_2)\|<\ep/3
\eneq
for all $f\in {\mathcal G},$ where $p_1, p_2$ are projections in $A$
with $[p_1]=[p_2]=0$ and $h_i': C(X)\to (1-p_i)A(1-p_i)$ is a
unital monomorphism with $[h_i']=[h_i]=[h_1]$ in $KL(C(X),A),$
$i=1,2,.$ Denote $\phi_i(f)=f(\xi_i)p_i$ for all $f\in C(X).$
Suppose that $\{x_1, x_2,...,x_m\}\subset C(X)$ such that
\beq\label{GHT3-2}
\max\{\|f(x_i)\|: i=1,2,...,m\}\ge (3/4)\|f\|\rforal f\in {\mathcal
G}.
\eneq
There are $m$ non-zero mutually orthogonal projections
$p_{i,1},p_{i,2},...,p_{i,m}\in p_iAp_i$ such that $[p_{i,1}]=0,$
$i=1,2.$ Define
\beq\label{GHT3-3}
\phi_i'(f)=\sum_{k=1}^mf(x_i)p_{i,k}\rforal f\in C(X),\,\,\,i=1,2.
\eneq
By \ref{LengthL} there is a unital \hm\, $H_i': C(X)\to
C([0,1],p_iAp_i)$ such that
\beq\label{GHT3-4}
\pi_0\circ H_i'=h_i'+\phi_i\andeqn \pi_1\circ H_i'=h_i'+\phi_i'
\eneq
$i=1,2.$ Moreover,
\beq\label{GHT3-5}
\overline{\text{Length}}(\{\pi_t\circ H_i'\})\le L(X,
\xi_i)+\ep/2,\,\,\,i=1,2.
\eneq
Note that both $h_i' +\phi_i$ are $1/2$-${\mathcal G}$-injective \hm s
and
\beq\label{GHT3-6}
[h_1'+\phi_1']=[h_1']=[h_1]=[h_2'+\phi_2'].
\eneq
Combining \ref{NEW} and \ref{uo}, we obtain a unitary $u\in A$
with $[u]=0$ in $K_1(A)$ such that
\beq\label{GHT3-7}
{\rm ad}\, u\circ (h_1'+\phi_1')\approx_{\ep/3}
h_2'+\phi_2'\,\,\,\text{on}\,\,\,{\mathcal G}.
\eneq
Let $\{u_t: t\in [0,1/2]\}$ be a continuous path of unitaries in
$A$ such that $u_0=1,$ $u_1=u$ and
$${\rm{Lenghth}}(\{u_t\})\le
\pi+\ep/8(1+L(X)).$$

Now define $H_1: C(X)\to C([0,1], A)$ as follows
$$
\pi_t\circ H_1=\begin{cases} \pi_{2t}\circ H_1' & \rforal t\in
[0,1/2]\\
{\rm ad}\,u_t\circ (h_1'+\phi_1') & \rforal t\in (1/2,1].
\end{cases}
$$
Define $H_2: C(X)\to C([0,1],A)$ by $\pi_t\circ H_2=\pi_{1-t}\circ
H_2'.$

We check (applying \ref{lnprop1}) that $H_1$ and $H_2$ meet the
requirements of the theorem.

\end{proof}

\begin{lem}\label{K-semi}
Let $X$ be a finite CW complex, let $\ep>0$ and let ${\mathcal
F}\subset C(X)$ be a finite subset. Let $\eta=\sigma_{X, \ep/16,
{\mathcal F}}.$
Let $\{x_1, x_2,...,x_m\}$ be an
$\eta/2$-dense subset of $X$ for which $O_i\cap O_j=\emptyset$ (if
$i\not=j$), where
$$
O_i=\{x\in X: \text{dist}(x, x_i)<\eta/2s\},\,\,\, i=1,2,...,m
$$
for some integer $s\ge 1.$ Let $0<\sigma<1/2s.$

There exists $\dt>0,$  a finite subset ${\mathcal G}\subset C(X)$ and
a finite subset ${\mathcal P}\subset \underline{K}(C(X))$ satisfying
the following:

Suppose that $A$ is a unital separable simple \CA\, with tracial
rank zero and suppose that $\psi: C(X)\to A$ is a unital
$\dt$-${\mathcal G}$-multiplicative \morp\, such that
\beq\label{ksemi1}
[\psi]|_{\mathcal P}=[h]|_{\mathcal P}\andeqn \mu_{\tau\circ \psi}(O_i)\ge
\sigma\cdot \eta,\,i=1,2,...,m
\eneq
for some \hm\, $h: C(X)\to A.$ Then, there exists a monomorphism
$\phi: C(X)\to A$ such that
\beq\label{ksemi2}
\psi\approx_{\ep} \phi\,\,\,\text{on}\,\,\, {\mathcal F}
\eneq

\end{lem}

\begin{proof}
To simplify the notation, by considering each component, without loss of generality, we may assume
$X$ is connected finite CW complex.
Let $\ep>0$ and ${\mathcal F}\subset C(X)$ be a finite subset. Let
$\eta>0,$ $\{x_1,x_2,..,x_m\},$ $\sigma>0,$ $s>0$ and $O_i$ be as
in the statement.

Let $\sigma_1=\sigma/2.$  Let $\gamma>0,$ ${\mathcal G}_1\subset
C(X),$ $\dt_1>0$ ( in place of $\dt$) and ${\mathcal P}\subset
\underline{K}(C(X))$ be as required by Theorem 4.6 of \cite{Lncd}
for $\ep/4$ (instead of $\ep$ ), ${\mathcal F},$ $\sigma_1$ (in place
of $\sigma$) and $\eta$ above.

We may also assume that
$$
[L_1]|_{\mathcal P}=[L_2]|_{\mathcal P}
$$
for any pair of $\dt_1$-${\mathcal G}_1$-multiplicative \morp s from $C(X)$ to
(any) unital \CA\, $B,$ provided that
$$
L_1\approx_{\dt_1} L_2\,\,\,\text{on}\,\,\,{\mathcal G}_1.
$$

Let $\dt_2>0$ and ${\mathcal G}\subset C(X)$ be a finite subset
required in \ref{ndig} for $\min\{\ep/4, \dt_1/2\}$ (in pace of
$\ep$) and ${\mathcal G}_1$ (in place of ${\mathcal F}$). We may
assume that $\dt_2<\min\{\dt_1, \ep/4\}$ and ${\mathcal G}_1\cup
{\mathcal F}\subset {\mathcal G}.$

Let $\dt=\dt_2/4.$ Let $\psi$ be as in the statement. Without loss
of generality, by applying Lemma \ref{ndig} (with $\psi=\phi$) we
may write that
\beq\label{ksem1}
\psi(f)\approx_{\dt_1/2} \psi_1\oplus h_0(f)\,\rforal \, f\in
{\mathcal G}_1,
\eneq
where $\psi_1: C(X)\to (1-p)A(1-p)$ is a unital
$\dt_1/2$-${\mathcal G}_1$-multiplicative \morp\, and $h_0:
C(X)\to  pAp$ is a \hm\, with finite dimensional range, $p$ is a
projection with $\tau(1-p)<\gamma/2$ for all $\tau\in T(A).$ We
may further assume that
\beq\label{ksm2}
\mu_{\tau\circ h_0}(O_i)>\sigma_1\cdot \eta/2=\sigma\cdot
\eta,\,\,\, i=1,2,...,m.
\eneq
Note that we assume that $X$ is connected. Fix $\xi_X\in X$  and
let $Y_X=X\setminus \{\xi_X\}.$ Since $h_0$ is a \hm\, with finite
dimensional range,
\beq\label{ksemi3-}
[h_0|_{C_0(Y_X)}]=0\,\,\,\text{in}\,\,\,KL(C_0(Y_X),A).
\eneq

It follows from \ref{shk} that there is a unital monomorphism
$h_1: C(X)\to (1-p)A(1-p)$ such that
\beq\label{ksm3}
[h_1|_{C_0(Y_X)}]=[h|_{C_0(Y_X)}].
\eneq
By (\ref{ksem1}), (\ref{ksm3}) and (\ref{ksemi3-}),
\beq\label{ksemi3+}
[h_1]|_{\mathcal P}=[\psi_1]|_{\mathcal P}\,\,\,\text{in}\,\,\, KK(C(X), A).
\eneq
 By applying
Theorem 4.6 of \cite{Lncd}, we obtain a unitary $u\in A$ such that
\beq\label{ksmi4}
{\rm ad}\, u\circ (h_1\oplus h_0)\approx_{\ep/2} \psi_1\oplus
h_0\,\,\,\text{on}\,\,\, {\mathcal F}
\eneq
Put $h={\rm ad}\, u\circ (h_1\oplus h_0).$ Then $h$ is a
monomorphism. We have
\beq\label{ksm5}
h\approx_{\ep}\psi\,\,\,\text{on}\,\,\,{\mathcal F}.
\eneq

\end{proof}

\begin{lem}\label{chdig}
Let $X$ be a finite CW complex and let $A$ be a unital separable
simple \CA\, with tracial rank zero. Suppose that $h: C(X)\to A$
is a unital \hm\, with the spectrum $F.$ Then, for any $\ep>0$ and
any finite subset ${\mathcal F}\subset C(X)$ and $\gamma>0,$ there is
a projection $p\in A$ with $\tau(1-p)<\gamma$ for all $\tau\in
T(A)$ and a unital \hm\, $h_0: C(Y)\to pAp$ with finite
dimensional range and a unital monomorphism  $\phi: C(Y)\to
(1-p)A(1-p)$ for a compact subset $Y$ which contains $F$ and which
is a finite CW complex such that
$$
h\approx_{\ep} h_0\circ s_1\oplus \phi\circ s_1\,\,\,{\rm on}\,\,\, {\mathcal F},
$$
where $s_1: C(X)\to C(Y)$ is the quotient map.

\end{lem}

\begin{proof}

Let $\ep>0,$ $\gamma>0$ and let ${\mathcal F}\subset C(X)$ be a
finite subset. Let $\eta>0$ be such that $|f(x)-f(x')|<\ep/32$
whenever ${\rm dist}(x,x')<2\eta$ and for all $f\in {\mathcal F}.$
In other words, $\eta<{\sigma_{X,\ep/32,{\mathcal F}}\over{2}}.$

There is a finite CW complex $Y$ such that $F\subset Y\subset X$
and
\beq\label{chdig-01}
Y\subset \{x\in X: {\rm dist}(x, F)\le \eta/4\}.
\eneq

Denote by $s_1$ the quotient map from $C(X)$ onto
$C(Y).$
Suppose that $Y=\sqcup_{i=1}^K Y_i,$ where each $Y_i$ is a
connected component of $Y.$  Let $E_1, E_2,..,E_K$ be mutually
orthogonal projections in $C(Y)$ corresponding to the components
$Y_1, Y_2,...,Y_K.$

Define  $h_Y: C(Y)\to A$ such that $h_Y\circ s_1=h.$ Let
$\{x_1,x_2,...,x_m\}$ be an $\eta/4$-dense subset of $Y.$ By
replacing $x_1,x_2,...,x_m$ by points in $F,$ we can still assume
that they are $\eta/2$-dense in $Y.$ There is an integer $s\ge 1$
for which $O_i\cap O_j=\emptyset$ ($i\not=j$), where
$$
O_i=\{x\in Y: {\rm dim}(x,x_i)<\eta/2s\}.
$$
Note that $2\eta\le \sigma_{Y, \ep/32, s_1({\mathcal F})}.$

 Put
$$
\sigma_0=\inf\{\mu_{\tau\circ h_Y}(O_i): 1\le i\le m,\,\,
\,\tau\in T(A)\}/\eta.
$$
Since $A$ is simple and $x_i\in F,$ $i=1,2,...,m,$ $\sigma_0>0.$  We may assume that
$m>2.$ Choose $\sigma>0$ so small that it is smaller than $\sigma_0$
and $(1-1/m)\sigma_0.$

 Let $\dt>0,$ ${\mathcal G},$ ${\mathcal P}$ and $\gamma_1>0$ be
as required by 4.6 of \cite{Lncd}  for $Y$  corresponding to the
above $\ep/4,$ $s_1({\mathcal F})$ and $\sigma.$ Let
$$
\gamma_2=\min\{\gamma,\gamma_1, \sigma_0\eta/m\}.
$$
We may assume that
$$
[L_1]|_{\mathcal P}=[L_2]|_{\mathcal P}
$$
for any pair of $\dt$-${\mathcal G}$-multiplicative \morp\, $L_1, L_2:
C(Y)\to A$ provided that
$$
L_1\approx_{\dt} L_2\,\,\,\text{on}\,\,\,{\mathcal G}.
$$
We may also assume that $E_1,E_2,..,E_K\subset {\mathcal G}.$

 Let ${\mathcal G}_1\subset C(Y)$ be a finite subset such that ${\mathcal G}\subset
 {\mathcal G}_1.$
 By applying \ref{L2dig}, we may write
\beq\label{chdig01}
h_Y\approx_{\dt/2} \psi\oplus h_0'\,\,\,{\rm on}\,\,\, {\mathcal G}_1,
\eneq
where $\psi: C(Y)\to (1-p)A(1-p)$ is a unital $\dt$-${\mathcal
G}_1$-multiplicative \morp\, and $h_0': C(Y)\to pAp$ is a unital
\hm\, with finite dimensional range and $p\in A$ is projection in
$A$ with $\tau(1-p)<\gamma_2$ for all $\tau\in T(A).$

With sufficiently large ${\mathcal G}_1$ and smaller $\dt,$ we may
assume that
$$
\mu_{\tau\circ h_0'}(O_i)\ge (1-1/m)\sigma_0\eta\ge \sigma\eta
$$
for all $\tau\in T(A),$ $i=1,2,...,m.$

By (\ref{chdig01}), we may also assume that
\beq\label{chdig02}
[h_Y]|_{\mathcal P}=[\psi\oplus h_0']|_{\mathcal P}.
\eneq

Suppose that $h_Y(E_i)=e_i$ and $h_0'(E_i)=e_i',$ $i=1,2,...,K.$

Fix a point $\xi_i\in Y_i$ and let $\Omega_i=Y_i\setminus
\{\xi_i\}.$ By \ref{shk},  there is a unital monomorphism
$h_{00}^{(i)}: C(Y_i)\to (e_i-e_i')A(e_i-e_i')$ such that
\beq\label{chdig02+}
[h_{00}^{(i)}|_{C_0(\Omega_i)}]=[h_Y|_{C_0(\Omega_i)}].
\eneq
Define $h_{00}: C(Y)\to (1-p)A(1-p)$ by $h_{00}(f)=\sum_{i=1}^K
h_{00}^{(i)}(f|_{Y_i})$ for $f\in C(Y).$  It follows from
(\ref{chdig02}) that
$$
[h_{00}]|_{\mathcal P}=[\psi]|_{\mathcal P}.
$$

It follows from 4.6  of \cite{Lncd} that there is a unitary $w\in
A$ such that
$$
h_0'\oplus \psi\approx_{\ep/4} {\rm ad}\, w\circ (h_0'\oplus h_{00}
) \,\,\,{\rm on}\,\,\,s_1({\mathcal F}).
$$
Define $h_0={\rm ad}\, w\circ h_0'$ and $\phi={\rm ad}\, w\circ
h_{00}.$

From the above, $h_0$ and $\phi$ satisfy the requirements.

\end{proof}





\begin{thm}\label{SHLn}
Let $X$ be a connected finite CW complex and let $A$ be a unital
separable simple \CA\, with tracial rank zero. Suppose that
$\psi_1, \psi_2: C(X)\to A$ are two unital monomorphisms such that
\beq\label{shln1}
[\psi_1]=[\psi_2]\,\,\,\text{in}\,\,\, KL(C(X), A).
\eneq
Then, for any $\ep>0$ and any finite subset ${\mathcal F}\subset
C(X),$
 there exists a unital \hm\, $H: C(X)\to C([0,1], A)$ such that
\beq\label{shln2}
\pi_0\circ H\approx_{\ep} \psi_1\andeqn \pi_1\circ H\approx_{\ep}
\psi_2\,\,\,\text{on}\,\,\,{\mathcal F}.
\eneq
Moreover, each $\pi_t\circ H$ is a monomorphism and
\beq\label{shln3}
\overline{\rm{Length}}(\{\pi_t\circ H\})\le
L(X)+\underline{L}_p(X)2\pi+\ep.
\eneq
\end{thm}

\begin{proof}
Let $\dt>0$ and ${\mathcal G}\subset C(X)$ be a finite subset
required by Theorem 3.3   of \cite{Lncd} (for $\psi_2$ being the
given unital monomorphism) and for $\ep/4$ and ${\mathcal F}$
given. By \ref{chdig} we may write that
\beq\label{shln4}
\psi_1\approx_{\min\{\ep/4, \dt/2\}} h_1\oplus
h_0\,\,\,\text{on}\,\,\, {\mathcal G}.
\eneq
where $h_1: C(X)\to (1-p)A(1-p)$ is a unital monomorphism,
$h_0: C(X)\to pAp$ is  a unital \hm\, with finite dimensional
range and $p\in A$ is a projection such that $\tau(1-p)<\dt/2$ for
all $\tau\in T(A).$ By \ref{mmm}, there is a unital \hm\, $h_{00}: C(X)\to pAp$
with finite dimensional range such that
\beq\label{shln5}
\sup\{|\tau\circ h_{00}(f)-\tau\circ \psi_2(f)|: f\in {\mathcal G},
\tau\in T(A)\}<\dt/2.
\eneq

Write
\beq\label{shln5+1}
h_0(f)=\sum_{i=1}^m f(x_i)p_i\andeqn h_{00}(f)=\sum_{j=1}^m f(x_j)e_j\rforal f\in C(X),
\eneq
where $x_i, y_j\in X$ and $\{p_1,p_2,...,p_n\}$ and
$\{e_1,e_2,...,e_m\}$ are two sets of mutually orthogonal
projections with $p=\sum_{i=1}^np_i=\sum_{j=1}^m e_j.$ It follows
from \ref{zhr2} that there is a unital commutative finite
dimensional \SCA\, $B$ of $pAp$ which contains $p_1,p_2,...,p_n$
and a set of mutually orthogonal projections
$\{e_1',e_2',...,e_m'\}$ such that $[e_j']=[e_j],$ $j=1,2,...,m.$
Define
\beq\label{sln5+2}
h_{00}'(f)=\sum_{j=1}^m f(x_j)e_j'\rforal f\in C(X).
\eneq

Note, by \ref{LengthL}, there is unital \hm\, $\Phi_1: C(X)\to
C([0,1],B)$ such that
$$
\pi_0\circ \Phi_1=h_0\rforal f\in C(X),\,\,\,\pi_1\circ
\Phi_1=h_{00}'\andeqn
$$
$$
\overline{\text{Length}}(\{\pi_t\circ \Phi_1\})\le
L(X)+\ep/2.
$$
 It follows that there is a \hm\, $H_1: C(X)\to
pAp$ such that
\beq\label{shln6}
\pi_0\circ H_1=h_1\oplus h_0\,\,\,{\rm
on}\,\,\,{\mathcal F},\,\,\,\pi_1\circ H_1=h_1\oplus h_{00}'\andeqn
\eneq
\beq\label{shln7}
\overline{\text{Length}}(\{\pi_t\circ H_1\})\le L(X)+\ep/2.
\eneq

Since $[e_j']=[e_j],$ $j=1,2,...,m,$ by (\ref{shln5}),
$$
\sup\{|\tau\circ (h_1\oplus h_{00}')(f)-\tau\circ \psi_2(f)|:\tau\in T(A)\}<\dt
$$
for all $f\in {\mathcal G}.$

By applying Theorem 3.3 of \cite{Lncd}, we obtain a unitary $u\in A$ such that
\beq\label{shln8}
{\rm ad}\, u\circ (h_1\oplus h_{00})\approx_{\ep/4}\psi_2\,\,\,\text{on}\,\,\,{\mathcal F}.
\eneq
It follows from  \ref{uo} that we may assume that
$u\in U_0(A).$ It follows that there is a continuous path of
unitaries $\{u_t: t\in [0,1]\}$ of $A$ for which
\beq\label{shln9}
u_0=1_A, \,\,\, u_1=u\andeqn \text{Length}(\{u_t\})\le
\pi+{\ep\over{4(1+L(X))}}.
\eneq
Connecting $\pi_t\circ H_1$ with ${\rm ad}\, u_t\circ (h_1\oplus
h_{00}),$ we obtain a \hm\, $H: C(X)\to C([0,1], A)$ such that
\beq\label{shln10}
\pi_0\circ H=h_1\oplus h_0\,\,\,{\rm on}\,\,\,{\mathcal
F}, \,\,\,\pi_1\circ H_1={\rm ad}\, u\circ (h_1\oplus
h_{00})\andeqn
\eneq
\beq\label{shln11}
\overline{\text{Length}}(\{\pi_t\circ H\})\le L(X)
+\underline{L}_p(X)2\pi+\ep.
\eneq
It is clear that $H$ meets the requirements.

\end{proof}

\begin{thm}\label{GHTh2}
Let $X$ be a connected finite CW complex and let $A$ be a unital
separable simple \CA\, with tracial rank zero. Suppose that $\psi_1,
\psi_2: C(X)\to A$ are unital \hm s such that
$$
[\psi_1]=[\psi_2]\,\,\,\text{in}\,\,\, KL(C(X),A).
$$
Then, for any $\ep>0$  and any finite subset ${\mathcal F}\subset
C(X),$ there exist  two \hm\, $H_i: C(X)\to C([0,1],A)$ {\rm
(}$i=1,2${\rm ) } such that
$$
\pi_0\circ H_1\approx_{\ep/3} \psi_1, \pi_1\circ
H_1\approx_{\ep/3} \pi_0\circ H_2 \andeqn \pi_1\circ
H_2\approx_{\ep/3} \psi_2\,\,\,\text{on}\,\,\,{\mathcal F}.
$$
Moreover,
$$
\overline{\rm{Length}}(\{\pi_t\circ H_1\})\le L(X)+2\pi
\underline{L}_p(X)+\ep/2 \andeqn
$$
$$
\overline{\rm{Length}}(\{\pi_t\circ H_2\})\le L(X)+\ep/2.
$$
\end{thm}

\begin{proof}
It follows from \ref{shk} that there is a unital monomorphism
$\phi: C(X)\to A$ such that
\beq\label{ghth2-1}
[\phi]=[\psi_1]=[\psi_2]\,\,\,\text{in}\,\,\,KL(C(X),A).
\eneq
Let $\dt>0$ and ${\mathcal G}\subset C(X)$ be a finite subset
required by Theorem 3.3   of \cite{Lncd} (for $\phi$ being the
given unital monomorphism) and for $\ep/8$ and ${\mathcal F}$
given. By \ref{chdig} we may write that
\beq\label{ghth2-2}
\psi_i\approx_{\min\{\ep/4, \dt/2\}} h_i\oplus
h_0^{(i)}\,\,\,\text{on}\,\,\, {\mathcal G}.
\eneq
where $h_i: C(X)\to (1-p_i)A(1-p_i)$ is a unital \hm,
$h_0^{(i)}: C(X)\to p_iAp_i$ is  a unital \hm\, with finite dimensional
range and $p_i\in A$ is a projection such that $\tau(1-p_i)<\dt/2$ for
all $\tau\in T(A),$ $i=1,2.$ By \ref{mmm}, there is a unital \hm\, $h_{00}^{(i)}: C(X)\to p_iAp_i$
with finite dimensional range such that
\beq\label{ghth2-3}
\sup\{|\tau\circ h_{00}^{(i)}(f)-\tau\circ \phi(f)|: f\in
{\mathcal G}, \tau\in T(A)\}<\dt/2.
\eneq
By  applying \ref{zhr2}, \ref{LengthL} and the proof of
\ref{SHLn}, without loss of generality, we may assume that there
is unital \hm\, $\Phi_i: C(X)\to C([0,1],p_iAp_i)$ such that
$$
\pi_0\circ \Phi_i(f)=h_0^{(i)},\,\,\,\pi_1\circ
\Phi_i=h_{00}^{(i)}\andeqn
$$
$$
\overline{\text{Length}}(\{\pi_t\circ \Phi_i\})\le
L(X)+\ep/4.
$$
It follows that there is a unital \hm\, $H_i': C(X)\to C([0,1],A)$ such that
\beq\label{ghth2-4-1}
\pi_0\circ H_i'=h_i\oplus h_0^{(i)},\,\,\,\pi_1\circ H_i'=h_i\oplus h_{00}^{(i)}
\andeqn
\eneq
\beq\label{ghth2-4-2}
\overline{\text{Length}}(\{\pi_t\circ H_i'\})\le L(X)+\ep/4.
\eneq

Note that
$$
\sup\{|\tau\circ (h_i\oplus h_{00}^{(i)})(f)-\tau\circ \phi(f)|: \tau\in T(A)\}<\dt
$$
for all $f\in {\mathcal G}.$
By applying Theorem 3.3 of \cite{Lncd}, we obtain a unitary $u\in A$ such that
\beq\label{ghth2-4}
{\rm ad}\, u\circ (h_1\oplus h_{00}^{(1)})\approx_{\ep/4}
h_2\oplus h_{00}^{(2)}\,\,\,\text{on}\,\,\,{\mathcal F}.
\eneq
It follows from \ref{uo} that we may assume that $u\in U_0(A).$ It
follows that there is a continuous path of unitaries $\{u_t: t\in
[0,1]\}$ of $A$ for which
\beq\label{ghth2-5}
u_0=1_A, \,\,\, u_1=u\andeqn \text{Length}(\{u_t\})\le
\pi+{\ep\over{8(1+L(X)}}.
\eneq
Connecting $\pi_t\circ H_1'$ with ${\rm ad}\, u_t\circ (h_1\oplus
h_{00}^{(1)}),$ we obtain a \hm\, $H_1: C(X)\to C([0,1], A)$ such
that
\beq\label{ghth2-6}
\pi_0\circ H_1=h_1\oplus h_0^{(1)}, \,\,\,\pi_1\circ H_1\approx_{\ep/4} h_2\oplus
h_{00}^{(2)}\,\,\,{\rm on}\,\,\,{\mathcal F}
 \andeqn
\eneq
\beq\label{ghth2-7}
\overline{\text{Length}}(\{\pi_t\circ H_1\})\le L(X)
+\underline{L}_p(X)2\pi+\ep.
\eneq
Define $H_2$ by $\pi_t\circ H_2'=\pi_{1-t}\circ H_2'.$
It is clear that $H_1$ and $H_2$ meets the requirements.

\end{proof}

\vspace{0.2in}

\section{Approximate homotopy for approximately multiplicative maps}

\begin{thm}\label{SHLpiinj}
Let $X$ be a compact metric space, let $\ep>0$ and let ${\mathcal
F}\subset C(X)$ be a finite subset. Then there exits $\dt>0,$ a
finite subset ${\mathcal G}\subset C(X)$ and a finite subset
${\mathcal P}\subset \underline{K}(C(X))$ satisfying the
following:

Suppose that $A$ is a unital separable purely infinite simple
\CA\, and suppose that $\psi_1, \psi_2: C(X)\to A$ are two
$\dt$-${\mathcal G}$-multiplicative and $1/2$-${\mathcal G}$-injective \morp s. If
\beq\label{piinj1}
[\psi_1]|_{\mathcal P}=[\psi_2]|_{\mathcal P},
\eneq
then there exists an $\ep$-${\mathcal F}$-multiplicative \morp s
$H_1: C(X)\to C([0,1], A)$ such that
$$
\pi_0\circ H\approx_{\ep} \psi_1,\,\,\, \pi_1\circ H \approx_{\ep}
\psi_2 \,\,\,{\rm on}\,\,\,{\mathcal F}.
$$
Moreover,
\beq\label{piinj2}
{\rm Length}(\{\pi_t\circ H\})\le 2\pi+\ep
\eneq
\beq\label{piinj2+}
\overline{{\rm Length}}(\{\pi_t\circ H\})\le 2\pi
\underline{L}_p(X)+\ep
\eneq
\end{thm}

\begin{proof}
Let $\ep>0$ and let ${\mathcal F}\subset C(X)$ be a finite subset. It
follows from \ref{NEW}, \ref{inj} and \ref{LP} that, for a choice
of $\dt,$ ${\mathcal G}$ and ${\mathcal P}$ as in the statement of the
theorem, there exists a uniatry $u\in A$ and there exists a
projection $p\in A$ such that
\beq\label{piinj-3}
{\rm ad}\, u\circ \psi_1(f)\approx_{\ep/2} \psi_2(f)\andeqn
\psi_2(f)\approx_{\ep/2} (1-p)\psi_2(f)(1-p)+f(\xi)p
\eneq
for $f\in{\mathcal F}$ and for some $\xi\in X$). As in the proof
\ref{uo}, one finds a unitary $w\in pAp$ with $[1-p+w]=[u]$ in
$K_1(A).$ Thus (by replacing $\ep/2$ by $\ep$), we may assume that
$u\in U_0(A).$ We then, by a result of N. C. Phillips, obtain a
continuous path of unitaries $\{u_t : t\in [0,1]\}$ such that
$$
u_0=1,\,\,\, u_1=u\andeqn \text{Length}(\{u_t\})\le \pi +{\ep\over{4(1+L(X))}}.
$$
Define $H: C(X)\to C([0,1])$ by $\pi_t\circ H={\rm ad}\, u_t\circ h_1$ for $t\in [0,1].$

\end{proof}

\begin{thm}\label{SHLpi}
Let $X$ be a  a connected finite CW complex,
let $\ep>0$ and let ${\mathcal F}\subset C(X)$ be a finite subset.
Then there exits $\dt>0,$ a finite subset ${\mathcal G}\subset C(X)$
and a finite subset ${\mathcal P}\subset \underline{K}(C(X))$
satisfying the following:

Suppose that $A$ is a unital separable purely infinite simple
\CA\, and suppose that $\psi_1, \psi_2: C(X)\to A$ are two
$\dt$-${\mathcal G}$-multiplicative \morp s. If
\beq\label{ship1}
[\psi_1]|_{\mathcal P}=[\psi_2]|_{\mathcal P},
\eneq
then there exist two $\ep$-${\mathcal F}$-multiplicative \morp s
$H_1, H_2: C(X)\to C([0,1], A)$ such that
$$
\pi_0\circ H_1\approx_{\ep/3} \psi_1,\,\,\, \pi_1\circ H_1
\approx_{\ep/3} \pi_0\circ H_2 \andeqn
$$
$$
\pi_1\circ H_2\approx_{\ep/3} \psi_2 \,\,\,{\rm on}\,\,\,{\mathcal F}.
$$
Moreover,
\beq\label{ship2}
\overline{{\rm Length}}(\{\pi_t\circ H_1\})\le {\bar L}_p(X)
(1+2\pi)+\ep\andeqn
\eneq
\beq\label{ship3}
{\overline{{\rm Length}}}(\{\pi_t\circ H_2\})\le {\bar L}_p(X)+\ep.
\eneq
\end{thm}

\begin{proof}
Fix $\ep>0$ and fix a finite subset ${\mathcal F}\subset C(X).$

Let $\dt_1>0$ and let  ${\mathcal G}_1\subset C(X)$ be a finite
subset and let  ${\mathcal P}\subset \underline{K}(C(X))$ be a
finite subset required by \ref{NEW} for $\ep/4$ and ${\mathcal
F}.$ We may assume that
$$
[L_1]|_{\mathcal P}=[L_2]|_{\mathcal P}
$$
for any pair of $\dt_1$-${\mathcal G}_1$-multiplicative \morp s $L_1, L_2: C(X)\to A,$ provided that
$$
L_1\approx_{\dt_1} L_2\,\,\,\text{on}\,\,\,{\mathcal G}_1.
$$
Let $\dt_2>0$ and a finite subset ${\mathcal G}_2\subset C(X)$
required by \ref{LP} corresponding to the finite subset ${\mathcal
G}_1$ and positive number $\min\{\ep/6, \dt_1/2\}.$

Choose $\dt=\min\{\dt_1, \dt_2, \ep/6\}$ and ${\mathcal
G}={\mathcal G}_1\cup{\mathcal G}_2\cup {\mathcal F}.$

Now suppose that $\psi_1$ and $\psi_2$ are as in the statement
with above $\dt,$ ${\mathcal G}$ and ${\mathcal P}.$ We may write
$$
\psi_i\approx_{\min\{\dt_1/2,\ep/6\}} \phi_i\oplus h_i \,\,\,{\rm
on}\,\,\, {\mathcal G}_1,
$$
where $\phi_i: C(X)\to (1-p_i)A(1-p_i)$ is a $\dt_1/2$-${\mathcal G}_1$-multiplicative \morp\, and $h_i:
C(X)\to p_iAp_i$ is defined by $h_i(f)=f(\xi_i)p_i$ for all $f\in C(X),$ where
$\xi_i\in X$ is a point and $p_i$ is a non-zero projection, $i=1,2.$
It is clear that we may assume that $[p_i]=0$ in $K_0(A).$
There are $x_1,x_2,...,x_m$ in $X$ such that
\beq\label{SHLpi-1}
\max\{\|g(x_i)\|: i=1,2,...,m\}\ge (3/4)\|g\|\rforal g\in {\mathcal G}_1.
\eneq
There are non-zero mutually orthogonal projections $p_{i,1},p_{i,2},...,p_{i,m}$ in $p_iAp_i$ such that
$[p_{i,k}]=0,$ $k=1,2,...,m,$ and $p_i=\sum_{k=1}^m p_{i,k},$  $i=1,2.$

Define  $h_0^{(i)}(f)=\sum_{k=1}^mf(x_k)p_{i,k}$ for $f\in C(X).$
Note that
$$
[h_0^{(i)}]=0\,\,\,\text{in}\,\,\, KL(C(X), A).
$$

By \ref{LengthL},
 there is \hm\, $H_i': C(X)\to
C([0,1],p_iAp_i)$ such that
\beq\label{ship11}
\pi_0\circ H_i'=h_i \andeqn \pi_1\circ H_i'=h_0^{(i)}, i=1,2.
\eneq
Moreover, we can
require that
\beq\label{ship12}
{\overline{\text{Length}}}(\{\pi_t\circ H_i'\})\le L(X,
\xi_i)+\ep/4,\,\,\, i=1,2.
\eneq
It follows from \ref{SHLpi-1} that $\phi_i\oplus h_0^{(i)}$ is
$\dt_1$-${\mathcal G}_1$-multiplicative and $1/2$-${\mathcal G}_1$-injective.

Thus, by \ref{NEW} there exists a unitary $u\in A$ such that
\beq\label{ship13}
{\rm ad}\, u\circ (\phi_1\oplus h_0^{(1)})\approx_{\ep/3} \phi_2\oplus
h_0^{(2)}\,\,\,\text{on}\,\,\, {\mathcal F}.
\eneq
It follows from \ref{uo} that we may assume that $u\in U_0(A).$
 Thus we obtain a continuous rectifiable
path $\{u_t: t\in [0,1]\}$ of $A$ such that
\beq\label{ship15}
u_0=1_A,\,\,\, u_1=u\andeqn
\eneq
$$
\text{Length}(\{u_t\})\le \pi+{\ep\over{4(1+L(X))}}.
$$
Define $H_1: C(X)\to C([0,1], A)$ by
$$
\pi_t\circ H_1=\begin{cases} \phi_1\oplus \pi_{2t}\circ H_1' &\text{if $t\in [0,1/2]$;}\\
            {\rm ad}\, u_{2(t-1/2)}\circ (\phi_1\oplus h_0^{(1)})\}) &\text{if $t\in (1/2,1]$.}
            \end{cases}
$$
Then
\beq\label{ship17}
\pi_0\circ H_1= \phi_1\oplus h_1,\,\,\,\pi_1\circ
H_1\approx_{\ep/3} \phi_2\oplus h_0^{(2)}\,\,\,\text{on}\,\,\, {\mathcal
F}.
\eneq
Moreover,
\beq\label{ship18}
{\overline{\text{Length}}}(\{\pi_t\circ H_1\})\le
L(X,\xi_1)(1+2\pi)+\ep.
\eneq

We then define $H_2: C(X)\to C([0,1],A)$ by $\pi_t\circ H_2=\pi_{1-t}\circ H_2''.$
We see that $H_1$ and $H_2$ meet the requirements of the theorem.
\end{proof}

\begin{thm}\label{SHL}
Let $X$ be a metric space which is a connected finite CW complex,
let $\ep>0$ and let ${\mathcal F}\subset C(X)$ be a finite subset.
Then there exists $\dt>0,$  a finite subset ${\mathcal G}\subset C(X)$
and a finite subset ${\mathcal P}\subset \underline{K}(C(X))$
satisfying the following:

Suppose that $A$ is a unital separable simple \CA\, with tracial
rank zero and suppose $\psi_1, \psi_2: C(X)\to A$  are two unital
$\dt$-${\mathcal G}$-multiplicative \morp s. If
\beq\label{shl1}
[\psi_1]|_{\mathcal P}=[\psi_2]|_{\mathcal P},
\eneq
then there exist two  $\ep$-${\mathcal F}$-multiplicative \morp s
$H_1, H_2:C(X)\to C([0,1], A)$ such that
\beq\label{shl2}
\pi_0\circ H_1\approx_{\ep/3} \psi_1, \pi_1\circ
H_1\approx_{\ep/3} \pi_1\circ H_2, \andeqn
\eneq
\beq\label{shl3}
\pi_0\circ H_2\approx_{\ep/3} \psi_2 \,\,\,\text{on}\,\,\, {\mathcal
F}.
\eneq
Moreover,
\beq\label{shl4}
{\overline{\rm{Length}}}(\{\pi_t\circ H_1\})\le L(X)+2\pi {\bar
L}_p(X)+\ep \andeqn
\eneq
\beq\label{shl5}
{\overline{\rm{Length}}}(\{\pi_t\circ H_2\})\le L(X)+\ep.
\eneq

\end{thm}

\begin{proof}
Fix $\ep>0$ and fix a finite subset ${\mathcal F}\subset C(X).$ Let
$\eta=\sigma_{X,\ep/32,{\mathcal F}}.$ Suppose that
$\{x_1,x_2,...,x_m\}$ is $\eta/2$-dense in $X.$ There is $s\ge 1$
such that $O_i\cap O_j=\emptyset,$ where
$$
O_i=\{x\in X: \di(x, x_i)<\eta/2s\},\,\,\,i=1,2,...,m.
$$
Put $\sigma={1\over{(2sm+1)}}.$

Let $\dt_1$ (in place of $\dt$), $\gamma>0,$ ${\mathcal G}_1$ (in
place of ${\mathcal G}$) and ${\mathcal P}\subset
\underline{K}(C(X))$ be as required by Theorem 4.6 of \cite{Lncd}
for $\ep/6$ (in place of $\ep$) and ${\mathcal F}$ and $\sigma.$
We also assume that $2m\sigma\eta<1-\gamma.$ We may assume that
$\dt<\ep$ and ${\mathcal F}\subset {\mathcal G}.$ Furthermore, we
assume that
$$
[L_1]|_{\mathcal P}=[L_2]|_{\mathcal P}
$$
for any pair of $\dt_1$-${\mathcal G}_2$-multiplicative \morp s $L_1, L_2: C(X)\to A,$ provided that
$$
L_1\approx_{\dt_1} L_2\,\,\,\text{on}\,\,\,{\mathcal G}_1.
$$

Let $\dt$ and ${\mathcal G}$ be as in \ref{ndig} corresponding
$\min\{\dt_1/2, \ep/3\}$ ( instead of $\ep$) and ${\mathcal G}_1$
(instead of ${\mathcal G}$) and $\gamma/2.$

Let $A$ be a unital separable simple \CA\, with tracial rank zero
and $\psi_1, \psi_2: C(X)\to A$ be as described. By \ref{ndig}, we
may write
\beq\label{shl10}
\psi_i\approx_{\min\{\dt_1/2, \ep/6\}} \phi_i\oplus
h_i\,\,\,\,\text{on}\,\,\, {\mathcal G}_1\,(i=1,2),
\eneq
where $\phi_i: C(X)\to (1-p_i)A(1-p_i)$ is a unital \morp\, and
$h_i: C(X)\to p_iAp_i$ is a unital \hm\, with finite dimensional
range.

 Note that we assume that $p_1$ and $p_2$ are unitarily equivalent
 and
\beq\label{shl11}
\tau(1-p_i)<\gamma/2\, \rforal\, \tau\in T(A),\,\,\,i=1,2.
\eneq

There are non-zero mutually orthogonal projections
$p_{i,1},p_{i,2},...,p_{i,m}\in p_iAp_i$ such that
$\sum_{k=1}^mp_{i,k}=p_i$ and
\beq\label{shl12-}
\tau(p_{i,k})\ge {1-\gamma/2\over{m+1}},\,\,\,k=1,2,...,m\andeqn
i=1,2.
\eneq
Since $[p_1]=[p_2],$ we also require that
\beq\label{shl12-2}
[p_{1,k}]=[p_{2,k}],\,\,\,k=1,2,...,m.
\eneq
Define
$h_0^{(i)}: C(X)\to p_iAp_i$ by
$$
h_0^{(i)}(f)=\sum_{k=1}^m f(x_k)p_{i,k}\rforal f\in C(X),\,\,\,i=1,2.
$$
By \ref{zhr2}, without loss of generality, we may assume that
there is a unital commutative finite dimensional \SCA\,
$B_i\subset p_iAp_i$ such that $h_i(C(X)), h_0^{(i)}(C(X))\subset
B_i,$ $i=1,2.$

It follows from \ref{LengthL} that there is \hm\, $H_i': C(X)\to p_iAp_i$
such that
\beq\label{shl12}
\pi_0\circ H_i'=h_i,\,\,\,
 \pi_1\circ H_i'=h_0^{(i)}\andeqn
\eneq
\beq\label{shl13}
{\overline{\text{Length}}}(\{\pi_t\circ H_i'\})\le
L(X)+\ep/4,\,\,\, i=1,2.
\eneq
Note that, since $X$ is path connected,
\beq\label{shl13+2}
[h_1]=[h_0^{(1)}]=[h_0^{(2)}]=[h_2]\,\,\,\text{in}\,\,\, KL(C(X),A).
\eneq
It follows from (\ref{shl1}), (\ref{shl10})  as well as
(\ref{shl13+2}) that
\beq\label{shl13+3}
[\psi_1\oplus h_0^{(1)}]|_{\mathcal P}=[\psi_2\oplus h_0^{(2)}]|_{\mathcal P}
\eneq
We also have that
\beq\label{shl13+4}
\mu_{\tau\circ (\psi_i\oplus h_0^{(i)})}(O_k)\ge \tau(p_{i,k})\ge \sigma\eta
\eneq
for all $\tau\in T(A),$ $k=1,2,...,m,$ $i=1,2.$ Moreover,
\beq\label{shl14-1}
|\tau(\phi_1\oplus h_0^{(1)})(f)-\tau(\phi_2\oplus h_0^{(2)})(f)|<\gamma
\eneq
for all $\tau\in T(A)$ and all $f\in {\mathcal G}_1.$

Thus, by Theorem 4.6 of \cite{Lncd} there exists a unitary $u\in A$
such that
\beq\label{shl18}
{\rm ad}\, u\circ (\phi_1\oplus h_0^{(1)})\approx_{\ep/6} \phi_2\circ
h_0^{(2)}\,\,\,\text{on}\,\,\, {\mathcal F}.
\eneq
By \ref{uo},  we may assume
that $u\in U_0(A).$ Thus we obtain a continuous rectifiable path
$\{u_t: t\in [0,1]\}$ of $A$ such that
\beq\label{shl20}
u_0=1_A,\,\,\, u_1=u\andeqn
\eneq
\beq\label{shl21}
\text{Length}(\{u_t\})\le \pi+{\ep\over{4(1+L(X))}}.
\eneq

Define $H_1: C(X)\to C([0,1], A)$ by
$$
\pi_t\circ H_1=\begin{cases} \phi_1\oplus \pi_{2t}\circ H_1' &\text{if $t\in [0,1/2]$;}\\
            {\rm ad}\, u_{2(t-1/2)}\circ (\phi_1\oplus h_0^{(1)})\}) &\text{if $t\in (1/2,1]$.}
            \end{cases}
$$
Then
\beq\label{shl21+}
\pi_0\circ H_1= \phi_1\oplus h_1,\,\,\,\pi_1\circ
H_1\approx_{\ep/3} \phi_2\oplus h_0^{(2)}\,\,\,\text{on}\,\,\, {\mathcal
F}.
\eneq
Moreover,
\beq\label{shl22}
{\overline{\text{Length}}}(\{\pi_t\circ H_1\})\le
L(X)+ 2\pi L(X,\xi_1)+\ep.
\eneq

We then define $H_2: C(X)\to C([0,1],A)$ by $\pi_t\circ H_2=\pi_{1-t}\circ H_2''.$
We see that $H_1$ and $H_2$ meet the requirements of the theorem.

\end{proof}

\begin{rem}

{\rm
It should be noted that results in this section are not generalization of those in the previous section.
It is important to know that $\pi_t\circ H_i$ in Theorem \ref{SHLpi} and those
in Theorem \ref{SHL} are not \hm s, while $\pi_t\circ H_i$ in \ref{GHTh2} and \ref{GHP3} are  unital \hm s.
Furthermore, $\pi_t\circ H$ in \ref{SHLn} and \ref{GHT2} are unital monomorphisms.


 }

\end{rem}

\chapter{Super Homotopy}

\setcounter{section}{14}
\section{Super Homotopy Lemma --- purely infinite case}


In this section  and the next we study the so-called Super
Homotopy Lemma of \cite{BEEK} for higher dimensional spaces.

\begin{thm}\label{SHLoldpi}
Let $X$ be a path connected metric space and let
 ${\mathcal F}\subset C(X)$ be a finite subset.

Then, for any $\ep>0,$ there exists $\dt>0,$ a finite subset ${\mathcal G}\subset
C(X)$ and a finite subset ${\mathcal P}\subset \underline{K}(C(X))$ satisfying the following:

Suppose that $A$ is a unital separable purely infinite simple
\CA\,  and suppose that $h_1, h_2: C(X)\to A$ are two unital
monomorphisms such that
\beq\label{shpip1}
[h_1]=[h_2]\,\,\,\text{in}\,\,\, KL(C(X), A).
 \eneq

 If there
are unitaries $u, v\in A$ such that

\beq\label{shpip2}
\|[h_1(f), u]\|<\dt,\,\,\,\|[h_2(f), v]\|<\dt\,\rforal\, f\in
{\mathcal G}\andeqn
\eneq
\beq\label{shpip2+}
 \rm{Bott}(h_1,u)|_{\mathcal P}=\rm{Bott}(h_2,v)|_{\mathcal P},
\eneq
then there exists  a unital  monomorphism  $H: C(X)\to C([0,1],
A)$ and there exists a continuous rectifiable path $\{u_t: t\in
[0,1]\}$
such that, for each $t\in [0,1],$ $\pi_t\circ H$ is a unital monomorphism,
\beq\label{shpip3}
u_0=u,\,\,\, u_1=v,\,\,\, \pi_0\circ
H\approx_{\ep} h_1\andeqn \pi_1\circ H\approx_{\ep} h_2\,\,\,\text{on}\,\,\,{\mathcal F}
\eneq
 and
\beq\label{shpip3+}
\|[u_t, \pi_t\circ H(f)]\|<\ep
\eneq
for all $t\in [0,1]$ and for all $f\in {\mathcal F}.$ Moreover,

\beq\label{shpip4}
\rm{Length}(\{\pi_t\circ H\}) &\le &
2\pi+{\ep\over{8(1+L(X)}},\\\label{shpip4+}
\overline{\rm{Length}}(\{\pi_t\circ H\}) &\le & 2\pi
\underline{L}_p(X)
+\ep\andeqn\\\
\label{shpip4++} {\rm{Length}}(\{u_t\}) &\le & 4\pi+\ep.
\eneq
\end{thm}

\begin{proof}
Let $\ep>0$ and let ${\mathcal F}\subset C(X)$ be a finite subset. We
may assume that $1_{C(X)}\in {\mathcal F}.$ Put
$$
\ep_1=\min\{\ep/8, 2\sin (\ep/16)\}.
$$
Let ${\mathcal G}_0\subset C(X\times S^1)$ be a finite subset, ${\mathcal
P}_1\subset \underline{K}(C(X\times S^1))$ be a finite subset,
$\dt_1>0$  be a positive number associated with $\ep_1/2,$ and ${\mathcal
F}\otimes S\subset C(X\times S^1)$ as required by \ref{NEW}, where
$S=\{1_{C(S^1)}, z\}\subset C(S^1).$ By choosing even smaller
$\dt_1,$ we may assume that ${\mathcal G}_0={\mathcal G}_1\otimes S.$ Let
${\mathcal P}_2=({\rm id}-\hat{{\boldsymbol{ \bt}}})({\mathcal P}_1)\cup \hat{{\boldsymbol{ \bt}}}({\mathcal
P}_1).$

We may also assume that $\dt_1<\ep_1$ and ${\mathcal F}\subset {\mathcal
G}_1.$ Moreover, we may further assume that
$$
[L_1]|_{{\mathcal P}_1\cup {\mathcal P}_2}=[L_2]|_{{\mathcal P}_1\cup {\mathcal
P}_2}
$$
for any pair of $\dt_1$-${\mathcal G}_0$-multiplicative unital \morp s $L_1,L_2: C(X\otimes S^1)\to A,$
provided that
$$
L_1\approx_{\dt_1} L_2\,\,\,\text{on}\,\,\,{\mathcal G}_0.
$$

Let $\dt_2>0$ and let ${\mathcal G}_2$ be a finite subset required by \ref{pihf2} corresponding to
${\mathcal G}_1$ (in place of ${\mathcal F}$) and $\dt_1/2$ (in place of $\ep$).
We may assume that $\dt_2<\dt_1/2$ and ${\mathcal G}_1\subset {\mathcal G}_2.$

Let $\dt>0$ and let ${\mathcal G}\subset C(X)$ (in place of ${\mathcal
F}_1$) be a finite subset required in \ref{appn} for $\dt_2/2$ (
in place of $\ep$) and ${\mathcal G}_2$ (in place of ${\mathcal F}$).

Let ${\mathcal P}\in \underline{K}(C(X))$ be a finite subset such that
${\boldsymbol{ \bt}}({\mathcal P})\supset \hat{{\boldsymbol{ \bt}}}({\mathcal P}_1).$  Now suppose that
$h_1,h_2$ and $u$ and $v$ are in the statement of the theorem with
the above $\dt,$ ${\mathcal G}$ and ${\mathcal P}.$

By \ref{appn}, there are $\dt_2/2$-${\mathcal G}_2\otimes S$-multiplicative unital \morp s $\psi_1, \psi_2: C(X)\to
A$ such that
\beq\label{suphpi-1}
\|\psi_1(f\otimes g)-h_1(f)g(u)\|
&<& \dt_2/2\,\,\,\andeqn\\\label{suphpi-2}
\|\psi_2(f\otimes g)-h_2(f)g(v)\|&<&\dt_2/2
\eneq
for all $f\in {\mathcal G}_2$ and $g\in S.$

By applying \ref{pihf2}, we obtain two continuous paths of
unitaries $\{w_{i,t}: t\in [0,1]\}$ in $A$ and two unital
monomorphisms $\Phi_i: C(X\times S^1) \to A$ such that
\beq\label{suphpi-3}
w_{1,0}=u,\,\,\, w_{2,0}=v,\,\,\, w_{i,1}=\Phi_i(1\otimes z),
\eneq
\beq\label{suphpi-4}
\|\Phi_i(f\otimes 1)-h_i(f)\|<\dt_1/2 \,\,\andeqn\|[h_i(f),
w_{i,t}]\|<\dt_1/2
\eneq
for all $ f\in {\mathcal G}_1$ and all $t\in [0,1],$ $i=1,2.$  Moreover,
\beq\label{suphpi-5}
\text{Length}(\{w_{i,t}\})\le \pi+\ep/8.
\eneq
It follows from (\ref{shpip1}), (\ref{shpip2+}), (\ref{suphpi-1}), (\ref{suphpi-2}) and \ref{Khp}
that
\beq\label{suphpi-6}
[\Phi_1]|_{{\mathcal P}_2}=[\Phi_2]|_{{\mathcal P}_2}.
\eneq
Since $\Phi_1$ and $\Phi_2$ are \hm s, we have
\beq\label{suphpi-6+}
[\Phi_1]|_{{\mathcal P}_1}=[\Phi_2]|_{{\mathcal P}_1}.
\eneq
 It follows from
\ref{NEW} that there is a unitary $U\in A$ such that
\beq\label{suphpi-7}
{\rm ad}\, U\circ \Phi_1\approx_{{\ep_1\over{2}}} \Phi_2\,\rforal
f\in {\mathcal F}\otimes S.
\eneq
By \ref{uo}, we may assume that $U\in U_0(A)$ ( by replacing $\ep_1/2$ by $\ep_1$).
It follows that there is a continuous path of unitaries
$\{U_t: t\in [0,1]\}$ such that
\beq\label{suphpi-8}
U_0=1,\,\,\, U_1=U\andeqn \text{Length}(\{U_t\})\le
\pi+{\ep\over{8(1+L(X))}}.
\eneq
By (\ref{suphpi-7}), there exists $a\in A_{s.a}$ such that
$\|a\|\le \ep/8$ such that
\beq\label{suphpi-9}
({\rm ad}\, (U\circ \Phi_1(1\otimes z))^*\Phi_2(1\otimes
z)=\exp(ia).
\eneq
Moreover
\beq\label{suphpi-10}
\|\exp(ita)-1\|<\ep/8\rforal t\in [0,1].
\eneq
Put $z_t={\rm ad}\,U\circ \Phi_1(1\otimes z)\exp(ita).$ Then
\beq\label{suphpi-11}
z_0={\rm ad}\,U\circ \Phi_1(1\otimes z),\,\,\, z_1=\Phi_2(1\otimes
z)\andeqn
\eneq
\beq\label{suphpi-12}
\text{Length}(\{z_t\})\le \ep/8.
\eneq

Define $u_t$ as follows
$$
u_t=\begin{cases} w_{1,3t} &\text{if $ t\in [0,1/3]$;}\\
                               U_{3(t-1/3)}^*\Phi_1(1\otimes z)U_{3(t-1/3)} &\text{if $t\in (1/3,2/3]$;}\\
                               z_{6(t-2/3)} &\text{if $t\in
                               (2/3,5/6]$;}\\
                               w_{2,3(2/3-t)} & \text{if $t\in (5/6,1]$}
                               \end{cases}
                               $$
     Define $H: C(X\otimes S^1)\to A$ by
     $$
        \pi_t\circ H=\begin{cases} h_1 &\text{if $t\in [0,1/3]$;}\\
                              {\rm ad}\, U_{3(t-1/3)}\circ h_1 & \text{if $t\in (1/3, 2/3]$;}\\
                                           {\rm ad}\, U_1\circ h_1 &\text{if $t\in (2/3,1]$}
                                           \end{cases}
                                           $$

We compute that
\beq\label{suphpi-13}
\text{Length}(\{u_t\})\le
2\pi+2\ep/8+2\pi+2\ep/8+\ep/8<4\pi+\ep\andeqn
\eneq
\beq\label{suphpi-14}
\text{Length}(\{\pi_t\circ H\})\le
2\pi+{\ep\over{4(1+L(X,\xi_X)}}.
\eneq
In particular, by \ref{Lnprop1},
\beq\label{suppi-14+}
\overline{\text{Length}}(\{\pi_t\circ H\})\le 2\pi
\underline{L}_p(X)+\ep/4.
\eneq

 It follows from (\ref{suphpi-4})
and (\ref{suphpi-10}) that
\beq\label{suphpi-15}
\|[\pi_t\circ H(f), u_t]\|<2(\dt_1/2)+\dt_1/2+2(\ep/8)<\ep\rforal
f\in {\mathcal F}.
\eneq

We see that $\{u_t:t\in [0,1]\}$ and $H$ meet the requirements.

\end{proof}

\begin{thm}\label{suppiT}
Let $X$ be a connected metric space and let
 ${\mathcal F}\subset C(X)$ be a finite subset.

Then, for any $\ep>0,$ there exists $\dt>0,$ a finite subset ${\mathcal G}\subset
C(X)$ and a finite subset ${\mathcal P}\subset \underline{K}(C(X))$ satisfying the following:

Suppose that $A$ is a unital separable purely infinite simple
\CA\,  and suppose that $h_1, h_2: C(X)\to A$ are two unital
\hm s such that
\beq\label{suppiT-1}
[h_1]=[h_2]\,\,\,\text{in}\,\,\, KL(C(X), A).
 \eneq

 If there
are unitaries $u, v\in A$ such that

\beq\label{suppiT-2}
\|[h_1(f), u]\|<\dt,\,\,\,\|[h_2(f), v]\|<\dt\,\tforal\, f\in
{\mathcal G}\andeqn
\eneq
\beq\label{suppiT-3}
 \rm{Bott}(h_1,u)|_{\mathcal P}=\rm{Bott}(h_2,v)|_{\mathcal P},
\eneq
then there exist  two  unital  \hm s $H_1, H_2: C(X)\to
C([0,1], A)$ and a continuous rectifiable path $\{u_t:
t\in [0,2]\}$
such that, for each $t\in [0,1],$ $\pi_t\circ H$ is a unital \hm,
\beq\label{suppiT-4}
&&\hspace{0.4in}u_0= u,\,\,\, u_2=v,\\
&& \hspace{-0.4in}\pi_0\circ
H_1\approx_{\ep} h_1,\,\,\, \pi_1\circ H_1\approx_{\ep} \pi_0\circ
H_2 \andeqn \pi_1\circ H_1\approx_{\ep}
h_2\,\,\,\text{on}\,\,\,{\mathcal F},
\eneq
\beq\label{suppiT-5}
\|[u_t, \pi_t\circ H_1(f)]\| &<&\ep \andeqn\\
 \|[u_{1+t}, \pi_t\circ H_2(f)]\| &<&\ep \rforal t\in [0,1]
\eneq
and for all $f\in {\mathcal F}.$ Moreover,

\beq\label{suppiT-6}\nonumber
\overline{\rm{Length}}(\{\pi_t\circ H_1: t\in [0,1/2]\}) &\le & {\bar L}_p(X)+\ep/2,\\
\rm{Length}(\{\pi_t\circ H_1: t\in [1/2,1]\}) &\le & 2\pi+{\ep\over{8(1+L(X))}},\\
\overline{\rm{Length}}(\{\pi_t\circ H_1\})  &\le& {\bar
L}_p(X)+2\pi \underline{L}_p(X)
+\ep,\\
\overline{\rm{Length}}(\{\pi_t\circ H_2\}) &\le & {\bar L}_p(X)+\ep \andeqn\\
\rm{Length}(\{u_t\}) &\le & 4\pi+\ep.
\eneq
\end{thm}

\begin{proof}
Let $\ep>0$ and let ${\mathcal F}\subset C(X)$ be a finite subset. We
may assume that $1_{C(X)}\in {\mathcal F}$ and ${\mathcal F}$ is in the unit ball of $C(X).$  Put
$$
\ep_1=\min\{\ep/8, \sin (\ep/16)\}.
$$
Let ${\mathcal G}_0\subset C(X\times S^1)$ be a finite subset, ${\mathcal
P}_1\subset \underline{K}(C(X\times S^1))$ be a finite subset,
$\dt_1>0$  be  a positive number associated with $\ep_1/2$ and ${\mathcal
F}\otimes S\subset C(X\times S^1)$ as required by \ref{NEW}, where
$S=\{1_{C(S^1)}, z\}\subset C(S^1).$ By choosing even smaller
$\dt_1,$ we may assume that ${\mathcal G}_0={\mathcal G}_1\otimes S.$ Let
${\mathcal P}_2=({\rm id}-\hat{{\boldsymbol{ \bt}}})({\mathcal P}_1)\cup \hat{{\boldsymbol{ \bt}}}({\mathcal
P}_1).$

We may also assume that $\dt_1<\ep_1$ and ${\mathcal F}\subset {\mathcal
G}_1.$ Moreover, we may further assume that
$$
[L_1]|_{{\mathcal P}_1\cup {\mathcal P}_2}=[L_2]|_{{\mathcal P}_1\cup {\mathcal
P}_2}
$$
for any pair of $\dt_1$-${\mathcal G}_0$-multiplicative unital \morp s
$L_1,L_2: C(X\otimes S^1)\to A,$
provided that
$$
L_1\approx_{\dt_1} L_2\,\,\,\text{on}\,\,\,{\mathcal G}_0.
$$

Let $\dt_2>0$ and let ${\mathcal G}_2\subset C(X)$ be a finite subset
required by \ref{pcldig} for $\dt_1/4$ (in place of $\ep$) and
${\mathcal G}_0.$ We may assume that $\dt_2<\dt_1/4$ and ${\mathcal
G}_0\subset {\mathcal G}_2.$

Let $\dt>0$ and let ${\mathcal G}\subset C(X)$ ( in place of ${\mathcal
F}$) be a finite subset required in \ref{appn} for $\dt_2/2$ ( in
place of $\ep$) and ${\mathcal G}_2$ (in place of ${\mathcal F}$).

Let ${\mathcal P}\in \underline{K}(C(X))$ be a finite subset such that
${\boldsymbol{ \bt}}({\mathcal P})\supset
\hat{{\boldsymbol{ \bt}}}({\mathcal P}_1).$

Now suppose that $h_1,h_2$ and $u$ and $v$ are in the statement of
the theorem with the above $\dt,$ ${\mathcal G}$ and ${\mathcal
P}.$ Suppose that the spectrum of $h_i$ is $F_i,$ $i=1,2.$ Let
$s_i:  C(X)\to C(F_i)$ be the quotient map, $i=1,2.$ Define ${\bar
h}_i: C(F_i)\to A$ such that ${\bar h}_i\circ s_i=h_i,$ $i=1,2.$

Let $\eta_i>0,$ let ${\mathcal G}_3^{(i)}\subset C(F_i)$ be a finite
subset and let ${\mathcal P}_2^{(i)}\subset \underline{K}(C(F_i))$ be
a finite subset required by \ref{presemi} corresponding to
$s_i({\mathcal G}_1)$ (in place of ${\mathcal F}$) and $\dt_1/4$ which
work for $F_i,$ $i=1,2.$ We may assume that $\eta_i<\dt_1/4$ and
$s_i({\mathcal G}_2)\subset {\mathcal G}_3^{(i)},$ $i=1,2.$

We may further assume that $
[L^{(i)}]|_{{\mathcal P}_2^{(i)}}$ is well defined for any
$\eta_i$-${\mathcal G}_3^{(i)}$-multiplicative \morp\, $L^{(i)}:
C(F_i)\to A.$
 Put $\eta=\min\{\eta_1, \eta_2\}$ and let ${\mathcal
G}_3\subset C(X)$ be a finite subset such that $s_i({\mathcal
G}_3)\supset {\mathcal G}_3^{(i)},$ $i=1,2.$

By \ref{appn},
there are $\dt_2/2$-${\mathcal G}_2\otimes S$-multiplicative unital
\morp s $\psi_1, \psi_2: C(X)\to
A$ such that
\beq\label{suppiT-01}
&&\|\psi_1(f\otimes g)-h_1(f)g(u)\|<\dt_2/2
\eneq
for all $f\in {\mathcal G}_2\andeqn g\in S,$ and
\beq
\label{suppiT-02}
&&\|\psi_2(f\otimes g)-h_2(f)g(v)\|<\dt_2/2
\eneq
for all $f\in {\mathcal G}_2\andeqn g\in S.$

By applying \ref{pcldig} (and see the remark \ref{Rdig}), there is a nonzero projection $p_i\in A$
and $\xi_i\in F_i\subset X$ and $t_i\in S^1$ ($i=1,2$) such that
\beq\label{suppiT-03}
\|h_1(f)g(u)-((1-p_1)h_1(f)g(u)(1-p_1)+f(\xi_1)g(t_1)p_1)\|<\dt_1/4\\
\|h_2(f)g(v)-((1-p_2)h_2(f)g(v)(1-p_2)+f(\xi_2)g(t_1)p_2)\|<\dt_1/4\label{suppiT-03+}
\eneq
for all $f\in {\mathcal G}_1$ and $ g\in S,$
\beq\label{suppiT-04}
\|h_i(f)-((1-p)h_i(f)(1-p)+f(\xi_i)p_i)\|<\eta/2\rforal f\in {\mathcal
G}_3,
\eneq
$i=1,2.$ Moreover, we may assume that $(1-p_i){\bar h}_i(1-p_i)$
is $\eta_i$-${\mathcal G}_3^{(i)}$-multiplicative, $i=1,2.$
Furthermore, there are unitaries $w_i\in (1-p_i)A(1-p_i)$
($i=1,2$) such that
\beq\label{suppiT-04+}
\|u-w_1\oplus t_1p_1\|<\dt_1/2\andeqn \|v-w_2\oplus
t_2p_2\|<\dt_1/2.
\eneq
(Note that $\eta_i$ and ${\mathcal G}_3^{(i)}$ can be chosen, which depend on $F_i$ among other things,
after $\dt$ and ${\mathcal G}$ are chosen, thanks to \ref{pcldig}.)

By replacing $p_i$ by a non-zero sub-projection, we may assume
that $[p_i]=0$ in $K_0(A)$ (see \ref{Rdig}) and $(1-p_i){\bar h}_i(1-p_i)$ is
$1/2$-${\mathcal G}_3^{(i)}$-injective.

It follows that
\beq\label{suppiT-05}
[(1-p){\bar h}_i(1-p)]|_{{\mathcal P}_2^{(i)}}=[{\bar h}_i]|_{{\mathcal
P}_2^{(i)}},\,\,\,i=1,2.
\eneq
It follows from \ref{presemi} that there is a unital \hm\, $h_i':
C(X)\to (1-p_i)A(1-p_i)$ (with spectrum $F_i$) such that
\beq\label{suppiT-06}
h_i'\approx_{\dt_1/4}(1-p_i)h_i(1-p_i)\,\,\,\text{on}\,\,\, {\mathcal
G}_1.
\eneq

Let $\{x_1,x_2,...,x_m\}\subset X$ and
$\{\zeta_1,\zeta_2,...,\zeta_k\}\subset S^1$ such that
\beq\label{suppiT-07}
\max\{\|f(x_j\times \zeta_l)\|: 1\le j\le m, 1\le l\le k\}\ge
(3/4) \|f\|\rforal f\in {\mathcal G}_0.
\eneq
There are mutually orthogonal non-zero projections
$\{e_{j,l}^{(i)}: 1\le j\le m \andeqn 1\le l\le k\}$ in $p_iAp_i$
such that
$$
\sum_{j,l}e_{j,l}^{(i)}=p_i\andeqn
[e_{j,l}^{(i)}]=0\,\,\,\text{in}\,\,\,K_0(A),
$$
$i=1,2.$ Define $h_{00}^{(i)}: C(X\times S^1)\to p_iAp_i$ by
\beq\label{suppit-08}
h_{00}^{(i)}(f)=\sum_{j,l} f(x_j\times
\zeta_l)e_{j,l}^{(i)}\rforal f\in C(X\times S^1).
\eneq
It follows from \ref{LengthL} that there is a unital \hm\,
$\Phi_i: C(X\times S^1)\to C([0,1], p_iAp_i)$ such that
\beq\label{suppit-09}
\pi_0\circ \Phi_i(f)=f(\xi_i\times t_i)p_i \andeqn \pi_1\circ
\Phi_i(f)= h_{00}^{(i)}(f)
\eneq
for all $f\in C(X\times S^1)$ and
\beq\label{suppit-10}
\overline{\text{Length}}(\{\pi_t\circ \Phi_i|_{C(X)\otimes
1}\})&\le& {\bar L}_p(X)+\ep/4\\
\text{Length}(\{\pi_t\circ \Phi_i(1\otimes z)\})&\le & \pi,
\eneq
$i=1,2.$

 Now $(1-p_1)\psi_1(1-p_1)\oplus h_{00}^{(1)}$ and $(1-p_2)\psi_2(1-p_2)\oplus h_{00}^{(2)}$ are both
 $\dt_1/2$-${\mathcal G}_1\otimes S$-multiplicative, by (\ref{suppiT-01}), (\ref{suppiT-02}), (\ref{suppiT-03})
 and (\ref{suppiT-03+}), and are both $1/2$-${\mathcal G}_0$-injective, by (\ref{suppiT-07}).
 Put $L_i'=(1-p_1)\psi_1(1-p_1)\oplus h_{00}^{(1)},$ $i=1,2.$

 Moreover, since $[p_i]=0$ in $K_0(A)$ and $[h_{00}^{(i)}]=0$ in $KL(C(X), A),$ by the assumption  of $\dt_1$ and
 ${\mathcal G}_0,$ and by the assumption (\ref{suppiT-1}) and (\ref{suppiT-3}),
 $$
 [L_1']|_{{\mathcal P}_1}=[L_2']|_{{\mathcal P}_1} .
 $$
 It follows from \ref{NEW} that there is a unitary $U\in A$ such that
 \beq\label{suppit-11}
 {\rm ad}\, U\circ L_1'\approx_{\ep_1/2} L_2'\,\,\,\text{on}\,\,\, {\mathcal F}\otimes S.
 \eneq
 By \ref{uo}, without loss of generality, we may assume that $U\in U_0(A),$ by replacing $\ep_1/2$ by
 $\ep_1$ above.
 It follows that we obtain a continuous path of unitaries \\
  $\{U_t: t\in [0,1]\}$ such that
 \beq\label{suppit-12}
 U_0=1,\,\,\, U_1=U\andeqn \text{Length}(\{U_t\})\le
\pi+{\ep\over{8(1+L(X))}}.
\eneq
By (\ref{suppit-11}), there exists $a\in A_{s.a}$ such that
$\|a\|\le \ep/8$ such that
\beq\label{suppit-13}
(U^*( w_1\oplus h_{00}^{(1)}(1\otimes z))U)^*(w_2\oplus
h_{00}^{(2)}(1\otimes z)) =\exp(ia).
\eneq
Moreover
\beq\label{suppit-14}
\|\exp(ita)-1\|<\ep/8\rforal t\in [0,1].
\eneq
Put $z_t=U^*(w_1\oplus h_{00}^{(1)}(1\otimes z))U\exp(ita).$ Then
\beq\label{suppit-15}
z_0=U^*(w_1\oplus (1\otimes z))U,\,\,\, z_1=w_2\oplus h_{00}^{(2)}(1\otimes
z)\andeqn
\eneq
\beq\label{suphpi-16}
\text{Length}(\{z_t\})\le \ep/8.
\eneq
Similarly, by (\ref{suppiT-04+}), we obtain two continuous paths
of unitaries $\{V_{1,t}:t\in  [0,1]\}$ of $A$ such that
\beq\label{suppit-17}
V_{1,0}=u,\,\,\, V_{2,0}=v,\,\,\,V_{i,1}=w_i\oplus t_i
p_i,\,\,\,\text{Length}(\{V_{i,t}\})\le \ep/8\andeqn
\eneq
\beq\label{suppit-18}
\|V_{i,t}-1\|<\ep/8.
\eneq

Define $u_t$ as follows
$$
u_t=\begin{cases}   V_{1,4t} &\text{if $t\in [0,1/4]$;}\\
w_1\oplus \pi_{4(t-1/4)}\Phi_1(1\otimes z) &\text{if $ t\in (1/4,1/2]$;}\\
                               U_{2(t-1/2)}^*(w_1\oplus \Phi_1(1\otimes z))U_{2(t-1/2)} &\text{if $t\in (1/2,1]$;}\\
                               z_{2(t-1)} &\text{if $t\in
                               (1,3/2]$;}\\
                               w_2\oplus \pi_{4(7/4-t)}\circ \Phi_2(1\otimes z)& \text{if $t\in (3/2, 7/4]$;}\\
                               V_{2,4(2-t)} &\text{if $t\in (7/4,
                               2]$}\,.
                               \end{cases}
                               $$
     Define $H_1: C(X)\to A$ by
     $$
        \pi_t\circ H_1=\begin{cases} h_1'\oplus \pi_0\circ \Phi_1|_{C(X)\otimes 1}  &\text{if $t\in [0,1/4]$;}\\
                              h_1 '\oplus \pi_{4(t-1/4)}\circ \Phi_1|_{C(X)\otimes 1}       & \text{if $t\in (1/4, 1/2]$;}\\
                              {\rm ad}\, U_{2(t-1/2)}\circ (h_1'\oplus \Phi_1|_{C(X)\otimes 1})   & \text{if $t\in (1/2,
                              1]$}\,.
                                 \end{cases}
                                           $$

and define $H_2: C(X)\to C([0,1], A)$ by
$$
\pi_t\circ H_2=\begin{cases} h_2'\oplus h_{00}^{(2)} & \text{if $t\in [1, 3/2]$;}\\
             h_2'\oplus \pi_{4(7/4-t)}\circ \Phi_2|_{C(X)\otimes 1} & \text{if $t\in (3/2, 7/4]$;}\\
             h_2'\oplus \pi_0\circ \Phi_2|_{C(X)\otimes 1} & \text{if $t\in
             (7/4,2]$}\,.
             \end{cases}
             $$

We compute that
\beq\label{suppit-19}
\text{Length}(\{u_t\})\le
\ep/8+\pi+2\pi+\ep/8+\pi+\ep/8<4\pi+\ep,
\eneq
\beq\label{suppit-20-}
\overline{\text{Length}}(\{\pi_t\circ H_1: t\in [0,1/2]\}) &\le& \ep/8+{\bar L}_p(X)+\ep/4\andeqn\\
\text{Length}(\{\pi_t\circ H_1: t\in [1/2,1]\}) &\le & 2\pi+{\ep\over{8(1+L(X)}}.
\eneq
Therefore (by \ref{Lnprop1})
\beq\label{suppit-20}
\overline{\text{Length}}(\{\pi_t\circ H_1\})\le {\bar L}_p(X)+2\pi \underline{L}_p(X)+\ep/2.
\eneq
Moreover,
\beq\label{suppit-21}
\overline{\text{Length}}(\{\pi_t\circ H_2\})\le {\bar L}_p(X)+\ep/4.
\eneq

 It follows from (\ref{suppit-14}),
(\ref{suppit-18}), (\ref{suppiT-04+})  and the construction that
\beq\label{suppi-22}
\|[\pi_t\circ H_1(f), u_t]\|<\ep\rforal
f\in {\mathcal F}\andeqn t\in [0,1],
\eneq
and
\beq\label{suppit-23}
\|[\pi_t\circ H_2, u_{1+t}]\|<\ep\rforal f\in {\mathcal F}
 \andeqn  t\in [0,1].
 \eneq
It is then easy to check  that $\{u_t:t\in [0,2]\},$ $H_1$ and $H_2$ meet the requirements.

\end{proof}

The following is an improved version of the original Super
Homotopy Lemma  of \cite{BEEK} for purely infinite simple \CA s.
One notes that length of $\{U_t\}$ and that of $\{V_t\}$ are
considerably shortened.

\begin{cor}\label{supbeek}
For any $\ep>0,$ there exists $\dt>0$ satisfying the following:
Suppose that $A$ is a unital purely infinite simple \CA\,and $u_1,
u_2, v_1,v_2\in U(A)$ are four  unitaries. Suppose that
$$
\|[u_1,\,v_1]\|<\dt,\,\,\, \|[u_2, v_2]\|<\dt,
$$
$$
[u_1]=[u_2], \,\,\,[v_1]=[v_2]\andeqn \rm{bott}_1(u_1,v_1)=
\rm{bott}_1(u_2, v_2).
$$
Then there exist two continuous paths of unitaries $\{U_t: t\in
[0,1]\}$ and $\{V_t: t\in [0,1]\}$ of $A$ such that
$$
U_0=u_1,\,\,\, U_1=u_2,\,\,\, V_0=v_1,\,\,\,V_1=v_2,
$$
$$
\|[U_t, V_t]\|<\ep\rforal t\in [0,1]\andeqn
$$
$$
\rm{Length}(\{U_t\})\le 4\pi +\ep \andeqn \rm{Length}(\{V_t\})\le
4\pi+\ep.
$$
Moreover, if $sp(u_1)=sp(u_2)=S^1,$ then the length of $\{U_t\}$ can be shortened to $2\pi+\ep.$
\end{cor}

\begin{proof}
When $sp(u_1)=sp(u_2)=S^1,$ the corollary follows immediately from
Theorem \ref{SHLoldpi}. Here $X=S^1.$ The subset ${\mathcal F}$
should contain $1_{C(X)}$ and $z.$ Note that the function $z$ has
the  Lipschitz constant ${\bar D}_z=1.$ Note also that if we take
$\{u_t\}$ for $\{V_t\}$ and $U_t=\pi_t\circ H(z\otimes 1),$ then
the pair will almost do the job. It is standard  to bridge the
ends of the path to $u_1$ and $u_2$ within an arbitrarily small
error (both in norm and length).

For the general case, we apply Theorem \ref{suppiT}.
Note again that it is standard to bridge the
 the arbitrarily
small gap between two paths of unitaries.
\end{proof}

\section{Super Homotopy Lemma --- finite case}

\begin{lem}\label{supdig}
Let $X$ be a finite CW complex  and ${\mathcal F}\subset C(X)$ be a
finite subset. Let $\ep>0$ and $\gamma>0,$  there exits $\dt>0$
and a finite subset ${\mathcal G}\subset C(X)$ satisfying the
following:

Suppose that $A$ is a unital separable simple \CA\, with tracial
rank zero and that $h: C(X)\to A$ is a unital \hm\, and $u\in A$
is a unitary such that
\beq\label{supdig-1}
\|[h(g), u]\|<\dt\rforal g\in {\mathcal G}.
\eneq
Then, there is a projection $p\in A,$ a unital \hm\,  $h_1: C(X)\to
(1-p)A(1-p),$ a unital \hm\,  $\Phi: C(X\otimes S^1)\to pAp$ with
finite dimensional range and a unitary $v\in (1-p)A(1-p)$ such
that
\beq\nonumber
&&\|u-v\oplus \Phi(1\otimes z)\|<\ep,\,\,\,\|[v,\,
h_1(f)]\|<\ep,\\\label{supdig-2}
&&\|h(f)-(h_1(f)\oplus
\Phi(f\otimes 1))\|<\ep
\eneq
for all $f\in {\mathcal F}$  and
\beq\label{supdig-3}
\tau(1-p)<\gamma\rforal \tau\in T(A).
\eneq

\end{lem}

\begin{proof}

Let $\ep>0,$
$\gamma>0$ and a finite subset ${\mathcal F}\subset C(X)$ be given. We
may assume that $1_{C(X)}\in {\mathcal F}.$

Let $\dt>0$ and a finite subset ${\mathcal G}\subset C(X)$ be as
required by \ref{CL2dig} corresponding to $\ep/4$ and $\gamma/2$
and ${\mathcal F}_1.$

We may assume that ${\mathcal F}\subset {\mathcal G}$ and $\dt<\ep/16.$

Let $\eta_1=\sigma_{X, \ep/32, {\mathcal F}}.$ We may assume that
$\eta_1<1.$ Now suppose that $h$ and $u$ are in the statement of
the lemma for the above $\dt$ and ${\mathcal G}.$ Suppose that
$F\subset X$ is the spectrum of $h.$

 Let $\{x_1,x_2,...,x_n\}\subset F$ be
an $\eta_1/4$-dense subset. Suppose that $O_i\cap O_j=\emptyset$
for $i\not=j,$ where
$$
O_i=\{x\in X: \di(x, x_i)<\eta_1/2s\},\,\,\,i=1,2,...,n,
$$ where $s\ge 1$ is an integer.

Let $g_i\in C(X)$ be such that $0\le g_i(x)\le 1,$ $g_i(x)=1$ if
$x\in O(x_i, \eta_1/4s)$ and $g_i(x)=0$ if $x\not\in O_i,$
$i=1,2,...,n.$

Put ${\mathcal F}_1={\mathcal F}\cup\{g_i: i=1,2,...,n\}$ and
put
$$
Q_i=O(x_i, \eta_1/4s),\,\,\,i=1,2,...,n.
$$
 Let
\beq\label{supdig+1}
\inf\{\mu_{\tau\circ h}(Q_i): \tau\in T(A), i=1,2,...,n\}\ge
\sigma_1\eta_1
\eneq
for some $\sigma_1<1/2s.$

Since $F$ is compact, there are finitely many $y_1,y_2,...,y_K\in
F$ such that
$$
\cup_{k=1}^K O(y_k,\eta_1/4)\supset F.
$$
Denote $\Omega_k=\{x\in X: \di(x, y_k)\le \eta_1/4\}$ and
$$
Y=\cup_{k=1}^K \Omega_k.
$$
Then $Y$ is a finite CW complex and $F\subset Y\subset X.$
 Moreover, $\{x_1,x_2,...,x_n\}$ is $\eta_1/2$ -dense in $Y.$
Let $s: C(X)\to C(Y)$ be the
quotient map and let ${\bar h}: C(Y)\to A$ be such that ${\bar
h}\circ s=h.$ Put
$$
{\bar O}_i=\{y\in Y: \di(x_i, y)<\eta_1/2s\},\,\,\,i=1,2,...,n.
$$

Choose an integer $L$ such that $1/L<\gamma/2$ and let
$\sigma_2={\sigma_1\over{2(L+1)}}.$

Let $\dt_1>0$  and ${\mathcal G}_1\subset C(Y)$ be a finite subset and
${\mathcal P}\subset \underline{K}(C(Y))$ be a finite subset required
by \ref{K-semi} corresponding to $\ep/2,$ $s({\mathcal F}),$ $\sigma_2$
(and $\eta_1$) above. Let ${\mathcal G}_2\subset C(X)$ be a finite
subset such that $s({\mathcal G}_2)\supset {\mathcal G}_1.$

Let $\dt_2=\min\{{\dt_1\over{2}}, {\sigma_2\cdot
\eta_1\over{2}}\}.$ We may assume that
$$
[L_1]|_{\mathcal P}=[L_2]|_{\mathcal P}
$$
for any pair of $\dt_2$-${\mathcal G}_1$-multiplicative \morp s $L_1,
L_2: C(Y)\to A$ provided that
$$
L_1\approx_{\dt_2} L_2\,\,\,\,\text{on}\,\,\,{\mathcal G}_1.
$$



 By applying \ref{CL2dig}, there is
a projection $q\in A$  and
  a unital \hm\, $\Psi: C(X\otimes S^1)\to qAq$ with finite dimensional range
  and a unitary $v\in (1-q)A(1-q)$ such that
\beq\label{supdig-6}
\|u-v_1\oplus \Psi(1\otimes z)\|<\ep/2,
\eneq
\beq\label{supdig-6+}
\|h(f)g(u)-((1-q)h(f)g(u)(1-q)\oplus \Psi(f\otimes g)))\|<\ep/2
\eneq
for all $f\in {\mathcal F}_1$  and $g\in S,$ where $S=\{1_{C(S^1)},
z\},$ and
\beq\label{supdig-7}
\tau(1-q)<\gamma/2\rforal \tau\in T(A).
\eneq
Moreover,
\beq\label{supdig-8}
\|h(f)-((1-q)h(f)(1-q)+\Psi(f\otimes  1))\|<\dt_2/2\andeqn
\eneq
\beq\label{supdig-9}
\|[(1-p),\,h(f)]\|<\dt_2/2\rforal f\in {\mathcal G}_2.
\eneq
In particular, $(1-q)h(1-q)$ is $\dt_2/2$-${\mathcal
G}_2$-multiplicative.

We may write
\beq\label{supdig-10}
\Psi(f\otimes g)=\sum_{k=1}^m
(\sum_{j=1}^{n(k)}f(\xi_k)g(t_{k,j})q_{k,j})
\eneq
for all $f\in C(X)\andeqn g\in C(S^1),$ where $\xi_k\in X,$
$t_{k,j}\in S^1$ and where $\{q_{k,j}:k,j\}$ is a set of mutually
orthogonal projections with
$\sum_{k=1}^m\sum_{j=1}^{n(k)}q_{k,j}=q.$

For each $j$ and $k,$ there are mutually orthogonal projections
$p_{k,j,l},$ $l=1,2,...,L+1$ in $A$ such that
\beq\label{supdig-11}
q_{k,j}=\sum_{l=1}^{L+1}p_{k,j,l},\,\,\, [p_{k,j, l}]=[p_{k,j,1}],
\eneq
$l=1,2,...,L$ and $[p_{k,j,L+1}]\le [p_{k,j,1}].$ This can be
easily arranged since $A$ has tracial rank zero but it also
follows from \ref{half}.

Put $p_{k,j}=\sum_{l=2}^{L+1}p_{k,j,l},$ $j=1,2,...,n(k)$ and
$k=1,2,...,m$ and put
$$
p=\sum_{k=1}^m\sum_{j=1}^{n(k)}p_{k,j}.
$$
Define $\Phi: C(X\times S^1)\to pAp$ by
\beq\label{supdig-12}
\Phi(f)=\sum_{k=1}^m
\sum_{j=1}^{n(k)}f(\xi_k)g(t_{k,j})p_{k,j}\rforal f\in C(X)\andeqn
g\in C(S^1)
\eneq
and define $\psi: C(X)\to (1-p)A(1-p)$ by
\beq\label{supdig-13}
\psi(f)=(1-q)h(f)(1-q)+\sum_{k=1}^nf(\xi_k)(\sum_{j=1}^{n(k)}p_{k,j,1})
\eneq
for all $f\in C(X).$ Note that $\psi$ is $\dt_2/2$-${\mathcal
G}_1$-multiplicative. Define $\psi_1: C(Y)\to (1-p)A(1-p)$ by
\beq\label{supdig-13+}
\psi_1(f)=(1-q){\bar
h}(f)(1-p)+\sum_{k=1}^nf(x_k)(\sum_{j=1}^{n(k)}p_{k,j,1})
\eneq
for all $f\in C(Y).$

Define $\phi: C(Y)\to pAp$ by
$$
\phi(f)=\sum_{k=1}^m f(x_k)(\sum_{j=1}^{n(k)}p_{k,j})
$$ for all $f\in C(Y).$
Write $Y=\sqcup_{k=1}^m Y_k,$ as a disjoint union of finitely many
path connected components, where each $Y_k$ is a connected finite
CW complex. Write $Z_k=Y_k\setminus \{y_k\},$ where $y_k\in Y_k$
is a point.

Note that
\beq\label{supdig-13++}
[\phi|_{C_0(\sqcup_k Z_k)}]=0\,\,\,\text{in}\,\,\,KL(C_0(\sqcup_k
Z_k), A).
\eneq
Let $E_k\in C(Y)$ be the projection corresponding to $Y_k,$
$k=1,2,...,m.$
 It follows \ref{shk} that there is a
unital \hm\, $h_{00}: C(Y)\to (1-p)A(1-p)$ such that
\beq\label{supdig-13+++}
[h_{00}|_{C_0(\sqcup_k
Z_k)}]=[h|_{C_0(\sqcup_kZ_k)}]\,\,\,\text{in}\,\,\, KL(C_0(Y_X),
A)
\eneq
and
\beq\label{supdig-13+1}
[h_{00}(E_k)]=[\psi_1(E_k)]\,\,\,\text{in}\,\,\, K_0(A),
k=1,2,...,m.
\eneq

 It follows from (\ref{supdig-8}) that
\beq\label{supdig-14-}
[h_{00}]|_{\mathcal P}=[\psi_1]|_{\mathcal P}.
\eneq

On the other hand, by (\ref{supdig-8}) again,
\beq\label{supdig-14}
\|h(f)-(\psi(f)-\Phi(f\otimes 1))\|<\dt_2/2\rforal f\in {\mathcal G}_2.
\eneq
We estimate, by (\ref{supdig-14}) and (\ref{supdig+1}) that
\beq\label{supdig-15}
\tau(\psi_1\circ s(g_i))&\ge & \tau(h(g_i))-\tau(\phi(s(g_i)))-\dt_2/2\\
&\ge &
\tau(h(g_i))-{L\tau(\Phi(g_i))\over{L+1}}-{\sigma_1\cdot\eta_1\over{4(L+1)}}\\
&\ge& {\tau(h(g_i))\over{L+1}}-\dt_2/2-{\sigma_1\cdot\eta_1\over{4(L+1)}}\\
&\ge &
{\sigma_1\cdot\eta_1\over{L+1}}-{\sigma_1\cdot\eta_1\over{2(L+1)}}\\
&=&{\sigma_1\cdot\eta_1\over{2(L+1)}}\ge \sigma_2\cdot \eta_1,
\eneq
$i=1,2,...,n.$
 It follows that
\beq\label{supdig-16}
\mu_{\tau\circ \psi_1}({\bar O}_i)\ge \sigma_2\cdot \eta_1\rforal
\tau\in T(A),
\eneq
$i=1,2,...,n.$ Therefore
$$
\mu_{t\circ\psi_1}({\bar O}_i)\ge \sigma_2\cdot\eta_1\rforal t\in
T((1-p)A(1-p)).
$$

By applying \ref{K-semi} to $\psi_1,$ we obtain a unital
monomorphism ${\bar h}_1: C(Y)\to (1-p)A(1-qp)$ such that
\beq\label{supdig-17}
\psi_1\approx_{\ep/2} {\bar h}_1\,\,\,\,\text{on}\,\,\,s({\mathcal F}).
\eneq
Define $h_1={\bar h}_1\circ s.$ Then (\ref{supdig-17})  becomes
\beq\label{supdig-17+}
\psi\approx_{\ep/2} h_1\,\,\,\text{on}\,\,\, {\mathcal F}.
\eneq
It follows from (\ref{supdig-17+}) and (\ref{supdig-8}) that
\beq\label{supdig-18-}
h\approx_{\ep} h_1\oplus \Phi\,\,\,\text{on}\,\,\,{\mathcal F}.
\eneq
Put $v=v_1\oplus
\sum_{k=1}^n\sum_{j=1}^{n(k)}t_{k,j}p_{k,j,1}.$ It is a unitary in
$(1-p)A(1-p).$ Moreover,
$$
\|u-v\oplus \Phi(1\otimes z)\|<\ep.
$$
The other inequality in (\ref{supdig-2})  also follows.
Finally,
\beq\label{supdig-18}
\tau(1-p)&=&\tau(1-q)+\tau(\sum_{k=1}^m\sum_{j=1}^{n(k)}p_{k,j,1})\\
&<& \gamma/2+{1\over{L}}<\gamma
\eneq
for all $\tau\in T(A).$

\end{proof}

\begin{thm}\label{SHLold}
Let $X$ be a connected finite CW complex, let $\ep>0$ and
${\mathcal F}\subset C(X)$ be a finite subset. Then, there exists
$\dt>0$ and a finite subset ${\mathcal G}\subset C(X)$ satisfying
the following:

Suppose that $A$ is a unital separable simple \CA\, with tracial
rank zero and suppose that $h_1, h_2: C(X)\to A$ are two unital
\hm s such that \beq\label{shlold1}
[h_1]=[h_2]\,\,\,\text{in}\,\,\, KL(C(X), A) \andeqn
 \eneq
 if there
are unitaries $u, v\in A$ such that
\beq\label{shlold2}
\|[h_1(f),\, u]\|&<&\dt,\,\,\,\|[h_2(f),\, v]\|<\dt\,\rforal\,
f\in
{\mathcal G}\andeqn\\
\rm{Bott}(h_1, u)&=&\rm{Bott}(h_2, v)
\eneq
then there exist a unital \hm\,  $H_i: C(X)\to
C([0,1], A)$ ($i=1,2$) and continuous rectifiable path $\{u_t:
t\in [0,2]\}$
such that
\beq\label{shiold3}
\hspace{-0.2in}u_0=u,\,\,\, u_2=v,\,\,\, \pi_0\circ
H_1\approx_{\ep}h_1,\,\pi_1\circ H_1\approx_{\ep}\pi_0\circ H_2
\andeqn \pi_1\circ H_2\approx_{\ep}h_2
\eneq
on ${\mathcal F},$
 and
\beq\label{shiold3+}
\|[u_t,\, \pi_t\circ H_1(f)]\|<\ep \andeqn \|[u_{t+1},\,
\pi_t\circ H_2(f)]\|<\ep
\eneq
for all $t\in [0,1]$ and for all $f\in {\mathcal F}.$ Moreover,
\beq\label{shiold4}
\overline{\rm{Length}}(\{\pi_t\circ H_1: t\in [0,2/3] \})&\le& L(X)+\ep/2\\
{\rm{Length}}(\{\pi_t\circ H_1: t\in [2/3, 1]\}) &\le & 2\pi +{\ep\over{2(1+L(X))}},\\
\overline{\rm{Length}}(\{\pi_t\circ H_1\})&\le& L(X)+2\pi \underline{L}_p(X)+\ep, \\
\overline{\rm{Length}}(\{\pi_t\circ H_2\})&\le& L(X)+\ep\andeqn\\
{\rm{Length}}(\{u_t\})&\le & 4\pi+\ep.
\eneq
\end{thm}

\begin{proof}
Let $\ep>0$ and ${\mathcal F}\subset C(X)$ be a finite subset. Let
$$
{\mathcal F}_1=\{f\otimes g\in C(X\times S^1):  f\in {\mathcal F}, g=1_{C(S^1)},
g=z\}.
$$
 We may
assume that ${\mathcal F}$ is in the unit ball of $C(X).$

Let $\eta_1=\sigma_{X,\ep/32, {\mathcal F}}$ and
$\eta_2=\min\{\sigma_{X\times S^1, \ep/32,{\mathcal F}_1}, \eta_1\}.$
 Let $\{x_1, x_2,...,x_m\}$ be an $\eta_1/4$-dense in $X.$
 Let $\{t_1,t_2,...,t_l\}$ be $l$ points on the unit circle
 which divides the unit circle into $l$ even arcs. Moreover,
 we assume that $8 m \pi/l<\min\{\ep/4,\eta_2/2\}.$
Choose $s\ge 1$ such that $O(x_i\times t_j)\cap O(x_{i'}\times
t_{j'})=\emptyset$ if $i\not=i'$ or $j\not=j',$ where
$$
O(x_i\times t_j)=\{ x\times t\in X\times S^1: {\rm
dist}(x,x_i)<\eta_1/2s \andeqn {\rm dist}(t,t_j)<\pi/sl\}.
$$

Choose $0<\sigma_1<1/2ml\eta_2$ and $\sigma_1<1/2s.$

Let ${\mathcal G}_1\subset C(X\times S^1)$ be a finite subset,
$\dt_1>0,$ $1>\gamma>0$ and ${\mathcal P}\subset
\underline{K}(C(X\times S^1))$ be as required by Theorem 4.6 of
\cite{Lncd} corresponding to $\ep/4,$ ${\mathcal F}_1$ and $\sigma_1$
above.

Let ${\mathcal P}_1\subset \underline{K}(C(X))$ be a finite subset
such that
$$
{\boldsymbol{ \bt}}({\mathcal P}_1)\supset \hat{{\boldsymbol{ \bt}}}({\mathcal P}).
$$
Let ${\mathcal P}_2={\mathcal P}\cup\{({\rm id}-\hat{{\boldsymbol{ \bt}}})({\mathcal P}),
\hat{{\boldsymbol{ \bt}}}({\mathcal P})\}.$

We may also assume that
$$
[L_1]|_{{\mathcal P}_2}=[L_2]|_{{\mathcal P}_2}
$$
for any pair of $\dt_1$-${\mathcal G}_1$-multiplicative \morp s $L_1,
L_2: C(X\times S^1)\to A,$ provided that
$$
L_1\approx_{\dt_1} L_2\,\,\,\text{on}\,\,\,{\mathcal G}_1.
$$

Without loss of generality, we may assume that $\dt_1<\sin
(\ep/4)$ and we may assume that ${\mathcal G}_1={\mathcal G}_1'\otimes S,$
where ${\mathcal G}_1'\subset C(X),$  where $S=\{1_{C(S^1)}, z\}.$
Moreover, we may assume that both sets are in the unit balls of
$C(X)$ and $C(S^1),$ respectively.

Without loss of generality, we may also assume that $\gamma<1/2.$
Let $\dt_2>0$ ( in place of $\dt$) and ${\mathcal G}_2\subset
C(X)$ (in place of ${\mathcal F}_1$) be a finite subset required
by \ref{appn} for $\dt_1/2$ (in place of $\ep$) and ${\mathcal
G}_1'$ (in place of ${\mathcal F}_0$). We may assume that
$\dt_2<\dt_1/2$ and ${\mathcal G}_1'\subset {\mathcal G}_2.$

 Let $\dt>0$ and ${\mathcal G}_3\subset C(X)$ be a finite
subset required in \ref{supdig} for $\dt_2/2,$ $\gamma/4$ and
${\mathcal G}_2.$ We may assume that ${\mathcal G}_2\subset {\mathcal G}_3$ and
$\dt_3<\dt_2/2.$

Let  ${\mathcal G}={\mathcal G}_3\cup {\mathcal F}.$

 Now suppose that

\beq\label{supf-1}
\|[u,h_1(g)]\|<\dt\andeqn \|[v, h_2(g)]\|<\dt
\eneq
for all $g\in {\mathcal G}$ and
\beq\label{supf-2}
[h_1]=[h_2]\,\,\,\text{in}\,\,\,KL(C(X),A)\andeqn
\text{Bott}(h_1,u)=\text{Bott}(h_2,v).
\eneq
In particular,
\beq\label{supf-3}
\text{Bott}(h_1,u)|_{{\mathcal P}_1}=\text{Bott}(h_2,v)|_{{\mathcal P}_1}.
\eneq

By applying \ref{supdig}, we obtain projections $p_i\in A$ and
unital \hm s $h_{1,i}: C(X)\to (1-p_i)A(1-p_i),$ unital
\hm s $\Phi_i: C(X\times S^1)\to p_iAp_i$ and unitaries $w_i'\in
(1-p_i)A(1-p_i)$ such that
\beq\label{supf-4}
\|u-w_1\oplus \Phi_1(1\otimes z)\|<\dt_2/2, \|v-w_2\oplus
\Phi_2(1\otimes z)\|<\dt_2/2 ,
\eneq
\beq\label{supf-5}
\|[w_i,\,h_i(f)]\|<\dt_2/2\andeqn h_i\approx_{\dt_2/2}
h_{1,i}\oplus \phi_i\,\,\,\text{on}\,\,\,{\mathcal G}_2,
\eneq
where $\phi_i(f)=\Phi_i(f\otimes 1)$ for $f\in C(X),$ $i=1,2.$
Moreover,
\beq\label{supf-6}
\tau(1-p_i)<\gamma/4\rforal \tau\in T(A),\,\,\,i=1,2.
\eneq
Using Zhang's Riesz interpolation as in the proof of \ref{ndig},
we may assume (by replacing $\gamma/4$ by $\gamma/2,$ if
necessary, in (\ref{supf-5})) that
$$
[p_1]=[p_2].
$$
Define $\Psi^{(0)}: C(X\times S^1)\to p_1Ap_1$ by
\beq\nonumber
\Psi^{(0)}(f)=\sum_{k=1}^m\sum_{j=1}^l
f(x_k)g(t_j)e_{k,j} \rforal f\in C(X)\andeqn g\in C(S^1),
\eneq
where $\{e_{k,j}: k,j\}$ is a set of mutually orthogonal
projections in $p_1Ap_1$ such that
\beq\label{supf-8}
\tau(e_{k,j})\ge
{\tau(1-p_1)\over{2ml}}\ge \sigma_1\cdot\eta_2\rforal \tau\in T(A).
\eneq
($i=1,2$).

It follows from \ref{zhr2} that there is unital commutative finite
dimensional \SCA\,  $B_i\subset  p_iAp_i$ which contains the image
of $\Phi_i(C(X\times S^1))$ and a set of mutually orthogonal
projections $e_{k,j}^{(i)}$ such that
\beq\label{supf-7-1}
[e_{k,j}^{(i)}]=[e_{k,j}],\,\,\,k=1,2,...,m, j=1,2,...,l,
\eneq
and $i=1,2.$

Define $\Psi_{0,i}: C(X\times S^1)\to B_i$  by
\beq\label{supf-7}
\Psi_{0,i}(f\otimes g)=\sum_{k=1}^m\sum_{j=1}^l
f(x_k)g(t_j)e_{k,j}^{(i)}
\eneq
for all  $f\in C(X)\andeqn g\in C(S^1).$

 By \ref{LengthL}, there is a unital \hm\, $H_i': C(X)\to
C([0,1],B_i)$ and a continuous path of unitaries $w_{i,t}:
t\in [0,1]$ in $p_iAp_i$ such that
\beq\label{supf-9}
w_{i,0}=\Phi_i(1\otimes z),\,\,\, w_{i,1}=\Psi_{0,i}(1\otimes z),
\eneq
\beq\label{supf-10}
\pi_0\circ H_i'(f)=\Phi_i(f\otimes 1), \pi_1\circ
H_i'(f)=\Psi_{0,i}(f\otimes 1)\andeqn
\eneq
\beq\label{supf-11}
[\pi_t\circ H_i'(f),\, w_{i,t}]=0\rforal t\in [0,1]\andeqn\rforal
f\in C(X),
\eneq
$i=1,2.$ Moreover,
\beq\label{supf-12}
\overline{\rm{Lenghth}}(\{\pi_t\circ H_i'\})\le L(X)+\ep/2\andeqn
\text{Length}(\{w_{i,t}\})\le \pi,
\eneq
$i=1,2.$

From (\ref{supf-4}) and (\ref{supf-5}), by applying \ref{appn}, we obtain a $\dt_1/2$-${\mathcal G}_1'$ -multiplicative \morp\,
$\psi_i': C(X\times S^1)\to A$  such that
\beq\label{supf-13}
\|\psi_i'(f\otimes g)-h_{1,i}(f)g(w_i)\|<\dt_1/2
\eneq
for all $f\in {\mathcal G}_1'$ and $g\in S,$ $i=1,2.$
Define $\psi_i: C(X\times S^1)\to A$ by
\beq\label{supf-13+}
\psi_i(f\otimes g)=\psi_i'(f\otimes g)\oplus \Phi_{0,i}(f\otimes g)
\eneq
for all $f\in C(X)\andeqn g\in C(S^1).$
 Moreover, by the fact that
$[p_1]=[p_2],$ (\ref{supf-2}), (\ref{supf-3}), (\ref{supf-4}),
(\ref{supf-5}) and (\ref{supf-7-1}),  we compute that
\beq\label{supf-14}
[\psi_1]|_{\mathcal P}=[\psi_2]|_{\mathcal P}.
\eneq
By \ref{supf-8}), (\ref{supf-7}) and (\ref{supf-7-1}),
\beq\label{supf-15}
\mu_{\tau\circ \phi_i}(O(x_k\times t_j))\ge
\sigma_2\cdot\eta_2\rforal \tau\in T(A),
\eneq
$k=1,2,..,m,$ $j=1,2,...,l$ and $i=1,2.$ Furthermore, since
$[e_{k,j}^{(1)}]=[e_{k,j}^{(2)}],$ $k=1,2,...,m,$ $j=1,2,...,l,$
we compute that
\beq\label{supf-16}
|\tau\circ \psi_1(f)-\tau\circ \psi_2(f)|<\gamma\rforal f\in {\mathcal
G}_1.
\eneq
Now, from the choice of $\dt_1,$ $\gamma,$ ${\mathcal G}_1,$ by
applying Theorem 4.6 of \cite{Lncd}, we obtain a unitary $U\in A$
such that
\beq\label{supf-17}
{\rm ad}\, U\circ \psi_1\approx_{\ep/4}
\psi_2\,\,\,\text{on}\,\,\, {\mathcal F}_1.
\eneq
There is a continuous path of unitaries $\{U_t:t\in [0,1]\}$ such
that
\beq\label{supf-18}
U_0=1,\,\,\,U_1=U\andeqn \text{Length}(\{U_t\})\le
\pi+{\ep\over{4(1+L(X))}}.
\eneq
Define $H_1: C(X)\to C([0,1],A)$ by, for $f\in C(X),$
\beq\label{supf-19}
\hspace{-0.2in}\pi_t\circ H_1(f)=\begin{cases} h_{1,1}(f)\oplus
\phi_1(f)
& \hspace{-0.13in}\text{if $t\in [0,1/3]$;}\\
 h_{1,1}(f)\oplus \pi_{3(t-1/3)}\circ H_1'(f)
&\hspace{-0.13in}\text{if $t\in (1/3,2/3]$;}\\
   U_{3(t-2/3})^* (h_{1,1}(f)\oplus \Psi_{0,1}(f\otimes
   1))U_{3(t-2/3)}
   &\hspace{-0.13in}\text{if $t\in (2/3,1]$}\,.
\end{cases}
\eneq
Define $H_2: C(X)\to C([0,1], A)$ by
\beq\label{supf-20-}
\pi_t\circ H_2=\begin{cases} h_{1,2}\oplus \pi_{2(1/2-t)}\circ
H_2'
& \text{if $t\in [0,1/2]$;}\\
               h_{1,2}\oplus \phi_2 & \text{if $t\in (1/2,1]$}.
               \end{cases}
\eneq
Then, by (\ref{supf-5}) and (\ref{supf-17}),
\beq\label{supf-20}
\pi_0\circ H_1\approx_{\ep} h_1,\,\,\, \pi_1\circ H_1\approx_{\ep}
\pi_0\circ H_2\andeqn \pi_1\circ H_2\approx_{\ep} h_2
\eneq
on ${\mathcal F}.$ Moreover, by (\ref{supf-12}) and (\ref{supf-18})
\beq\label{supf-21}
\overline{\text{Length}}(\{\pi_t\circ H_1:t\in [0,2/3]\})&\le& L(X)+\ep/2,\\
\text{Length}(\{\pi_t\circ H_1: t\in [2/3,1]\} &\le &
2\pi+{\ep\over{2(1+L(X))}}
\eneq
\andeqn\beq\label{supf-21+} \overline{\text{Length}}(\{\pi_t\circ
H_2\})\le L(X)+\ep.
\eneq
So, by \ref{Lnprop1},
\beq\label{supf-21+1}
\overline{\text{Length}}(\{\pi_t\circ H_1\})\le L(X)+2\pi
\underline{L}_p(X)+\ep.
\eneq

By (\ref{supf-4}), since $\dt_1<\sin (\ep/4),$ there is $a_i\in
A_{s.a}$ with $\|a_i\|<\ep/4$ such that
$$
u^*(w_1\oplus \Phi_1(1\otimes z))=\exp(ia_1)\andeqn v^*(w_2\oplus
\Phi_2(1\otimes z))=\exp(i a_2).
$$

Define
\beq\label{supf-22}
u_t=\begin{cases} u\exp(i 3ta_1) &\text{if $t\in
[0,{1\over{3}}]$;}\\
w_1\oplus w_{1,3(t-1/3)} & \text{if $t\in
({1\over{3}}, {2\over{3}}]$;}\\
   U_{3(t-1/3)}^*\circ (w_1\oplus w_{1,1})U_{3(t-1/3)} & \text{if $t\in
   ({2\over{3}}, 1]$;}\\
   w_1\oplus w_{2,2(3/2-t)} &\text{if $t\in (1,{3\over{2}}]$;}\\
  v\exp(i2(2-t) a_2) & \text{if $t\in ({3\over{2}}, 2]$}\,.
   \end{cases}
   \eneq
Then
$$
u_0=u,\,\,\, u_2=v\andeqn
$$
$$
\text{Length}(\{u_t\})\le \ep/4+\pi+2\pi+\pi +\ep/4=4\pi +\ep/2.
$$

 We also have
$$
\|[u_t, \pi_t\circ H_1(f)]\|<\ep\rforal f\in {\mathcal F}
$$
and
$$
\|[u_{t+1}, \pi_t\circ H_2(f)]\|<\ep\rforal f\in {\mathcal F}
$$
for  $t\in [0,1].$

\end{proof}

We obtain the following improvement of original Super Homotopy
Lemma in the case that $A$ is a unital separable simple \CA\, of
tracial rank zero.

\vspace{0.1in}

\begin{cor}\label{subeek2}
For any  $\ep>0,$ there exists $\dt>0$ satisfying the following:

Suppose that $A$ is a unital separable simple \CA\, with tracial
rank zero and $u_1, u_2, v_1,v_2\in U(A)$ are four  unitaries. Suppose that
$$
\|[u_1,\,v_1]\|<\dt,\,\,\, \|[u_2, v_2]\|<\dt,
$$
$$
[u_1]=[u_2], \,\,\,[v_1]=[v_2]\andeqn \rm{bott}_1(u_1,v_1)=
\rm{bott}_1(u_2, v_2).
$$
Then there exist two continuous paths of unitaries $\{U_t: t\in
[0,1]\}$ and $\{V_t: t\in [0,1]\}$ of $A$ such that
$$
U_0=u_1,\,\,\, U_1=u_2,\,\,\, V_0=v_1,\,\,\,V_1=v_2,
$$
$$
\|[U_t, V_t]\|<\ep\rforal t\in [0,1]\andeqn
$$
$$
\rm{Length}(\{U_t\})\le 4\pi +\ep \andeqn \rm{Length}(\{V_t\})\le
4\pi+\ep.
$$
\end{cor}

\vspace{0.1in}

\begin{rem}\label{Rsup}
{\rm As stated in  \cite{BEEK}, the Basic Homotopy Lemma is not a
special case of the Super Homotopy Lemma. It mentioned not only
the estimate on the length of the paths but more importantly, the
presence of two paths of unitaries in the Super Homotopy Lemma.
But in the Basic Homotopy Lemma it is essential to keep one
unitary fixed. In our cases also, in the Basic Homotopy Lemma, the
map $h$ remains fixed in the homotopy while \ref{suppiT} and
\ref{SHLold} allow a path from $h_1$ to $h_2.$ We would point out
that beside the obvious difference in appearance of paths and the
estimates of the length of paths, one additional essential
difference appears: In the finite case, the constant $\dt$ in the
Super Homotopy Lemma \ref{SHLold} is universal while in \ref{MT3}
$\dt$ depends on a measure distribution which as shown in
\ref{CEX} can not be universal. This phenomenon does not occur in
the case that $X$ is one-dimensional. This again shows the
additional difficulties involved in  the Basic Homotopy Lemma for
higher dimensional spaces $X.$ Moreover, in the Basic Homotopy
Lemma, $h$ is assumed to a monomorphism and it fails when $h$ is
not assumed to be monomorphism whenever ${\rm dim} X\ge 2.$
However, in both \ref{suppiT} and \ref{SHLold}, $h_1$ and $h_2$
are only assumed to be \hm s.

On the other hand, there are additional burden to connect two \hm
s $h_1$ and $h_2$ as well as two unitaries while the Basic
Homotopy Lemma considers only one unitary and one monomorphism. It
should be noted  that the Super Homotopy Lemma does not follow
from \ref{SHLpi} or \ref{SHL}, since $\pi_t\circ H_i$ there is not
a \hm\,

It is also important to note that even in the case that
$L(X)=\infty,$ Theorem \ref{GHT1}, \ref{GHT2}, \ref{GHP3},
\ref{SHLn}m \ref{GHTh2}, \ref{SHLpiinj}, \ref{SHLpi}, \ref{SHL},
\ref{SHLoldpi}, \ref{suppiT} and \ref{SHLold} remain valid.
However, the lengths of these homotopy may be infinite.

}\end{rem}

\chapter{Postlude}
\setcounter{section}{16}

\section{Non-commutative cases}

\begin{df}\label{NCD1}

{\rm

(i)\,\,\, Let $X$ be a compact metric space, let $r$ be a positive
integer and let $A=M_r(C(X))=C(X,M_r).$  Denote by $tr$ the
standard
  normalized trace on $M_r.$  Let
$\tau\in T(A).$ Then there is a probability Borel measure $\nu$
such that
$$
\tau(f)=\int_X tr(f) d\nu\rforal f\in C(X,M_r).
$$
Define $\mu_{\tau}=r\nu.$

Suppose that $B$ is a unital \CA\, and $h: A\to B$ is a unital
\hm. Let $\{e_{i,j}\}$ be a system of matrix unit for $M_r.$
Identify $e_{1,1}$ with the (constant) rank one projection in $A.$
Then $e_{1,1}Ae_{1,1}\cong C(X).$  Suppose that $t\in T(B).$  Let
$\tau=t\circ h.$ Define $h_1=h|_{e_{1,1}Ae_{1,1}}.$ Define
$\tau_1(f)=r\tau\circ h_1(f)$ for $f\in C(X).$  Then $\tau_1$ is a
tracial state of $C(X).$ The associated measure is $\mu_{\tau}.$

(ii) \,\,\, Suppose that $X$ is a connected finite CW complex with
covering dimension $d.$ Let $A=PM_m(C(X))P,$ where $P$ is a rank
$k$ ($\le m$) projection in $M_m(C(X)).$ Let $\tau\in T(A).$
Denote by $tr$ the standard trace on $M_k.$ Then there exists a
Borel probability measure $\nu$ on $X$ such that
$$
\tau(f)=\int_X tr(f) d\nu\rforal f\in A.
$$
Define $\mu_{\tau}=k\nu.$

Let $B$ be a unital \CA, and let $h: A\to B$ be a unital \hm.
There is a projection $e\in M_{d+1}(A)$ such that $e+1_A$ is a
trivial projection on $X$ (see for example 8.12 of \cite{Hu}).
Therefore $(e+1_A)M_{d+1}(A)(e+1_A)\cong M_r(C(X))$ for some
integer $r\ge 1.$ Let $t=({r\over{k(d+1)}})(\tau\otimes Tr),$ $Tr$
is the standard trace on $M_{d+1}.$ Then $t$ gives a tracial state
on $M_r(C(X)).$ Using notation in (i), the measure $\mu_t$ on $X$
is same as $\mu_{\tau}$ defined above.

(iii)\,\,\, Let  $X$ be (not necessary connected) finite CW
complex with finite covering dimension and let $A=PM_m(C(X))P,$
where $P$ is a projection in $M_m(C(X)).$ Then one may write
$X=\sqcup_{i=1}^k X_i,$ where each $X_i$ is a connected. So we may
also write $A=\oplus_{i=1}^k P_iM_{m(k)}(C(X_i))P_i,$ where
$P_i\in M_{m(i)}(C(X_i))$ is a rank $r(i)$ projection. Put $A_i=
P_iM_{m(k)}(C(X_i))P_i,$ $i=1,2,...,k.$ Suppose that $\tau\in
T(A).$  Put $\af_i=\tau(P_i)$ and $\tau_i=(1/\af_i)\tau|_{A_i}.$
By applying (ii), one then obtains a probability Borel measure
$\mu_{\tau_i}$ on $X_i.$ For any Borel subset $S\subset X,$ write
$S=\sqcup_{i=1}^k S_i,$ where each $S_i\in X_i$ is a Borel subset.
Define $\mu_{\tau}$ by
$$
\mu_{\tau}(S)=\sum_{i=1}^k \af_i\mu_{\tau_i}(S_i)
$$
for all Borel subset $S\subset X.$

 }

\end{df}

\begin{lem}\label{NCL1}
Let $X$ be a compact metric space and let $A=M_r(C(X)),$ where
$r\ge 1$ is an integer. Then, for any $\ep>0,$  any finite subset
${\mathcal F}\subset A$ and map ${\varDelta}: (0,1)\to (0,1),$ there
exists $\dt>0,$ a finite subset ${\mathcal G}\subset A$  and a finite
subset ${\mathcal P}\subset \underline{K}(A)$ satisfying the
following:

Suppose that $B$ is a unital separable simple \CA\, with tracial
rank zero, $h: A\to B$ is a unital monomorphism and $u\in B$ is a
unitary
 such that $\mu_{\tau\circ h}$ is $\varDelta$-distributed for all $\tau\in T(B),$
\beq\label{ncl1-1}
\|[h(a), u]\|<\dt\tforal\, a\in {\mathcal G} \andeqn
\rm{Bott}(h,u)|_{\mathcal P}=0.
\eneq
Then, there exists a continuous rectifiable path of unitaries
$\{u_t:t\in [0,1]\}$ of $B$ such that
\beq\label{ncl1-2}
 u_0=u,\,\,\,u_1=1_B\andeqn \|[h(a),u_t]\|<\ep\tforal  a\in {\mathcal F}\andeqn
t\in [0,1].
\eneq
Moreover,
\beq\label{ncl1-3}
\rm{Length}(\{u_t\})\le 2\pi+\ep.
\eneq

\end{lem}

\begin{proof}
Let $\ep>0$ and ${\mathcal F}\subset M_r(C(X))$ be given. There is
a finite subset ${\mathcal F}_1\subset C(X)$ such that
$$
{\mathcal F}\subset \{(g_{i,j}): g_{i,j}\in {\mathcal F}_1\}.
$$
Let $\dt_1>0,$  ${\mathcal G}_1\subset C(X)$ be a finite subset  and ${\mathcal P}_1\subset
\underline{K}(C(X))= \underline{K}(A)$ be a finite subset as required by \ref{MT2}
associated with $\ep/2,$ ${\mathcal F}_1$ and $\varDelta.$
We may assume that
$$
[L_1]|_{\mathcal P}=[L_2]|_{\mathcal P}
$$
for any pair of $\dt_1$-${\mathcal G}_1$-multiplicative \morp s $L_1,
L_2: C(X)\to A,$ provided that
$$
L_1\approx_{\dt_1} L_2\,\,\,\text{on}\,\,\,{\mathcal G}_1.
$$

 We may assume that $\dt_1<\ep/2.$

Let $\dt_2>0$ and ${\mathcal G}_2\subset A$ be a finite subset required by \ref{appn} for
$\dt_1/2$ and ${\mathcal G}_1.$ We may assume that $\dt_2<\min\{\dt_1/2, \ep/2\}$ and ${\mathcal G}_2\supset {\mathcal G}_1\cup {\mathcal F}_1.$

Let $\{e_{i,j}: i,j=1,2,...,r\}$ be a system of matrix unit for $M_r.$

It is easy to see that there is
$\dt>0$ and a finite subset ${\mathcal G}\subset A$  satisfying the following: if
$B$ is a unital \CA, if $h': A\to B$ is a unital \hm\,  and if  $u\in
B$ is a  unitary with
$$
\|[h'(f), u]\|<\dt\rforal f\in {\mathcal G},
$$
then, there exists a unitary $v\in B$ such that
$$
\|u-v\|<\dt_2/2, \,\,\,h'(e_{1,1})v=vh'(e_{1,1})\andeqn \|[h'(g),
v]\|<\dt_2/2
$$
for all $f\in {\mathcal G}_3,$ where
$$
{\mathcal G}_3=\{(g_{ij})\in M_k(C(X)): g_{i,j}\in {\mathcal G}_2\}.
$$
Suppose that $h: A\to B$ is a monomorphism and $u\in A$ is a
unitary satisfying the conditions in the lemma with $\dt$ and
${\mathcal G}$ as above. Define $h_1: C(X)\to h(e_{1,1})Bh(e_{1,1})$
by $h_1=h|_{e_{1,1}Ae_{1,1}}.$ Then, \ref{MT3} applies to $h_1$
and $e_{1,1}ve_{1,1}.$

One sees that we reduce the general
case to the case that $r=1.$

\end{proof}

\begin{lem}\label{NCLpi1}
Let $X$ be a compact metric space and let $A=M_r(C(X)),$ where
$r\ge 1$ is an integer. Then, for any $\ep>0,$  any finite subset
${\mathcal F}\subset A,$ there
exists $\dt>0,$ a finite subset ${\mathcal G}\subset A$  and a finite
subset ${\mathcal P}\subset \underline{K}(A)$ satisfying the
following:

Suppose that $B$ is a unital purely infinite  simple \CA, $h: A\to B$ is a unital monomorphism and $u\in B$ is a
unitary
 such that
\beq\label{nclp1-1}
\|[h(a), u]\|<\dt\tforal\, a\in {\mathcal G} \andeqn
\rm{Bott}(h,u)|_{\mathcal P}=0.
\eneq
Then, there exists a continuous rectifiable path of unitaries
$\{u_t:t\in [0,1]\}$ of $B$ such that
\beq\label{nclp1-2}
u_0=u,\,\,\,u_1=1_B\andeqn \|[h(a),u_t]\|<\ep\tforal\, a\in
{\mathcal F} \andeqn t\in [0,1].
\eneq
Moreover,
\beq\label{ncl1p-3}
\rm{Length}(\{u_t\})\le 2\pi+\ep
\eneq

\end{lem}

\begin{proof}
The proof of this lemma is almost identical to that of \ref{NCL1}. We will apply \ref{Tpi}.
\end{proof}


%
%
%
%

The essential significance in the following lemma is the estimate
of the length of $\{v_t\}.$ It should be noted that the proof
would be much shorter if we allow the length to be bounded by
$2L+\ep.$

\vspace{0.1in}

\begin{lem}\label{Pert}

\rm{(1)}\,\,\,There is a positive number $\dt$ (with $1/2>\dt>0$)
and  there is a function $\eta: (0,1/2)\to \R_+$ with $\lim_{t\to
0}\eta(t)=0$ satisfying the following:

Let $B$ be a unital \CA\, and $\{u_t: t\in [0,1]\}$ be a path of
unitaries with $u_1=1$ such that
$$
\text{Length}(\{u_t\})\le L.
$$
Suppose that there is a non-zero projection $p\in A$ such that
$$
\|[u_t, p]\|<\dt\rforal t\in [0,1].
$$
Then there exists a path of unitaries $\{v_t: t\in [0,1]\}\subset
pAp$  and  there exist partitions
$$
0=t_0<t_1\cdots <t_m=1\andeqn 0=s_0<s_1\cdots <s_m=1
$$
such that
$$
\|pu_{t_i}p-v_{s_i}\|<{\dt\over{\sqrt{1-\dt}}},\,\,\,i=1,2,...,m,\,\,v_1=p\andeqn
$$
$$
\text{Length}(\{v_t\})\le L+(7L\eta(\dt)),
$$
where $m=[{6L\over{\pi}}]+1.$

\rm{(2)}\,\,\, For any $\ep>0$ and any $M>0.$ There exists $\dt>0$
satisfying the following:

Let  $B$ be a unital
\CA\, and $\{u_t: t\in [0,1]\}$ be a path of unitaries with
$u_1=1$ such that
$$
\text{Length}(\{u_t\})\le L.
$$
Suppose that there is a non-zero projection $p\in A$  and a
self-adjoint subset ${\mathcal G}\subset pAp$ with $p\in {\mathcal G}$ and
$\|a\|\le M$ for all $a\in {\mathcal G}$ such that
$$
\|[u_t, a]\|<\dt\rforal t\in [0,1].
$$
Then there exists a path of unitaries $\{v_t: t\in [0,1]\}\subset
pAp$ such that
$$
\|pu_0p-v_0\|<\ep,\,\,\,v_1=p,
$$
$$
\|[v_t, a]\|<\ep\rforal a\in {\mathcal G}\andeqn
$$
$$
\text{Lengh}(\{v_t\})\le L+\ep,
$$
\end{lem}

\begin{proof}
The proof is somewhat similar to that of \ref{eqexp}. Fix $M>0.$
Put
$$
S=\{ e^{i\pi t}: t\in [-1/2, 1/2]\}.
$$
Denote by $F: S\to [-1/2, 1/2]$ by $F(e^{i\pi t})=t.$ $F$ is a
continuous function. Let
$$
S_1=\{e^{i\pi t}: t\in [-1/6,1/6]\}.
$$

If $U_1$ and $U_2$ are  two unitaries in a unital \CA\, $A$  and
$p\in A$ is a projection such that
\beq\label{Pert-01}
\|[p,U_i]\|<\dt,\,\,\,i=1,2,
\eneq
then
\beq\label{Pert-02}
\|pU_ip|pU_ip|^{-1}-pU_ip\|<{\dt\over{\sqrt{1-\dt}}},\,\,\,i=1,2.
\eneq
Here the inverse $|pU_ip|^{-1}$ is taken in $pAp.$ Note that
$pU_ip|pU_ip|^{-1}$ is a unitary in $pAp.$  If
$\text{sp}(U_1U_2)\subset S_1,$ then
\beq\label{Pert-021}
\|1-U_1U\|=\max\{|1-e^{i\pi\theta}|: \theta\in [-1/6,1/6]\} \le
1/2.
\eneq
We estimate that
\beq\label{Pert-022}
&&\hspace{-0.5in}
\|p-pU_1p|pU_1p|^{-1}pU_2p|pU_2p|^{-1}\|\\\nonumber
&\le &
\|p-pU_1U_2p\|+\|pU_1U_2p-pU_1p|pU_1p|^{-1}pU_2p|pU_2p|^{-1}\|\\\nonumber
&< & 1/2+\|pU_1U_2p-pU_1pU_2p\|+\|pU_1pU_2p-pU_1p|pU_1p|^{-1}pU_2p|pU_2p|^{-1}\|
\\
&<& 1/2+\dt+2{\dt\over{1-\dt}}.
\eneq
There is $\dt_0>0$ such that if $\dt<\dt_0,$
$$
1/2+\dt+2{\dt\over{1-\dt}}<1.
$$
Thus  when $\dt<\dt_0,$ $\text{sp}(U_1U_2)\subset S_1$ and  when
(\ref{Pert-01}) holds,
\beq\label{Pert-03}
\text{sp}(pU_1p|pU_1p|^{-1}pU_2p|pU_2p|^{-1})\subset S.
\eneq
Fix $M>0.$

 Moreover,  by choosing smaller $\dt_0>0,$ when $0<\dt<\dt_0$ and
when (\ref{Pert-01}) holds, there is a positive number $\eta(\dt)$ (for each $\dt$)
with
\beq\label{Pert-05}
\lim_{t\to 0} \eta(t)=0
\eneq
such that
\beq\label{Pert-04}
\|\exp(i \pi F(pU_1p|pU_1p|^{-1}pU_2p|pU_2p|^{-1}))-p\exp(i \pi
F(U_1U_2))p\|<\eta(\dt)
\eneq
and at the same time
\beq\label{Pert-04+}
\|F(pU_1p|pU_1p|^{-1}pU_2p|pU_2p|^{-1})-pF(U_1U_2)p\|<\eta(\dt)
\eneq
and
\beq\label{pert-041}
\|[a, F(pU_1p|pU_1p|^{-1}pU_2p|pU_2p|^{-1})]\|<\eta(\dt),
\eneq
if $\|[a, U_i]\|<\dt,$ $a\in pAp$ and $\|a\|\le M$ (cf. 2.5.11 of
\cite{Lnbk} and 2.6.10 of \cite{Lnbk}). It should be noted that
\beq\label{Pert-04++}
\|F(pU_1pU_2p|pU_1pU_2p|^{-1})\|\le \|F(U_1U_2)\|+\eta(\dt).
\eneq

We will use these facts in the rest of the proof and assume that
$0<\dt<\dt_0.$

\vspace{0.1in}

 Let $0=t_0<t_1<\cdots <t_m=b$ be a partition so that
\beq\label{Pert-1}
\pi/7\le \text{Length}(\{u_t: t\in [t_{i-1}, t_i]\})\le
\pi/6,\,\,\,i=1,2,...,m,
\eneq
where $m=[{6L\over{\pi}}]+1.$

Put $u_i=u_{t_i}$ and $z_i=pu_ip|pu_ip|^{-1},$ where the inverse
is taken in $pAp.$ Then $\text{sp}(u_{i-1}^*u_i)\subset S_1.$
Define ${\tilde w}_i=z_{i-1}^*z_i,$ $i=1,2,...,m.$

 If
\beq\label{Pert-2}
\|[u_t, p]\|<\dt,
\eneq
then (by (\ref{Pert-02})
\beq\label{Pert-3}
\|pu_ip|pu_ip|^{-1}-pu_ip\|<{\dt\over{\sqrt{1-\dt}}},
\eneq
 $i=1,2,...,m.$  It follows from (\ref{Pert-03})that
 $$sp({\tilde w_i})\subset S.$$

Moreover
 by (\ref{Pert-2}), (\ref{Pert-03}) and (\ref{Pert-04}),
\beq\label{Pert-4}\nonumber
&&\hspace{-0.5in}\|z_{i-1}\exp(\sqrt{-1}\pi F({\tilde
w_i}))-pu_ip\|\\\nonumber
&=&\|z_{i-1}\exp(\sqrt{-1}\pi F({\tilde
w_i}))-pu_{i-1}u_{i-1}^*u_ip\|\\\nonumber
&<&\|z_{i-1}\exp(\sqrt{-1}\pi F({\tilde
w_i}))-pu_{i-1}pu_{i-1}^*u_ip\|+\dt\\\nonumber
&<&\|z_{i-1}\exp(\sqrt{-1}\pi F({\tilde w_i}))-pu_{i-1}p
\exp(\sqrt{-1}\pi F({\tilde w_i}))\|+\eta(\dt)+\dt\\
&<&{\dt\over{\sqrt{1-\dt}}}+\eta(\dt)+\dt,\,\,\,i=1,2,...,m.
\eneq

Suppose that ${\mathcal G}\subset pAp$ is self-adjoint such that
$\|a\|\le M$ for all $a\in {\mathcal G}.$
Then, if $\dt<\dt_0$ and
if
$$
\|[a, u_t]\|<\dt\rforal a\in {\mathcal G}\andeqn t\in [0,1],
$$
then, by (\ref{pert-041}),
\beq\label{Pert-5-1}
\|[a, F({\tilde{w_i}})]\|\le \eta(\dt)\rforal a\in {\mathcal G}.
\eneq
Thus, by (\ref{Pert-4}), (\ref{Pert-05}) and (\ref{Pert-3}), a
positive number $\dt$ depending on $\ep,$ $M$ and $F$ only, such
that, if
$$
\|[a, u_t]\|<\dt\rforal a\in {\mathcal G}\andeqn t\in [0,1],
$$
then,
\beq\label{Pert-5}
\|[a, z_{i-1}\exp(\sqrt{-1}\,\pi t\,F({\tilde w_i}))]\|<\ep
\eneq
for all $t\in [0,1]$  and for all $a\in {\mathcal G}.$  Let
$$
l_i=\text{Length}(\{u_t:t\in [t_{i-1}, t_i]\}), i=1,2,...,m.
$$

 Now, define
$$
v_0=u_0,\,\,\, v_t=z_{i-1}(\exp(\sqrt{-1}\pi
{t-s_{i-1}\over{s_i-s_{i-1}}}F({\tilde w_i}))) \rforal \, s\in
[s_{i-1},s_i],
$$
where $s_0=0,$ $s_i=\sum_{j=1}^il_j/L,$ $i=1,2,...,m.$ It follows
that
\beq\label{Pert-5+}
\|[a,v_t]\|<\ep\rforal t\in [0,1]
\eneq
and for all $a\in {\mathcal G}.$

Clearly, if $t,t'\in [s_{i-1},s_i],$ then (also using
(\ref{Pert-04++}) and (\ref{Pert-1}))
\beq\label{Pert-6}\nonumber
\|v_t-v_{t'}\|&=&\|\exp(\sqrt{-1}\pi
{t-s_{i-1}\over{s_i-s_{i-1}}}F({\tilde w_i})) -\exp(\sqrt{-1}\pi
{t'-s_{i-1}\over{s_i-s_{i-1}}}F({\tilde w_i})) \|\\
\nonumber
 &\le & \|\pi F({\tilde w_i})\| {|t-t'|\over{l_i/L}}
 \le
L{(l_i+\eta(\dt)\cdot \pi)\over{l_i}}|t-t'|\\
&=& (L+7L\eta(\dt))|t-t'|.
\eneq

One then computes that
\beq\label{Pert-7}
\|v_t-v_{t'}\|\le (L +7L\eta(\dt))|t-t'|\rforal t,t'\in [0,1]
\eneq
The lemma follows from the fact that $\lim_{\dt\to 0}\eta(\dt)=0.$

\end{proof}

\begin{lem}\label{NCL2}
Let $X$ be a finite CW complex, let $F\subset X$ be a compact subset
and let $A_1=QM_m(C(X))Q,$ where $m\ge 1$ is an integer and $Q\in
M_m(C(X))$ is a projection.
Let $A=PM_m(C(F))P,$ where $P=s(Q)$ and $s: A_1\to A$ is the surjective
\hm\, induced by the quotient map $C(X)\to C(F).$

Then, for any $\ep>0,$  any finite subset ${\mathcal F}\subset A$
and a non-decreasing map ${\varDelta}: (0,1)\to (0,1),$ there
exists $\dt>0,$ a finite subset ${\mathcal G}\subset A$  and a
finite subset ${\mathcal P}\subset \underline{K}(A)$ satisfying
the following:

Suppose that $B$ is a unital separable simple \CA\, with tracial
rank zero, $h: A\to B$ is a unital monomorphism and $u\in B$ is a
unitary
 such that $\mu_{\tau\circ h}$ is $\varDelta$-distributed for all $\tau\in T(B),$
\beq\label{ncl2-1}
\|[h(a), u]\|<\dt\tforal\, a\in {\mathcal G} \andeqn
\rm{Bott}(h,u)|_{\mathcal P}=0.
\eneq
Then, there exists a continuous rectifiable path of unitaries
$\{u_t:t\in [0,1]\}$ of $B$ such that
\beq\label{ncl2-2}
u_0=u,\,\,\,u_1=1_B\andeqn \|[h(a),u_t]\|<\ep\tforal\, a\in
{\mathcal F}\andeqn t\in [0,1].
\eneq
Moreover,
\beq\label{ncl2-3}
\rm{Length}(\{u_t\})\le 2\pi+\ep.
\eneq

\end{lem}

\begin{proof}
It is clear that we may reduce the general case to the case that
$X$ is connected. So we now assume that $X$ is connected. Let $d$
be the dimension of $X.$ By 8.15 of \cite{Hu}, there is a
projection $e'\in M_{d+1}(QM_m(C(X))Q)$ such that $e'+Q$ is a trivial
projection. So
$$
(e'+Q)M_{d+1}(C(X))(e'+Q) \cong M_r(C(X))
$$
for some integer $r\ge 1.$
Put $e=s(e').$ It follows that
$$
(e+P)M_{d+1}(C(F))(e+P)\cong M_r(C(F)).
$$

Note that $P=1_A.$
Put $C=(e+1_A)M_{d+1}(C(F))(e+1_A).$ Note also
that
$$\underline{K}(PCP)=\underline{K}(A).$$

If $h: A\to B$ is a unital \hm, then it extends to a unital \hm\,
$h_1: M_{d+1}(A)\to M_{d+1}(B).$ Let  $h'=(h_1)|_C.$
For any $\dt_1>0$ and a finite subset ${\mathcal G}_1\subset A,$ there exists
$\dt_2>0$ and a finite subset ${\mathcal G}_2$ satisfying the following:
if $u'\in B$ is a unitary such that
$$
\|[h(f), u']\|<\dt_2\rforal f\in {\mathcal G}_2,
$$
then
$$
\|[h_1(g), U']\|<\dt_1\rforal {\mathcal G}_1'
$$
where $U'={\rm diag}(u',u',...,u')$ and where ${\mathcal G}'_1=\{(g_{i,j}): g\in {\mathcal G}_1\}.$
In particular, with larger ${\mathcal G}_1,$
$$
\|[1_C, U']\|<\dt_1.
$$
Thus there is a unitary $V\in h_1(1_C)Bh_1(1_C)$ such that
$$
\|V-h_1(1_C)U'h_1(1_C))\|<2\dt_1\andeqn \|[h'(g), V]\|<4\dt_1\rforal g\in {\mathcal G}_1'\cap C.
$$

Let $\{e_{i,j}\}$ be a system of matrix for $M_r$ and be elements in $C\cong M_r(C(F)).$

Therefore, by applying \ref{NCL1} and its proof,
for any $\ep_1>0$ and any ${\mathcal F}_1\subset C(F),$ there is $\dt>0$ and a finite subset ${\mathcal G}\subset
A$ and a finite subset ${\mathcal P}\subset  \underline{K}(A)$ satisfying the following:
suppose that $h$ and $u$ satisfy the assumption of the lemma for the above $\dt$ and ${\mathcal G}$ and ${\mathcal P},$
there exists a path of unitaries $\{w_t: t\in [0,1]\}$ in $h'(e_{1,1})Bh'(e_{1,1})$ such that
\beq\label{ncl2-4}
\|w_0-w\|<\ep_1/(4r)^2,\,\,\,w_1=e_{1,1}, \,\,\,
\|[h'(g), W_t]\|<\ep_1
\eneq
for all $t\in [0,1]$ and $g\in {\mathcal F}_1',$ and
\beq\label{ncl2-5}
\text{Length}(\{w_t\})\le 2\pi+\ep_1/2,
\eneq
where $w$ is a unitary in $e_{1,1}Ce_{1,1}$ with
\beq\label{ncl2-6}
\|w-e_{1,1}U'e_{1,1}\|<\ep_1/(2r)^2,
\eneq
$W_t={\rm diag}(w_t,w_t,...,w_t): t\in [0,1]$ and
$$
{\mathcal F}_1'=\{(g_{i,j})\in M_r(C(F))\cong C: g_{i,j}\in {\mathcal F}_1\}.
$$
We may assume that $P=1_A\in {\mathcal F}_1.$ For any finite subset
${\mathcal F}\subset A,$ we may assume that ${\mathcal F}\subset {\mathcal F}_1',$ if
${\mathcal F}_1$ is sufficiently large.
Note that $h'(g)P=Ph'(g)$ for all $g\in C$ (since $h$ commutes with $P$).
Note also that
\beq\label{ncl2-7}
\text{Length}(\{W_t\})\le 2\pi+\ep_1/4.
\eneq
With sufficiently small $\ep_1$ and sufficiently large ${\mathcal F}_1,$ by applying \ref{Pert},
we obtain a path of unitaries
$\{v_t: t\in [0,1]\}$ of $B$ such that
\beq\label{ncl2-8}
&&\|v_0-u\| < \ep/8,\,\,\, v_1=1_B,\,\,\,\|[h(f),v_t]\|<\ep/2\andeqn\\
&&\text{Length}(\{v_t\}) \le   2\pi+\ep/2.
\eneq

The lemma then follows by connecting $v_0$ to $u$ with
the length no more than $\ep/2.$

\end{proof}

\begin{lem}\label{NCL3}
Let $X$ be a finite CW complex, let $F\subset X$ be a compact subset
and let $A_1=QM_m(C(X))Q,$ where $m\ge 1$ is an integer and $Q\in
M_m(C(X))$ is a projection.
Let $A=PM_m(C(F))P,$ where $P=s(Q)$ and $s: A_1\to A$ is the surjective
\hm\, induced by the quotient map $C(X)\to C(F).$

Then, for any $\ep>0,$  any finite
subset ${\mathcal F}\subset A,$ there exists $\dt>0,$ a finite
subset ${\mathcal G}\subset A$  and a
finite subset ${\mathcal P}\subset \underline{K}(A)$ satisfying the
following:

Suppose that $B$ is a unital purely infinite  simple \CA,
$h: A\to B$ is a unital monomorphism and $u\in B$ is a
unitary such that
\beq\label{3ncl2-1}
\|[h(a), u]\|<\dt\tforal\, a\in {\mathcal G} \andeqn
\rm{Bott}(h,u)|_{\mathcal P}=0.
\eneq
Then, there exists a continuous rectifiable path of unitaries
$\{u_t:t\in [0,1]\}$ of $B$ such that
\beq\label{3ncl2-2}
u_0=u,\,\,\,u_1=1_B\andeqn \|[h(a),u_t]\|<\ep\tforal\, a\in
{\mathcal F}\andeqn  t\in [0,1].
\eneq
Moreover,
\beq\label{3ncl2-3}
\rm{Length}(\{u_t\})\le 2\pi+\ep\pi.
\eneq

\end{lem}

\begin{proof}
The proof is virtually identical to that of \ref{NCL2} but apply \ref{NCLpi1} instead of \ref{NCL1}.
\end{proof}

\begin{df}\label{Ddst2}
{\rm Let $A$ be a unital AH-algebra. Then $A$ has the   form
$A=\lim_{n\to\infty} (A_n, \phi_n),$ where
$A_n=\oplus_{i=1}^{r(n)}P_{i,n}M_{k(i,n)}(C(X_{i,n}))P_{i,n}$ and
$P_{i,n}\in M_{k(i)}(C(X_{i,n}))$ is a projection and $X_{i,n}$ is
a path connected finite CW complex. For the convenience, without
loss of generality, we may assume that each $\phi_n$ is unital. By
replacing $X_{i,n}$ by a compact subset of $X_{i,n},$ we may
assume  that $\phi_n$ is injective for each $n.$

Let $B$ be a unital \CA\, with a tracial state $\tau\in T(B)$ and
$\psi: A\to B$ be a positive linear map.
 Define $\tau_n$ by $\tau_n=\tau\circ \psi\circ
\phi_n.$ Define  $\mu_{\tau_n}$ as in (iii) of \ref{NCD1}. Let
$\varDelta_n: (0,1)\to (0,1)$ be an increasing map. We say that $\tau\circ \psi$
is $\{\varDelta_n\}$ distributed if $\mu_{\tau_n}$ is
$\varDelta_n$ -distributed.


}

\end{df}

\begin{thm}\label{NCT1}
Let $A$ be a unital AH-algebra,  let $\ep>0$ and ${\mathcal F}\subset
A$ be a finite subset. Let $\varDelta_n: (0,1)\to (0,1)$ be a
sequence of increasing maps.

Then there exists $\dt>0,$ a finite subset ${\mathcal G}\subset A$ and a finite
subset ${\mathcal P}\subset \underline{K}(A)$ satisfying the following: Suppose that
$B$ is a unital separable simple \CA\, with tracial rank zero, $h: A\to B$ is a unital monomorphism, and suppose that
there is a unitary $u\in B$ such that
\beq\label{NCT1-1}
\|[h(a), u]\|<\dt\rforal f\in {\mathcal G},
\rm{Bott}(h,u)|_{\mathcal P}=0\andeqn
\eneq
$\mu_{\tau\circ h}$ is $\{\varDelta_n\}$-distributed for all $\tau\in T(B).$  Then there exists a
continuous path of unitaries $\{u_t:t\in [0,1]\}$ such that
\beq\label{NCT1-2}
u_0=u,\,\,\, u_1=u,\,\,\, \|[h(a), v_t]\|<\ep\tforal f\in
{\mathcal F}\andeqn t\in [0,1].
\eneq
Moreover,
\beq\label{NCT1-3}
&&\|u_t-u_{t'}\|\le (2\pi+\ep)|t-t'|\tforal t,t'\in [0,1]\\\label{NCT1-4}
&&\andeqn
\rm{Length}(\{u_t\})\le 2\pi+\ep.
\eneq

\end{thm}

\begin{proof}
Let $\ep>0$ and ${\mathcal F}\subset A$ be given.
Write $A=\lim_{n\to\infty}(A_n, \phi_n),$ where
each $A_n=\oplus_{i=1}^{r(n)}P_{i,n}M_{k(i)}(C(X_{i,n}))P_{i,n},$ where each $X_{i,n}$ is a
connected finite CW complex. By replacing $X_{i,n}$ by one of its compact subset
$F_{i,n},$ we may write $A=\overline{\cup_n A_n}.$
Thus, without loss of generality, we may assume that ${\mathcal F}\subset A_n$ for some
integer $n\ge 1.$ We then apply \ref{NCL2}.

Note to get \ref{NCT1-3}, we can either apply \ref{eqexp} or state an improved version of
\ref{NCT2} by
applying \ref{MT3+Lip} instead of \ref{MT3}.
\end{proof}



%

\begin{cor}\label{NCT1C}
Let $A$ be a unital AH-algebra, let $\ep>0$ and ${\mathcal F}\subset
A$ be a finite subset. Suppose that
$B$ is a unital separable simple \CA\, with tracial rank zero and $h: A\to B$ is a unital monomorphism
Then there exists $\dt>0,$ a finite subset ${\mathcal G}\subset A$ and a finite
subset ${\mathcal P}\subset \underline{K}(A)$ satisfying the following: Suppose that
there is a unitary $u\in B$ such that
\beq\label{NCT1C-1}
\|[h(a), u]\|<\dt\tforal f\in {\mathcal G}\andeqn
\rm{Bott}(h,u)|_{{\mathcal P}}=0.
\eneq
 Then there exists a
continuous path of unitaries $\{u_t:t\in [0,1]\}$ such that
$$
u_0=u,\,\,\, u_1=u, \|[h(a), v_t]\|<\ep\rforal f\in {\mathcal
F}\andeqn t\in [0,1],
$$
$$
\|u_t-u_{t'}\|\le (2\pi+\ep)|t-t'|\tforal t,t'\in [0,1] \andeqn
$$
$$
\rm{Length}(\{u_t\})\le 2\pi+\ep.
$$

\end{cor}

\begin{thm}\label{NCT2}
Let $A$ be a unital AH-algebra,  let $\ep>0$ and ${\mathcal
F}\subset A$ be a finite subset. Then there exists $\dt>0,$ a
finite subset ${\mathcal G}\subset A$ and a finite subset
${\mathcal P}\subset \underline{K}(A)$ satisfying the following:
Suppose that $B$ is a unital purely infinite  simple \CA\,  $h:
A\to B$ is a unital monomorphism and suppose that there is a
unitary $u\in B$ such that
\beq\label{NCT2-1}
\|[h(a), u]\|<\dt\rforal f\in {\mathcal G}, \rm{Bott}(h,u)|_{\mathcal P}=0.
\eneq
Then there exists a
continuous path of unitaries $\{u_t:t\in [0,1]\}$ such that
$$
u_0=u,\,\,\, u_1=u, \|[h(a), v_t]\|<\ep\tforal f\in {\mathcal F}
\andeqn \tforal t\in [0,1]
$$
and
$$
\|u_t-u_{t'}\|\le (2\pi+\ep)|t-t'|\tforal t,t'\in [0,1].
$$
Consequently,
$$
{\rm{Length}}(\{u_t\})\le 2\pi+\ep.
$$

\end{thm}

\begin{proof}
The proof follows the same argument used in \ref{NCT1} but applying \ref{NCL3}.
\end{proof}

The additional significance  of the following theorem is that,
when $A$ is also assumed to be simple, the ``measure-distribution"
is not needed in the assumption.

\begin{thm}\label{NCT3}
Let $A$ be a unital simple AH-algebra. For any $\ep>0$ and any
finite subset ${\mathcal F}\subset A,$ there exists $\dt>0,$  a
finite subset ${\mathcal G}\subset A$ and a finite subset ${\mathcal P}\subset \underline{K}(A)$ satisfying the following:
Suppose that $B$ is a unital simple \CA\, of tracial rank zero, or
$A$ is a purely infinite simple \CA, suppose that $h: A\to B$ is a
unital \hm\, and suppose that there is a unitary $u\in B$ such
that
$$
\|[h(a), u]\|<\dt\andeqn \rm{Bott}(h, u)|_{\mathcal P}=0.
$$
Then there exists a rectifiable continuous path of unitaries
$\{u_t: t\in [0,1]\}$ such that
$$
u_0=u,\,\,\, u_1=1_B, \|[h(a), u_t]\|<\ep\rforal f\in {\mathcal F}
$$
and
$$
\|u_t-u_{t'}\|\le (2\pi+\ep)|t-t'|\tforal t, t'\in [0,1].
$$
Consequently
$$
{\rm{Length}}(\{u_t\})\le 2\pi+\ep.
$$

\end{thm}

\begin{proof}
We will prove the case that $B$ is a unital separable simple \CA\, with tracial rank zero.
The case that $B$ is a unital purely infinite simple \CA\, follows from \ref{NCT2}.

In the definition \ref{Ddst2}, let $B=A$ and $\psi={\rm id}_A.$
Then since $A$ is a unital simple \CA\, there is a sequence $\varDelta_n: (0,1)\to (0,1)$ such that
$\tau\circ \psi$ is $\{\varDelta_n\}$-distributed for all $\tau\in T(A).$

For any unital \CA\, $B$ and any $t\in T(B)$ and any unital \hm\, $h: A\to B,$
$t\circ h=\tau$ for some $\tau\in T(A).$ Thus
$t\circ h$ is $\{\varDelta_n\}$-distributed. In other words, once $A$ is known, $\{\varDelta_n\}$
can be determined independent of $B$  and $h.$ Thus the Theorem follows \ref{NCT1}.

\end{proof}

\section{Concluding remarks}

\begin{NN}

{\rm  Suppose that $X=\sqcup_{i=1}^n X_i,$ where each $X_i$ is a
connected compact metric space. We can define
\beq\nonumber
\underline{L}_{p}(X)&=&\max_{1\le i\le n}
\underline{L}_p(X_i),\\\nonumber \overline{L}_p(X)&=&\max_{1\le
i\le n} \overline{L}_p(X_i)\andeqn\\\nonumber
 L(X)&=&\max_{1\le i\le n}
L(X_i)
\eneq

A modification of the proofs of \ref{GHP3} shows that the
following version of Theorem \ref{GHP3}:

\begin{thm}\label{NGHP3}
Let $X$ be a compact metric space with finitely many path
connected components and let $A$ be a unital purely infinite
simple \CA. Suppose that $h_1, h_2: C(X)\to A$ are two unital \hm
s such that
$$
[h_1]=[h_2]\,\,\,\text{in}\,\,\, KL(C(X),A).
$$
Then, for any $\ep>0$ and any finite subset ${\mathcal F}\subset
C(X),$ there exist two unital \hm s $H_1, H_2: C(X)\to C([0,1],
A)$ such that
$$
\pi_0\circ H_1\approx_{\ep/3} h_1,\pi_1\circ H_1\approx_{\ep/3}
\pi_0\circ H_2\andeqn \pi_1\circ H_2\approx_{\ep/3}
h_2\,\,\,\text{on}\,\,\,{\mathcal F}.
$$
Moreover,
\beq\nonumber
\overline{\rm{Length}}(\{\pi_t\circ H_1\})&\le& L(X)(1+2\pi)+\ep/2
\andeqn\\
\overline{\rm{Length}}(\{\pi_t\circ H_2\}) &\le & L(X)+\ep/2.
\eneq
\end{thm}

Similar statements also hold for Theorem \ref{SHLn}, Theorem
\ref{GHTh2}, Theorem \ref{SHLpiinj}, Theorem \ref{SHLpi}, Theorem
\ref{SHL}, Theorem \ref{SHLoldpi}, Theorem \ref{suppiT} and
Theorem \ref{SHLold}.

For example, \ref{SHLn} holds for general finite CW complex $X:$

\begin{thm}\label{SHLnF}
Let $X$ be a  finite CW complex and let $A$ be a unital separable
simple \CA\, with tracial rank zero. Suppose that $\psi_1, \psi_2:
C(X)\to A$ are two unital monomorphisms such that
\beq\label{Fshln1}
[\psi_1]=[\psi_2]\,\,\,\text{in}\,\,\, KL(C(X), A).
\eneq
Then, for any $\ep>0$ and any finite subset ${\mathcal F}\subset
C(X),$
 there exists a unital \hm\, $H: C(X)\to C([0,1], A)$ such that
\beq\label{Fshln2}
\pi_0\circ H\approx_{\ep} \psi_1\andeqn \pi_1\circ H\approx_{\ep}
\psi_2\,\,\,\text{on}\,\,\,{\mathcal F}.
\eneq
Moreover, each $\pi_t\circ H$ is a monomorphism and
\beq\label{Fshln3}
\overline{\rm{Length}}(\{\pi_t\circ H\})\le
L(X)+\underline{L}_p(X)2\pi+\ep.
\eneq
\end{thm}

}

\end{NN}

\begin{NN}

{\rm One may also define $L(X)$ for general compact metric spaces.
There are at least two possibilities. One may consider the
following constant:
$$
\sup\{{\rm{Length}}(\{\gamma\}): \text{all rectifiable continuous
paths in}\,\,\, X\}.
$$

If $X=\lim_{\leftarrow\,\,n }X_n,$ where each $X_n$ is a finite CW
complex, then one may consider the constant:
$$
\inf\sup\{L(X_n): n\in \N\},
$$
where the infimum is taken among all possible reverse limits.

It seems that the second constant is easy to handle for our
purpose. However, in general, one should not expect such constant
to be  finite.

}
\end{NN}

\begin{NN}
{\rm The last section presents some versions of the Basic Homotopy
for AH-algebras. It seems possible that more general situation
could be discussed. More precisely, for example, one may ask if
\ref{NCT1C} holds for general separable amenable \CA\, $A.$
However, it is beyond the scope of this research to answer this question.
While we know
several cases are  valid, in general, further
study  is required.

There are several valid versions of  Super Homotopy Lemmas for
non-commutative cases and could be presented here.  However, given
the length of this work, we choose not include them.

}
\end{NN}

\begin{NN}
{\rm Applications of the results in the work have been mentioned
in the introduction. It is certainly desirable to present some of
them here. Nevertheless, the length of this work again prevent us
from including them here without further prolonging the current
work. We would like however, to mention that, for example, Theorem
\ref{NCT1} and Theorem \ref{NCT1C} are used in \cite{Lnemb3},
where we show crossed products of certain AH-algebras by $\Z$ can
be embedded into a unital simple AF-algebras. It also plays an
essential role in the study of  asymptotic unitary equivalence of
monomorphisms(\cite{Lnasym}).

}
\end{NN}








\backmatter
\bibliographystyle{amsalpha}
\bibliography{}

\printindex

\end{document}